%% file: SimSze9.tex
\documentclass{book}
\usepackage{amsthm,amsmath,amssymb}

\numberwithin{equation}{section}

\usepackage{graphics,epsfig,boxedminipage,wrapfig}
\usepackage[bookmarks]{hyperref}
\usepackage[utf8]{inputenc}

\usepackage{color}

\makeatletter
\renewcommand{\@seccntformat}[1]{\@nameuse{the#1}.\quad}
\makeatother

\gdef\silent#1\par{}

\renewcommand\thesection{\arabic{section}}

\newtheorem{theorem}{\sc Theorem}[section]
\newtheorem{corollary}[theorem]{\sc Corollary}
\newtheorem{conjecture}[theorem]{\sc Conjecture}

\newtheorem{problem}[theorem]{\sc Problem}

\theoremstyle{definition}

\newtheorem{definition}[theorem]{\sc Definition}
\newtheorem{remark}[theorem]{\sc Remark}
\newtheorem{remarks}[theorem]{\sc Remarks}

\gdef\Proclaim#1. #2\par{\goodbreak\medskip\noindent{\sc #1.}~{\sl #2\par}\medskip}

\begin{document}

\newtheorem{megjegy}[section]{Remarks}

\gdef\megjegyname{\sc Remarks-X}
\renewcommand{\themegjegy}{\thesection.\arabic{remark}.}%

\gdef\proof. {\medbreak\noindent{\sc Proof. \enspace}}

\gdef\girth{{\bf girth}}
\gdef\circ{{\bf circ}}
\def\pcup{\buildrel{\bullet}\over {\bigcup}}

\renewcommand{\theproblem}{\thesection.\arabic{problem}.}%

\newtheorem{construction}[theorem]{\sc Construction}
\newtheorem{MetaTheorem}[theorem]{\sc Meta-Theorem}
\newtheorem{MetaConj}[theorem]{\sc Meta-Conjecture}

\newtheorem{defn}[theorem]{\sc Definition}
\newtheorem{examp}[theorem]{{\sc Example}}

\newtheorem{questio}[theorem]{\sc Question}
\newtheorem{HistRem}[theorem]{\sc Historical Remarks}

\def\Rd{\color[rgb]{0.7,0,0}}

\def\Vo#1 {{\bf #1}}
\def\vo#1(#2){{\bf #1} (#2)}
\def\voy#1(#2) (#3){{\bf #1 (#2)} (#3)}

\def\Arx#1 {{\Bl\margofont ~ #1.}}

\def\Arx#1 {}

\gdef\girth{{\bf girth}}
\gdef\circ{{\bf circ}}

% --------- Silent-ek --------

\font\margoFont=cmr6
\font\margofont=cmbx7
\newcount\margonum
\newcount\margonumma

\def\SilentFig#1\par{\marginpar{\fbox{Hidden figure}}\par}

\gdef\dense{\itemsep=0pt\parskip=0pt}

\def\chapname{}

\makeatletter
\renewcommand{\theenumi}{(\alph{enumi})}
\renewcommand{\labelenumi}{\theenumi}
\renewcommand{\theenumii}{(\roman{enumii})}
\renewcommand{\labelenumii}{\theenumii}
\renewcommand{\theenumiii}{[\arabic{enumiii}]}
\renewcommand{\labelenumiii}{\theenumiii}
\makeatother

\gdef\densebox#1{{\baselineskip=5pt\tiny #1\par}}
\gdef\densebox#1{\medskip ~\hskip-9mm\begin{boxedminipage}{\the\hsize}\baselineskip=8.5pt\margofont #1\end{boxedminipage}\medskip}

\gdef\DenseBox#1{\par}

\gdef\DenseBox#1{\medskip \hrule \vskip2pt {\baselineskip=8.5pt\margofont\Bl #1 \par}\vskip2pt\hrule \medskip}

\gdef\xx{{\mathbf x}}

\gdef\ev#1){#1),}

%% LovaszSurvey vege

\def\beq#1{\begin{equation}\label{#1}}
\def\eeq{\end{equation}}

\def\beqn#1{\begin{eqnarray}\label{#1}}
\def\eeqn{\end{eqnarray}}

\def\beqns{\begin{eqnarray*}}
\def\eeqns{\end{eqnarray*}}

\def\text#1{~\mbox{#1}~}
\def\Text#1{\qquad\mbox{#1}\qquad}

\hyphenation{Della-monica Raz-bo-rov}

%% ================== Szegmentalas, section ==========

\gdef\PaperName{Simonovits-Szemerédi} % ezt az elejere kell beirni, a macro-k behivasa ele.

\gdef\Section#1#2{\section{#1}\label{#2}\renewcommand\rightmark{\thesection.\qquad\jobname~\quad #1}}

\silent

\gdef\Section#1#2{
 \pagestyle{myheadings}
\section{#1}
\setcounter{theorem}{0}
\label{#2} 
 \markboth{\PaperName}{%%\PaperName: 
\thesection.~ #1,\qquad %%\Ugrato
%% \hfill
{~~\quad ~{~}} %%\newpage 
}}

% ----------- vegyes kepletek -------

\gdef\eqref#1{{\rm(\ref{#1})}}

\gdef\loglog{\log\log }
% ------- Grafelmelet ------

\gdef\cA{{\mathcal A}}
\def\cB{{\mathcal B}}
\def\cC{{\mathcal C}}
\gdef\cD{{\mathcal D}}
\def\cE{{\mathcal E}}

\def\bF{{\mathbb F}}
\def\cF{{\mathcal F}}
\def\cG{{\mathcal G}}
\def\cK{{\mathcal K}}

\gdef\bH{{\mathbb H}}
\gdef\cH{{\mathcal H}}
\gdef\cI{{\mathcal I}}
\gdef\calI{{\mathcal I}}

\gdef\cP{{\mathcal P}}
\gdef\cQ{{\mathcal Q}}
\gdef\cS{{\mathcal S}}
\gdef\cT{{\mathcal T}}

\gdef\EE{{\mathbb E}}
\gdef\LL{{\mathbb L}}
\gdef\NN{{\mathbb N}}
\gdef\PP{{\mathbb P}}
\gdef\RR{{\mathbb R}}

\def\VV{{\Rd\mathbb V}}
\def\WW{{\mathbb W}}
\def\ZZ{{\mathbb Z}}

\gdef\LiL{{L\in\mathcal L}}

\gdef\Cc#1{{C_{#1}}}

 \gdef\P#1{{P_{#1}}}

\gdef\cL{{\mathcal L}}
\gdef\cM{{\mathcal M}}

\def\Hm{{H_m}}
\def\Hk{{H_k}}
\def\Sm{{S_m}}

\gdef\Tn#1,#2{{T_{#1,#2}}}

\def\Pk{{P_k}}
\def\Kp{{K_p}}
\gdef\Tk{{T_k}}

\def\Kfm{{K_4^-}}

\def\De{{\Delta}}

\def\de{\delta }
\def\eps{{\varepsilon}}
\def\la{\lambda }
\def\teta{{\vartheta}}

\def\be{\beta }
\def\ga{{\gamma}}
\def\ka{\kappa }

\gdef\dext{{\bf dex}}

\gdef\EXT{{\bf EX}}
\gdef\RT{{\bf RT}}
\gdef\turtwo#1{\lfloor {#1^2\over 4~}\rfloor}
\gdef\Tur#1#2{\left(1-{1\over #2}\right){#1\choose2}}
\gdef\tur#1#2{\left(1-{1\over #2}\right){#1}}

\gdef\K#1{{K_{#1}}}
\def\KK#1#2{{K^{(#1)}_{#2}}}

\gdef\Chyp#1#2{{\mathcal C}_{#1}^{#2}}
\gdef\Cthyp#1#2{\widetilde{\mathcal C}_{#1}^{#2}}
\gdef\Fhyp#1#2{{\mathbb F_{#1}^{(#2)}}}
\gdef\Hhyp#1#2{{\mathbb H_{#1}^{(#2)}}}
\gdef\tHhyp#1#2{{\widetilde{\mathbb H}}_{#1}^{(#2)}}
\gdef\Khyp#1#2{{\mathbb K_{#1}^{(#2)}}}
\gdef\Lhyp#1#2{{\mathbb L_{#1}^{(#2)}}}
\gdef\Krr{{K_r^{(r)}}}

\def\down{\hfill\break}
\def\ddori{\hfill\break \hbox{~~~~~~~~~~} }
\def\dDori{\hfill\break \hbox{~~~} }
\def\dorix{\hfill\break \hbox{~~~~~~~~~~}--\quad }
\def\dori{\hfill\break\phantom{~}\hskip13pt } 
\def\dorib{\hfill\break\hbox{\quad} $\bullet$\hbox{~~} }
\def\dorim#1 #2 {\hfill\break\phantom{~}\hskip-#1pt#2\quad}

\gdef\abrax#1#2{\epsfig{file=#1.eps,width=#2mm}}

\gdef\BalAbra#1#2{\begin{wrapfigure}{L}{#2mm}\epsfig{file=#1.eps,width=#2mm}
\end{wrapfigure}}

\gdef\BalAbraCap#1#2{\setlength\intextsep{4pt}  
% a label ugyana, mint az eps-nev
\begin{wrapfigure}{l}{0.2\textwidth}
\includegraphics[width=0.2\textwidth]{#1.eps}
\caption{#2~}
\label{#1}
 \end{wrapfigure}}

\gdef\BalAbraCapMedium#1#2{\setlength\intextsep{4pt}  
% The name in the label is the same as the epsname
\begin{wrapfigure}{l}{0.3\textwidth}
\includegraphics[width=0.3\textwidth]{#1.eps}
\caption{#2~}
\label{#1}
 \end{wrapfigure}}

\gdef\BalAbraCapMed#1#2#3{\setlength\intextsep{4pt}  
% a label ugyana, mint az eps-nev
\begin{wrapfigure}{l}{#3mm}
\includegraphics[width=#3mm]{#1.eps}
\caption{#2~}
\label{#1}
 \end{wrapfigure}}

\gdef\BalAbraX#1#2#3\par{\noindent\epsfig{file=#1.eps,width=30mm}~\vskip-#2mm\hfill \begin{minipage}{100mm} #3 \end{minipage} \medskip
}

\gdef\BalAbraD#1#2#3{\begin{wrapfigure}{L}{#2mm}
\epsfig{file=#1.eps,width=#3mm}
\epsfig{file=#2.eps,width=#3mm}
\end{wrapfigure}}

\gdef\ok{\,\fbox{\tiny $\surd$}}
\gdef\qu{\,\fbox{\small $?$}}

\def\emC{\sc}
\def\Sc#1 {{\sc #1} }

\def\Style{\global\advance \margonum by 1
\, \fbox{\margofont\the\margonum}
\marginpar{\epsfig{file=Rdanger.eps,height=6mm} {
 \fbox{\tiny\the\margonum}\\ \margofont Style}}}

\def\Danger{\global\advance \margonum by 1
\, \fbox{\margofont\the\margonum}
\marginpar{~~~\epsfig{file=Rdanger.eps,height=6mm} {
 \fbox{\tiny\the\margonum} \\\phantom{~} ~\margofont Sharp curve}}} 

\def\subsub#1//{\medskip\noindent {\bf#1~}}

\def\IttSep#1{
\noindent~~~~~~\rule{50mm}{0.3mm}\fbox{#1}\rule{50mm}{0.3mm}
\marginpar{~\hskip10mm$\longleftarrow$\fbox{\bf Itt tartok}\,\fbox{#1}\label{Itt#1}}}

\def\bL{{\Rd\mathbf L}}

\newcount\sepcikk

\def\th{^{\rm th}}

\def\La{{\Prp \Lambda}}

\def\RR{{\mathbb R}}
\def\XX{{\mathbb X}}

\def\ti{\tilde}
\def\Ti{\widetilde}

\def\mindeg{d_{\sf min}}
\def\maxdeg{d_{\sf max}}

\def\margo#1{\marginpar{{\sf {\small \raggedright #1\\ }}}}

\def\Unfinished{
\global\advance \margonum by 1 
{\fbox{\tiny \the\margonum}
\epsfig{file=Rdanger.eps,height=4mm} }
 \marginpar{\epsfig{file=Rdanger.eps,height=6mm} 
{\margofont
\fbox{\the\margonum}\\ Unfinished}}}

\def\reci#1 {{1\over #1}}
\def\Reci#1{{1\over #1}}
\def\half{\frac{1}{2}}

\newcount\margonum

\def\sepa{\smallskip\hrule\smallskip}

\def\minisepa{\medskip\centerline{--- $\cdot$ ---}\smallskip}

\gdef\Code#1#2{\marginpar{\margofont Label: \\ {\sc #1}\\ File: \\ {\bf #2}}}

\def\Ref#1{\ref{#1} \fbox{#1}}
\gdef\oms{\omega^*}
\gdef\YY{{\mathbb Y}}

\gdef\G{{G}}

\def\UiUj{{(U_i,U_j)}}

\def\ViVj{{(V_i,V_j)}}

\def\proofend{~\ifhmode \unskip\nobreak\hfill\vrule height5pt width4pt depth2pt\medskip\fi \ifmmode{\vrule height5pt width4pt depth2pt}\fi}

\def\mqed{\hfill{\Box}\medskip}
\def\qed{\hfill$\Box$\medskip}
\def\Qed{\hfill$\Box$\medskip}

\gdef\Bibitem#1#2\par{\bibitem{#1}{\baselineskip=8pt\Orange \margofont #2\par}\par}

\gdef\Bibitem#1#2\par{}

\newcount\megv
\gdef\Megvan#1{\global\advance\megv by 1 \marginpar{\margofont \the\megv:\\ #1}}

\gdef\Megvan#1{\global\advance \megv by 1~\fbox{\the\megv:~\Bl\margofont #1}}

\gdef\Megvan#1{} %% Kikapcs

%% Ez a Sarkozy javitasokhoz van visszateve.

\def\Margo#1{{
\global\advance \margonum by 1 ${}^{[\the\margonum]}$
\marginpar{~\hskip10mm\begin{boxedminipage}{35mm}
\epsfig{file=maple.eps,height=8mm}
\vskip-4mm\hskip8mm{\tiny~~(\the\margonum)}\\
\Rd  {\margofont\baselineskip=10pt \raggedright #1\\}
\end{boxedminipage}}}}

\gdef\calH#1#2{{\mathcal H}_{#2}^{(#1)}}
\def\minco{\delta_2}

\gdef\cR{{\Rd\mathbb R}}

\gdef\ckH#1{{\mathbb H^{(k)}_{#1}}}

\def\codeg{{\delta}_2}

\newcount\Bib
\def\cFa{{{\mathbb F}_7^{(3)}}}

\def\onep{(1+\eta)}

\def\Ckk{{C_{2k}}}

%% %% \gdef\unif#1{^{(#1)}}
%% %% \gdef\runif#1{_#1^{(#1)}}

%% \def\Exercise #1. #2\par{\goodbreak\smallskip\noindent{\bf Exercise\enspace}{#2}\hfill\fbox{{#1}}\par\smallskip}

\def\Imply{\quad\Longrightarrow\quad }

\newcount\thx
\newcount\prx

\def\sz{\global\advance\thx by 1\hfill\fbox{\the\thx}~}

 \def\psz{\global\advance\prx by 1\hfill\fbox{\tiny\rmPr\the\prx}~}

\def\ext{{\bf ex}}

\def\cLs{{\Rd{\mathcal L}^*}}

\def\Tm{{T_m}}

\def\rL{{L}}

\def\Xe{{\sfX}}

\def\Badly1{{\Rd\sc Badly}} % ???

\def\TMC1{{TMC}}

\newenvironment{BegsThm}[1]{\medskip \noindent{\Rd\sc Theorem #1 of {\bf [BEGS]}.} \sf\em\Bl } {\medskip}

\renewenvironment{quote}[1]{{
\parindent=0pt
\parskip=7pt
\leftskip=12pt\rightskip=12pt\baselineskip=7pt
%% \medskip
\small #1\par } \leftskip=0pt\rightskip=0pt
\parindent=15pt
\parskip=0pt
\par\medskip }

\def\fele#1{\left\lfloor{#1\over 2}\right\rfloor}
\def\felen#1{\left\lfloor\displaystyle{{#1\over 2}}\right\rfloor}

\def\chiSI{{\chi_{SI}}}
\def\chiS{{\chi_S}}

%% \def\X{{\Rd\chi}}
%% \def\Xe{{\sfX}}

%% \def\BC1{{\Bl Badly Coloured}}
%% \def\Well1{{\Rd\sc Well}}

%% \def\embed{{\hookrightarrow\, }}

%% \linenumbers
%% \nolinenumbers\leftskip=-1cm\rightskip=-26mm

\def\separule{\bigskip\hrule\vskip2pt\hrule\bigskip}

\def\boxit#1{\medskip \noindent\begin{boxedminipage}{350pt}#1\end{boxedminipage}\medskip}

\gdef\boxitS#1{} %  legalabbis atmenetileg silent
%% {\marginpar{~~~~~~~~~~~~~~~~~~~\epsfig{file=rbox.eps,height=8mm}}}

\def\NarrowBox#1#2#3{\medskip\noindent\begin{boxedminipage}{#1mm}{\DBl\fbox{\margofont #2}}\vskip1mm \Rd\margofont\baselineskip=11pt #3\end{boxedminipage}\medskip}

\def\NarrowBoxDan#1#2{\global\advance\dcount by 1
\medskip\noindent
\begin{boxedminipage}{\the\textwidth}\vskip1mm 
\DBl\fbox{\small\the\dcount:~\margofont #1}
\Rd\margofont\baselineskip=11pt #2\end{boxedminipage}\marginpar{~~~\epsfig{file=Rdanger.eps,height=13mm}}\medskip}

\def\narrowbox#1#2{\medskip 
\noindent\begin{boxedminipage}{#1mm}\Rd\margofont\baselineskip=11pt #2\end{boxedminipage}\medskip}

\def\Narrowbox#1{\medskip \dcount=0
\noindent\begin{boxedminipage}{300pt}#1\end{boxedminipage}\medskip}

\title{\bf Embedding graphs into larger graphs: results, methods, and problems}
\author{Miklós Simonovits and Endre Szemerédi}
\date{\small Alfréd Rényi Institute of Mathematics,\\ Hungarian Academy of Sciences, \\ ~\\ email: miki@renyi.hu, szemered@rutgers.edu }

\renewcommand{\thesection}{\arabic{section}}

\renewcommand{\theproblem}{\thesection.\arabic{problem}.}%

\renewcommand\leftmark{Miklós Simonovits and Endre Szemerédi: ~\qquad\ March29, 2019~}

\gdef\Final{\DBl\sf }

\maketitle

\tableofcontents

\newpage 

\hfill{Dedicated to Laci Lovász}

\hfill{on his 70th birthday}

\subsub Keywords:// Extremal graphs, Stability, Regularity Lemma, Semi-random methods, Embedding. {\em Subject Classification:} Primary 05C35; Secondary 05D10, 05D40, 05C45, 05C65, 05C80.

\Section{Introduction}{Intro}

%%\ResetCounters

%% \renewcommand\rightmark{\thesection: Introduction: x}

%% \newtheorem{remark}[remarks]{\sc Remarks}

%%\renewcommand{\theremarks}{\thesection.\arabic{remark}.}%

We dedicate this survey paper to Laci Lovász, our close friend, on the occasion of his 70th birthday. Beside learning mathematics from our professors, we
learned a lot of mathematics from each other, and we emphasise that we learned a lot of beautiful mathematics from Laci.

Extremal Graph Theory is a very deep and wide area of modern combinatorics. It is very fast developing, and in this long but relatively short survey we select some of those results which either we feel very important in this field or which
are {\em new breakthrough} results, or which -- for some other reasons -- are very close to us. Some results discussed here got stronger emphasis, since they are connected to Lovász (and sometimes to us).\footnote{We shall indicate the given names mostly in case of ambiguity, in cases where there are two mathematicians with the same family name, (often, but not always, father and son). We shall ignore this ``convention'' for Erdős, Lovász and Turán.} The same
time, we shall have to skip several very important results; often we just start describing some areas but then we refer the reader to other surveys or research papers, or to some survey-like parts of research papers. (Fortunately,
nowadays many research papers start with excellent survey-like introduction.)

\begin{quote}{
 Extremal graph theory became a very large and important part of Graph Theory, and there are so many excellent surveys on parts of it that we could say that this one is a ``survey of surveys''. Of course, we shall not try to cover the whole area, that would require a much longer survey or a book.

Also we could say that there are many subareas, ``rooms'' in this area, and occasionally we just enter a ``door'', or ``open a window'' on a new area, point out a few theorems/phenomena/problems, explain their essence, refer to some more detailed surveys, and move on, to the next ``door''.

It is like being in our favourite museum, having a very limited time, where we must skip many outstanding paintings. The big difference is that here we shall see many very new works as well.}
\end{quote}

One interesting feature of Extremal Graph Theory is its very strong connection and interaction with several other parts of Discrete Mathematics, and more generally, with other fields of Mathematics. It is connected to Algebra, Commutative Algebra, Eigenvalues, Geometry, Finite Geometries, Graph
Limits (and through this to Mathematical Logic, e.g., to Ultra Product, to Undecidability), to Probability Theory, application of Probabilistic Methods, to the evolution of Random Structures, and to many other topics. It is also strongly connected to Theoretical Computer Science through its methods,
(e.g., using random and pseudo-random structures), expanders, property testing, and also it is strongly connected to algorithmic questions. This connection is fruitful for both sides.

One reason of this fascinating, strong interaction is that in Extremal Graph Theory often seemingly simple problems required the invention of very deep new methods (or their improvements). Another one can be that combinatorial methods start becoming more and more important in other areas of
Mathematics. A third reason is, perhaps, -- as Turán thought and used to emphasise, -- that Ramsey's theorem and {\em his theorem} are applicable because they are generalizations of the Pigeon Hole Principle. Erdős wrote several papers on how to apply these theorems in Combinatorial Number Theory,
(we shall discuss \cite{Erd38Tomsk}, but see also, e.g., \cite{Erd1968-09,Erd1969-14,Erd1970-20,Erd1971-13,Erd1980-39,Erd1981-29}, Erdős, A. Sárközy, and T. Sós \cite{ErdSarkoSos95Prod}, Alon and Erdős \cite{AlonErd85Sidon}).\footnote{Sometimes we list papers in their time-order, in some other cases
  in alphabetical order.} Beside putting emphasis on the results we shall emphasise the development of methods also very strongly. We shall skip discussing our results on the Erdős-Sós and Loebl-Komlós-Sós conjectures on tree embedding, since they are well described in \cite{FureSim13Degen},
\footnote{The proof of the approximate Loebl-Komlós-Sós conjecture was first attacked in \cite{AjtKomSzem92Loebl} and then proved in a sequence of papers, from Yi Zhao \cite{YiZhao11Loebl}, Cooley \cite{Cooley09LoeblKS}, \dots Hladk\'y, KomPiguet, Simonovits, M. Stein, \cite{HladKomPigu17-1}-\cite{HladKomPigu17-4}.}  however, we shall
discuss some other tree-embedding results. We are writing up the proof of the Erdős-Sós Conjecture \cite{AjtKomSimSzem}.\footnote{The first result on the Loebl Conjecture was an ``approximate'' solution of Ajtai, Komlós, and Szemerédi \cite{AjtKomSzem92Loebl}.}

Since there are several excellent books and surveys, like Bollobás \cite{Bollo78ExtreBook,Bollo95Hand}, Füredi \cite{Fure91Sur,Fure94Zurich}, Füredi-Simonovits \cite{FureSim13Degen}, Lovász \cite{Lov79CombExerc,Lov12LimitBook}, Simonovits \cite{Sim77TurInflu,Sim97ErdInfluPrague,Sim99Stirin,Simo13ErdInfluPragueB}, Simonovits and Sós \cite{SimSosV01RT}, about Extremal Graph Theory, or some parts of it, here we shall often shift the discussion into those directions which are
less covered by the ``standard'' sources.\footnote{Essential parts of this survey are connected to Regularity Lemmas, Blow-up lemmas, applications of Absorbing techniques, where again, there are several very important and nice surveys,
  covering those parts, e.g., Alon \cite{Alon95HB}, Gerke and Steger \cite{GerkeSteg05SparseRegu}, Komlós and Simonovits \cite{KomSim96SzemRegu}, Kühn and Osthus \cite{KuhnOsthus09Embedding,KuhnOsthus14ICM}, Rödl and Ruciński
  \cite{RodlRuc10HyperSurv}, Steger \cite{Steger14DeterRandom}, and many others.} Also, we shall emphasise/discuss the surprising connections of Extremal Graph Theory and other areas of Mathematics. We had to leave out several important new
methods, e.g,
\begin{itemize}\dense
\item We shall mostly leave out results on Random Graphs;
\item Razborov's Flag Algebras \cite{Razbor07Flag}, (see also 
Razborov \cite{Razbor08K3}, 
Pikhurko and Razborov \cite{PikhuRazbo17Struct}, 
Grzesik \cite{Grzesik12}, Hatami, Hladký, Král', Norine, and Razborov, \cite{HatamiHladKralNorineRazbor12}, 
 and many others);
\item the Hypergraph Container Method, ``recently'' developed, independently, by J. Balogh, R. Morris, and W. Samotij \cite{BalMorrSamo15Indep}, and by D. Saxton and A.~Thomason \cite{SaxtonThomason15Container}. Very recently Balogh,
  Morris, and Samotij published a survey on the Container method in the ICM volumes \cite{BalMorrSamo18ICM_Contain}; \item the method of Dependent Random Choice, see e.g. Alon, Krivelevich, and Sudakov \cite{AlonKriveSudakov03Dependent} or,
  for a survey on this topic see J. Fox and B. Sudakov
\cite{FoxSudak11Dependent}, or many other sources\dots 
\item and we have to emphasise, we had to leave out among others almost everything connected to the Universe of Integers, e.g., Sum-free sets, the 
 Cameron-Erdős conjecture, the Sum-Product problems, \dots

\end{itemize} 

Further, we skipped several areas, referring the readers to more authentic sources. Thus, e.g., one of the latest developments in Extremal Graph Theory is the surprising, strong development in the area of Graph Limits, coming from several sources. A group of researchers meeting originally at
Microsoft Research, started investigating problems connected to graph limits, for various reasons. We mention here only Laci Lovász and Balázs Szegedy \cite{LovSzege07Anal}, Borgs, Chayes, Lovász, Sós, and Vesztergombi \cite{BorgsChayLovSosVesz06Count}, \cite{BorgsChayLovSosVesz06Testing},
\cite{BorgsChayLovSosVesz08Conv}, \cite{BorgsChayLovSosVesz10Uniq}, \cite{BorgsChayLovSosVesz12ConvII}, and M. Freedman, Lovász and Schrijver \cite{FreedLovSchrij}. Lovász has published a 500pp book \cite{Lov12LimitBook} about this area.

\subsub How do we refer to papers?// We felt that in a survey like this it is impossible to refer to all the good papers. Also, we tried to ``introduce'' many authors. So if a paper is missing from this survey that does not mean that we felt it was not worth including it. Further, wherever we
referred to a paper, we mentioned (all) its authors, unless the paper and its authors had been mentioned a few lines earlier. (We never use {\it at al!}) In the ``References'' we mentioned the given names as well, mostly twice: (i) first occurrence and (ii) first occurrence as the first author.

\minisepa

Extremal Graph Theory could have started from a number theoretical question of Erdős \cite{Erd38Tomsk}, however, he missed to observe that this is a starting point of a whole new theory. Next came Turán's theorem \cite{TuranML}, with his questions, which somewhat later triggered a fast development
of this area. In between, starting from a topological question, Erdős and Stone \cite{ErdStone46} proved their theorem (here Theorem~\ref{ErdStoneTh}), which later turned out to be very important in this area. Among others, this easily implies the Erdős-Simonovits Limit Theorem \cite{ErdSim66Lim},
strengthened to the Erdős-Simonovits Stability theory \cite{Erd66Tihany, Sim66Tihany}, and much later led to Szemerédi's Regularity Lemma \cite{Szem78Regu}. All these things will be explained in more details.

We shall also discuss several important methods: Stability, Regularity Lemma, Blow-up Lemma, Semi-random Methods, Absorbing Lemma, and also some further, sporadic methods.

In this survey we mention hypergraph extremal problems only when this does not become {\em too technical} for most of the readers: thus we shall not discuss, among others, the very important Hypergraph Regularity Lemmas.\footnote{For a
  concise ``description'' of this topic, see Rödl, Nagle, Skokan, Schacht, and Kohayakawa \cite{RodlNagleSkokSchKoha05PNAS}, and also Solymosi \cite{Solymo05PNAS}.} %%\NagyZ{Van-e R-?}

\subsub Repetitions.// This paper ``covers'' a huge area, with a very involved structure. So we shall occasionally repeat certain assertions, to make the
paper more readable.

\subsub Notation.// Below, for a while we shall consider simple graphs, (hypergraphs, loops and multiple edges are excluded) and for graphs (or hypergraphs) the first subscript almost always denotes the number of vertices: $G_n$, $S_n$, $T_{n,p}$
\dots are graphs on $n$ vertices. There will be just two exceptions: $K_d(n_1,\dots,n_d)$ means a $d$-partite complete graph with $n_i$ vertices in its $i\th$ class; if we list a family of graphs $\cL=\{L_1,\dots,L_t\}$, then again, the
subscript is not necessarily the number of vertices. The maximum and minimum degrees will be denoted by $\maxdeg(G)$ and $\mindeg(G)$, respectively. (In the hypergraph section $\de_1(\Hhyp{}k)$ is also the minimum degree.)  $C_\ell$,
$\P\ell$, $K_\ell$ denote the cycle, path, and the complete graph of $\ell$ vertices, respectively. Further, $e(G)$, $v(G)$, $\chi(G)$ and $\alpha(G)$ denote the number of edges, vertices, the chromatic number, and the independence
number,\footnote{often called stability number.}  respectively. For a graph $G$, $V(G)$ and $E(G)$ denote the sets of vertices and edges. If $U\subseteq V(G)$ then $G[U]$ is the subgraph induced by $U$.  \medskip

\BalAbraCap{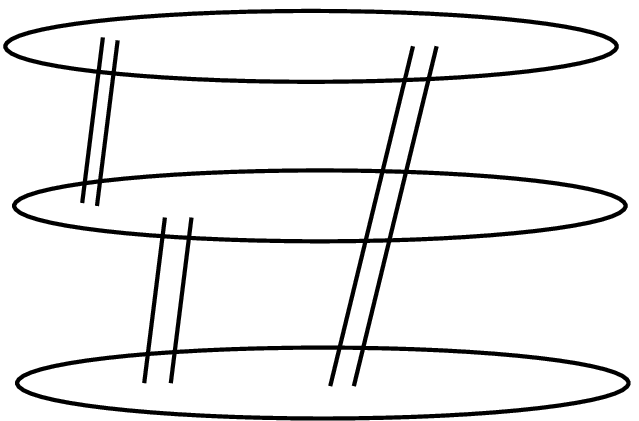}{Turán graph} Given a family $\cL$ of excluded graphs, $\ext(n,\cL)$ denotes the maximum of $e(G_n)$ for a graph $G_n$ not containing excluded subgraphs, and $\EXT(n,\cL)$ denotes the family of extremal graphs for
$\cL$: the graphs not containing excluded subgraphs but attaining the maximum number of edges.  Let $T_{n,p}$ be the Turán graph: the graph obtained by partitioning $n$ vertices into $p$ classes as equally as possible and joining two
vertices if and only if they are in distinct classes.

Given some graphs $L_1,\dots,L_r$, $R(L_1,\dots,L_r)$ is the Ramsey number: the smallest integer $R$ for which all $r$-edge-coloured $K_R$ have a
(monochromatic) subgraph $L_i$ in some colour $i$.
%% \footnote{We use both the lower Ramsey number, $R^-(\dots)=R(\dots)-1$, and this ``upper'' one.} 
If all
these graphs are complete graphs, $L_i=K_{p_i}$, then we use the abbreviation $R(p_1,\dots,p_r)$ for this Ramsey number.

\begin{remark}[Graph Sequences]
  Speaking of $o(n^2)$ edges, or $o(n)$ vertices, \dots, we cannot speak of an individual graph, only about a sequence of graphs.  In this paper, using $o(.)$, we always assume that $n\to\infty$. Also, since $A(n)+o(n)$ and $A(n)-o(n)$ mathematically are exactly the same, we shall be cautious with
  formulating some of our theorems.
\end{remark}

\Section{The beginnings}{Kezdet}

\subsection{Very early results} 

The first, simplest extremal graph result goes back to Mantel.

\begin{theorem}[Mantel, \ev(1907) \cite{Mantel07}]
  If a graph $G_n$ does not contain a $\K3$ then it has at most $\turtwo n$ edges.
\end{theorem} 

Turán, not knowing of this theorem,\footnote{The last paragraph of Turán's original paper is as follows: ``\dots Further on, I learned from the kind communication of Mr. József Krausz that the value of $d_k(n)$ given on p438 for $k=3$ was found already in 1907 by W. Mantel (Wiskundige Opgaven, vol
  10, pp. 60-61). I know this paper only from the reference Fortschritte d. Math. vol 38, p. 270.''} proved the following generalization:

~

\begin{theorem}[Turán, \ev(1941) \cite{TuranML,Tur54Col,Tur89Coll}\footnote{Turán's papers originally written in Hungarian were translated into English after his death, thus \cite{Tur89Coll} contains the English translation of \cite{TuranML}.}]\label{TurTh}
 
  If a graph $G_n$ does not contain a $\K{p+1}$ then $e(G_n)\le e(T_{n,p}),$ and the equality is achieved only for $G_n=T_{n,p}$.
\end{theorem}

The above theorems are sharp, since $T_{n,p}$ is $p$--chromatic and therefore it does not contain $\K{p+1}$. These theorems have very simple proofs and there was also an earlier extremal graph theorem discovered and proved by Erdős.  He set out from the following ``combinatorial number theory''
problem.

\begin{problem}[\sc Erdős \ev(1938) \cite{Erd38Tomsk}]\label{MultipSidon} 
  Assume that $\cA:=\{a_1,\dots,a_m\}\subseteq[1,n]$ satisfies the ``multiplicative Sidon property'' that all the pairwise products are different. More precisely, assume that \beq{MulSid}\Text{if} a_ia_j= a_ka_\ell\Text{then}\{i,j\}=\{k,\ell\}.\eeq How large can $|\cA|$ be?
\end{problem}

The primes in $[2,n]$ satisfy \eqref{MulSid}. Can one find many more integers with this ``multiplicative Sidon property''? To solve Problem \ref{MultipSidon}, Erdős proved

\begin{theorem}[Erdős, \cite{Erd38Tomsk}]\label{MultiC4}
If $G\subseteq K(n,n)$ and $C_4\not\subseteq\G$ then $e(G)\le3n\sqrt{n}.$
\end{theorem} 

This may have been the first non-trivial extremal graph result. It is interesting to remark that this paper contained the first {\em finite geometric construction}, due to Eszter~Klein, to prove the lower bound $\ext(n,C_4)\ge cn\sqrt n$.  Let $\pi(n)$ be the number of primes in $[2,n]$.  Using Theorem~\ref{MultiC4}, Erdős proved

\begin{theorem}\label{MultipSidonTh} 
 Under condition \eqref{MulSid},
$$|\cA|\le \pi(n)+O(n^{3/4}).$$
\end{theorem}

As to the sharpness of this theorem, Erdős wrote:

\begin{quote}{
 ``Now we prove that the error term cannot be better than $O(n^{3/4}/(\log n)^{3/2})$. First we prove the following lemma communicated to me by Miss E. Klein:

{\bf Lemma. }{\it Given $p(p+1)+1$ elements, ($p$ a prime), we can construct $p(p+1)+1$ combinations taken $p+1$ at a time, having no two elements in common.}''}
\end{quote}

Though the language is somewhat archaic, it states the existence of a finite geometry on $p^2+p+1$ points: this seems to be the first application of finite geometric constructions in this area. Later Erdős applied graph theory in number theory several times. Among others, he returned to
the above question in \cite{Erd1968-09} and proved that the lower bound is sharp.

Much later a whole theory developed around these types of questions. Here we mention only some results of A. Sárközy, Erdős, and Sós \cite{ErdSarkoSos95Prod}, the conjecture of Cameron and Erdős \cite{CamerErd99Matra}, and the papers of Ben Green \cite{Green04CamerErd}, and of Alon, Balogh, Morris, and Samotij \cite{AlonBalMorrSamo14Cameron}.

\begin{remark}
  (a) One could ask, what is the connection between $\ext(n,C_4)$ and Theorem \ref{MultipSidonTh}. Erdős wrote each non-prime integer $a_i\in\cA$ as $a_i=b_{j(i)}d_{j(i)}$, where $b_j$ are the primes in $[n^{2/3},n]$ and all the integers in $[1,n^{2/3}]$, and $d_j$ are the integers in
  $[1,n^{2/3}]$ . So the non-prime numbers defined a bipartite graph $G[\cB,\cD]$, and each $a_i$ defined an edge in it.  If this graph contained a $C_4$, the corresponding four integers would have violated \eqref{MulSid}.  Consider those integers for which $a_{j(i)}\le10\sqrt{n}$ and
  $b_{j(i)}\le10\sqrt{n}$. Not having $C_4$ in the corresponding subgraph, we can have at most $c\sqrt n^{3/2}=cn^{3/4}$ such integers. This does not cover all the cases, however, we can partition all the integers into a ``few'' similar subclasses.\footnote{E.g. we can put $a_i$ into $\cA_i$ if the
    smallest prime divisor of $a_t$ is in $(2^t,2^{t+1}]$ and use a slight generalization of Theorem \ref{MultiC4} to $K(m,n)$.}  This proves that $$|\cA|\le\pi(n)+O(n^{3/4}).$$

  (b) There is another problem solved in \cite{Erd38Tomsk}, also by reducing it to an extremal graph problem. It is much simpler, and we skip it.

(c) The \Sc ``additive~Sidon'' condition assumes that $a_i+a_j=a_k+a_\ell$ implies $\{i,j\}=\{k,\ell\}$. In case of \Sc Sum-free sets we exclude
$a_i+a_j=a_k$.  See e.g. Cameron \cite{Camer87BCC}, Bilu \cite{Bilu98} and many other papers on integers, or groups.

(d) Some related results can be found, e.g., in Chan, Győri, and A. Sárközy \cite{ChanGyoriSarko10Prod}.
\end{remark} 

\begin{remarks}
(a) Many extremal problems formulated for integers automatically extend to finite Abelian groups, or sometimes to any finite group.

Thus, e.g., a {\em Sidon sequence} can be (and is) defined in any Abelian group as a subset for which $a_i+a_j\neq a_k+a_\ell$, unless $\{i,j\}=\{k,\ell\}$. A paper of Erdős and Turán \cite{ErdTur41Sidon} estimates the maximum size Sidon subset of $[1,n]$. Several papers of Erdős investigate Sidon
problems, e.g., \cite{Erd44Sidon, Erd54Sidon}. This problem was extended to groups, first by Babai \cite{Babai85Sid}, Babai and Sós \cite{BabaiSos85Sidon}. (Surprisingly, for integers there are sharp differences in analysing the density of finite and of infinite Sidon sequences, see Subsection
\ref{InfiSidon}.) Problems on ``sum-free sets'' were generalized to groups by Babai, Sós, and then by Gowers \cite{GowersT08QuasiGr}, by Balogh, Morris, and Samotij \cite{BalMorrSamo14Abel}, and by Alon, Balogh, Morris and Samotij \cite{AlonBalMorrSamo14Cameron}.\footnote{A longer annotated
  bibliography of O'Bryant can be downloaded from the Electronic Journal of Combinatorics \cite{Bryant04Sidon} on Sidon sets.}

(b) A few more citations from this area are:
Alon \cite{Alon91Sumfree},
Green and Ruzsa \cite{GreenRuzsa05Abelian},
Lev, Łuczak, and Schoen \cite{LevLuczSchoen01Abel},
Sapozhenko \cite{Sapo02Abel,Sapozhenko02Sumfree,Sapo09CamerErd}, \dots
\end{remarks}

\minisepa 

The next extremal graph result was motivated by topology. It is

\begin{theorem}[Erdős-Stone Theorem \ev(1946) \cite{ErdStone46}]\label{ErdStoneTh} For any fixed $t\ge1$,
$$\ext(n,K_{p+1}(t,\dots,t))=\left(1-{1\over p}\right){n\choose 2}+o(n^2).$$
\end{theorem}

We mention here one more ``early extremal graph theorem'', strongly connected to the Erdős-Stone Theorem:

\begin{theorem}[Erdős, Kővári-T. Sós-Turán \ev(1954) \cite{KovSosTur54}]\label{KovSosTurTh}\footnote{Mostly we call this result the Kővári-T. Sós-Turán theorem.
    Here we added the name of Erdős, since \cite{KovSosTur54} starts with a footnote according to which ``As we learned after giving the manuscript to the Redaction, from a letter of P. Erdős, he has found most of the results of this
    paper.'' Erdős himself quoted this result as Kővári-T. Sós-Turán theorem.} \footnote{For Hungarian authors we shall mostly use the Hungarian spelling of their names, though occasionally this may differ from the way their name
    was printed in the actual publications.} \beq{KSTFo}\ext(n,K(a,b))\le \half \root a\of {b-1}\cdot n^{2-{1\over a}}+O(n).\eeq
\end{theorem}

Below let $a\le b$. A very important conjecture is the sharpness of \eqref{KSTFo}:

\begin{conjecture}[KST is sharp]\label{KSTConj}
For every pair of positive integers $a\le b$, there exists a constant $c_{a,b}>0$ for which
\beq{KSTLower}\ext(n,K(a,b))> c_{a,b} \,{n^{2-\reci a }}.\eeq
\end{conjecture}

We can replace $c_{a,b}$ by $c_{a,a}$. The conjecture follows from Erdős \cite{Erd38Tomsk} if $a=2$. A sharp construction can be found in \cite{KovSosTur54} for the bipartite case\footnote{Here ``sharp'' means that not only the exponent $2-{1\over p}$ but the value of $c_{a,b}$ is also sharp.}, (and
later an even sharper one from Reiman \cite{Reim58Zara}), and a sharp construction is given by Erdős, Rényi, and Sós \cite{ErdRenyiSos} and by Brown for $C_4=K(2,2)$ and not necessarily bipartite $G_n$. Brown also found a sharp construction \cite{Brown66Thomsen}, for $K(3,3)$\footnote{The sharpness
  of the multiplicative constant followed from a later result of Füredi.}. Much later, Conjecture \ref{KSTConj} was proved for $a\ge4$ and $b>a!$ by a Kollár-Rónyai-Szabó construction \cite{KollRonyaiSzab96}. This was slightly improved by Alon-Rónyai-Szabó \cite{AlonRonyaiSzab99}: the weaker
condition $b>(a-1)!$ is also enough for \eqref{KSTLower}.

\minisepa

Theorems of Erdős and Stone and of Kővári, T. Sós, and Turán have a very important consequence.

\begin{corollary}
 $\ext(n,\cL)=o(n^2)$ if and only if $\cL$ contains a bipartite graph.
\end{corollary}

Actually, the theorem of Kővári, T. Sós, and Turán and the fact that $T_{n,2}$ is bipartite imply the stronger dichotomy:

\begin{corollary}\label{DegenJump}
 If $\LiL$ is bipartite, then $\ext(n,\cL)=O(n^{2-(2/v(L))})=o(n^2)$ and if $\cL$ contains no bipartite graphs, then $\ext(n,\cL)\ge\turtwo n$.
\end{corollary}

The case when $\cL$ contains at least one bipartite graph will be called \Sc Degenerate. It is one of the most important and fascinating areas of Extremal Graph Theory, and Füredi and Simonovits have a long survey \cite{FureSim13Degen} on this. Here we shall deal only very shortly with \Sc Degenerate extremal graph problems.

\minisepa

The main contribution of Turán was that he had not stopped at proving Theorem~\ref{TurTh} but continued, asking the ``right'' questions: 

\begin{quote}{
 What is the maximum of $e(G_n)$ if instead of excluding $\K{p+1}$ we exclude some arbitrary other subgraph $L$?\footnote{Later this question was generalized to excluding an arbitrary family of subgraphs, however, that was only a small
 extension.}}
\end{quote} 

\begin{figure}[h]
\begin{center}
\begin{tabular}{cccc}
\epsfig{file=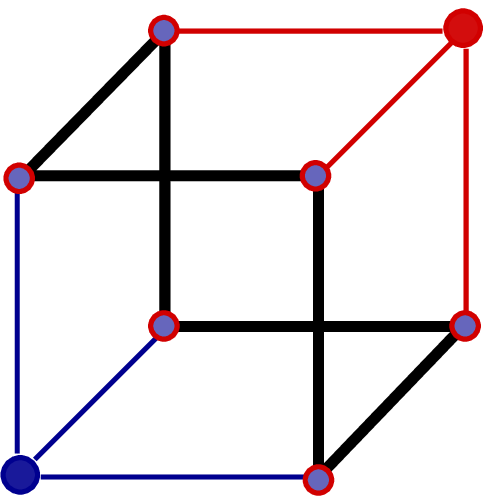,height=17mm,width=17mm}&
\epsfig{file=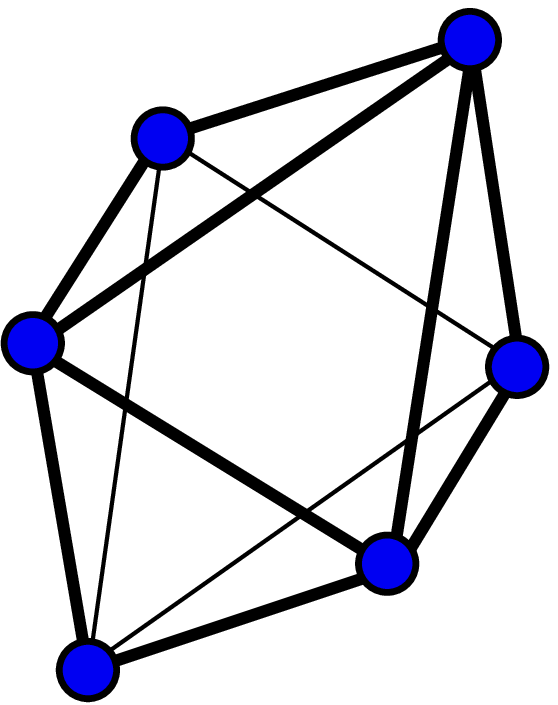,height=17mm,width=17mm}&
\epsfig{file=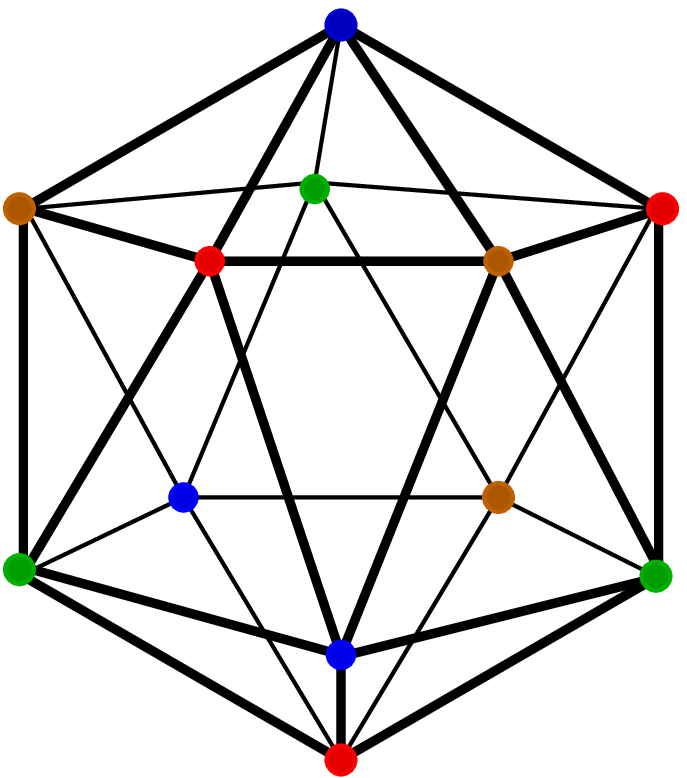,height=17mm,width=17mm}&
\epsfig{file=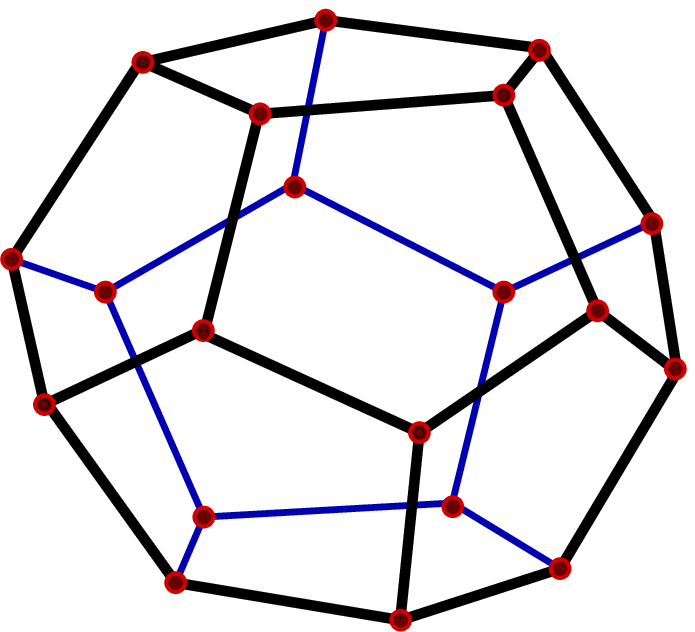,height=17mm,width=17mm}\\
Cube &Octahedron& Icosahedron, & Dodecahedron \\
\end{tabular}
\caption{Excluded platonic graphs}
\end{center}
\label{PlatonicFigs}
\end{figure}

To provide a starting point, Turán asked for determining the extremal numbers and graphs for the Platonic bodies: Cube, Octahedron, Dodecahedron, and Icosahedron,\footnote{The tetrahedron is $K_4$, covered by Turán theorem.} for paths, and for lassos \cite{Erd65Smole},\footnote{Lasso is a graph
  where we attach a path to a cycle. Perhaps nobody considered the lasso-problem carefully, however, very recently Sidorenko solved a very similar problem of the keyrings \cite{Sidore18KeyRings}.} see Figure~\ref{PlatonicFigs}. Here the extremal graph problem of the Cube was (partly) solved by Erdős
and Simonovits \cite{ErdSim69Cube}. We return to this problem in Subsection \ref{DegenS}. The problems of the Dodecahedron and Icosahedron were solved by Simonovits, \cite{Sim74Symm},\cite{Sim74Ico}, using the Stability method, see Subsection \ref{Stability}.  
\bigskip

\noindent
The general question can be formulated as follows:

\begin{quote}{
Given a family $\cL$ of forbidden graphs, what is the {maximum of $e(G_n)$} if $G_n$ does not contain subgraphs $\LiL$?}
\end{quote}

\subsection{Constructions}

Mostly in an extremal graph problem first we try to find out how do the extremal structures look like. In the nice cases this is equivalent to finding a construction providing the lower bound in our extremal problems, and then we try to find the matching upper bound.  Here we shall not go into
more details, however, we shall return to this question in Section \ref{HyperLargeExcludedS}, on Matchings, 1-factors, and the Hamiltonicity of Hypergraphs.

\gdef\mytabline#1//#2//#3//{&&\\
#1&\begin{minipage}{7cm}#2\end{minipage}&\begin{minipage}{14mm}#3\end{minipage}\\
&&\\
 \hline}

\medskip

\subsection{Some historical remarks}\label{Historical}

Here we make two remarks.

(A) When Turán died in 1976, several papers appeared in his memory, analysing, among others, his influence on Mathematics. Erdős himself wrote several such papers, e.g., \cite{Erd80OnTuranActa,Erd80OnTuranMath}. In \cite{Erd78OnTuranJGT} he wrote, on Turán's influence on Graph Theory:

\begin{quote}{
`` In this short note, I will restrict myself to Turán's work in Graph Theory, even though his main work was in analytic number theory and various other branches of real and complex analysis. Turán had the remarkable ability to write
 perhaps only one paper or to state one problem in various fields distant from his own; later others would pursue his idea and a new subject would be born. In this way Turán initiated the field of Extremal Graph Theory. \dots \\ \dots
 Turán also formulated several other extremal problems on graphs, some of which were solved by Gallai and myself \cite{ErdGallai59Path}. I began a systematic study of extremal problems in Graph Theory in 1958 on the boat from Athens to
 Haifa and have worked on it since then. The subject has grown enormously and has a very large literature;\dots ''} 
\end{quote}

Observe that Erdős implicitly stated here that until the early 60's most of the results in this area were sporadic.

(B) Here we write about Extremal Graph Theory at length, still, if one wants to tell what Extremal Graph Theory is, and what it is not, that is rather difficult. We shall avoid answering this question, however, we remark that since the
Goodman paper \cite{Goodman59} and the Moon-Moser paper \cite{MoonMoser62Turan} an alternative answer was the following. Consider some ``excluded subgraphs'' $L_1,\dots,L_t$, count the multiplicities of their copies, $m(L_i,G_n)$, in $G_n$,
and Extremal Graph Theory consists of results asserting some inequalities among them. Since the emergence of Graph Limits this approach became stronger and stronger. One early ``counting'' example is

\begin{theorem}[Moon and Moser \ev(1962) \cite{MoonMoser62Turan}] 
If $t_k$ is the number of $K_k$ in $G_n$, then 
$$k(k-2)t_k\geq t_{k-1}\left(\frac{(k-1)^2t_{k-1}}{t_{k-2}}-n\right). $$
\end{theorem}

\subsection{Early results}\label{EarlyS} 

If we restrict ourselves to simple graphs, {\em some central theorems} assert that for ordinary graphs the general situation is {almost the same as} for ${K_{p+1}}$: the extremal graphs $S_n$ and the {\em almost extremal} graphs $G_n$ are {\em very similar} to $T_{n,p}$. The similarity of two graph
sequences $(G_n)$ and $(H_n)$ means that one can delete {$o(n^2)$} edges of $G_n$ and add $o(n^2)$ edges to obtain $H_n$.

The general asymptotics of $\ext(n,\cL)$ and the asymptotic structure of the extremal graphs are described by

\begin{theorem} [Erdős-Simonovits, \ev(1966) \cite{Erd67Rome},\cite{Erd66Tihany},\cite{Sim66Tihany}]\label{GenExtre} 
 Let
 \beq{subchrom}p:=\min_{\LiL}\chi(L)-1.\eeq
 If $S_n\in\EXT(n,\cL)$, then one can delete from and add to $S_n$ $o(n^2)$ edges to obtain $T_{n,p}$.
\end{theorem}

A much weaker form of this result immediately follows from Erdős-Stone Theorem.

\begin{theorem}[Erdős-Simonovits \ev(1966) \cite{ErdSim66Lim}]\label{ErdSimLimTh}

Defining $p$ by \eqref{subchrom},
$$\ext(n,\cL)=\left(1-{1\over p}\right){n\choose 2}+o(n^2).
$$\end{theorem}

One important message of these theorems is that for simple graphs the extremal number and the extremal structure are determined up to $o(n^2)$ by the minimum chromatic number of the excluded subgraphs. In some sense this is a great luck:
there are many generalizations of the Turán type extremal graph problems, but we almost never get so nice answers for the natural questions in other areas.

Speaking about the origins of Extremal Graph Theory, we have to mention the {\em dichotomy} that occasionally we have very nice extremal structures but in the \Sc Degenerate case the extremal graphs seem to have much more complicated
structures (unless $\cL$ contains a tree or a forest). Actually this may explain why Erdős missed to observe the importance of his Theorem~\ref{MultiC4}, on $\ext(n,C_4)$.

To solve some extremal graph problems, Simonovits defined the {\em Decomposition Family} (see, e.g., \cite{Sim83Fra}).

\begin{definition}[Decomposition, $\cM=\cM(\cL)$]
Given a family $\cL$ with $p:=\min\{\chi(L):\LiL\}-1$, then
$M\in\cM$ if $L\subseteq M\otimes K_{p-1}(t,\dots,t)$ for $t=v(M)$, where $L\otimes H$ denotes the graph obtained by joining each vertex of $L$ to each vertex of $H$.
\end{definition}

The meaning of this is that $M\in\cM$ if we cannot embed $M$ into one class of a $T_{n,p}$ without obtaining an excluded $\LiL$. Thus, e.g, if we put a $C_4$ into the first class of $T_{n,p}$, then the resulting graph contains a
$K_{p+1}(2,2,\dots,2)$. Therefore $C_4$ is in the decomposition class of $K_{p+1}(2,\dots,2)$.

\dots Given a family $\cL$ of excluded subgraphs, the decomposition family $\cM=\cM(\cL)$ determines (in some sense) the error terms and the finer structure of the $\cL$-extremal graphs. Namely, the error terms depend on $\ext(n,\cM)$, see \cite{Sim66Tihany}.

The next theorem ``explains'', why is $T_{n,p}$ extremal for $K_{p+1}$.

\begin{definition}[Colour-critical edge]
 The edge $e\in E(L)$ is called critical, if $\chi(L-e)<\chi(L)$.
\end{definition}

\begin{figure}[h]
\begin{center} 
\begin{tabular}{ccc}
\epsfig{file=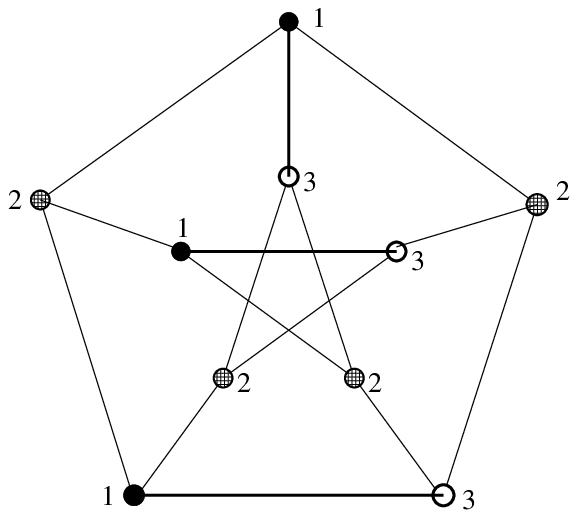,height=21mm}&
\epsfig{file=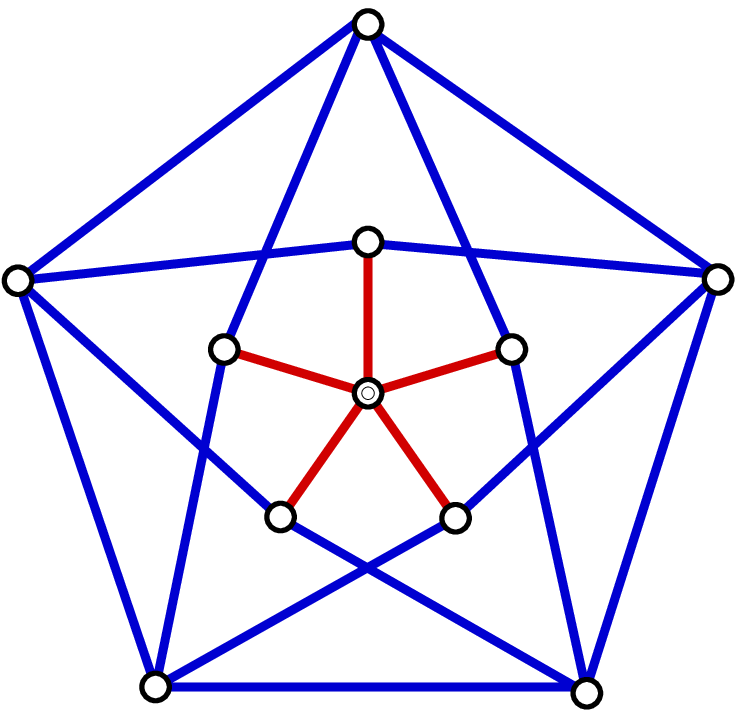,height=21mm}& 
\epsfig{file=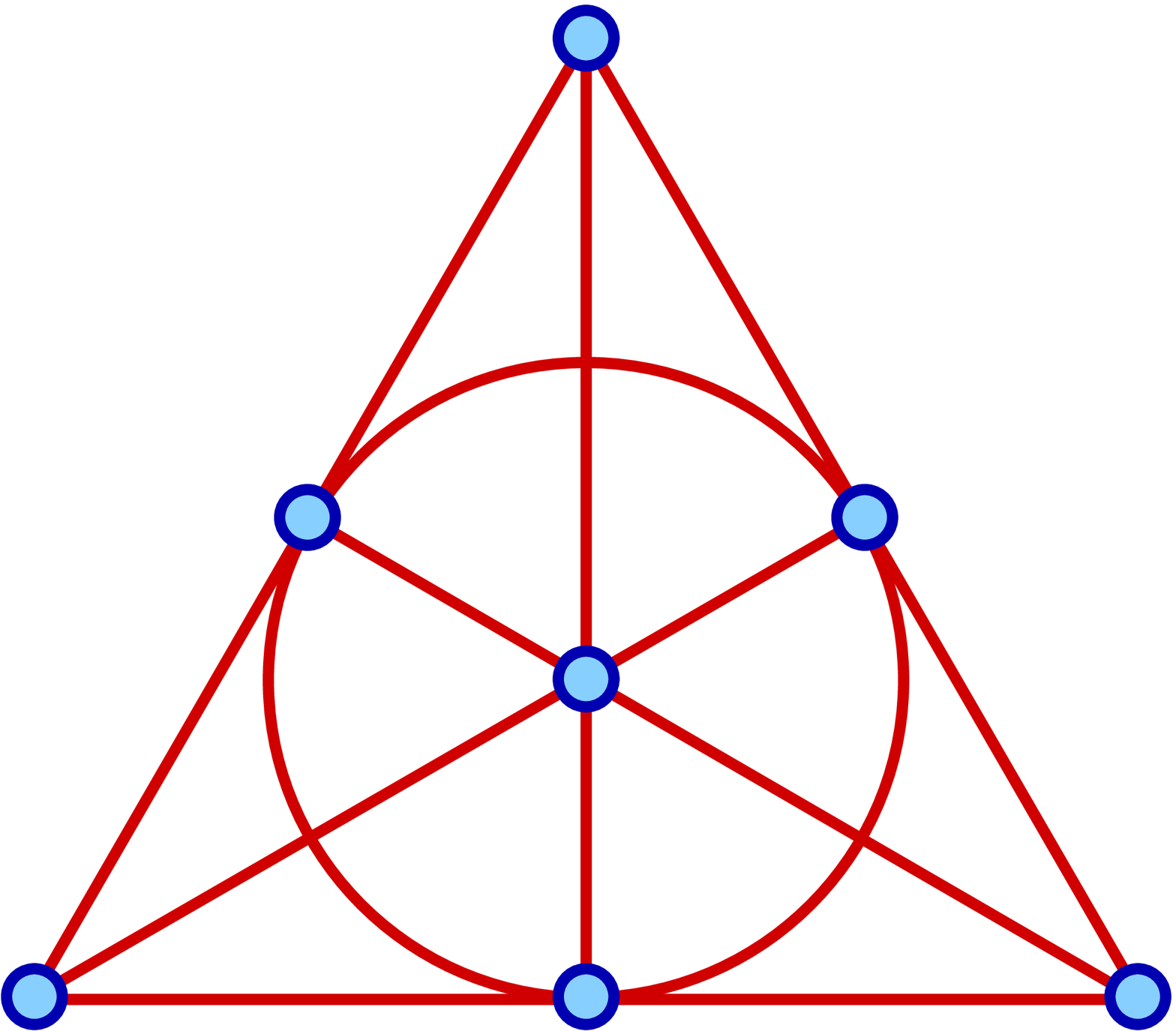,height=20mm} \\
Petersen graph,&Grötzsch graph,& Fano hypergraph\\
\end{tabular}
%% \caption{}
\end{center}
\label{Grotzsch-Peter}
\end{figure}

Of course, in such cases, $\chi(L-e)=\chi(L)-1$. Each edge of an odd cycle is critical. In Figure \ref{Grotzsch-Peter}, e.g., one can see the Grötzsch graph (often incorrectly called Mycielski graph). It is 4-chromatic but all its edges are critical. On the other hand, in the Petersen, or the
Dodecahedron graphs, there are no critical edges.\footnote{If the automorphism group of $G$ is edge-transitive, then either all the edges are critical, or none of them. By the way, in \cite{Sim99Stirin}, Simonovits discusses these questions in more details, among others, the extremal problems of
  generalized Petersen graphs.}  \footnote{The definition applies to hypergraphs as well, the triples of the Fano hypergraph are also critical.}  The next theorem solves all cases when $L$ has a critical edge.

\medskip

\begin{theorem}[Simonovits]\footnote{For $p=2$ this was also known (at least, implicitly) by Erdős.}
  Define $p=p(\cL)$ by \eqref{subchrom}.  The following statements are equivalent:

(a) Some $(p+1)$-chromatic $\LiL$ has a critical edge.

(b) There exists an $n_0$ such that for $n>n_0$, $T_{n,p}$ is extremal for $\cL$.

(c) There exists an $n_1$ such that for $n>n_1$, $T_{n,p}$ is the {\em only} extremal graph for $\cL$.\footnote{Here 
(c)$\to$(b) 
is trivial, and one can prove that (b) implies (c) with $n_1=n_0+3p$. 
}
\end{theorem}

Here (a) is equivalent to that adding an edge to $K_p(t,\dots,t)$, for $t=v(L)$, we get a graph containing $L$, and equivalent to that $K_2\in\cM$.  The extremal results on the Dodecahedron $D_{20}$ and Icosahedron $I_{12}$ follow from
\cite{Sim74Symm}, and \cite{Sim74Ico}. We skip the details referring the reader to the survey of Simonovits \cite{Sim99Stirin}.

\begin{MetaTheorem}[Simonovits]\label{MetaCritical}
``Whatever'' we can prove for $L=K_{p+1}$, with high probability, we can also prove it for any $L$ having a critical edge.
\end{MetaTheorem}

\subsection{Which Universe?}\label{UniverseS}

Extremal problems exist in a much more general setting: Theorem \ref{MultipSidonTh} is, e.g., an extremal theorem on sets of integers. In general, we fix the family of some objects, e.g., integers, graphs, hypergraphs, $r$-multigraphs, where some $r$ is fixed and the edge-multiplicity is bounded by
$r$. We exclude some substructures, and try to optimize some (natural) parameters. More generally, putting some bounds on the number of one type of substructures, we try to maximize (or minimize) the number of some other substructures. This approach can be found in the paper of Moon and Moser
\cite{MoonMoser62Turan}, or in Lovász and Simonovits, \cite{LovSim83Birk}.\footnote{Bollobás \cite{Bollo76Aberd} also contains similar, strongly related results.}  The paper of Alon and Shikhelman \cite{AlonShikhel16} is also about this question, in a more general setting. (We should also mention
here the famous conjecture of Erdős, on the number of pentagons in a triangle-free graph, well approximated by Győri \cite{Gyori89C5}, (and by Füredi, unpublished) and then solved by Grzesik \cite{Grzesik12} and by Hamed Hatami, Hladký, Král, Norine and Razborov \cite{HatamiHladKralNorineRazbor12}.)
With the development of the theory of graph limits this view point became more and more important. Below we list some most common Universes, and some related papers/surveys.

\begin{enumerate}\dense
\item \Sc Integers, as we have seen above, in Problem \ref{MultipSidon} or Theorem \ref{MultipSidonTh}. Among many other references, here we should mention the book of Tao and Vu on Additive Combinatorics \cite{TaoVu06Addit}, the
  Geroldinger-Ruzsa book \cite{GeroldRuzsa09Doc}, and also a new book of Bajnok \cite{Bajnok19Additive}.
\item \Sc Abelian~Groups, see e.g. Babai-Sós, \cite{BabaiSos85Sidon} Gowers \cite{GowersT06-3Unif,GowersT08QuasiGr} and also \cite{GeroldRuzsa09Doc}, a survey of Tao and Vu \cite{TaoVu17SumFree} and \cite{Bajnok19Additive}.\footnote{There
    are several earlier results on similar questions, e.g., Yap, \cite{Yap69SumFree,Yap75SumFree}, Diananda and Yap \cite{DianYap69}, yet they are slightly different, or several papers of A. Street, see \cite{Street72SumFree}.}
\item \Sc Graphs: this is the main topic of this survey;
\item \Sc Digraphs and \Sc multigraphs, with bounded arc/edge multiplicity,\footnote{One has to assume that the edge-multiplicity is bounded, otherwise even for the excluded $K_3$ in the Universe of multigraphs we would get arbitrary many edges. As an exception, in the Füredi-Kündgen theorem
    \cite{FureKund02Weight} no such bound is assumed.}  see Brown-Harary \cite{BrownHarary70}, Brown-Erdős-Simonovits \cite{BrownErdSim72JCTB,BrownErdSim85Algo}; Sidorenko, \cite{Sidor93Multi}, for longer surveys see Brown and Simonovits \cite{BrownSim84DM} and \cite{BrownSim99Erd}, Bermond and
  Thomassen \cite{BermondThomassen81Digra}, Thomassen \cite{Thomassen79BCC}, Bang-Jensen and Gutin \cite{BangGutin09Book}, Jackson and Ordaz \cite{JacksonOrdaz90Survey};

\item \Sc Hypergraphs, see e.g. de Caen \cite{Caen91Surv}, Füredi \cite{Fure91Sur}, Sidorenko \cite{Sidor95Know}, Keevash \cite{Keev11HypergBCC};\footnote{and many others, see e.g. Rödl and Rucinski \cite{RodlRuc10HyperSurv}, or the much earlier Bermond, Germa, Heydemann, and Sotteau
    \cite{BermondGermaHeydSotte76Orsay}, and the corresponding Sections \ref{HyperSmallExcludedS} and \ref{HyperLargeExcludedS}.}

\item \label{ExtreSub} {\sc Extremal subgraphs of random graphs}: Babai, Simonovits, Spencer \cite{BabaiSimSpenc90}, Brightwell, Panagiotou, and Steger \cite{BrightPanaSteger12}, Rödl-Schacht \cite{RodlSchacht13RandExt}, Rödl \cite{Rodl14Seoul}, Schacht \cite{Schacht16ExtremalRandom}, \dots DeMarco
  and J. Kahn \cite{DeMarco12Thes}, \cite{DeMarcoKahn15Man}, \cite{DeMarcoKahn15ArxTur}, or

\item {\sc Extremal subgraphs of Pseudo-random graphs}, see e.g., Krivelevich and Sudakov \cite{KriveSudak06Pseudo} Aigner-Horev, Hàn, and Schacht \cite{AigHanSchacht14Odd}, Conlon and Gowers \cite{ConGow16Sparse,Con14Congress}, or Conlon, Fox and Zhao \cite{ConFoxZhao14SparsePseudo}, or e.g., Allen, Böttcher, Kohayakawa, Person, \cite{AllenBottHanKohaPers14,AllenBottHanKohaPers17Powers}.

\end{enumerate}

Perhaps one of those who first tried to compare various universes and analyze their connections was Vera T. Sós \cite{Sos76Lincei,Sos91Additive}. In \cite{Sos76Lincei} she considered the connections between Graph Theory, Finite Geometries, and Block Designs. The emphasis in these papers was on the
fact that basically the same problems occur in these areas in various settings, and these areas are in very strong connection, interaction, with each other.\footnote{Vera Sós did not call these areas Universes.}  Most of the above universes we shall skip here, to keep this survey relatively short, however,
below we consider some extremal problems on integers.

\subsubsection{Integers}\label{IntegerS} 

We do not intend to describe this very wide area in details, yet we start with some typical, important questions in this area, i.e., in the theory of extremal problems on subsets of integers. As we have stated in Problem \ref{MultipSidon}, Erdős considered the multiplicative Sidon problem in
\cite{Erd38Tomsk}. Even earlier Erdős and Turán formulated the following conjecture for subsets of integers:

\begin{conjecture}[Erdős-Turán \ev(1936) \cite{ErdTur36r3n}]\label{ErdTuranRkn}\footnote{Actually, here they formulated this for $r_3(n)$.}
  If $\cA\subseteq[1,n]$ does not contain a $k$-term arithmetic progression, then $|\cA|=o(n)$ (as $n\to\infty$).\footnote{Speaking of arithmetic progressions we always assume that its terms are distinct.}
\end{conjecture}

The proof seemed those days very difficult. Even the simplest case $k=3$ is highly non-trivial: it was first proved by K.F. Roth, in 1953 \cite{Roth53London}, and for $k=4$ those days the conjecture seemed even more difficult.\footnote{The conjectures on $r_k(n)$ were not always correct. Vera Sós
  wrote a paper \cite{SosV02Turbulent} on the letters between Erdős and Turán during the war, where one can read that Szekeres e.g. conjectured that for $n=\half(3^\ell+1)$ $r_k(n)\le 2^\ell$. This was later disproved by Behrend \cite{Behrend46AP3}. (This conjecture is also mentioned in
  \cite{ErdTur36r3n}.)}  The conjecture was proved by Szemerédi, first for $k=4$ \cite{Szemer69AP4}, and then for any $k$:

\begin{theorem}[Szemerédi \ev(1975) \cite{Szem75APk}]\label{Szem75APkT}
  Let $k$ be a fixed integer. If $r_k(n)$ is the maximum number of integers, $a_1,\dots,a_m\in[1,n]$ not containing a $k$-term arithmetic progression, then $r_k(n)=o(n)$.
\end{theorem}

\subsub Ergodic theory and Szemerédi Theorem.// Not much later that Szemerédi proved this theorem, Fürstenberg gave an alternative proof, in 1977, using ergodic theory \cite{Fursten77Szemer}.  This again is an example where seemingly simple combinatorial problems led to very deep theories. One advantage of Fürstenberg's approach was that it made possible for him and Katznelson and their school to prove
several important generalizations, e.g., the high dimensional version \cite{FurstenKatz78Szemer}, Bergelson and Leibman \cite{BergelLeib96Poly} proved some polynomial versions of the original theorem, and later the density version of
Hales-Jewett theorem \cite{FurstenKatz91DensityHalesJewett}.

\subsub Polymath on Hales-Jewett theorem.// As we just stated, one of the important generalizations of Theorem \ref{Szem75APkT} is the {\em density version of} Hales-Jewett theorem, obtained by Fürstenberg and Katznelson
\cite{FurstenKatz91DensityHalesJewett} which earlier seemed hopeless.

\begin{quote}{There is a big difference between the Hales-Jewett Theorem and Szemerédi's theorem: just to explain the meaning of the Hales-Jewett Theorem or its density version, is more difficult than to explain the earlier ones. It is one
    of the most important Ramsey type theorems, asserting, that -- fixing the parameters appropriately -- a high dimensional $r$-coloured structure will contain a (small) monochromatic substructure, a so called ``combinatorial line''.  We
    remark that there was a similar result of Graham and Rothschild (on $n$-parameter sets) earlier, \cite{GrahRoth71Parameter}.}\end{quote}

A simplified version of this was given by Austin \cite{Austin11DensityHalesJewett}.\footnote{See also Austin \cite{Austin10Remo}, Beigleböck \cite{Beigl08HL}, (?) Bergelson and Leibman \cite{BergelLeib99HJ} Gowers,
  \cite{Gowers10PolyMathHJ}, Polymath \cite{Polymath},\dots } The Polymath project \cite{Polymath12} provided a completely elementary proof of this theorem.  A nice description of this is the MathSciNet description of the ``PolyMath'' proof
of this (see MathSciNet MR2912706, or the original paper, \cite{Polymath12}).\footnote{Similarly to the proof of $r_3(n)=o(n)$ from the Ruzsa-Szemerédi Triangle Removal Lemma, (see Theorem \ref{RuzsaSzemTh}) Rödl, Schacht, Tengen and
  Tokushige proved $r_k(n)=o(n)$ and several of its generalizations in \cite{RodlSchachtTengToku06} ``elementarily'', i.e. not using ergodic theoretical tools. On the other hand, they remarked that those days no elementary proof was known
  on the Density Hales-Jewett theorem.}

%% \\Irni kellene itt néhány mondatot Gowers idevonatkozó eredmenyeirol!}

\subsub Erdős conjecture on the sum of reciprocals.// One of the central questions in this area was if there are arbitrary long arithmetic progressions consisting of primes. In 1993, Erdős wrote a paper on his favourite theorems
\cite{Erd93Favo}, where he wrote that in those days the longest arithmetic progression of primes had 17 integers (and was obtained with the help of computers). Now we all know the celebrated result:

\begin{theorem}[Green and Tao \ev (2008) \cite{GreenTao08P-AP}]\label{GreenTaoTh}
 The set of primes contains arbitrary long arithmetic progressions.
\end{theorem}

We close this part with the related famous open problem of Erdős which would imply Theorem \ref{GreenTaoTh}. Then we list some estimates on $r_k(n)$.

\begin{problem}[Erdős \cite{Erd93Favo}] 
  Let $\cA=\{a_1,\dots,a_n,\dots\}\subseteq\ZZ$ be a set of positive integers.  Is it true that if $\sum \reci a_i =\infty$, then for each integer $k>2$, $\cA$ contains a $k$-term arithmetic progression?
\end{problem}

\subsub Estimates on $r_k(n)$.// First Roth proved that $r_3(n)=O({n\over\log\log n})$.\footnote{Using $\loglog n$ we always assume that $n\ge100$, and therefore $\loglog n>3/2$.} Many researchers worked on improving the estimates on
$r_3(n)$, or more generally, on $r_k(n)$. Roth's estimate was followed by the works of Heath-Brown \cite{HeathBrown87r3n} and Szemerédi \cite{Szem90r3n}, and then Bourgain \cite{Bourg99r3n}. One of the last breakthroughs
was

\begin{theorem}[Sanders, \cite{Sanders11r3n}] Suppose that $\cA\subseteq\{1,\dots,N\}$ contains no 3-term arithmetic progressions. Then
$$ |\cA|=O\left(\frac{(\log\log N)^5}{\log N}N\right). $$
\end{theorem}

The exponent of $\loglog n$ was brought down to 4 by T. Bloom \cite{Bloom16AP3}. 

\begin{remark}[Lower bounds] 
  Clearly, for $k\ge3$, $r_k(n)\ge r_3(n)$. Behrend \cite{Behrend46AP3} proved that there exists a $c>0$ for which \beq{BehrendLower} r_3(n)\ge {n\over e^{c\sqrt{\log n}}}.\eeq This was improved by Elkin \cite{Elkin11AP}, and then, in a
  much more compact way, by Green and Wolf \cite{GreenWolf10Behrend}.
\end{remark}

\minisepa

One could ask, what do we know about $r_4(n)$. A major breakthrough was due to Gowers \cite{Gowers98AP4}, according to which for every $k\ge3$ there exists a $c_k>0$ for which
$$r_k(n)=O\left({n\over{\loglog^{c_k} n}}\right).\qquad \text{(Actually,} c_k=2^{-(2^{k+9})} \text{~~works.)}$$
Green and Tao improved this for $k=4$ to 
$$r_4(n)=O\left({n\over{\log^c n}}\right) \Text{for some constant} c>0.$$
See also Green and Tao \cite{GreenTao09NewBounds}, \cite{GreenTaoT17LowR4n} and the survey of Sanders \cite{Sanders14ICM}.

\subsub Other important problems on Integers:// We start with the following remark. A property $\cP$ is always a family of subsets of some fixed set. It is {\em monotone decreasing} if $X\in\cP$ and $Y\subseteq X$ implies that
$Y\in\cP$. (Two examples of this are (i) the sets of integers not containing solutions of some given equations, and (ii) the family of graphs not containing an $L$.)  When we fix a Universe and a ``monotone'' property $\cP$, then, beside
asking for the size of the extremal sets $X$ for $\cP$, we may also ask, e.g., how many $X\in\cP_n$ are there, where $\cP_n\subseteq\cP$ is defined by some parameter $n$ of these objects. We may also ask, what is their typical
structure. For graphs these are the Erdős-Kleitman-Rothschild type problems \cite{ErdKleiRoth76Lincei}, discussed in Section \ref{ErdKleitRothS}. The same questions can also be asked for extremal problems on integers, and we shall not
return to them later, therefore we list some of them here, together with ``their extremal problems''.

\dorib We have started with the multiplicative Sidon problem (Thm \ref{MultipSidonTh}), and there is also the problem of additive Sidon sets.  \dorib We wrote about excluding the $k$-term arithmetic progressions, in this
subsection.\footnote{Subsubsections will also be called Subsections.}  \dorib Another important area is the problem of sum-free sets, see, e.g., Cameron and Erdős \cite{CamerErd99Matra}, Alon, Balogh, Morris, and Samotij
\cite{AlonBalMorrSamo14Count}, Łuczak and Schoen \cite{LuczSchoen00SumFree,LuczSchoen01SumFree}, Sapozhenko \cite{Sapo02Abel},\dots Balogh, Liu, Sharifzadeh, and Treglown \cite{BalLiuSharTreg15}, Balogh, Morris, and Samotij
\cite{BalMorrSamo14Abel}. \dorib Very active research characterises the Sum-Product problems, introduced by Erdős and Szemerédi \cite{ErdSzem83SumProd}, see also Gy. Elekes \cite{Elekes97SumProd}, Bourgain, Gambourd, and Sarnak
\cite{BourgGambSarnak10SumProd}, Solymosi \cite{Solymo05SumProd} among the very many related papers.  \dorib Important and deep questions can be listed in connection with Freiman-Ruzsa type results (see e.g. \cite{Freiman73Struct}).

\subsub Recommended surveys, papers:// 
Ruzsa \cite{Ruzsa93Linear,Ruzsa95Linear}, Solymosi \cite{Solymo05SumProd}, Granville and Solymosi \cite{GranvSolymo16Surv}, Pomerance and A. Sárközy \cite{PomerSarko95HB}, \dots and the survey of Shkredov \cite{Shkredov06rkn}.

\subsubsection{Groups}

We have considered some problems about the Universe of Integers. It is, of course, natural to ask the analogous questions for groups. There are
many results of this type. Here we mention only a few papers on groups. Most of the definitions immediately generalize from integers to any Abelian
group. We have already mentioned that Sidon sets in groups were investigated by Babai and Sós \cite{BabaiSos85Sidon}. From among the many-many further
similar extensions we mention only Alon, Balogh, Morris and Samotij \cite{AlonBalMorrSamo14Count}, \cite{BalMorrSamo14Abel}, Gowers \cite{GowersT08QuasiGr},
on quasi-random groups, Green \cite{Green05Abelian}, Green and Ruzsa \cite{GreenRuzsa05Abelian}, Lev, Łuczak and Schoen \cite{LevLuczSchoen01Abel},
Sapozhenko \cite{Sapo02Abel,Sapo09CamerErd}, and B. Szegedy on Gowers norms and groups \cite{SzegeGowers}.

There are also results on non-Abelian groups, e.g., Sanders \cite{Sanders10NonAbel}, and the paper of Babai and Sós \cite{BabaiSos85Sidon} considers both Abelian and non-Abelian groups (and try to determine the maximum size of a Sidon set in them), and we refer the reader to these papers.

\begin{remark} 
For each property $\cP$, one can also investigate the $\cP$-maximal, or the $\cP$-minimal structures. Here, e.g., one can try to count the maximal subsets of property $\cP$, in $[1,n]$, or in a group $\cG$, \dots Thus e.g.,
Balogh, Liu, Sharifzadeh and Treglown \cite{BalLiuSharTreg15} count the maximal sum-free subsets, while Balogh, Bushaw, Collares, Liu, Morris, Sharifzadeh \cite{BalBushCollMorrShar} describe the typical structure of graphs with no large cliques\footnote{The description of the typical structure is a stronger result than just counting them.}.
\end{remark}

\begin{remark}
  There are several results on subsets of $\cG$ without non-trivial arithmetic progressions where $\cG$ is a group, or a linear vector space. Here we mention only the paper of Croot, Lev, and Péter P. Pach \cite{CrootLevPachP17AP} on the
  linear vector space $\Bbb Z^n_4$. J. Wolf provides a very clear and detailed description of this paper and related results, in the MathSciNet-MR3583357. See also the related post of T. C. Tao, and the paper of Ellenberg and Gijswijt
  \cite{EllenbergGijswijt17}, building on \cite{CrootLevPachP17AP}.
\end{remark}

\subsection{Ramsey or density?}

One difference between Ramsey and Turán theory is that in the Turán case we have density statements, while in the Ramsey case the densities are not enough
to ensure the occurrence of a monochromatic structure. A trivial example of this is the problem of $R(4,4)$, yet, instead of this we consider another trivial
example: the connectedness. If we RED-BLUE-edge-colour a $K_n$, then either we have a RED connected spanning subgraph, or a BLUE one. However, we may have
$\approx\half{n\choose2}$ edges in both colours, not enough for a connected spanning subgraph.

More generally, if we $r$-colour the edges of a graph $G$, and consider its subgraphs $G^{i}$ defined by the $i\th$ colour, and we assert that -- under some conditions, -- $G$ has a monochromatic $L$, because it has at least $\reci r e(G)$
edges, that is a Turán type, density theorem.

The Erdős-Turán conjecture, and its proof, the Szemerédi theorem, came from the van der Waerden theorem, \cite{Waerden27} according to which, if $r$ and $k$ are fixed, and we $r$-colour the integers in an arbitrary way, then there will
be a monochromatic arithmetic progression of length $k$. The Erdős-Turán conjecture was the corresponding density conjecture: any infinite sequence of integers of positive lower density contains an arithmetic progression of $k$ terms.

Many Ramsey problems are very different from density problems, however, in some other cases a Ramsey problem may be basically a density problem. 
Sometimes a density theorem generalizes a Ramsey type results in a very non-trivial way. In this survey mostly we are interested in density problems. 

\begin{examp} 
For any tree $\Tk$, trivially, $\ext(n,\Tk)<(k-2)n$. This implies that
$R(\Tk,\Tk)<4k$, and for $r$ colours $R_r(\Tk,\dots,\Tk)\le 2kr$.
\end{examp}

So, for trees the Ramsey problem is a density problem, up to a constant. For more details, see e.g. the paper of Faudree and Simonovits \cite{FaudSim92Ramsey}.

\subsection{Why are the extremal problems interesting?}

Extremal graph problems are interesting on their own, they emerge in several branches of Discrete Mathematics, e.g., in some parts of Graph Theory not directly connected to Extremal Graph Theory, in Combinatorial Number Theory, and also
they are strongly connected to Ramsey Theory.

Erdős wrote several papers on how can Graph Theory be applied in Combinatorial Number Theory, or in Geometry, see e.g. \cite{Erd1968-09}. András Sárközy
returned to the investigation and generalization of Erdős' results discussed in Problem~\ref{MultipSidon}: the next step was to analyze the case when 
no product of six distinct numbers from $\cA$ was a square\footnote{If \eqref{MulSid} is violated then $a_ia_ja_ka_\ell$ is a square.}.  The corresponding graph
theoretical lemmas were connected to $\ext(n,m,C_6)$\footnote{$\ext(n,m,\cL)$ is the maximum number of edges an $\cL$-free graph $G\subset K(n,m)$ can
  have. This problem may produce surprising phenomena when $n=o(m)$. } and were established by Erdős, A. Sárközy, and Sós, \cite{ErdSarkoSos95Prod} by
G.N. Sárközy \cite{SarkoG95Cycles} and by E. Győri \cite{Gyori97C6}. (Similar cycle-extremal results and similar methods were used also in a paper of
Dietmann, Elsholz, Gyarmati and Simonovits \cite{DietElshGyarmSim05}, but for somewhat different problems.)

Extremal graph theory is strongly connected to many other parts of Mathematics, among others, to Number Theory, Geometry, the theory of Finite Geometries, Random Graphs, Quasi-Randomness, Linear Algebra, Coding Theory.

The application of constructions based on finite geometry became important and interesting research problems, we mention here just a few, like Reiman \cite{Reim58Zara}, Hoffman-Singleton \cite{HoffmanSingle60}, Benson \cite{Benson66},
W.G. Brown \cite{Brown66Thomsen}. Erdős, Rényi, T. Sós \cite{ErdRenyiSos} \dots and refer the reader again to the surveys of Vera Sós \cite{Sos76Lincei}, Füredi and Simonovits \cite{FureSim13Degen} or the papers of Lazebnik, Ustimenko,
Woldar, e.g., \cite{LazebUstimWold94} and others.

These constructions are connected to the construction of expander graphs (Ramanujan graphs) by Margulis \cite{Margul84Arith,Margul87Sochi,Margul88Explic},
and Lubotzky, Phillips and Sarnak \cite{LuboPhilSarnak88CCA}, which use highly non-trivial mathematics, and in some sense are strongly connected to Extremal
Graph Theory.\footnote{The Margulis-Lubotzky-Phillips-Sarnak papers are eigenvalue-extremal, however, as Alon pointed out, (see the last pages of
 \cite{LuboPhilSarnak88CCA}), these constructions are ``extremal'' for many other graph problems as well. } Here we recommend the survey of Alon in the
Handbook of Combinatorics \cite{Alon95HB} and several of his results on the eigenvalues of graphs, e.g. \cite{Alon86Eigen} or the Alon-Milman paper
\cite{AlonMilman}. Another surprise was that in \cite{KollRonyaiSzab96} deeper results from algebra also turned out to be very useful. In some other cases
(e.g., Bukh and Conlon \cite{Bukh15Random,BukhConlon18Rati}) randomly chosen polynomial equations were used for constructions in extremal graph theory. We
return to this question in Section \ref{DegenS}.

\begin{remark}  
  There are cases, when important methods came from that part of Discrete Mathematics which is not directly Extremal Graph Theory, however, is very strongly connected to it. One example is (perhaps) the Lovász Local Lemma
  \cite{ErdLov75Local}, originally invented for problems very strongly connected to Extremal Graph Theory.\footnote{The Lovász Local Lemma is one of the most important tools in Probabilistic Combinatorics (including the application of
    probabilistic methods). Its proof is very short, and it is described, among others, in the Alon-Spencer book \cite{AlonSpencer16}, in Spencer \cite{Spencer87TenLectures}, or in the original paper, available at the ``Erdős homepage''
    \cite{ErdHomepage}.}
\end{remark}

\subsection{Ramsey Theory and the birth of the Random Graph Method}

There are many cases when some Ramsey type theorem is very near to a density theorem. In graph theory perhaps one of the first such results was that of Chvátal \cite{Chvat77Tree}. Faudree, Schelp and others also proved many results
on the Ramsey topic, where in one colour we exclude a large tree. Faudree and Simonovits discussed in \cite{FaudSim92Ramsey} this connection.

\subsub Erdős Magic.// Turán thought that two-colouring $K_n$, say, in RED and BLUE, we shall always have a monochromatic $K_m$ with $m>\sqrt{n}$. The reason he thought this was (most probably) that for $n=m^2$, $\Tn n,m$ yields a
2-colouring of $K_n$ (where the edges of $\Tn n,m$ are RED, the others are BLUE), and thus $K_n$ does not contain any RED $K_{m+1}$, neither a BLUE $K_{m+1}$ and this construction seemed to be very nice. So Turán thought this maybe the
best. When Turán asked Erdős about this, right after the war, Erdős answered that in a random colouring of $K_n$ the largest monochromatic $K_p$ has order at most $(2+o(1))\log_2 n$ (for sharper results see e.g., Bollobás and Erdős
\cite{BolloErd76RandClique}). In some sense, this was the beginning of the Theory of Random graphs. Joel Spencer calls this ``the Erdős Magic'' and discusses this story in details, e.g., in \cite{Spencer11-80Years}, or in
\cite{SpencerJ13Magic}, where he describes this and also the whole story of $R(3,k)$, its estimate by Erdős \cite{Erd47Random}, by Ajtai, Komlós and Szemerédi \cite{AjtKomSzem80Ramsey}, the application of the Lovász Local Lemma
\cite{ErdLov75Local} by Spencer\cite{Spencer75R3k}, and finally the matching deep result of Jeong Han Kim \cite{Kim95RamseyK3}, using the Rödl nibble and many other deep tools.

So we see that most of the graphs $G_n$ are counterexamples to this conjecture of Turán, however, we cannot construct graph sequences $(G_n)$ without complete graphs on $\lfloor c\log n\rfloor$ vertices and independent set of vertices of
size $\lfloor c\log n\rfloor$. Actually, to construct such graphs is a famous open problem of Paul Erdős, weakly approximated, but still unsolved.\footnote{One problem with this sentence is that the notion of ``construction'' is not well
  defined, one of us witnessed a discussion between Erdős and another excellent mathematician about this, but they strongly disagreed. As to the constructions, we mention the Frankl-Wilson construction of Ramsey graphs \cite{FranklWilson},
  or some papers of Barak, Rao, Shaltiel, and Wigderson \cite{BarakRaoWigder12} and others.}

One of the beautiful conjectures is 

\begin{conjecture}[Vera Sós]
A Ramsey graph is quasi-random.
\end{conjecture}

Of course, here we should know what is a Ramsey graph and when is a graph quasi-random. We formulate this only in the simplest case. Given an integer $m$, let $N=R(m,m)$ be the smallest integer for which 2-colouring $K_N$ we must have a
monochromatic $K_m$. A Ramsey graph is a graph on $R(m,m)-1$ vertices not containing $K_m$, nor $m$ independent vertices. The notion of quasi-randomness came originally, in a slightly hidden form from the works of Andrew Thomason
\cite{Thomason87Poznan} (connected to some Ramsey problems). Next it was formulated in a more streamlined form by Chung, Graham, and Wilson \cite{ChungGrahWilson89} and here, without going into details, we ``define'' it as follows.

\begin{definition}
 For $p>0$ fixed, a sequence $(G_n)$ of graphs is $p$-quasi-random if $e(G_n)=p{n\choose2}+o(n^2)$, and the number of $C_4$'s in $G_n$ is
 $6{n\choose4}p^4+o(n^4)$, as in the random Binomial graph $R_{n,p}$, with edge-probability~$p$.
\end{definition}

The following beautiful theorem also supports the Sós Conjecture.

\begin{theorem}[Prömel and Rödl \cite{PromRodl99non-Ramsey}]
For any $c>0$ there exists a $c^*>0$ such that if neither $G_n$ nor its complementary graph contains a $K_{[c\log n]}$ then $G_n$ contains all the graphs $H_\ell$ of $\ell=[c^*\log n]$ vertices.
\end{theorem}

\minisepa

\Proclaim Lovász Meta-Theorem. Many years ago Lovász formulated the principle that the easier is to obtain a ``construction'' for a problem by Random Methods, the more complicated it is to obtain it by ``real construction''.

Supporting examples were those days, among others, \dorib the above ``missing'' Ramsey Construction, \dorib the Expander Graph problem, and \dorib also some good codes from Information Theory.

The random graphs have good expander properties. Expander graphs are important in several areas, among others, in Theoretical Computer Science. 
The fact that the random graphs are expanders were there in several early papers implicitly or explicitly, see e.g. Erdős and Rényi, \cite{ErdRenyi62Connect}, Pósa \cite{Posa76RandomHamil}.

\begin{remark} 
  Pinsker \cite{Pinsker73Conc} proved the existence of bounded degree expanders, using Random Methods, and the construction of Margulis \cite{Margul73Expan} was a breakthrough in this area. They were used e.g. in the AKS Sorting networks
  \cite{AjtKomSzem83SortingNetw}.
\end{remark}

\subsection{Dichotomy, randomness and matrix graphs}\label{DichotomyS}

In the areas considered here, there are two extreme cases: (a) sometimes for some small constant $\nu$ we partition the $n$ vertices into $\nu$ classes $U_1,\dots,U_\nu$ and join the vertices according to the partition classes they
belong to: if some vertices $x\in U_i$ and $y\in U_j$, $x\ne y$ are joined then all the pairs $x',y'$ are joined for which $x'\in U_i$ and $y'\in U_j$, $x'\ne y'$. These graphs can be described by a $\nu\times\nu$ matrix, and therefore can
be called \Sc matrix-graphs. \footnote{A generalization of these graphs is the generalized random graph, where we join the two vertices with probability $p_{ij}$, independently.} (Similar approach can be used in connection with
edge-coloured graphs, multigraphs and digraphs.)

Often such structures are the extremal ones, in some others the random graphs. We could say that a dichotomy can be observed: sometimes the
extremal structures are very simple, in some other cases they are very complicated, randomlike, fuzzy. (One very important feature of the random graphs is
that they are expanders. This is why in a random graph much fewer edges ensure Hamiltonicity than in an arbitrary graph.)

\subsection{Ramsey problems similar to extremal problems}

In some cases the Ramsey graphs are chaotic, see above, in some other (mainly off-diagonal) cases they are very similar to $\Tn n,k$. Below we shall discuss
only those cases of off-diagonal Ramsey Numbers that are strongly connected to Extremal Graph Theory.\footnote{The Ramsey numbers $R(L,M)$ form a twice
 infinite matrix whose rows and columns are indexed by the graphs $L$ and $M$. If $L\neq M$, then $R(L,M)$ is called ``off-diagonal''.} The area
where the required monochromatic subgraphs are not complete graphs started with a paper of Gerencsér and Gyárfás \cite{GerenGyarf67}. Given two graphs, $L$
and $M$, the Ramsey number $N=R(L,M)$ is the minimum integer $N$ for which any RED-BLUE-colouring of $K_N$ contains a RED $L$ or a BLUE $M$. Gyárfás and
Gerencsér started investigating these problems, Bondy and Erdős \cite{BondyErd73Rams} discussed the case when both $L$ and $M$ are cycles. Chvátal \cite{Chvat77Tree} proved that

\begin{theorem} If $\Tm$ is any fixed $m$-vertex tree, then
$$R(T_m,K_\ell)=(\ell-1)(m-1)+1.$$
\end{theorem}

A construction yielding a lower bound is obvious: consider the Turán graph $\Tn (\ell-1)(m-1),{m-1}$. Colour its edges in RED and the complementary graph in BLUE. Observe that the construction yielding a lower bound in this theorem is the
Turán graph $K_{\ell-1}(m-1,\dots,m-1)$. Faudree, Schelp and others proved many results on these types of off-diagonal Ramsey problems.

\begin{remark}[Simple Ramsey extremal structures]
  In several cases the Ramsey-extremal, or at least the Ramsey-almost extremal structures can be obtained from partitioning $n$ vertices into a bounded number $\cC_1,\dots,\cC_m$ of classes of vertices and colouring an edge $xy$ according
  to the classes of their endpoints: all the edges where $x\in \cC_i$ and $y\in \cC_j$ have the same colour, for all $1\le i,j\le m$.  This applies, e.g., to the path Ramsey numbers described by Gerencsér and Gyárfás~\cite{GerenGyarf67}.
\end{remark}

Occasionally some slight perturbation of such structures also provides Ramsey-extremal colourings. Such examples occur in connection with the cycle Ramsey
numbers, e.g., the ones in the Bondy-Erdős conjecture on the Ramsey numbers on odd cycles \cite{BondyErd73Rams}, and in many other cases. 

\subsection{Applications in Continuous Mathematics}

Toward the end of his life Turán wrote a series of papers, starting perhaps with \cite{Tur69Calga}, the last ones with Erdős, Meir, and Sós, \cite{ErdMeirSosTur1,ErdMeirSosTur2,ErdMeirSosTur3,ErdMeirSosTurCorr} on the application of his
theorem in estimating the number of short distances in various metric spaces, or in estimating some integrals, potentials.\footnote{We mention just a few related papers, for a more detailed description of this area see the remarks of
  Simonovits in \cite{Tur89Coll}, and the surveys of Katona \cite{Kat83ProbPoznan,Katon14Gruyter}.} He also liked mentioning a similar result of G.O.H. Katona \cite{Kat69XiEta} where Katona applied Turán's theorem to distributions of random
variables.  Perhaps the first result of Katona in this area was

\begin{theorem}
Let $a_1,\dots,a_n$ be $d$-dimensional vectors, ($d\ge1$), with $|a_i|\ge1$ for $i=1,\dots,n$. Then the number of pairs $(a_i,a_j)$ ($i\ne j$) satisfying $|a_i+a_j|\ge1$
is at least
$$\begin{cases} t(t-1)& \Text{if} $n=2t$ ~~~~~~~~~\text{(even)}\\ t^2& \Text{if} $n=2t+1$ \Text{(odd).} \end{cases} $$
\end{theorem}

Somewhat later A. Sidorenko (under the influence of Katona) also joined this research \cite{Sidor82,Sidor90}. They proved continuous versions of discrete (extremal graph) theorems, mostly to apply it in analysis and probability theory.

\begin{remark}
  Sidorenko also reformulated the Erdős-Simonovits conjecture \cite{Sim84Wat} in the language of integrals, \cite{Sidor91DMA,Sidor93CorrelGC}. The original conjecture had various forms, but all these forms asserted that if $L$ is bipartite,
  and $E$ is noticeably larger than $\ext(n,L)$, then among all the graphs $G_n$ with $E$ edges, the Random Graph has the fewest copies of $L$. The weakest form of this conjecture is that 

\begin{conjecture}[Erdős-Simonovits \cite{Sim84Wat}]
For any bipartite $L$, there exist two constants,
  $C=C_L>0$ and $\gamma=\gamma_L>0$, such that if $e(G_n)>C\ext(n,L)$, then $G_n$ contains at least
$$\ga\cdot n^v\left({E\over n^2}\right)^e$$ copies of $L$, for $e=e(L)$ and $v=v(L)$. 
\end{conjecture}

These forms were primarily referring to the sparse case, when $E$ is slightly above $\ext(n,L)$. On the other hand, Sidorenko's form becomes meaningful only for dense graph sequences.
For some more details on this, see Füredi and Simonovits, \cite{FureSim13Degen} or Simonovits \cite{Sim84Wat}, or Sidorenko \cite{Sidor91DMA}.
\end{remark}

\subsection{The Stability method}\label{Stability}

In this section we shall describe the Stability method in a somewhat abstract form, but not in its most general form. Stability in these cases mostly means that for a property $\cP$ we conjecture that the optimal objects have some simple
structure, and the almost optimal structures are very similar to the (conjectured) optimal ones, in some mathematically well defined sense. There are various forms of the stability methods, here we restrict ourselves to one of them. A
``property'' below is always a subset of the Universe. Generally we have two properties, $\cP$ and a much simpler property/subset $\cQ\subseteq\cP$. If a family $\cL$ of excluded graphs is given, then $\cP_n:=\cP(n,\cL)$ is the family of
$n$-vertex $\cL$-free graphs.

\begin{quote}{ The \Sc Stability~method means 
that the optimization is easy on $\cQ_n$ and we reduce the ``optimization on $\cP_n$'' to ``optimization on $\cQ_n$'', e.g., -- when we try to maximize the number of edges, -- by considering a conjectured extremal graph $S_n$ and showing that if $G_n\in\cP_n-\cQ_n$, 
then $e(G_n)<e(S_n)$. So, since the maximum is at least $e(S_n)$ it must be attained in $\cQ_n$.}
\end{quote}

We start with three examples. We wish to maximize some function $e(G)$ on the $n$-element objects of $\cP$, denoted by $\cP_n$.

\subsub Examples:// (where $e(G_n)$ is the number of edges).
\begin{enumerate}\dense
\item $\cP$ means that $K_{p+1}\not\subseteq\G$ and $\cQ$ is the family of $p$-chromatic graphs. It is easy to maximize $e(G_n)$ for $\le p$-chromatic graphs.
\item $\cP$ is the family of Dodecahedron-free graphs: $$\cP:=\{G~:~D_{12}\not\subseteq\G\},$$ $\cP_n:=\cP(n,D_{12})$, and $\cQ_n=\cQ(n,p,s)$ is the family of $n$-vertex graphs from which one can delete $\le s-1$ vertices to get a $\le
  p$-chromatic graph. It is easy to prove that $\cQ(n,6,2)\subseteq\cP$.\footnote{This is equivalent to that deleting any 5 vertices of $D_{12}$ one gets $\ge3$-chromatic graphs.} %%\NagyZ{szerint ezt duplan forditottam, de nem latom.}
  \dori Simonovits conjectured that the extremum is attained by a graph $H(n,2,6)$, where $H(n,p,s)$ is the generalization of $T_{n,p}$: the $n$-vertex graph having the maximum number of edges in $\cQ(n,p,s)$. Next he proved that if
  $G_n\in\cP-\cQ(n,2,6)$ then $e(G_n)<e(H(n,2,6))-\half n+O(1)$. So, to maximize $e(G_n)$ was reduced to maximizing it in $\cQ(n,2,6)$, which is easy.

\item\label{Octa} $\cP$ is the family of Octahedron-free graphs and $\cQ_n$ is the family of those graphs $G_n$ where $V(G_n)$ can be partitioned into $V_1$ and $V_2$ so that $G[V_1]$ does not contain $C_4$ and $G[V_2]$ does not contain $\P3$. Again,
  it is not too difficult to prove that $\cQ_n\subseteq\cP$. Using this, Erdős and Simonovits, applying a stability argument \cite{ErdSim71Octa}, determined the exact extremal graphs for large $n$. (Actually, they proved a much more general
  theorem on $\EXT(n,K_{p+1}(a_1,\dots,a_p))$.)

\end{enumerate}

It is worth mentioning the simplest case of \ref{Octa}:

\begin{theorem}[Erdős-Simonovits \cite{ErdSim71Octa}]
There exists an $n_0$ such that if $n>n_0$ and $S_n$ is extremal for the Octahedron graph $K(2,2,2)$, then $V(S_n)$ can be partitioned into into two parts, $A$ and $B$ so that $A$ spans a $C_4$-extremal graph in $S_n$, $B$ spans a $P_3$-extremal graph and each $x\in A$ is joined to each $y\in B$.
\end{theorem}

The product conjecture asserts that

\begin{conjecture}[Simonovits, see \cite{Sim83ProdBirk}]
If the decomposition class $\cM$ of a finite $\cL$ does not contain trees or forests, then each $S_n\in\EXT(n,\cL)$ is a product of $p$ subgraphs of $(n/p)+o(n)$ vertices, where $p$ is defined by \eqref{subchrom}.
\end{conjecture}

\minisepa

Let us fix a Universe, for the sake of simplicity, the universe of graphs or hypergraphs. Now we repeat what we said above, in a slightly more detailed form. The method of stability means that
\begin{enumerate}\dense
\item We consider an extremal graph problem, where some property $\cP$, and two parameters $n$ and $e$ are given (mostly $n$ is the number of vertices and $e$ is the number of the edges) and we try to optimize, say maximize $e$ for fixed
 $n$, on $\cP_n\subseteq\cP$, where $\cP_n$ is the family of objects in $\cP$ having the parameter $n$.
\item We have a property $\cQ\subseteq\cP$ ``strongly'' connected to the considered extremal problem. $\cQ_n$ is the corresponding subfamily of $\cQ$ with parameter $n$. We assume that the maximization is difficult on $\cP_n$ but easy on
  $\cQ_n$.
\item We prove that the maximum is smaller on $\cP_n-\cQ_n$ than on $\cP_n$, therefore the extremal objects in $\cP_n$, (i.e. the ones achieving the maximum) must be also in $\cQ_n$, where it is easy to find them.
\end{enumerate}

This approach, introduced by Simonovits \cite{Sim66Tihany}, turned out to be very fruitful for many problems, e.g., for the extremal problems of the icosahedron, dodecahedron, and octahedron, and for several other graph problems, and in
several hard hypergraph problems, e.g. in case of the Fano hypergraph extremal problem \cite{FureSim05Fano,KeevSudak05Fano}, or the results of Füredi, Pikhurko, and Simonovits \cite{FurePikhSim03,FurePikhSim05}, (and many similar hypergraph
results) see Subsection \ref{FanoExtreS}.  We can say that in the last twenty-thirty years it became widely used. Below we list some papers connected to the stability method, from a much longer list. See e.g. Balogh, Mousset, Skokan,
\cite{BalMousSkok17Cov}, D. Ellis, \cite{Ellis11Permu}, E. Friedgut \cite{Friedgut08Stab}, Füredi, Kostochka and Luo \cite{FureKostoLuo17}, Füredi, Kostochka, Luo and Verstraëte \cite{FureKostoVerstr16ErdGall,FureKostoLuoVerstr18ErdGall}
 W.T. Gowers and O. Hatami \cite{GowersTHatami15Groups}, Gyárfás, G.N. Sárközy, and Szemerédi, \cite{GyarfSarkoSzemer09Stabi}, Keevash \cite{Keev08Shadows}, Keevash and Mubayi \cite{KeevMubayi04StabCancel} Nikiforov and Schelp
\cite{NikifSchelp08CycleStab}, Mubayi \cite{Mubayi07Trianglefree}, Patkós
\cite{Patkos15Poset}, Samotij \cite{Samotij14Stabi}, Tyomkyn and Uzzell \cite{TyomkynUzzell15}.

Stability method is used, e.g., in Subsection \ref{PosaSeyRevisit}, more precisely, in the corresponding paper \cite{HajnalHerdadeSzeme18A-Seymour} of P. Hajnal, S. Herdade, and Szemerédi, -- however, in a much more complicated form, -- to
provide a new proof of the Pósa-Seymour conjecture, without using the Regularity Lemma, or the Blow-up Lemma.

These papers were selected from many-many more and below we add to them some on the stability of the Erdős-Ko-Rado \cite{ErdKoRado61} which
started with the paper of A.J. Hilton and E.C. Milner \cite{HiltonMilner67EKR}, 
Balogh, Bollobás, and Narayanan \cite{BalBolloNara15Transfer}
and continued with several further works, like 
Das and Tran \cite{DasTran16EKR}, 
Bollobás, Narayanan and A. Raigorodskii \cite{BolloNaraRai16EKR}, 
Devlin and Jeff Kahn \cite{DevlinKahn16EKR}, D. Ellis, N. Keller, and N. Lifshitz, \cite{EllisKellLif16Arx}, \dots 

\begin{remark}  
  In \cite{Sim66Tihany}, where Simonovits introduced this Stability Method, another stability proof method was also introduced, the \Sc Progressive~Induction.  That meant that the extremal graphs became more and more similar to the
  conjectured extremal graphs as $n$ increased, and finally they coincided. This approach was useful when the conjectured theorem could have been proved easily by induction, but it was difficult to prove the Induction Basis.
\end{remark}

\subsection{The ``typical structure''}\label{ErdKleitRothS} 

The results considered here are related to the situation described in the previous section, on the $\cP-\cQ$-stability, and have the following form: we have
two properties, a complicated one, $\cP$, and a simpler one, $\cQ\subset\cP$ and we assert:

\begin{quote}{Almost all $n$-vertex $\cP$-graphs have property $\cQ$.}
\end{quote} 

Among the simplest ones we have already mentioned or will discuss the following ones.

\begin{enumerate}\dense
\item \label{Kolai}
Almost all $K_{p+1}$-free graphs are $p$-chromatic \cite{ErdKleiRoth76Lincei,KolaPromRoth87}. 
\item Almost all Berge graphs are perfect \cite{PromStege92Berge} (cf Remark \ref{FanoTypi}).
\item \label{BalBolloSim} Almost all $K(2,2,2)$-free graphs have a vertex-partition, where the first class is $C_4$-free and the second one is $P_3$-free,
  \cite{BalBolloSim11Octa}.
\end{enumerate}

\minisepa

Erdős conjectured that, given a family $\cL$ of forbidden graphs, it may happen that a large part of the $\cL$-free graphs are subgraphs of some extremal graphs $S_n\in\EXT(n,\cL)$, in the following sense. Denote by $\cP(n,L)$ the family of
$n$-vertex $\cL$-free graphs.  The subgraphs of an $\cL$-extremal graph $S_n$ provide $2^{\ext(n,\cL)}$ $\cL$-free graphs.

\begin{conjecture}[Erdős] If $\cL$ contains no bipartite $L$, then $$|\cP(n,\cL)|=O(2^{(1+o(1))\ext(n,\cL)}).$$ \end{conjecture}

The first such result, by Erdős, Kleitman, and Rothschild \cite{ErdKleiRoth76Lincei} asserted that for
$\cL=\{K_{p+1}\}$, $$\log_2|\cP(n,\cL)|=(1+o(1))\ext(n,\cL).$$ Later Erdős, Frankl, and Rödl \cite{ErdFranklRodl86Asym} proved Erdős' conjecture, for any
non-degenerate case.  Kolaitis, H.J. Prömel, B.L. Rothschild, \cite{KolaPromRoth87} extended some related results to the case of $L$ with critical edges,
see Meta-Theorem \ref{MetaCritical}. An important related result is that of Prömel and Steger \cite{PromSteg92ColCrit}.  Slowly a whole theory was built up
around this question. Here we mention in details just a few results, and then list a few related papers.

\begin{theorem}[Erdős, Frankl, and Rödl \cite{ErdFranklRodl86Asym}] If $\cL$ does not contain bipartite graphs, and $\cP(n,\cL)$ denotes the family of $n$-vertex $\cL$-free graphs, then
$$|\cP(n,\cL)|<2^{\ext(n,\cL)+o(n^2)}.$$
\end{theorem}

For sharper results, see Balogh, Bollobás, and Simonovits \cite{BalBolloSim04,BalBolloSim09Typi,BalBolloSim11Octa}. The cases \ref{Kolai} and \ref{BalBolloSim} mentioned above are also connected to this Erdős-Frankl-Rödl theory.

If one counts the number of $L$-free graphs $G_n$ for a bipartite $L$, then one faces several difficulties. We recommend the papers of Kleitman and Winston \cite{KleitmanWinston82}, Kleitman and D. Wilson \cite{KleitmanWilson97} and of
Morris and Saxton \cite{MorrisSaxton16}.

We get another ``theory'' if we exclude {\em induced subgraphs}, see e.g. Prömel and Steger \cite{ProemSteg91Quad,PromSteg92ColCrit}, Alekseev \cite{Aleks92}, Bollobás and Thomason \cite{BolloThom95Proje,BolloThom97HeredMonot}. The theory
is similar, however, the minimum chromatic number of $\LiL$ must be replaced by another, similar colouring number.  For some further, related results see Alon, Balogh, Bollobás, and Morris \cite{AlonBalBolloMorris11}.  \dots

A more general and sharper question is when the considered family of graphs is $\cP$, and the property $\cQ$ is strongly connected to $\cP$, then one can ask: is it true that almost all graphs $G_n\in\cP$ are also in $\cQ$. This often
holds, e.g., almost all $K_3$-free graphs are bipartite.  Again, a finer result, explained below, is

\begin{theorem}[Osthus, Prömel and Taraz \cite{OsthusPromTaraz03}]\label{OstPromTarTh} Let $\cT_p(n,\Gamma)$ denote the family of $K_p$-free graphs with $\Gamma$ edges.
  If $$t_3=t_3(n):={\sqrt3\over 4}n^{3/2}\sqrt{\log n},$$ then for any fixed $\eps>0$, the probability that a random $K_3$-free graph on $n$ vertices and $\Gamma$ edges is bipartite,
  \beq{OstTarF}\PP(G_n\in\cT_3(n,\Gamma)\Longrightarrow\chi(G_n)=2)\to
\begin{cases} 1~\text{if}\Gamma=o(n),\cr 0~\text{if}\half n\le\Gamma\le (1-\eps)t_3(n);\cr 1~\text{if} \Gamma\ge (1+\eps)t_3(n). \end{cases} \eeq
\end{theorem}

The most important line of \eqref{OstTarF} is the third line.
Actually, this ``story'' started with a result of Prömel and Steger \cite{PromSteg96TriaFree},
where the threshold-estimate was around $n^{7/4}$ for the third line of \eqref{OstTarF}. Prömel and Steger conjectured that the right exponent is $3/2$. 
Łuczak \cite{Lucz00TriaFree} proved a slightly weaker related result, where the exponent was $3/2$, however, instead of asking for bipartite graphs, he asked only for ``almost bipartite'' graphs.

The meaning of this theorem is that for very small $\Gamma=e(G_n)$ most of the graphs will have no cycles, therefore they will be bipartite. For slightly
larger $e(G_n)$ odd cycles will (also) appear, so there will be a ``slightly irregular'' interval, and then, somewhat above $t_3(n)$ everything becomes nice: almost all triangle-free graphs are bipartite.

\begin{remark} (a) One could ask why $cn\sqrt{n\log n}$ is the threshold for our problem. As \cite{OsthusPromTaraz03} explains, this is connected to the fact that this is the threshold where the diameter of a random graph becomes 2.

(b) A nice result of this paper extends the theorems from $\K3$-free graphs to $C_{2h+1}$-free random graphs.

(c) Another important generalization of this result is due to Balogh, Morris, Samotij, and Warnke \cite{BalMorrSamoWarn16Typic} to any complete graph $K_p$.
\end{remark}

\Sc Further~information can be found on these questions in the paper of Balogh, Morris, Samotij, and Warnke \cite{BalMorrSamoWarn16Typic}, which, besides formulating the main results of \cite{BalMorrSamoWarn16Typic}, i.e., extending Theorem \ref{OstPromTarTh} to any $K_p$,\footnote{The Master Thesis of Warnke contained results on $K_4$.}  provides
an excellent survey of this area and its connection to several other areas, among them to the problems on extremal subgraphs of Random Graphs, investigated also by Conlon and Gowers \cite{ConGow16Sparse}, \dots, see Section~\ref{UniverseS}/\S\ref{ExtreSub}.

As we have mentioned, here the situation for the degenerate (bipartite) case
is completely different. Related results are, e.g., Balogh and Samotij \cite{BaloSamo10,BalSamo11Kmm,BalSamo11Kst}, or Morris and Saxton \cite{MorrisSaxton16}.

\begin{HistRem}
This whole story started (perhaps) outside of Graph Theory, with some works of Kleitman and Rothschild, see \cite{KleitRoth70,KleitRoth75,KleitRoth76}.

A hypergraph analog of these results was proved by Nagle and Rödl \cite{NagleRodl01}.

Among the newer results we have mentioned or should mention several results of
Prömel and Steger, e.g., \cite{PromSteg92General,PromStege92Berge},
Alon, Balogh, Bollobás, and Morris \cite{AlonBalBolloMorris11}
and, on hypergraphs,
Person and Schacht \cite{PersonSchacht09Fano}, 
Balogh and Mubayi \cite{BalMubayi11Semi,BalMubayi12Cancel}, \dots
\end{HistRem}

\subsection{Supersaturated Graphs}\label{SuperSatS} 

When Turán proved his theorem, Rademacher immediately improved it:

\begin{theorem}[Rademacher, unpublished]
If $e(G_n)>\turtwo n $ then $G_n$ contains at least $\fele n$ copies of $K_3$.
\end{theorem}

This is sharp: putting an edge into a larger class of $\Tn n,2$ we get $\fele n$ triangles. More generally, putting $k$ edges into the larger class of $T_{n,p}$ we get $\approx k({n\over p})^{p-1}$ copies of $K_{p+1}$, and in particular,
for $p=2$ we get $k\fele n $ triangles. So Erdős generalized Rademacher Theorem:

\begin{theorem}[Erdős \ev(1962) \cite{Erd62Radem}]
 There exists a $c>0$ such that for any $0<k<cn$, if $e(G_n)\ge\turtwo n+k $ then $G_n$ contains at least $k\fele n$ copies of $K_3$.
\end{theorem}

Erdős conjectured that his result holds for any $c\le\half$.\footnote{Again, there is some difference between the cases of even and odd $n$.} He also generalized his result to $K_{p+1}$ in \cite{Erd69Casop}.\footnote{Erdős' paper contains many further interesting and important results. } Lovász and Simonovits proved the Erdős conjecture, in
\cite{LovSim76Aber}, and a much more general theorem in \cite{LovSim83Birk}.

Let $F(n,L,E)$ be the minimum number of copies of $L\subseteq G_n$ with $e(G_n)=E>\ext(n,L)$ edges. Lovász and Simonovits determined $F(n,K_{p+1},E)$, for
$e(T_{n,p}) <E<e(T_{n,p})+c_pn^2$, for an appropriately small $c_p>0$, using the stability method, and, more generally, for any $q\ge p$, and $e(\Tn n,q) \le
E<e(\Tn n,q)+c_qn^2$. In Subsection~\ref{LovSimStabS} we formulate the related, widely applicable Lovász-Simonovits Stability Theorem.

\begin{remark}  
  The Lovász-Simonovits method did not work in the general case, farther away from the Turán numbers. Their Supersaturated Graph result, on the number of complete subgraphs, was extended by Fisher and Ryan \cite{Fisher89}, by Razborov
  \cite{Razbor08K3}, then by Nikiforov \cite{Nikif11Cliques}, and finally, by Reiher \cite{Reiher16CliqueDensi}. For a related structural stability theorem see also the paper of Pikhurko and Razborov~\cite{PikhuRazbo17Struct}.
\end{remark}

We complete this section with a $C_5$-Supersaturated theorem:

\begin{theorem}[Erdős \ev(1969) \cite{Erd69Casop}]
If $e(G_{2n})=n^2+1$ then $G_{2n}$ contains at least $n(n-1)(n-2)$ pentagons.
\end{theorem}

As Erdős remarks, a $K(n,n)$ with an extra edge added shows that his theorem is sharp. For some generalizations see Mubayi \cite{Mubayi10Critical}.
For early results on supersaturated graph results see e.g. the survey of Simonovits \cite{Sim84Wat}, explaining how the proofs of some extremal theorems depend on supersaturated graph results, the papers of Blakley and Roy \cite{BlakleyRoy65}, of Erdős
\cite{Erd62Radem,Erd69Casop}, Erdős and Simonovits \cite{ErdSim83Sup,ErdSim84Water}, Brown and Simonovits \cite{BrownSim84DM}, the next subsection, and many further results.

\subsection{Lovász-Simonovits Stability theorem}\label{LovSimStabS}

To prove and generalize Erdős' conjecture on $\K3$-supersaturated graphs, Lovász and Simonovits proved a ``sieve'', the simplest form of which is the following:

\begin{theorem} 
 For any constant $C>0$, there exists an $\eps>0$ such that if $|k|<\eps n^2$ and $G_n$ has $\turtwo n+k$ edges and fewer than $C |k| n$
 triangles $K_3$, then one can change $O(|k|)$ edges in $G_n$ to get a bipartite graph.
\end{theorem}

Here mostly we use $k>0$, but if we have the theorem for $k>0$, that immediately implies its extension for $k\le0$ as well. Let $m(L,G)$ denote the number of labeled copies of $L$ in $G$. The more general form is related to any $K_p$ and
not only around $\ext(n,K_p)$, but more generally, when we wish to estimate $m(K_p,\G)$ and $e(G_n)$ is around any $\ext(n,K_q)$, for $q\ge p$.

In the next, more general theorem $t$ and $d$ are defined by
$$e(G_n)=\left(1-{1\over t}\right){n^2\over 2~} \Text{and} d=\lfloor t\rfloor .$$

\begin{theorem}[Lovász--Simonovits \cite{LovSim83Birk}]\label{LovSimStab}
Let $C\ge0$ be an arbitrary constant. There exist positive constants $\de>0$ and a $C'>0$ such that if $-\de n^2<k<\de n^2$ and $G_n$ is a graph with 
$$ e(G_n)=e(T_{n,p})+k$$ edges and
$$ m(K_p,G_n)<{t\choose p}\left({n\over p}\right)^p+Ckn^{p-2},$$
then there exists a $K_d(n_1,\dots,n_d)$ such that $\sum n_i=n$,
$|n_i-{n\over d}|<C'\sqrt{k}$ and $G_n$ can be obtained from $K_d(n_1,\dots,n_d)$ by changing at most $C'k$ edges.
\end{theorem}

Here $t$ can be regarded as a ``fractional Turán-class-number''. To explain the meaning of this theorem, remember that if one puts $k$ edges into the first
class of a $T_{n,p}$, that creates $\approx c_1k n^{p-2}$ copies of $K_p$. This theorem asserts that in a graph $G_n$ with $e(T_{n,p})+k$ edges, either we
get much more copies of $K_p$, or $G_n$ must be very similar in structure to $T_{n,p}$.\footnote{This theorem may remind us of the Removal Lemma, (see
  Subsection \ref{RemovalS1}) yet, it is different in several aspects. Both they assert that either we have many copies of $L$ in $G_n$, or we can get an $L$-free graph
  from $G_n$ by deleting a few edges. However, the Removal Lemma has no condition on $e(G_n)$ and the Lovász-Simonovits theorem provides a much stricter
  structure.\dori This result can also be used for negative values of $k$, (and sometimes we need this), however, then we should replace $k$ by $|k|$ in
  some of the formulas.}

Theorem \ref{LovSimStab} can be used in many cases, e.g., it provides a clean and simple proof of the Erdős-Simonovits Stability Theorem.

\begin{remark}
  Assume that $p\ge3$ and $k=\ga n^2$, for some constant $\ga>0$. If one knows Theorem \ref{LovSimStab} for $K_p$, then one has it for any $p$-chromatic
  $L$, by applying the Erdős Hypergraph Theorem \ref{ErdHyperKST} to the $v$-uniform hypergraph on $V(G_n)$, for $v:=v(L)$, whose hyperedges are the
  vertices of the copies of $L$ in $G_n$. (For the details see, e.g., Brown and Simonovits \cite{BrownSim84DM}.) Hence Theorem~\ref{LovSimStab} is more
  general than the Erdős-Simonovits Stability theorem, since it does not completely exclude $L\subset G_n$, only assumes that $G_n$ does not contain too
  many copies of $L$.
\end{remark}

\begin{remark}  
  One could ask why do we call Theorem \ref{LovSimStab} a ``sieve''. Without answering this question, we make two remarks. 

(a) The methods used here could be considered in some sense ``primitive'' predecessors of what today is called Razborov's Flag algebras. 

(b) This whole story started with a ``survey'' paper of Lovász \cite{Lov73Sieve}, written in Hungarian, the title of which was ``Sieve methods''.
\end{remark}

\begin{remark} 
Often the Lovász-Simonovits sieve can be replaced by the Removal Lemma.
\end{remark}

\begin{remark}  
  The original proofs of the Erdős-Simonovits Limit theorem could have used the Regularity Lemma, (described in Section \ref{RegularityS}), or this theorem, but they had not: actually they were proved earlier. An alternative approach to prove the Erdős-Simonovits Stability theorem is to use the Regularity Lemma, however, then one needs the stability theorem itself for complete graphs. A simple, beautiful proof of that the stability holds for $K_p$ was found by Füredi \cite{Fured15SimStab}, who used for this purpose the Zykov symmetrization \cite{Zykov49}.
\end{remark} 

\subsection{Degenerate vs Non-degenerate problems}\label{DegenS}

We remind the reader that an extremal graph problem is \Sc Degenerate if $\ext(n,\cL)=o(n^2)$, or in case of $r$-uniform hypergraphs, $\ext(n,\cL)=o(n^r)$. Another way to describe a \Sc Degenerate extremal problem is that $\cL$ contains a bipartite $L$ (and for hypergraphs, $\cL$ contains an $L$ with strong chromatic number $r$, see Claim~\ref{HyperDegen}). The survey of Füredi and Simonovits \cite{FureSim13Degen} describes the details. Here we shall be very brief, describe just a few results and formulate three conjectures and describe some connections to Geometry, Finite Geometry, and Commutative Algebra.

So let us restrict ourselves to simple graphs. The simplest questions are when $L=K(a,b)$ and when $L=C_{2k}$. (The extremal problem of the paths $P_k$ is a theorem of Erdős and Gallai
\cite{ErdGallai59Path}.)

For $a=2$ and $a=3$ the sharp lower bounds came from finite geometric constructions of Erdős, Rényi, V. T. Sós, \cite{ErdRenyiSos} and W.~G.~Brown \cite{Brown66Thomsen}. The random methods gave only much weaker lower bounds. By \cite{ErdRenyi60Evol},\footnote{Actually, the First Moment Method yields a little better estimate, but far from being satisfactory. One can also see that as $a$ is fixed and $b$ gets larger, the ``random construction'' exponents converge to the optimal one. This motivates, among others, \cite{KollRonyaiSzab96,AlonRonyaiSzab99}}
$$\ext(n,K(a,b))>c_an^{2-{1\over a}-{1\over b}}.$$
The Kollár-Rónyai-Szabó construction \cite{KollRonyaiSzab96} and its improvement, the Alon-Rónyai-Szabó \cite{AlonRonyaiSzab99} construction, (proving the sharpness of Theorem \ref{KovSosTurTh} when $a$ is much smaller than $b$) used Commutative Algebra, and Lazebnik, Ustimenko, and Woldar have several more involved algebraic constructions.

\begin{remark}
  One of Erdős' favourite geometry problem was the following: Given $n$ points in $\EE^d$, how many equal (e.g., unit) distances can occur among them. Among others, he observed that if in $\EE^3$ we join two points iff their distance is 1, then the resulting graph does not contain $K(3,3)$. In
  Brown's construction this is turned around: the vertices of a graph $G_n$ are the $n=p^3$ points of a finite 3-dimensional affine space $AG(p,3)$. W.G. Brown joined two points $(x,y,z)$ and $(x',y',z')$ if their ``distance'' was ``appropriate'': $$(x-x')^2+(y-y')^2+(z-z')^2=\alpha,\qquad ({\rm mod}~p).$$ Then Brown proved that (for some primes $p$ and some $\alpha$) the resulting graph contains no $K(3,3)$ and has $\approx\half n^{2-(1/3)}$ edges. (Surprisingly, as Füredi proved in \cite{Fure96Zaran}, -- as to the multiplicative constant $\half$, -- for $K(3,3)$ the Brown construction
  is the sharp one, not the upper bound of Theorem \ref{KovSosTurTh}.) The Commutative Algebra constructions \cite{KollRonyaiSzab96,AlonRonyaiSzab99} can be regarded as extensions of Brown's construction, however, with deeper mathematics in the background.
\end{remark}

\begin{questio}\label{WhyNonComm}
Sometimes in our constructions we use commutative structures, sometimes non-commutative ones.
 One could ask, what is the advantage of using non-commutative structures.
\end{questio}

The same question can also be asked in connection with the so called Ramanujan graphs, see e.g., \cite{LuboPhilSarnak86STOC,LuboPhilSarnak88CCA}, or \cite{Margul84Arith}. The answer is simple: the Cayley graphs of commutative groups are
full of short even cycles.  We illustrate this through the girth problem.

\begin{theorem}[Bondy-Simonovits, Even Cycle: $C_{2k}$ \cite{BondySim74C2k}]\label{BondySimTh}
\beq{BonSimFo}\ext(n,C_{2k})\le c_1 k n^{1+(1/k)}.\eeq
\end{theorem}

\begin{conjecture}[Sharpness]\label{C2kSharp}
 The exponent $1+(1/k)$ is sharp, i.e.
$$\ext(n,C_{2k})\ge c_k n^{1+(1/k)} \Text{for some} c_k>0 .$$
 \end{conjecture}

 The first unknown case is $k=4$, $\ext(n,C_8)$. Theorem \ref{BondySimTh} is sharp for $C_4$, $C_6$, $C_{10}$, see \cite{Erd38Tomsk},
 \cite{ErdRenyiSos}, \cite{Brown66Thomsen} and \cite{Benson66}.
 For some related constructions see also Wenger \cite{Wenger91Constr}.

 Now, to answer Question \ref{WhyNonComm}, observe, that we often use in our algebraic constructions Cayley graphs, where in the commutative cases we have many ``coincidences'', leading to many $C_{2k}\subset G_n$, which can be avoided
 in the non-commutative cases. A very elegant (and important) example of this is:

\begin{construction}[Margulis, \cite{Margul82CCA}]
There exist infinite Cayley graph sequences $(G_{n,d})$ of degree $d=2\ell$ with girth greater than $c{\log n\over \log (d-1)}$.
\end{construction} 

\begin{remark} 
  We have emphasized that Extremal Graph Theory is connected to many other areas in Mathematics. The Margulis constructions \cite{Margul82CCA} are connected to Coding Theory. In somewhat different ways, several papers of Füredi and Ruszinkó are also extremal hypergraph results strongly connected to (or motivated by) Coding Theory, e.g., \cite{FureRuszi13Grid}.
\end{remark}

A more detailed analysis of these questions can be found, e.g., in Alon's survey \cite{Alon95HB}, or in the Füredi-Simonovits survey \cite{FureSim13Degen}.

\begin{remark}
  (a) It was a longstanding open question if one can improve the coefficient of $n^{1+(1/k)}$ in the Bondy-Simonovits theorem, from $ck$ to $o(k)$. After several ``constant''-improvements, Boris Bukh and Zilin Jiang \cite{BukhJiang17}\footnote{The original version claimed a slightly better estimate.} proved that \beq{BonSimBZ}\ext(n,C_{2k})\le 80 \sqrt k\log k\cdot n^{1+(1/k)}+O(n).\eeq According to \cite{BukhJiang17}, Bukh thinks that Conjecture \ref{C2kSharp} does not hold: he conjectures that for sufficiently large, but fixed $k$,
$$\ext(n,C_{2k})=o(n^{1+(1/k)}).$$
It is very ``annoying'' that we cannot decide this, not even for $C_8$.

(b) Related constructions were provided by Lazebnik, Ustimenko, and Woldar \cite{LazebUstimWold94,LazeUstiWold95Girth,LazebWold01} and by Imrich \cite{Imrich84Girth}.

(c) For some ordered versions of the $C_{2k}$ problem see, e.g., the (very new) results of Győri, Korándi, Tomon, Tompkins, and Vizer \cite{GyoriKoraVizer18C6}.
\end{remark}

\begin{HistRem}
(a) Actually, before the Erdős-Simonovits paper \cite{ErdSim69Cube} Erdős conjectured that the exponents can be only either $(2-\reci k )$ or $(1+\reci k )$. This conjecture was ``killed'' in \cite{ErdSim69Cube}, by some ``blow-up'' of the cube.

(b) Later Erdős and Simonovits conjectured that (i) for any rational exponent $\alpha\in[1,2)$ there exist degenerate extremal problems $\ext(n,\cL)$ for
which \beq{RatiExpo}\ext(n,\cL)/n^\alpha \to c_\cL>0,\eeq and (ii) for any degenerate problem $\ext(n,\cL)$ there exists a rational $\alpha\in[1,2)$ for which
\eqref{RatiExpo} holds. Recently, Bukh and Conlon \cite{BukhConlon18Rati} proved (i). (For hypergraphs Frankl \cite{Frankl85Rationals} has some earlier, corresponding results, see also Fitch \cite{Fitch16A-Rati}).
\end{HistRem}

\subsection{Dirac theorem: introduction}
 
Difficult problems always played a central role in the development of Graph Theory. We shall mention here two important problems, the Dirac theorem on Hamiltonian cycles (which is not so difficult) and the Hajnal-Szemerédi theorem on
equitable partitions.

\begin{theorem}[Dirac \ev(1952) \cite{Dirac52Abstr}]\label{DiracX} If $n\ge3$ and $\mindeg(G_n)\ge n/2$, then $G_n$ contains a Hamiltonian cycle.
\end{theorem}

If $n=2h$, then $K(h-1,h+1)$ has no Hamiltonian cycle, showing that Theorem~\ref{DiracX} is sharp.\footnote{The union of two complete graphs of $n/2$ vertices having at most one common vertex in common also shows the sharpness. For $=2\ell-1$ one can use $K(\ell,\ell-1)$ for sharpness.} One beautiful feature of this theorem is that as soon as we can guarantee a 1-factor, we get a Hamiltonian cycle.

 This theorem triggered a wide research, see, e.g., Ore's theorem \cite{Ore60Hamil}, Pósa's theorem on Hamiltonian graphs \cite{Posa62Hamil}, and related
 results. We shall return to discuss generalizations of Dirac's Theorem, above all, the Pósa-Seymour conjecture and the hypergraph generalizations in
 Sections \ref{PosaSeyRevisit} and \ref{HyperLargeExcludedS}.

\subsection{Equitable Partition}\label{HajSzemS}

We close this introductory part with a famous conjecture of Erdős, proved by András Hajnal and Szemerédi.

\begin{definition}[Equitable colouring]
 A proper vertex-colouring of a graph $G$ is \Sc Equitable if the sizes of any two colour classes differ by at most one.
\end{definition}

\begin{theorem}[Hajnal-Szemerédi \ev(1969) \cite{HajSzem69Corrad}]\label{HajSzemCorrTh}
For every positive integer $r$, every graph with maximum degree at most $r$ has an equitable colouring with $r+1$ colours.
\end{theorem}

The theorem is often quoted in its complementary form.
The sharpness is shown by the complementary graph of an almost-Turán graph, i.e. the union of complete graphs $K_r$ and $K_{r+1}$.

The theorem was proved first only for $K_3$, by K. Corrádi and A. Hajnal \cite{CorradiHajnal}. Then came the proof of Hajnal and Szemerédi. Much shorter and simpler proofs of Theorem \ref{HajSzemCorrTh} were found independently by Kierstead
and Kostochka, \cite{KierKosto08Equit} and Mydlarz and Szemerédi \cite{MydlarzSzemer07Manu}. The paper of Kierstead, Kostochka, Mydlarz
and Szemerédi \cite{KierKostoMydlSzemer10Equ} provides a faster (polynomial) algorithm to obtain the equitable colouring.\footnote{See the introduction of \cite{KierKostoMydlSzemer10Equ}. They also point out the applications of this theorem,
  e.g., in \cite{AlonFure92Spanning}, \cite{JansonRuc02Tail} for ``deviation bounds'' for sums of weekly dependent random variables, and in the Rödl-Ruciński proof of the Blow-up Lemma \cite{RodlRuc99Blow}.}

\subsub Multipartite case.// This theorem was generalized in several different ways, also, considered for multipartite graphs $G_n$ by Martin and Szemerédi \cite{MartinSzem08QuadHSz}, Csaba and Mydlarz \cite{CsabaMyd12Multi}, Extensions
towards multipartite hypergraphs can be found, e.g., in Lo-Markström \cite{LoMarks13MultiHajSzem}.  It is applied to prove some other graph theorems, e.g., in Komlós-Sárközy-Szemerédi
\cite{KomSarkoSzem98SeymApprox,KomSarkoSzem98SeymAnnals,KomSarkoSzemer01AlonYus}, and many other cases. We shall return to these questions in Section~\ref{ColouringProblems}.

The theorem was extended also to directed graphs, by Czygrinow, DeBiasio, Kierstead, and Molla \cite{CzygBiasKierMolla15}.

\subsub Random Graphs.// The problem of equitable partitions in Random Graphs was also discussed (and in some sense solved) in works of Bohman, Frieze,
Ruszinkó, and Thoma \cite{BohFrieRuszThoma01}. Their result was improved by Johansson, Jeff Kahn, and Van Vu \cite{JohKahnVu08Factors}.

\subsection{Packing, Covering, Tiling, $L$-factors}\label{TilingS}

When speaking of ``packing'', sometimes we mean edge-disjoint embedding of just two graphs (this is connected to some complexity questions from Theoretical Computer Science) and sometimes we try to cover the whole graph with some
vertex-disjoint copies of a graph.\footnote{In the Gyárfás Conjecture we try to pack many different trees into a complete graph.} 

There was a period, when -- because of Theoretical Computer Science, -- packing a graph into the complementary graph of another (i.e. the above problem for two {\em edge-disjoint} graphs) was a very actively investigated topic. This was connected to \Sc Evasiveness, i.e. to the problem, how many ``questions'' are needed to decide if a graph $G_n$ has property~$\cP$.

The whole area is described in a separate chapter of Bollobás' ``Extremal Graph Theory'' \cite{Bollo78ExtreBook}. For some further related details see also
the papers of Bollobás and Eldridge, e.g., \cite{BolloEldri78Packing} with the title ``Packing of Graphs and Applications to Computational Complexity''. Here
we mention also the result of P.~Hajnal \cite{HajnalP91Evasive}, (improving some important earlier results). He proves that the randomized decision tree
complexity of any nontrivial monotone graph property of a graph with $n$ vertices is $\Omega(n^{4/3})$. See also Bollobás \cite{Bollo76Elus}, and
\cite{King88ACM,King90CCA}. This is again a nice example where Combinatorics and Theoretical Computer Science are in a very strong interaction.

\minisepa

Given a graph $G_n$, and a sample graph $L$, a \Sc perfect~tiling is a covering of $V(G_n)$ with {\em vertex-independent} copies of $L$, and an \Sc almost-tiling is covering at least $n-o(n)$ vertices of it. Tiling is sometimes a tool, a method, in other cases it is the aim.
The Corrádi-Hajnal and the Hajnal-Szemerédi theorems, in Section \ref{HajSzemS}, were also about packing=tiling of graphs.

Perhaps the first case when tiling was used in a proof was the Rademacher-Turán theorem of Erdős \cite{Erd62Radem}, Here Erdős covered as many vertices of $G_n$ by vertex-independent triangles as he could, to prove the theorem. 

We can consider the problem of \Sc Tilings in a more general way. Here we have a (small) sample graph $L$ and wish to embed into $G_n$ as many vertex-independent copies\footnote{or, in other cases edge-disjoint copies. Here ``vertex-independent'' and ``vertex-disjoint'' are the same.} of $L$ as possible. The question is, given an integer $t$, which conditions on $G_n$ ensure, e.g, $t$ vertex-independent copies of $L$. We shall return to this question later, here we mention just a few related results,
as illustrations. As we mentioned, Theorem \ref{HajSzemCorrTh} is an example of such results. Below we formulate some more general results.

\begin{definition}
 Given two graphs $L$ and $G$, where $v(L)$ divides $v(G)$, an $L$-factor of $G$ is a family of $v(G)/v(L)$ vertex-disjoint copies of $L$,
\end{definition}

\begin{theorem}[Alon and Yuster \ev(1996) \cite{AlonYuster96Factor}]
For every $\eps>0$, and for every positive integer $h$, there exists an $n_0=n_0(\eps,h)$, such that for every ``sample'' graph $L_h$, 
every ``host'' graph $G_n$ with $n=h\ell>n_0$ vertices and
minimum degree
\beq{TilingFo}\mindeg(G_n)\ge \left(1-{1\over \chi(L)}+\eps\right)n \eeq
contains an $L$-factor.
\end{theorem}

As to the usage of names, tiling, packing, and perfect $L$-factor are almost the same: given a graph $G$ and a sample graph $L$, we wish to embed into $G$ as many vertex-independent copies of $L$ as possible, and if they (almost) cover $V(G)$, then we speak about an (almost) perfect tiling/packing.

Komlós extended the notion of $L$-factor by saying that $G$ has an $L$-factor, if it contains $\lfloor v(G)/v(L)\rfloor$ vertex-independent copies of $L$. Alon and Yuster conjectured and Komlós, G.N. Sárközy and Szemerédi
\cite{KomSarkoSzemer01AlonYus} proved the following\footnote{Actually, they formulated two ``similar'' conjectures, we consider only one of them.}

\begin{theorem}[Komlós, Sárközy, and Szemerédi]\label{KomSarkoSzem01AlonYusTil}

For every $L$ there is a constant $K=K_L$ such that if
$$\mindeg(G_n)\ge \left(1-{1\over \chi(L)}\right)n+K_L, $$
then $G_n$ has an $L$-factor.
\end{theorem}

Komlós surveyed the tiling situation in \cite{Koml00Tiling} in a more general way. He considered degree conditions for finding many disjoint copies of a
fixed graph $L$ in a large graph $G$. Let $\tau(n,L,M)$ be the minimum $m$ for which if $\mindeg(G_n) \geq m$, then there is an $L$-matching covering at
least $M$ vertices in $G$. For any fixed $x\in(0,1)$, Komlós determined $$f_{L}(x) = \lim_{n\to\infty} \reci n \tau(n,L,xn).$$

Thus, e.g, Theorem~8 of \cite{Koml00Tiling} determines (for a fixed but arbitrary $L$) a sharp min-degree condition on $G_n$ to enable us to cover $\approx
xn$ vertices of $G_n$ by copies of $L$. Among others, Komlós analyzed the strange and surprising differences between the cases when we try to cover $G_n$ with
vertex-independent copies of $L$ completely, and when we try only to cover it almost completely, or when we try to cover only a large part, $\lfloor
xn\rfloor$ vertices of it. Several further details and a careful analysis of this situation can be found in the paper of Komlós \cite{Koml00Tiling}.
%%We shall not repeat that, just mention a few newer results.

\minisepa

\def\hcr{\chi_{cr}}
\def\hcf{{\bf hc}_{\chi}}
\def\hcs{{\chi^*}}

Kawarabayashi conjectured that Theorem \ref{KomSarkoSzem01AlonYusTil} can be improved by taking into account the particular ``colouring structure'' of $L$.
This was proved by Cooley, Kühn, and Osthus \cite{CooleyKuhnOst07Pack}, and then by Kühn, and Osthus \cite{KuhnOst09Pack}. Komlós, in \cite{Koml00Tiling}
defined a chromatic number $\hcr(L)$, improving the earlier results, by using this $\hcr$. Kühn and Osthus defined another ``colouring parameter'' $\hcf$ depending on the sizes
of the colour-classes in the optimal colourings of $L$. Using $\hcf$, they defined a $\hcs(L)\in[\chi(L)-1,\chi(L)]$ and proved

\begin{theorem}[Kühn and Osthus \ev(2007) \cite{KuhnOst09Pack}]
If $\de(L,n)$ is the smallest integer $m$ for which every $G_n$ with
$\mindeg(G_n)\ge m$ has a perfect $L$-packing then
$$\de(L,n)=\tur{n}{\hcs(L)}+O(1). $$
\end{theorem}

\subsub Bipartite Packing. // The packing problems, as many similar problems, have also bipartite versions. Hladký and Schacht -- extending some results of Yi Zhao \cite{YiZhao09Tiling} -- proved 

\begin{theorem}[Hladký and Schacht \ev(2010) \cite{HladSchacht10Tiling}]
Let $1\le s<t$, $n=k(s+t)$. If 
$$\Phi_2(n, L):=\begin{cases}\half n+s-1 ~~~~~~~ \text{if}k\text{is even} \cr \half(n+t+s)-1 \text{if} k \text{is odd,}\end{cases}$$
then each subgraph $G\subseteq K(n,n)$ with minimum degree at least $\Phi_2(n,L)$ contains a $K(s,t)$-factor, and this is sharp, except if $t\ge 2s+1$ and $k$ is odd. 
\end{theorem}

For the ``missing case'' (when $t\ge 2s+1$ and $k$ is odd) see Czygrinow and DeBiasio~\cite{CzygBias11Tiling}.

\begin{remark}
Here we have to be careful with using the word ``factor'': Lovász, in his early papers \cite{Lov70Calgary,Lov70Prescrib} called a subgraph $H_n\subseteq G_n$ an $f$-factor if 
$$f:V(G_n)\to\NN,\Text{and}\deg_H(x)=f(x)\Text{for each}x\in V(G_n).$$ 
\end{remark}

\subsub Directed graphs.// There are many further results connected to this area. We shall return to some of them later, e.g., in Sections \ref{RT-Match} and \ref{TilingHyperS}.

Here we close this area with two analogous results on \Sc oriented~graphs: directed graphs without loops, where between any two vertices there is at most one arc. $\de^0(G)$ is the minimum of all the in- and out-degrees.

\begin{theorem}[Keevash and Sudakov \cite{KeevSudak09Triangle}]
There exist a fixed $c>0$, and an $n_0$ such that if $G_n$ is an oriented graph on $n>n_0$ vertices and $\de^0(G_n)>(\half-c)n$, then $G_n$ contains a ``cyclic triangle'' tiling which leaves out at most 3 vertices. This is sharp.
\end{theorem} 

\BalAbraCapMed{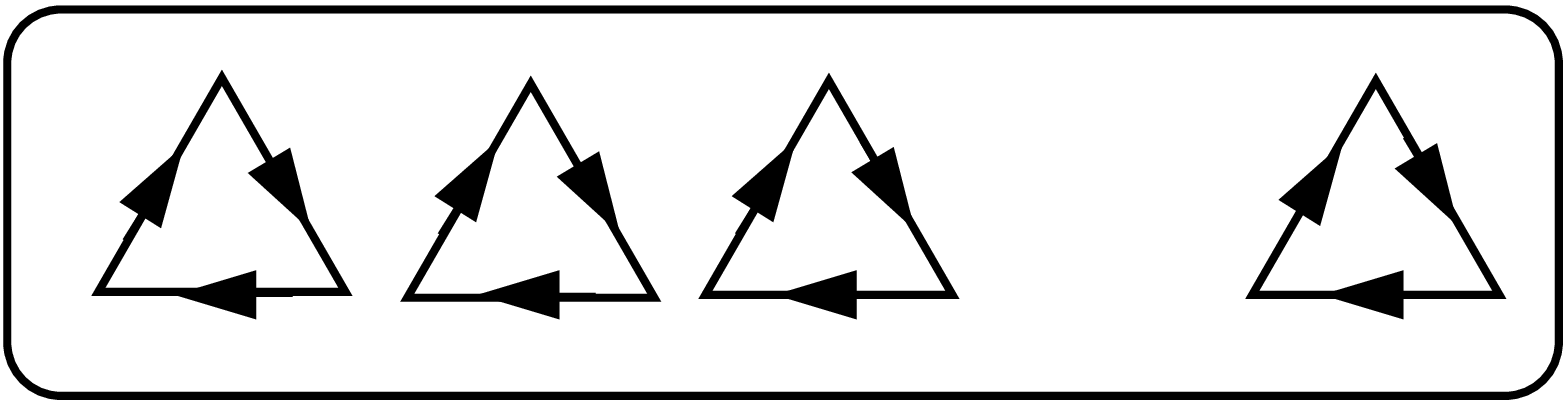}{Cyclic triangles}{40} Actually, Keevash and Sudakov \cite{KeevSudak09Triangle} describe the history of this theorem in details, explain several related results, and prove that the theorem is sharp in the sense that
if $n\equiv 3~ ({\rm mod}~ 18)$ then one cannot guarantee covering the whole graph with cyclic triangles, even under a stronger degree-condition. Further, they prove some other packing theorems where the lengths of some cycles are
prescribed but they need not be triangles. We close with

\begin{theorem}[Balogh, Lo, and Molla \cite{BalLoMol17Til}]
 There exists an $n_0$ such that for every $n\equiv0~ ({\rm mod }\, 3)$, if $n>n_0$, then any oriented graph $G_n$ on $n$ vertices with $\de^0(G_n)\ge{7n\over18}$ contains a $TT_3$-tiling, where $TT_3$ is the transitively oriented triangle.
\end{theorem}

\BalAbraCapMed{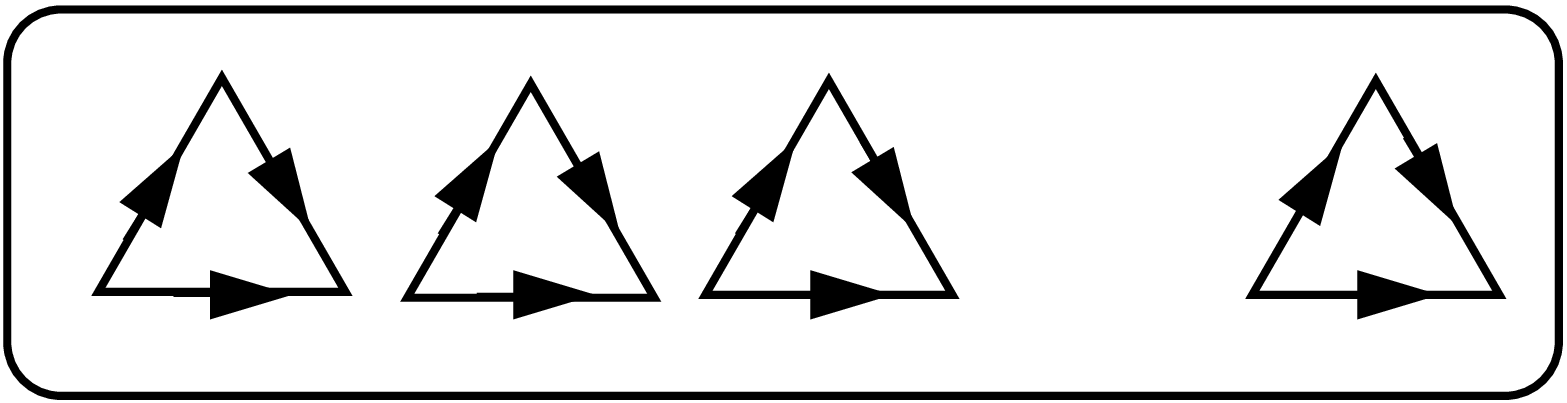}{Transitive triangles}{45} This area is fairly active nowadays, we refer to several papers on equitable colouring of Kostochka and others, e.g., \cite{KierKost08Ore,KostNakpraPemma05,KostoPelsWest}, and we mention
just some extensions to bipartite tiling, e.g., Czygrinow and DeBiasio \cite{CzygBias11Tiling}, and to oriented graphs, by Czygrinow, DeBiasio, Kierstead, and Molla \cite{CzygKierMolla14}, \cite{CzygBiasKierMolla15}, Yuster on tournaments
\cite{Yuster03Tiling} or \cite{CzygBias11Tiling}, 
\cite{CzygBiasNagle14Tiling}, or to hypergraphs, e.g., Pikhurko \cite{Pikhu08HyperTiling} or Czygrinow, DeBiasio, and Nagle \cite{CzygBiasNagle14Tiling}, and several papers of Yi Zhao and Jie Han, e.g., \cite{HanZhao15TilingK43}, and
many-many others.  Some of these results are good examples of the application of the Absorbing methods, discussed in Subsection \ref{Absorb}. This is, e.g., the case with the Balogh-Lo-Molla theorem or the Czygrinow-DeBiasio-Nagle result
just mentioned.

\Section{``Classical methods''}{Detours}

Here, (giving some preference to results connected to Laci Lovász) we still skip the area of graph limits, we skipped the applications of Lovász Local Lemma (though Lovász Local Lemma is among the most important tools in this area and it
came from research strongly connected to extremal graph problems \cite{ErdLov75Local}), but we mention two other, less known results of Lovász, strongly connected to extremal graph theory.

\subsection{Detour I: Induction?}\label{Linear-AS}

Speaking of methods in Extremal Graph theory, we mostly avoided speaking of hypergraph results, partly because they are often much more involved than the corresponding results on ordinary graphs.  One is the construction of graphs with high girth and high chromatic number. Erdős used Random Graphs to prove the existence of such graphs
\cite{Erd59GraphProb,Erd61GraphProbII} and there was a long-standing challenge to {\em construct} them.\footnote{We must repeat that the meaning of ``to construct them'' is not quite well defined. Let us agree for now that the primary aim was to eliminate the randomness.} Lovász often solved ``such problems'' by trying to use induction, and when
this did not work directly, to look for and find a stronger/more general assertion where the induction was already applicable. In {\em this} case Lovász generalized the problem to hypergraphs and used induction \cite{Lov68Chromatic}. His proof was a tour de force, rather involved but worked. It was the first ``construction'' to get graphs of high
girth and high chromatic number.

How is this problem connected to Extremal Graph Theory? 

\begin{theorem}
If $\cL$ is a finite class of excluded graphs, then
$\ext(n,\cL)=O(n)$ iff $\cL$ contains some tree (or forest).
\end{theorem}

To prove this, one needs the Erdős result \cite{Erd61GraphProbII} about high girth graphs $G_n$ with $e(G_n)>n^{1+\alpha}$ edges, or the Erdős-Sachs theorems \cite{ErdSachs63}, \dots, or some corresponding Lubotzky-Phillips-Sarnak-Margulis graphs, see e.g. \cite{LuboPhilSarnak86STOC,LuboPhilSarnak88CCA}.  (Lovász' tour de force construction was a big breakthrough into this direction though it was not quantitative, which is needed above).

\begin{remark}  
 Since that several alternative constructions were found. We mention here the construction of Nešetřil and Rödl \cite{NesetRodl79Amalg}. Perhaps one of the
 best is the construction of Ramanujan graphs by Lubotzky, Phillips, and Sarnak \cite{LuboPhilSarnak88CCA} and Margulis \cite{Margul88Explic}: it is very
 direct and elegant. It has only one ``disadvantage'': to verify that it is a good construction, one has to use deep Number Theory. (For more detailed
 description of the topic, see \cite{LuboPhilSarnak88CCA}, or the books of Lubotzky \cite{Lubo94Book} and of Davidov, Sarnak, and Valette \cite{DaviSarnakValeBook}.)
\end{remark}

\subsection{Detour II: Applications of Linear Algebra}\label{LinearS}

We start with repeating a definition.

\begin{definition}
$H$ is a colour-critical graph\footnote{more precisely, edge-colour-critical
graph.} if for any edge $e$ of $H$, $\chi(H-e)<\chi(H)$. It is
$k$-colour-critical if, in addition, it is $k$-chromatic.
\end{definition}

\BalAbraCap{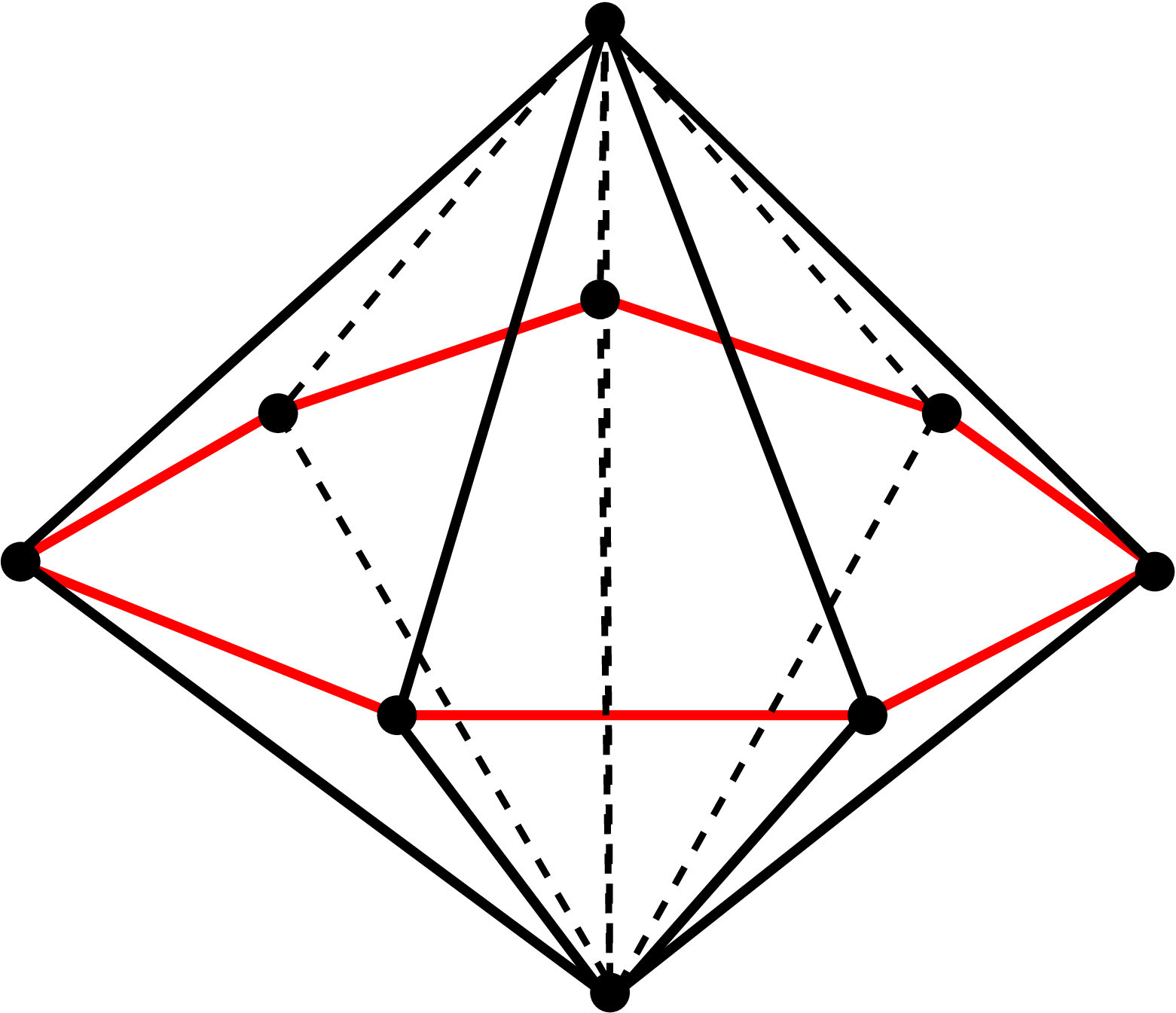}{Double Cone} The theory of colour-critical graphs is a fascinating area. Erdős writes about its beginnings, e.g., in his article in the memory of Gabriel Dirac \cite{Erd89Dir}.  Gallai also had several interesting
conjectures on colour critical graphs. One of them, on the independence number of a 4-chromatic colour-critical graph, was disproved by a construction of Brown and Moon \cite{BrownMoon69}, and then by Simonovits \cite{Sim72Criti} and Toft
\cite{Toft72Criti}. The disproof was strongly connected to a hypergraph extremal problem discussed also by Brown, Erdős, and Sós \cite{BrownErdSos73a}. Lovász improved the corresponding results, proving the following sharp and much more
general result.

\begin{theorem}[Lovász, \ev(1973) \cite{Lov73Criti,Lov76ChromHypLin}]\label{Lov73CritiTh}
  Let $\alpha_k(n)$ denote the maximum number of independent vertices in a $k$-critical graph on $n$ vertices.  Then
$$n-2kn^{1/(k-2)}\le\alpha_k(n)\le n-(k/6)n^{1/(k-2)}.$$
\end{theorem}

The lower bound is based on generalizing the Brown-Moon construction and the upper bound improves the result of Simonovits, $\alpha_4(n)\le n-c_2n^{2/5}$. Simonovits in \cite{Sim72Criti} and \cite{Sim72Hyper} used a hypergraph extremal problem, where the excluded hypergraphs were 3-uniform double-cones.\footnote{An $r$-cone is a 3-uniform hypergraph obtained from a cycle $x_1,\dots,x_k$ by adding $r$ new vertices $y_1,\dots,y_r$ and all the triples $y_jx_ix_{i+1}$ (where $x_{k+1}=x_1$).} One of his results was basically equivalent to a results of Brown, Erdős, Sós \cite{BrownErdSos73a}, where the excluded graphs were all the triangulations of a 3-dimensional sphere (the double cones are among these hypergraphs). Actually, Simonovits proved

\begin{theorem}[Simonovits \cite{Sim72Criti}] If $\Chyp r3$ denotes the (infinite) family of 3-uniform $r$-cones, then
$$\ext_3(n,\Chyp r3)=O(n^{3-\reci r }).$$
\end{theorem}

(Watch out, here the subscript $r$ is not the number of vertices!)
Simonovits, and then Toft, and Lovász, reduced the general case of the Gallai problem to the case when each $x\in X$ had degree $d(x)=k-1$, and these
vertices $x$ defined a $(k-1)$-uniform hypergraph $\cH_m$ on the remaining vertex-set $\cM:=V(G_n)-X$. Simonovits proved (for $k=4$) that this hypergraph $\cH_m$
does not contain any double-cone. Therefore $e(\cH_m)=O(m^{5/2})$. This led to his estimate. Lovász -- starting with the same approach -- excluded many more
$(k-1)$-uniform hypergraphs. To make his argument easier to follow, we restrict ourselves to the case $k=4$. Let $\cM:=\{u_1,\dots,u_m\}$. The neighbourhoods of $u_i$'s
defined a 3-uniform hypergraph on $[1,m]$, and Lovász attached to each $u_i$ a 0-1 vector of length ${m\choose 2}$, where, for a neighbourhood
$N(u_i)=\{a,b,c\}$, Lovász put 1's into the places $(a,b)$, $(b,c)$, and $(c,a)$, thus obtaining a 0-1 vector $\xx_i$ of length ${m\choose2}$, with three
coordinates equal to 1's (and all the others were equal to 0). The 4-criticality implied that these vectors were linearly independent.\footnote{In graph-theoretical language, Lovász excluded all the 3-graphs for
  which the Sperner Lemma holds: for which each pair was contained by an even number of triples.} Therefore, their number was at most ${m\choose2}$, i.e.,
we obtained that $e(\cH_m)=|X|\le {m\choose2}$. This gave the upper bound in Theorem \ref{Lov73CritiTh}.

 So he used the Linear Algebra method, basically unknown those days, to get the sharp result in his extremal graph problem. That gave sharp result also in
 the Gallai problem.

\begin{remark}\label{Criticality} For any monotone graph property $\cP$ we may define the critical structures: $G_n$ belonging to $\cP$ but after the deletion of any edge (or, in other cases, any vertex) we get a graph outside of $\cP$. If we have a graph-parameter on graphs, then criticality mostly means increasing/decreasing this parameter by deleting any edge. 
  Criticality was discussed for the stability number, chromatic numbers, diameter, and can be investigated for many other monotone properties.  Among criticality problems colour-criticality seems to be one of the most interesting ones.

  The fascinating area of colour critical graphs was started by G. Dirac, see e.g., Dirac \cite{Dirac52Criti,Dirac74Criti} and Erdős \cite{Erd89Dir}. There are several results on it in the Lovász book \cite{Lov79CombExerc}. We skip here the topic of colour-critical hypergraphs, listing just a few fascinating results on them: e.g., Abbott and Liu
  \cite{AbbottLiu78} and results of Krivelevich \cite{Krive97Crit}. Deuber, Kostochka, Sachs \cite{DeubKostoSachs96Criti}, Rödl and Siggers \cite{RodlSiggers06} and Anstee, Fleming, Füredi and Sali \cite{AnsteeFlemFureSali05,FureSali12}. Toft \cite{Toft72Criti}, and Simonovits \cite{Sim72Criti}, Rödl and Tuza \cite{RodlTuza85Crit}, Sachs and
  Stiebitz \cite{SachsStieb81Criti,SachsStieb85Criti,SachsStieb89Constr} Stiebitz, \cite{Stieb97Criti} Kostochka and Stiebitz \cite{KostoStieb02ListCrit}, Stiebitz, Tuza, and Voigt \cite{StiebTuzaVoigt09Criti}.

See also the survey of Sachs and Stiebitz, \cite{SachsStieb89Constr}.
\end{remark}

We conclude this part with two open problems.

\begin{problem}[Erdős]
Is it true that if $(G_n)$ is a sequence of 4-critical graphs, then $\mindeg(G_n)=o(n)$?
\end{problem}

Simonovits \cite{Sim72Criti} and Toft \cite{Toft72Criti} constructed 4-colour-critical graph sequences $(G_n)$ with $\mindeg(G_n)>c\root3\of n$.

\begin{problem}[Linial]
  How many triples can a 3-uniform $\cH_n$ have without containing a triangulation of the torus.\footnote{Observe that this is motivated by \cite{BrownErdSos73a}, and that we formulated it in its simplest case, however, we (more precisely, Nati Linial) meant a whole family of problems. He spoke about them in his talk in the Lovász Conference, 2018.}
\end{problem}

\begin{remark}
In some sense this geometric/linear algebraic approach helped Lovász to solve the famous Shannon conjecture on the capacity of $C_5$, in \cite{Lov79Shannon}.
\end{remark}

\Section{Methods: Randomness and the Semi-random method}{Methods}

Of course, the title of this chapter may seem too pretentious. We shall skip several important methods and related results here, or just touch on
them, primarily since they are well described in several other places, and partly since the bounded length of this survey does not allow us to
go into details.

Thus we shall skip most of the results {\em directly} connected to random graphs or random graph constructions, while we shall touch on the pseudo-random graphs. For random methods, the readers are referred to the books of Alon and
Spencer \cite{AlonSpencer16}, Bollobás \cite{Bollob01RandomBook}, or that of Janson, \L{u}czak, and Ruciński \cite{JansonLucRucRand}, the survey of Molloy \cite{Molloy98ProbMeth}, or the book of Molloy and Reed \cite{MolloyReed02Book}.

Here of course most of these sources are well known, like e.g., the books \cite{AlonSpencer16} or \cite{JansonLucRucRand}; we mention only Molloy's excellent chapter \cite{Molloy98ProbMeth}, which is a survey on this topic, and perhaps got less attention than it deserves. It contains many important details and explanations, and perhaps would fit the best to our topic, with the exception that it concentrates more on colouring and we concentrate more on independent sets in particular graphs.\footnote{The same applies to the book of Molloy and Reed \cite{MolloyReed02Book}.}

Listing the methods left out here, we should mention the extremely important works on Hypergraph Regularity Lemma, primarily that of Frankl and Rödl \cite{FranklRodl92Unif} on 3-uniform hypergraphs, (This was among the first ones). Of
course, we should mention the whole school of Rödl, among others, e.g., the papers of Rödl, Nagle, Skokan, Schacht, and Kohayakawa, \cite{RodlNagleSkokSchKoha05PNAS,KohaNagRodlSchacht10Linear}\footnote{This is a PNAS ``survey'', with an
  accompanying paper of Solymosi \cite{Solymo05PNAS}.}  Rödl and Schacht \cite{RodlSchacht07RegPartHypRL}, Rödl and Skokan \cite{RodlSkokan04Regu-k,RodlSkokan06RSA}, and Kohayakawa, and many further results, and, on the other hand, related
to this, several works of Tim Gowers \cite{GowersT06-3Unif,GowersT07HyperReg}, Green and Tao \cite{GreenTao10Arith} Terry Tao \cite{Tao06Revis} \dots

\subsection{Various ways to use randomness in Extremal Graph Theory}\label{VariousS}

 Random graphs are used in the area discussed here in several different ways.

\begin{enumerate}\dense

\item There are many cases where constructions seem unavailable but random graphs may be used to replace them. We have mentioned the Erdős Magic \cite{Erd47Random}. Another similar direct approach of Erdős and Rényi
  \cite{ErdRenyi60Evol} was that \beq{RandomLowerF} \ext(n,L)\ge c_L n^{2-v(L)/e(L)}.\eeq

\item In other cases we use modified random graphs: \eqref{RandomLowerF} is useless for cycles, but the modified random graph (in the simplest case the
  First Moment Method) worked in the papers of Erdős \cite{Erd59GraphProb} and \cite{Erd61GraphProbII}, and in many similar cases. We mention here the
  Erdős-Spencer book \cite{ErdSpencBook}, and its follower, the very popular book of Alon and Spencer \cite{AlonSpencer16}.

\item Since the very important paper of Erdős and Rényi \cite{ErdRenyi60Evol} -- on the evolution of Random Graphs, -- the investigation of the changes/phase-transitions in the structure of Random Graphs, as the number of edges is slowly
  increasing -- became a central topic of Combinatorics. Probably the first profound book on this was that of Bollobás \cite{Bollob01RandomBook}. Also we should definitely mention here the newer book of Janson, Łuczak and Ruciński
  \cite{JansonLucRucRand}, and the Molloy-Reed book \cite{MolloyReed02Book}.

\item Extremal graph theory optimizes on a Universe, and this Universe may be the family of Random Graphs.\footnote{Again, this case differs from the others: if we try to optimize some parameter on all $n$-vertex graphs, or on the subgraphs if the $d$-dimensional cube,\dots, that problem is well defined for individual graphs, while the assertions on the subgraphs of a random graphs make sense only in some asymptotic sense, the assertions always contain the expression ``almost surely as $n\to\infty$''.} Since a question of Erdős was answered by Babai, Simonovits, and Spencer \cite{BabaiSimSpenc90}, (here Theorem \ref{BabaiSimSpencTh}) this also became a very interesting and popular topic. We shall return to this question in Section~\ref{RandomUnivers}.

\item The Semi-random method was introduced by Ajtai, Komlós, and Szemerédi, for graphs in \cite{AjtKomSzem81Sidon} (to be applied in Combinatorial Number Theory, to Sidon Sequences), and later Rödl, in his famous solution of Erdős-Hanani
  problem \cite{ErdHanani63} about block designs developed the absolutely important method, now called Rödl Nibble \cite{Rodl85Nib}.

\end{enumerate}

In this short part we describe the Semi-Random method and the Rödl Nibble very superficially. According to their importance we should provide here a longer description, but the Rödl Nibble has an excellent description, a whole chapter,
in the Alon-Spencer book \cite{AlonSpencer16}, and it is more complicated and technical than that we could easily provide a sufficiently detailed description of it. (Beside referring to the Rödl Nibble, described in \cite{AlonSpencer16}, we
also mention the original Komlós-Pintz-Szemerédi paper \cite{KomlPintzSzem82HeilLower}, to see the origins, and also the Pippenger-Spencer paper \cite{PippeSpenc89}, and the Jeff Kahn paper \cite{Kahn91Colo} proving the asymptotic weakening
of the very famous Erdős-Faber-Lovász conjecture. (The Jeff Kahn paper also contains a fairly detailed description of the method.) The more recent paper of Alon, Kim, and Spencer \cite{AlonKimSpenc97} is also related to the previous topic.

\subsection{The semi-random method}\label{Nibble}

The semi-random method was introduced for graphs by Ajtai, Komlós, and Szemerédi \cite{AjtKomSzem81Sidon}, in connection with the infinite Sidon sequences. Later it was extended by Komlós, Pintz, and Szemerédi to 3-uniform hypergraphs in \cite{KomlPintzSzem82HeilLower}, to disprove a famous conjecture of Heilbronn, discussed in Subsections \ref{HeilbronnS}-\ref{HeilUpperS}. The method was further extended by Ajtai, Komlós, Pintz, Spencer, and Szemerédi in \cite{AjtKomPintzSpencSzem82} and by Duke, Lefmann, and Rödl~\cite{DukeLefRodl95Uncrowded}.

Here we are concerned with three topics, strongly connected to each other and to the estimates of the independence number of $K_3$-free
graphs, or analogous hypergraph results. The semi-random method had the form where the independence number of a graph or a hypergraph had to be estimated from
below, under the condition that the graph had no short cycles. The topics are

\begin{enumerate}\dense
\item {\bf Sidon's problem on infinite sequences.} A Sidon sequence is a sequence of integers in which no two (distinct pairs) have the same sum. What is the maximum density of such a sequence?
\item {\bf Heilbronn problem} for the ``minimum areas'' in geometric situations.
\item {\bf Ramsey problem} $R(3,k)$. This will be obtained as a byproduct, for free.
\end{enumerate}

\subsection{Independent sets in uncrowded graphs}\label{AKS-Indep}

\begin{figure}[h]
\begin{center}
\begin{tabular}{ccc}
\abrax{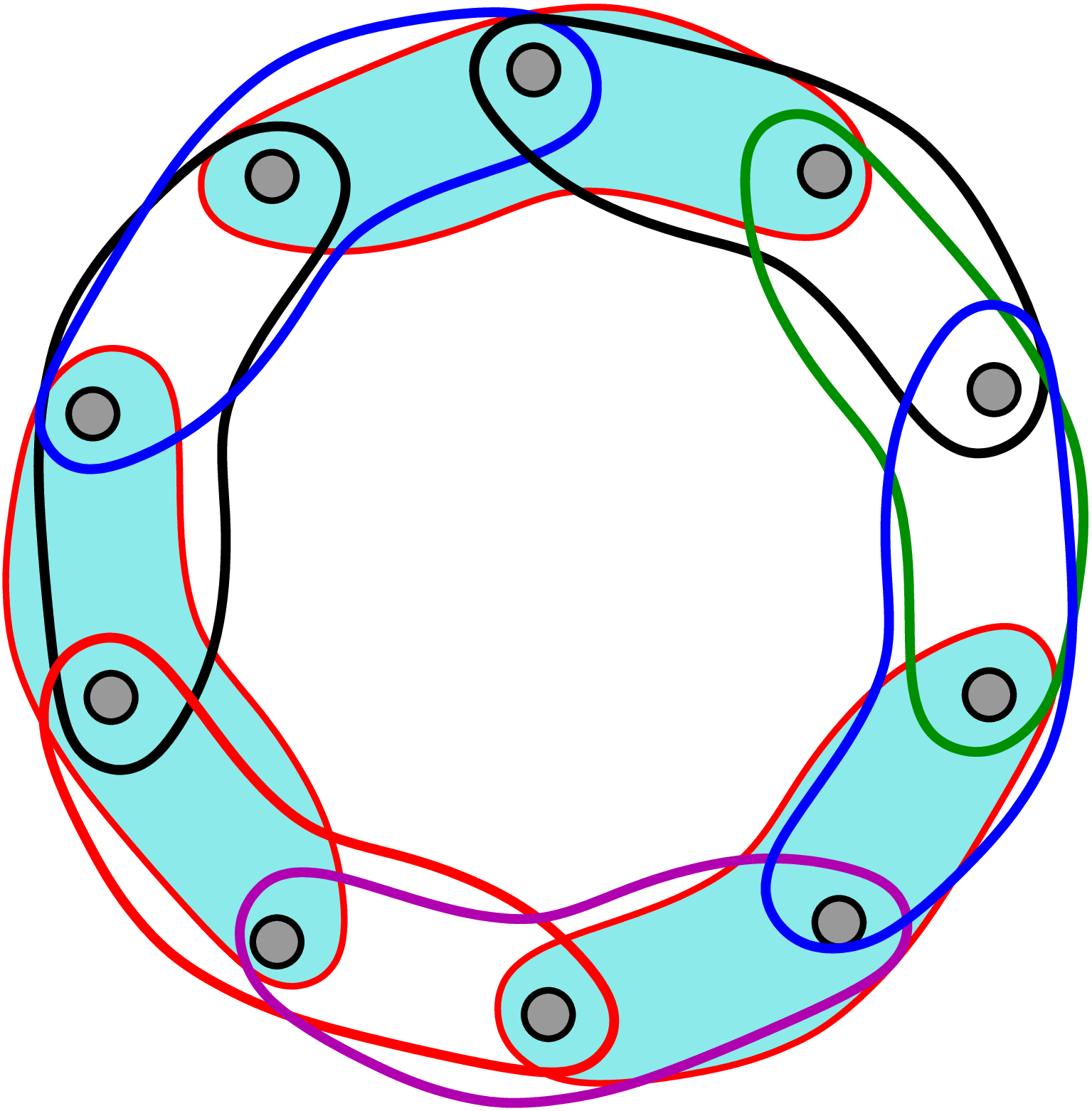}{22}&
\abrax{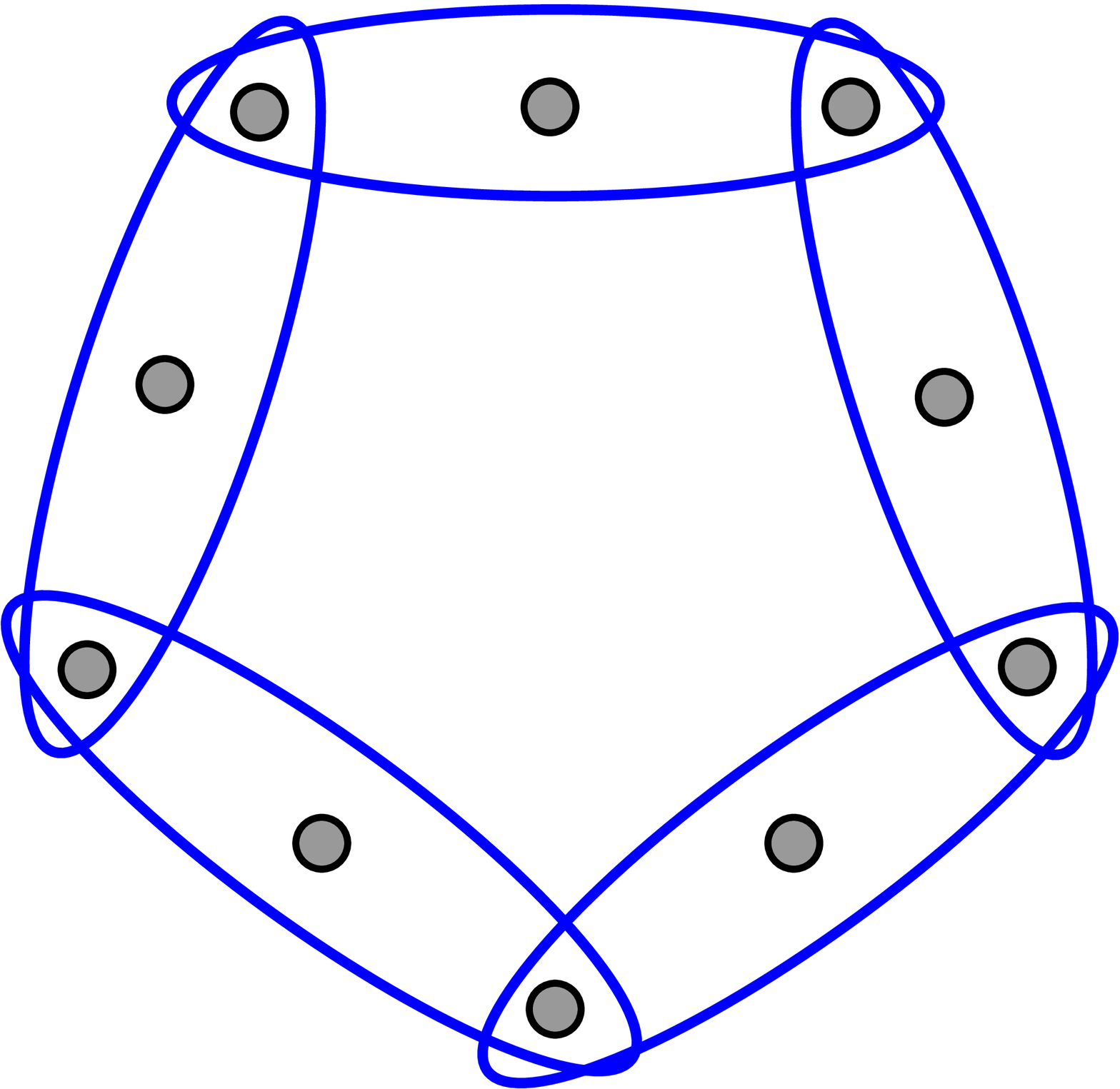}{22}&
\abrax{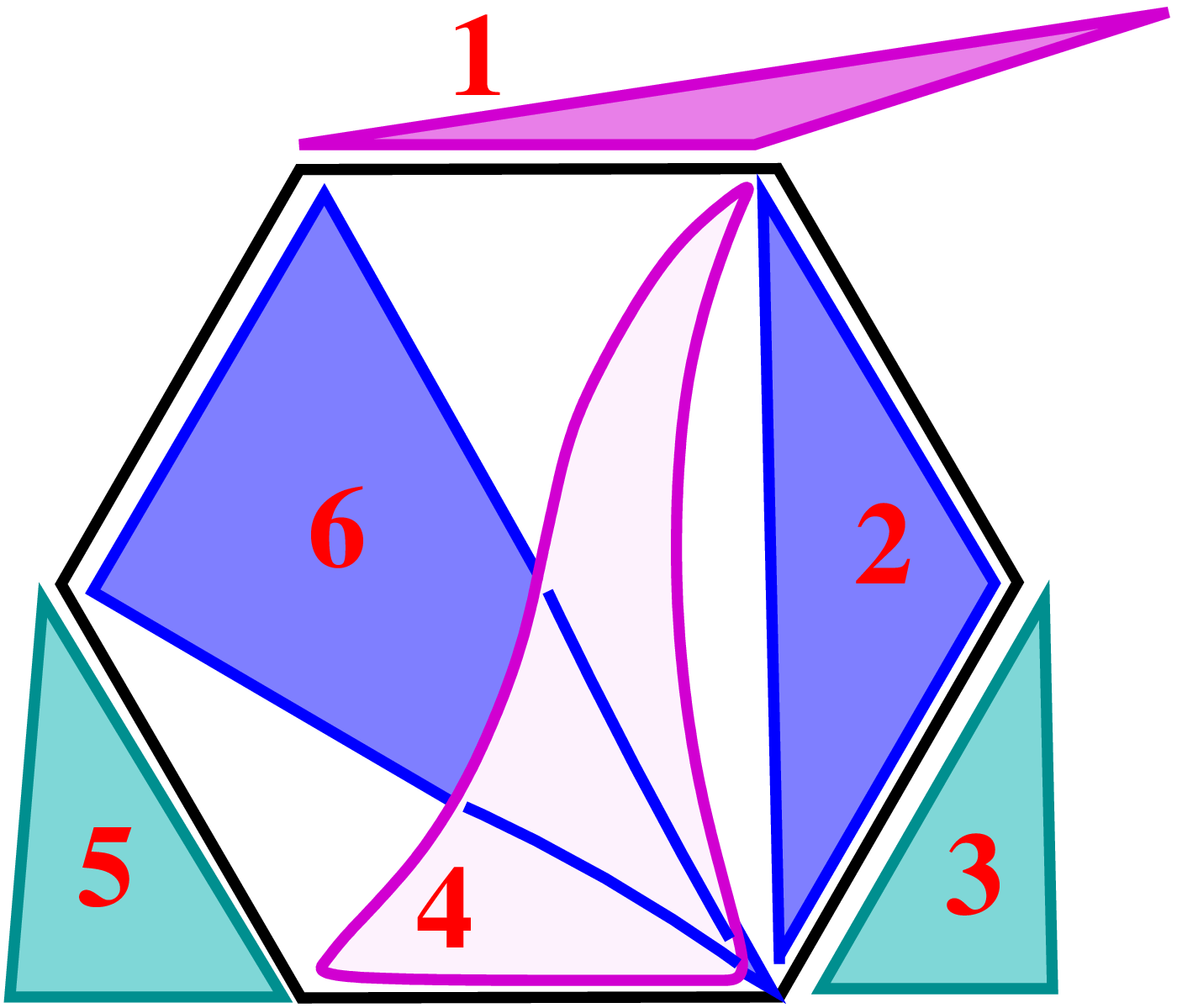}{22}\\
(a) Tight cycle& 
(b) Loose cycle&
(c) Berge cycle 
\end{tabular}
\end{center}
\caption{Various hypergraph cycles}
\label{CyclesV}
\end{figure}

A graph, hypergraph $\cH$ is called \Sc Uncrowded if it does not contain short cycles. For graphs we excluded triangles, for hypergraphs cycles of length 2, 3, or 4. In a hypergraph a cycle can be defined in several ways. Here a $k$-cycle ($k\geq 2$)\footnote{often called a Berge $k$-cycle: in Figure \ref{CyclesV}(b) the edges of a $C_6$ are covered by 3-tuples.} is a sequence of $k$ different vertices: $x_1,\ldots,x_{k-1},x_k=x_0$, and a sequence of $k$ different edges: $E_1,\ldots,E_k$ such that $x_{i-1},x_i\in E_i$ for $i=1,2\dots,k$. The cycle above is called \Sc simple if $E_i\cap\left(\cup_{j\not= i}E_j\right)=\{x_{i-1},x_i\}$ for $i=1,2\ldots,k$. In a hypergraph $\cH$, a 2-cycle is a pair of two hyperedges intersecting in at least two vertices; a vertex set $A\subset\cH$ is \Sc independent if it does not contain hyperedges. The maximum size of an independent set in $\cH$ is denoted by $\alpha(\cH)$. There are several lower bounds concerning independent sets in $k$-uniform uncrowded hypergraphs, mostly having the following form:

\begin{quote}{
``Hypergraphs having no short cycles have large independent sets.''}
\end{quote}

For ordinary graphs the following theorem, connected to infinite Sidon sequences, was the starting point.

\begin{theorem}[Ajtai, Komlós, Szemerédi \ev(1980) \cite{AjtKomSzem80Ramsey,AjtKomSzem81Sidon}]\label{AKSSidonTh} If in a triangle-free graph $G_n$ the average degree is $t:={2e(G_n)\over n}$, then the independence number
  \beq{AKS80IndepF} \alpha(G_n)\ge \reci100 \cdot{n\over t}\log t.\eeq
\end{theorem}

Turán's theorem, or the greedy algorithm would give only $n/t$ and in the random graphs we have the extra $\log t$ factor. The meaning of this theorem is, perhaps, that excluding $\K3$ forces a much larger independent set and a randomlike behaviour.

Ajtai, Erdős, Komlós, and Szemerédi also proved

\begin{theorem}[\cite{AjtErdKomSzem81Sparse}] If $G_n$ is $K_p$-free then $\alpha(G_n)\ge c'(n/t)\ln A$ where $A=(\ln t)/p$ and $c'$ is an absolute constant. 
\end{theorem}

See also Shearer \cite{Shearer83Indep,Shearer91Indep,Shearer95IndepE,Shearer95SparseInde} and Denley \cite{Denley94Indep}.

\subsubsection{Proof-Sketch}

Originally Theorem \ref{AKSSidonTh} had two different proofs. One of them was an induction on the number of vertices, \cite{AjtKomSzem80Ramsey}, (and a similar, perhaps more elegant proof -- also using induction on $n$ -- was given by
Shearer for its sharpening, see \cite{Shearer83Indep}). The other proof, from \cite{AjtKomSzem81Sidon} used an iterated random construction which later developed into the Rödl Nibble.
This approach turned out to be fairly important, so we ``sketch'' its main idea, suppressing most of the technical details, and following the description from p10 of \cite{AjtKomSzem81Sidon}.

\begin{enumerate}\dense

\item Since $G_n$ is triangle-free, $\alpha(G_n)>\maxdeg(G_n)\ge t$. So
we may assume that the degrees are smaller than $\reci100 \cdot{n\over t}\log t$. Similarly, we may assume that $t\ge \sqrt{n\log n}$, since otherwise 
$t>\reci100 {\log t\over t}n$, implying \eqref{AKS80IndepF}:\ $F(t)=t^{-2}\log t$ is monotone decreasing; for $t=\sqrt{n\log n}$ we have $F(t)\approx \reci 2n .$  This proves the assertion.
\item We select a subset $\cK\subset V(G_n)$ of $\reci110 (n/t)$ independent vertices.
\item Next we consider a vertex-set $\cM\subseteq V(G_n)-\cK$ of $\approx n/2$ vertices not joined to $\cK$: we need a lemma
about the existence of such an $\cM$.
\item Another lemma assures us that the crucial quantity $n/t$ does not drop too much when we move from $G_n$ to $G_m=G[\cM]$.
(It drops only by a factor $\teta:=1-\reci t -c_{10}\sqrt{t/n}> 1-t^{1/3}$.)
\item If we are lucky, then we can iterate this step $\approx \half\log t$ times:
we gain a $\half\log t$ factor and get Theorem \ref{AKSSidonTh}.
\item On the other hand, if we are ``unlucky'' and get stuck in the $r\th$ step, then for the corresponding $t_r$ we have that it is too large: $t_r^{-1/3}>1/\log t$. But then we get in this last step alone enough independent
  vertices:
$$ {n_r\over {t_r+1}}>(\log t)^{-3}{n\over 2^r}>{n(\log t)^{-3}\over\sqrt t}> n{\log t \over t}.$$
\end{enumerate}

\subsub Summarizing:// In the typical case we can choose small independent subsets in $V(G_n)$ basically $\log t$ times to gain a $\log t$ factor.  It is important that discarding these small vertex sets, we do not ruin the structure of the remaining part.

\subsection{Uncrowded hypergraphs}

Most probably, the earliest result on hypergraphs using the semi-random method was the following one.

\begin{theorem}[Komlós, Pintz, and Szemerédi \cite{KomlPintzSzem82HeilLower}, Lemma 1]\label{KPSzHypergTh} There exists a threshold $t_0$ and a constant $c>0$ such that if $\Hhyp n3$ is a $3$-uniform hypergraph on $n$ vertices, with average degree $t^{(3)}(\Hhyp n3)$, and not containing simple cycles of length at most $4$, then for $t:=\sqrt{t^{(3)}(\Hhyp n3)}>t_0$, $t=O(n^{1/10})$, we have
$$\alpha(\Hhyp n3)> c \frac{n}{t}\sqrt{\log t}.$$
\end{theorem}

The problems we discuss here were reduced to finding lower bounds on the independence number $\alpha(\bH)$ of a graph or hypergraph $\bH$ under the assumption that $\bH$ has no short cycles.\footnote{Actually, in\cite{KomlPintzSzem82HeilLower} one needs to exclude only cycles of length 2,3, and 4, where a cycle of length 2 is a pair of hyperedges intersecting in at least two vertices. Even this is improved in the next theorem.}  The above theorem
and its versions were enough for the early applications, in 1980's, however, as it turned out in \cite{DukeLefRodl95Uncrowded}, only hypergraph cycles of length 2 had to be excluded: the following much newer generalization can be very useful in some new applications.

\begin{theorem}[Duke, Lefmann, and Rödl \cite{DukeLefRodl95Uncrowded},
 Theorem 2]\label{DLR95Uncrow}
 Let $\bH$ be a $k$-uniform hypergraph on $n$ vertices. Let $\Delta$ be the maximum degree of $\bH$. Assume that $\Delta\le t^{k-1}$
 and $ t>t_0$. If $\bH$ doesn't contain $2$-cycles (two hyperedges with at least two common vertices), then
$$ 
\alpha(\cH)=\Omega\left(\frac{n}{t}(\log t)^{\frac{1}{k-1}}\right).
$$ 
\end{theorem}

Theorem \ref{DLR95Uncrow} for $k=3$ implies Theorem~\ref{KPSzHypergTh}. One advantage of it is that in our applications we may have many simple cycles of
length $3$ and $4$, but Theorem \ref{DLR95Uncrow} still can be applied.

\minisepa 

There are many results in this field. We mention here only (i) Duke, Lefmann, and Rödl \cite{DukeLefRodl95Uncrowded} on Uncrowded Hypergraphs, (ii) C. Bertram-Kretzberg, T. Hofmeister and H. Lefmann \cite{BertHofLef00Heil,BertLef99Uncrowd}, some generalizations and results on the algorithmic aspects of the Heilbronn Problem: finding efficiently the large independent set in an uncrowded hypergraph, and (iii) Lefmann and Schmitt \cite{LefmSchmitt02} and Lefmann \cite{Lefm03Heilbronn}, on the higher dimensional Heilbronn problem.

In \cite{AjtKomSzem81Sidon} it is remarked that it is enough to assume that the number of triangles is small, instead of assuming that there are no $K_3$'s at all. Shearer \cite{Shearer83Indep,Shearer91Indep} improved the constant in
Theorem \ref{AKSSidonTh} in an ingenious way:

\begin{theorem}[Shearer \cite{Shearer83Indep}] 
Let $f(t)={{t\log t -t+1}\over {(t-1)^2}}.$ Then for any triangle-free $G_n$,
with average degree $t$,
\beq{ShearFo} \alpha(G_n)\ge f(t)\cdot n.\eeq
\end{theorem}

This improves \eqref{AKS80IndepF}, and (in some sense) this is sharp.  Related extensions were found by T.~Denley \cite{Denley94Indep}, and by Shearer in cases when we assume that the odd girth of the graph is large \cite{Shearer95IndepE}, and also when we wish to use finer information on the degree distribution. Kostochka, Mubayi, and Verstraëte \cite{KostoMubVerstr13Kim} proved some hypergraph versions of this theorem, giving lower bounds on the stability number under the condition that certain cycles are excluded.

\subsection{Ramsey estimates}

Observe that -- as a byproduct, -- Theorem \ref{AKSSidonTh} immediately yields that

\begin{theorem}[Ajtai, Komlós, Szemerédi, on Ramsey theorem \cite{AjtKomSzem80Ramsey}]\label{AKS82RamsT} There exists a constant $c>0$ such that
\beq{AKSRamsFo} R(3,m)\le {m^2\over c\log m}.\eeq
\end{theorem}

\proof. Indeed, if $n>{m^2\over c\log m}$, and $K_3\not\subseteq G_n$, then either $G_n$ has a vertex $x$ of degree $\de_G(x)\ge m$, yielding an independent $m$-set $N_G(x)$, or by Theorem~\ref{AKSSidonTh}, 
$$\alpha(G_n)> c{{m^2\over c\log m}\over m}\log m=m $$ 
proving Theorem \ref{AKS82RamsT}.\qed

Theorem \ref{AKS82RamsT} improves Erdős' old result \cite{Erd59GraphProb}, by a $\log m$ factor.\footnote{Of course, Shearer's improvement yields an improvement of $c$ in \eqref{AKSRamsFo}.} For many years it was open if \eqref{AKSRamsFo} is just an improvement of the Erdős result or it is a breakthrough. Jeong Han Kim \cite{Kim95RamseyK3} proved much later, (using among others the Rödl Nibble) that this bound is sharp.

\begin{theorem}[J.H. Kim, Ramsey \ev(1995) \cite{Kim95RamseyK3}]\footnote{Actually, the proof works with $\ti c=\reci162 -o(1) $.} 
 \beq{KimRamseyFo} R(3,m)\ge \ti c {m^2\over\log m}.\eeq
\end{theorem}

So the $r(3,m)$-problem is one of the very few nontrivial infinite cases on Ramsey numbers where the order of magnitude is determined. Bohman and Keevash \cite{BohmanKeev13} and Fiz Pontiveros, Griffiths, and Morris
\cite{FizPontiGriffMorr13Rams} independently proved that $R(3,k)\ge(\reci4 +o(1)) {k^2\over{\log k}}$. Recently Shearer's estimate was ``strengthened'' as follows.

\begin{theorem}[Davies, Jenssen, Perkins, and Roberts \cite{DaviesJenssPerkRob}] If $G_n$ is a triangle-free graph with maximum degree $t$, and $\cI(G_n)$ is the family of independent sets in $G_n$,
then $${1\over {|\cI(G_n)|}}\sum_{I\in\cI(G_n)} |I| \ge (1+o(1)){\log t\over t}n.$$
\end{theorem}

This is a strengthening in the sense that here the average size is large, but it is a weakening: it uses the maximum degree, instead of the average degree.
For further details see, e.g., the introduction of \cite{DaviesJenssPerkRob}.

\begin{remark}
The above questions are connected to another important question: under some condition, what can be said about the number of independent sets in a graph or a hypergraph? Without going into details, we remark that these questions are connected to the container method, (for references see the Introduction).
\end{remark}

\begin{remark} Further related results can be found, e.g., in papers of Cooper and Mubayi, and Dutta \cite{CoopMubayi14Indep,CoopMubayi17Sparse,CoopDutMubayi14Indep} which count the number of independent sets in triangle-free graphs and
  hypergraphs.
\end{remark}

\subsection{Infinite Sidon sequences}\label{InfiSidon}

The finite Sidon sequences are well understood, the maximum size of a Sidon subset of $[1,n]$ is around $\sqrt{n}$, \cite{ErdTur41Sidon,Chowla44Sidon}, however, the infinite Sidon sequences seem much more involved. The greedy algorithm
provides an infinite Sidon sequence $(a_n)$ with $a_n>cn^{1/3}$. This was slightly improved by using Theorem~\ref{AKSSidonTh}, but only by $\sqrt[3]{\log n}$:

\begin{theorem}[Ajtai, Komlós, and Szemerédi \ev(1981) \cite{AjtKomSzem81Sidon}]\label{AKS-Sidon}
There exists an infinite Sidon sequence $B$ for which, if $B(n)$ denotes the number of elements of it in $[1,n]$, then
$$B(n)\ge c(n\log n)^{1/3}.$$
\end{theorem}

As it is remarked in \cite{AjtKomSzem81Sidon}, Erdős conjectured that $B(n)>n^{(1/2)-\eps}$ is possible.
As to Sidon sets, later Theorem \ref{AKS-Sidon} was improved ``dramatically'':

\begin{theorem}[Ruzsa \ev(1998) \cite{Ruzsa98Sidon}]
There exists an infinite Sidon sequence $B$ for which, if $B(n)$ denotes the number of elements in $[1,n]$, then
$$B(n)\ge n^{\sqrt2-1+o(1)}.$$
\end{theorem}

So the importance of this Ajtai-Komlós-Szemerédi result \cite{AjtKomSzem81Sidon} was that this was the beginning of the Semi-Random method.  The following generalization was proved by Cilleruelo \cite{Cilleruel14Sidon}. Call a sequence
$\cA=\{a_{i} \}_{i=1}^{\infty}$ $h$-Sidon if all the sums $a_{i_1}+ \cdots + a_{i_h}$ are distinct for $a_{i_1} \le \cdots \le a_{i_h}$. 
 
\begin{theorem}[Cilleruelo, \ev(2014) \cite{Cilleruel14Sidon}] For any $h \geq 2$ there exists an infinite $h$-Sidon sequence $\cA$ with $$ A(n) = n^{\sqrt{(h-1)^2+1}-(h-1) +o(1)}, $$ where $A(n)$ counts the number of elements of $\cA$ not exceeding $n$.
\end{theorem}

For $h=2$ Cilleruelo provides an explicit construction. We remark also that, by a ``random construction'' of Erdős and Rényi \cite{ErdRenyi60Additive}, for any $\de>0$, there exits an infinite sequence $Q:=(a_1,\dots,a_n,\dots)$ for which
the number of solutions of $a_i+a_j=h$ is bounded, for all $h$, and $a_k=O(k^{2+\de})$.

\subsection{The Heilbronn problem, old results}\label{HeilbronnS}

\begin{problem}[Heilbronn's problem on the area of small triangles] Consider $n$ points in the unit square (or in the unit disk) no three of which are collinear
What is the maximum of the minimum area of triangles, defined by these points where the maximum is taken for all point-sets?\footnote{If three of them are collinear that provides 0.}. 
\end{problem}

This maximum of the minimum will be denoted by $\Delta_n$. Erdős gave a simple construction where this minimum area was at least $\reci2n^2 $: for a prime
$p\approx n$, consider all the points $(\reci p (i,f(i))$, where $f(i)=i^2 ~(\mbox{mod}~p)$. So
$$\reci 2n^2 <\Delta_n \le\reci n-1 . $$
Heilbronn conjectured that $\Delta=O(\reci n^2 )$. This was disproved by

\begin{theorem}[Komlós, Pintz, and Szemerédi \ev(1981) \cite{KomlPintzSzem81Heil}]\label{AKSHeilLowerTh} $\Delta_n=\Omega({\log n\over n^2})$: For some constant $c>0$, for infinitely many $n$, there exist $n$ points in the unit square
 for which the minimum area is at least $c{\log n\over n^2}.$
\end{theorem} 

The proof of Theorem \ref{AKSHeilLowerTh} is based on Theorem \ref{KPSzHypergTh}.

\subsection{Generalizations of Heilbronn's problem, new results}\label{HeilbronnNewS}

Péter Hajnal and Szemerédi \cite{HajnalPSzem18Geom} used the Duke-Lefmann-Rödl lower bound (Theorem \ref{DLR95Uncrow}) to prove two new geometrical results. The first one \cite{HajnalPSzem18Geom} is closely related to Heilbronn's triangle problem, discussed in \cite{Roth51Heil,Schmidt72Heilb,Roth72Heil,Roth73Heilb,Roth76Heil,KomlPintzSzem81Heil}.

Consider an $n$-element point set $\PP\subset\EE^2$. Instead of triangles we can take $k$-tuples from $\PP$ and consider the area of the
convex hull of the $k$ chosen points. Denote the minimum area by $H_k(\PP)$, its maximum for the $n$-element sets $\PP$ by $H_k(n)$. So
$\Delta_n=H_3(n)$. The best lower bound on $H_3(n)$ from \cite{KomlPintzSzem82HeilLower},
and some trivial observations are summarized in the next line: there exists a constant $c>0$ such that 
$$
c\frac{\sqrt{\log n}}{n^2}\le \Delta_n=H_3 (n)\le H_4(n)\le H_5(n)\le \ldots= O\left(\frac{1}{n}\right).
$$

We mention two major open problems: 

\begin{problem}
Is it true that for any $\eps>0$, $H_3 (n)=O(1/n^{2-\eps})$?\dori
Is it true that $H_4(n)=o(1/n)$?
\end{problem}

One is also interested in finding good lower bounds on $H_4(n)$. Schmidt \cite{Schmidt72Heilb} proved that $H_4(n)=\Omega(n^{-3/2})$. The proof is a
construction of a point set $\PP$ by a simple greedy algorithm. In \cite{BertHofLef00Heil}, Bertram-Kretzberg, Hofmeister and Lefmann provided a new proof,
and extensions of this result. They also asked whether Schmidt's bound can be improved by a logarithmic factor. Using the semi-random method, Péter Hajnal
and Szemerédi improved Schmidt's bound and settled the problem of \cite{BertHofLef00Heil}.

\begin{theorem}[P. Hajnal and E. Szemerédi \cite{HajnalPSzem18Geom}] \label{heilbronn} For some appropriate constant $c>0$, for any $n>3$, there exist $n$ points in the unit square for which the convex hull of any 4 points has area at least $cn^{-3/2}(\log n)^{1/2}$.
\end{theorem}

\subsection{The Heilbronn problem, an upper bound}\label{HeilUpperS}

The first upper bound was Roth's fundamental result that $\Delta_n\ll\reci n\sqrt{\loglog~n} $.\footnote{Here $f\ll g$ is the same as $f\le cg$, for some
  absolute constant $c>0$.} Schmidt \cite{Schmidt72Heilb} improved this to $\De_n\ll \reci n\sqrt{\log~n} $. Roth returned to the problem and improved the
earlier results to $\Delta_n\ll 1/n^{1.117}$. Not only his bound was much better, his method was also groundbreaking. He combined methods from analysis,
geometry, and functional analysis. On these results see the survey of Roth \cite{Roth76Heil}. Roth's result was improved by Pintz, Komlós, and Szemerédi:

\begin{theorem}[Pintz, Komlós, Szemerédi \cite{KomlPintzSzem81Heil}] $\Delta_n\ll e^{c\sqrt{\log n}} n^{-8/7}.$
\end{theorem}

\subsection{The Gowers problem}\label{GowersCollin}

P. Hajnal and E. Szemerédi considered the following related but ``projective'' question:\footnote{Actually, Hajnal and Szemerédi found this problem on Gowers' homepage, but, as it turned out, from \cite{PayneWood13}, originally the problem was asked by Paul Erdős, \cite{Erd88-32}.}

\begin{problem}[Gowers \cite{GowersT16GeomQue}] Given a planar point 
 set $\PP$, what is the minimum size of $\PP$ that guarantees that one can find $n$ collinear points (points on a line) or $n$ independent
 points (no three on a line) in $\PP$?
\end{problem}

Gowers noted that in the grid at least $\Omega(n^2)$ points are needed, and if we have $2n^3$ points without $n$ points on a line, then a simple greedy algorithm finds $n$ independent points. Payne and Wood \cite{PayneWood13} improved the upper bound $O(n^3)$ to $O(n^2\log n)$. They considered an arbitrary point set with much fewer than $n^3$
points, and without $n$ points on a line, but they replaced the greedy algorithm by a Spencer lemma, based on a simple probabilistic sparsification.\footnote{Actually, above we spoke about the ``diagonal case'' but \cite{PayneWood13} covers some off-diagonal cases too.}

Péter Hajnal and Szemerédi improved the Payne-Wood upper bounds by improving the methods. They also started with a random sparsification, but after some additional preparation (getting rid of $2$-cycles) they could use the semi-random method (see \cite{DukeLefRodl95Uncrowded}) to find a large independent set.

\begin{theorem}[P. Hajnal and E. Szemerédi \cite{HajnalPSzem18Geom}]\label{HajnalSzemer18Gowers} There exists a constant $C>0$ such that in any planar point set $\PP$ of size $C\cdot\frac{n^2\log n}{\log\log n}$, there are $n$ points that are collinear or independent.
\end{theorem}

\subsection{Pippenger-Spencer theorem}\label{PippeSpencS}

In the theory developing around the semi-random methods, one should mention the papers of Rödl \cite{Rodl85Nib}, and of Frankl and Rödl \cite{FranklRodl85Perfect}, \dots

One important step was the Pippenger-Spencer result \cite{PippeSpenc89}, asserting that if the degrees in a $k$-uniform hypergraph are large and the \Sc codegrees are relatively small, then the hypergraph has an almost-1-factor.

\begin{definition}
A \Sc matching of a hypergraph $\cH$ is a collection of pairwise disjoint hyperedges. The \Sc chromatic~index $\chi'(\cH)$ of $\cH$ is the least number of matchings whose union covers the edge set of $\cH$. 
The \Sc codegree $\de_\ell(X)$ is the number of hyperedges containing the $\ell$-tuple $X\subseteq V(\cH)$, and $\de_\ell(\cH)$ is the minimum of this, taken over all the $\ell$-tuples of vertices in $\cH$.
\end{definition}

We formulate the theorem, in a slightly simplified form.

\begin{theorem}[Pippenger and Spencer \cite{PippeSpenc89}]
For every $k\ge2$ and $\de>0$ there exists a $\de'>0$ and an $n_0$ such that if $\cH_n$ is a $k$-uniform hypergraph with $n>n_0$ vertices, and
$$\mindeg(\cH_n)>(1-\de)\maxdeg(\cH_n),$$
and the codegrees are small:
$$\codeg(\cH_n)<\de'\mindeg(\cH_n)$$
then the chromatic index $$\chi'(\cH_n)\le(1+\de)\maxdeg(\cH_n) .$$
\end{theorem}

The meaning of the conclusion is that the set of hyperedges can be partitioned into packings (or matchings), almost all of which are almost perfect. Also the edges can be partitioned into coverings, almost all of which are almost
perfect. This theorem strengthens and generalizes a result of Frankl and Rödl \cite{FranklRodl85Perfect}.

\subsection{Erdős-Faber-Lovász conjecture}\label{ErdFabLovS}

We close this part with the beautiful result of Jeff Kahn on the famous Erdős-Faber-Lovász conjecture. Faber writes in \cite{Faber10ErdFabLov}:\footnote{In citations we use our numbers, not the original ones.}

\begin{quote}{ In 1972, Paul Erdős, László Lovász and I got together at a tea party in my apartment in Boulder, Colorado. This was a continuation of the discussions we had had a few weeks before in Columbus, Ohio, at a conference on hypergraphs. We talked about various conjectures for linear hypergraphs analogous to Vizing’s theorem for graphs (see \cite{FaberLov72EFL}). Finding tight bounds in general seemed difficult so we created an elementary conjecture that we thought would be easy to prove. We called this the $n$ sets problem: given $n$ sets, no two of which meet more than once and each with $n$ elements, color the elements with $n$ colors so that each set contains all the colors. In fact, we agreed to meet the next day to write down the solution. Thirty-eight years later, this problem is still unsolved in general. (See \cite{Romero07} for a survey of what is known.)}
\end{quote}

The original conjecture says: 

\begin{conjecture}
If $G_n$ is the union of $k$ complete graphs $K_k$, any two of which has at most one common vertex, then $\chi(G_n)\le k$. 
\end{conjecture}

As we stated, originally the conjecture was formulated using Linear Hypergraphs.\footnote{Perhaps the expression ``linear hypergraph'' was unknown those days.} A weakened asymptotic form of this was proved by Jeff Kahn: 

\begin{theorem}[Jeff Kahn \ev(1991) \cite{Kahn91Colo}]
If $A_1,\dots,A_k\subseteq [n]$ are nearly disjoint, in the sense that $|A_i\cap A_j|\le 1$ for all pairs $i\ne j$, then $\chi'(H)\le(1+o(1))n$.
\end{theorem}

Jeff Kahn gave an elegant proof of this assertion, and also a clear description of the sketch of his proof, which is also a nice description of the Semi-Random method. The proof is based on a technical generalization of the
Pippenger-Spencer Theorem \cite{PippeSpenc89}.

\begin{remark}
(a) Originally the Erdős-Faber-Lovász conjecture had a slightly different form, see above, or, e.g., Erdős \cite{Erd76Aberdeen}. Jeff Kahn refers to Hindman \cite{Hindman81EFL} who rephrased it in this form. Kahn also remarks that Komlós suggested to prove an asymptotic weakening.

(b) We also remark that the fractional form of this problem was solved by Kahn and Seymour \cite{KahnSeymour92ErdFabLov}.

(c) Vance Faber proved \cite{Faber10ErdFabLov} that if there are some counter-examples to the Erdős-Faber-Lovász conjecture, they should be in some sense in the ``middle range'', according to their densities.
\end{remark}

\begin{remark}[On Keevash' existence-proof of Block designs]\label{KeevashDesign} One of the recent results that is considered a very important breakthrough is that of
  Peter Keevash \cite{Keev14ExistDesign,Keev18ExistDesign}, according to which, if we fix some parameters for some block-designs, and the corresponding trivial divisibility conditions are also satisfied, then the corresponding block designs do exist, assumed that the set of elements is sufficiently large. Important but simpler results in the field were obtained by Richard Wilson \cite{Wilson72DesignI,Wilson72DesignII,Wilson75DesignIII,Wilson75Decomp}, Vojta Rödl \cite{Rodl85Nib}, using -- among others -- the methods described above, above all, the Rödl Nibble. The excellent paper of Gil Kalai \cite{Kalai16Keevash} explains this area: what one tries to prove and how the semi-random methods are used. As a very important contribution, approach, see also the papers of Glock, Kühn, Lo, and Osthus \cite{GlockKuhnLoOsthus16A-DesignIter,GlockKuhnLoOsthus17A-FDesign}. We do not go into details but again, refer the interested readers to the paper of Kalai written for the general audience, or, at least, for most of the combinatorists.
\end{remark}

\Section{Methods: Regularity Lemma for graphs}{RegularityS}

As we have already stated, Regularity Lemma is applicable in many areas of Discrete Mathematics. A weaker, ``bipartite'' version was used in the proof of Szemerédi Theorem \cite{Szem75APk} according to which $r_k(n)=o(n)$. Also weaker versions were used in the first applications in Graph Theory \cite{RuzsaSzem78,Szem72K4}. The first case when its
standard form (Theorem \ref{SzemReguTh}) was needed was the Chvátal-Szemerédi theorem \cite{ChvatSzem81ErdStone} on the parametrized form of the Erdős-Stone theorem\footnote{The question was that if $e(G_n)=e(T_{n,p})+\eps n^2$, how large $K_{p+1}(m,\dots,m)$ can be guaranteed in $G_n$? This maximum $m=m(p,\eps)$ had a very weak estimate in \cite{ErdStone46}. This was improved to $c(p,\eps)\log n$ by Bollobás and Erdős \cite{BolloErd73}, which was improved by Bollobás, Erdős, and Simonovits \cite{BolloErdSim76}. Chvátal and Szemerédi needed the Regularity Lemma to get the ``final'' result, sharp up to a multiplicative absolute constant.}. The Regularity Lemma is so successful, perhaps because of the
following.

To embed a graph $H$ into $G$ is much easier if $G$ is a random graph than if it is an arbitrary graph. The Regularity Lemma asserts that every graph $G$
can be approximated by a ``generalized random graph'', more precisely, by a ``generalized quasi-random graph''. But then we may embed $H$ into an almost
random graph which is much easier. (Also, many other things are easier for Random Graphs.)

 Lovász and Szegedy \cite{LovSzege07Anal} wrote a beautiful paper on the Regularity Lemma, where they wrote\footnote{This was formulated by many researchers.}:

``Roughly speaking, the Szemerédi Regularity Lemma says that the node set of every (large) graph can be partitioned into a small number of parts so that the subgraphs between the parts are random-like. There are several ways to make this
 precise, some equivalent to the original version, some not.\dots''

To formulate the Regularity Lemma, we define (i) the $\eps$-regular pairs of vertex-sets in a graph, (ii) the generalized random graphs, and (iii) the generalized quasi-random graphs.

Given a graph $G$, with the disjoint vertex sets $X,Y\subseteq V(\G)$, the edge-density between $X$ and $Y$ is
 $$d(X,Y):={{e(X,Y)}\over {|X||Y|}}.$$
 Regular pairs are highly uniform bipartite graphs, namely ones in which the density of any ``large'' induced subgraph is about the same as the overall density of the whole graph.

\begin{definition}[Regular pairs]\label{SzemReguDef}
 Let $\eps>0$. Given a graph $G$ and two disjoint vertex sets $A\subseteq V(G)$, $B\subseteq V(G)$, we say that the pair $(A,B)$ is $\eps$-regular if for every $X\subseteq A$ and $Y\subseteq B$ satisfying
$$|X|>\eps|A|\Text{and}|Y|>\eps|B|$$
we have \beq{OriginalUnif}|d(X,Y)-d(A,B)|<\eps.~\footnote{In random graphs this holds for sufficiently large disjoint vertex sets.}
\eeq
\end{definition}

%% One way to describe the Regularity Lemma is to say that any graph can be approximated by the union of a bounded number of randomlike bipartite graphs. 
We can also describe the Regularity Lemma as a statement asserting that any graph $G_n$ can be approximated well by generalized random graphs. However, first we have to define the generalized Erdős-Rényi random graph, then the generalized quasi-random graph sequences.

\begin{definition}[Generalized Random Graphs] 
Given a matrix of probabilities, $A:=(p_{ij})_{r\times r}$ and integers $n_1,\dots,n_r$.
We choose the subsets $U_1,\dots,U_r$ with $|U_i|=n_i$ and join $x\in U_i$ to $y\in U_j$ with probability $p_{ij}$, {\em independently}.
\end{definition}

The \Sc generalized~quasi-random graphs are obtained when we also fix an $\eps>0$ and instead of joining $U_i$ to $U_j$ using random independent edges, we join them with $\eps$-regular bipartite graphs of the given density $p_{ij}$.

 Regularity Lemma asserts that the graphs can be approximated by generalized quasi-random graphs well.

\def\mun{M_0}

\begin{theorem}[Szemerédi, 1978 \cite{Szem78Regu}]\label{SzemReguTh} 
 For every $\eps>0$ and $\mun$ there are two constants $M(\eps,\mun)$ and $N(\eps,\mun)$ with the following property: for every graph $G_n$ with $n\geq N(\eps,\mun)$ vertices there is a partition of the vertex set into {$\nu$} classes $$V=V_1\cup V_2\cup \ldots\cup V_\nu$$
 such that
\begin{itemize}\dense
\item $\mun\le \nu\le M(\eps,\mun)$,
\item $||V_i|-|V_j||\le1$, ($1\le i<j\le \nu$)
\item all but at most $\eps \nu^2$ of the pairs $(V_i,V_j)$ are $\eps$-regular.
\end{itemize}
\end{theorem}

The lower bound on the number of classes is needed mostly to make the number of edges within the clusters small. If e.g. $M_0>1/\eps$, then $\sum e(V_i)<\half\eps n^2$. So mostly we can choose $M_0=1/\eps$. We do not really need that $||V_i|-|V_j||\le1$, however, we do need that the $V_i$'s are not too large.
In the applications we very often use the {\em Reduced graph} or {\em Cluster graph} $H_\nu$ defined as follows:\footnote{Perhaps the name ``Reduced Graph'' comes from Simonovits, the ``Cluster Graph'' from Komlós, and the theorem itself was originally called ``Uniformity Lemma'': the name ``Regularity Lemma'' became popular only later.}

\begin{definition}[Cluster Graph $H_\nu=H_\nu(G_n)$]
Given a graph $G_n$, and two constants $\eps,\tau>0$, the corresponding Cluster Graphs are obtained as follows. We apply the Regularity Lemma with $\eps$, obtaining the partition $V_1,\dots,V_\nu$. The vertices of $H_\nu$ are the vertex-sets $V_i$ (called Clusters) and we connect the Cluster $V_i$ to $V_j$ if $\ViVj$ is $\eps$-regular in $G_n$ and $d\ViVj>\tau$.
\end{definition}

\begin{remark}\label{SzemReguScheme} 
{\em Mostly} we use the cluster graph as follows: %% \danger 
(i) We set out with a graph sequence $(G_n)$, satisfying some conditions $\cP$, (ii) take the cluster graphs $H_\nu$, (iii) derive that $(H_\nu)$ must have some properties $\cP^*$ (because of the combinatorial conditions on $G_n$), (iv)
therefore we can apply some ``classical'' graph theorem to $H_\nu$, (v) this helps to describe $H_\nu$\footnote{estimate $e(H_\nu)$ or prove some structural property of $H_\nu$.} and (vi) having this information on $H_\nu$ helps us to prove
what we wanted.

\dots Or, in case (v-vi), instead of, say, estimating $e(G_n)$, we embed some given graph $U_\mu$ into $H_\nu$, and using this, we embed a (much larger) $W_m$ into $G_n$.
\end{remark}

\begin{remark}
Given a graph sequence $(G_n)$, it may happen that we have very different regular partitions, e.g., in one of them the densities are around 0 and 1, in the other they are around $\half$.

The natural question if $(\eps,\tau)$ and $G_n$ determine $H_\nu$ is considered in the paper of Alon, Shapira, and Stav \cite{AlonShapiraStav08DM}. The answer is ``No, but under some conditions YES''. However, here the reader should be
cautions, we have not defined when do we consider two regular partitions near to each other.

\end{remark}

So the cluster graph is not determined by these parameters, (neither the $\eps$-regular partition), and $\tau$ is much larger than $\eps$, however, in the applications $\tau\to0$ as $\eps\to0$. (If e.g. we try to embed a relatively small
graph (?) with bounded degree $\Delta$, then $\tau\approx\eps^{1/(2\De)}$ is mostly enough for our proofs, see also, e.g., G.N. Sárközy, \cite{Sarko14QuantiBlow}.).

There are many surveys on Regularity Lemmas and its applications, e.g., Komlós and Simonovits \cite{KomSim96SzemRegu}, Komlós, Shokoufandeh, Simonovits, and Szemerédi \cite{KomShokoSimSzem00}; the survey of Komlós \cite{Kom99BlowSurvey} (formally on the Blow-up Lemma) is also very informative, and we also recommend a more recent survey of Szemerédi \cite{Szemer15ArithComm}.

There are several results on sparse graphs, see e.g. surveys of Kohayakawa and Rödl \cite{KohayaRodl03Quasi} or of Gerke and Steger on the Sparse Regularity
Lemma \cite{GerkeSteg05SparseRegu}.

There are several works on some other aspects of the Regularity Lemma (or Regularity Lemmas), e.g., whether one needs exceptional classes, or how many
clusters are needed (see e.g. Gowers \cite{Gowers97Tower}), or how can it be viewed from other points of views, e.g., Tao \cite{Tao06Revis}, however, we skip most of them. 

\subsub Regularity Lemma and Parameters.//
The Regularity Lemma is very inefficient in the sense that it works only for very-very large graphs.

There are three natural questions concerning this: (a) do we need the exceptional cluster-pairs, and (b) How large the threshold $n_0(M_0,\eps)$ must be,
(c) how many clusters are needed. The answer to the first question was given by the ``half-graph'' where one must use exceptional pairs. It is not quite
clear how many exceptional pairs are needed in general.

Several results are known of that if we fix $\eps>0$ and $M_0$, then we have to use many-many clusters. The first such result was due to Gowers \cite{Gowers97Tower}. (See also Moshkovitz and Shapira, \cite{MoshShapi16Gowers}.)
For sharper results we refer the reader to Fox, László Miklós Lovász, and Yufei Zhao \cite{FoxLovYufei17Regu}, \dots

\subsection{Ramsey problems, cycles}

As to the Ramsey Theory, we shall not introduce the standard notation here. For an excellent source, see the book of Graham, Rothschild, and Spencer \cite{GrahRothSpenBook}.
 
The Ramsey problems were extended to arbitrary graphs first by Gerencsér and Gyárfás \cite{GerenGyarf67}. Several results were proved in this area by
Faudree and Schelp \cite{FaudSchelp75Bolyai,FaudSchelp75Path,FaudSchelp75PathPath}, (and some parallel to Faudree and Schelp, by Vera Rosta \cite{Rosta73})
and more generally, by Erdős, Faudree, Rousseau, and Schelp, and others. Bondy and Erdős \cite{BondyErd73Rams} formulated several important conjectures for
the case when the excluded graphs are paths, or cycles. In case of cycles, it turned out that the parity of the length of cycles is also very
important. Again, we shall mention only a few related results, and then mention a few papers: this area is too large to describe it here in more details.

We have to start with the remark that in most cases considered below the (conjectured) extremal structures come from some ``matrix-graph-sequence'': $n$ vertices of a $K_n$ are partitioned into a few classes and then we
colour $E(K_n)$ so that the colours depend only on the classes of the endvertices. These structures are very simple and nice, not chaotic/randomlike at all, so the proofs are also similar to the proofs of some extremal problems. Often we have some stability of the Ramsey-extremal structures.

\begin{wrapfigure}{L}{27mm}
\abrax{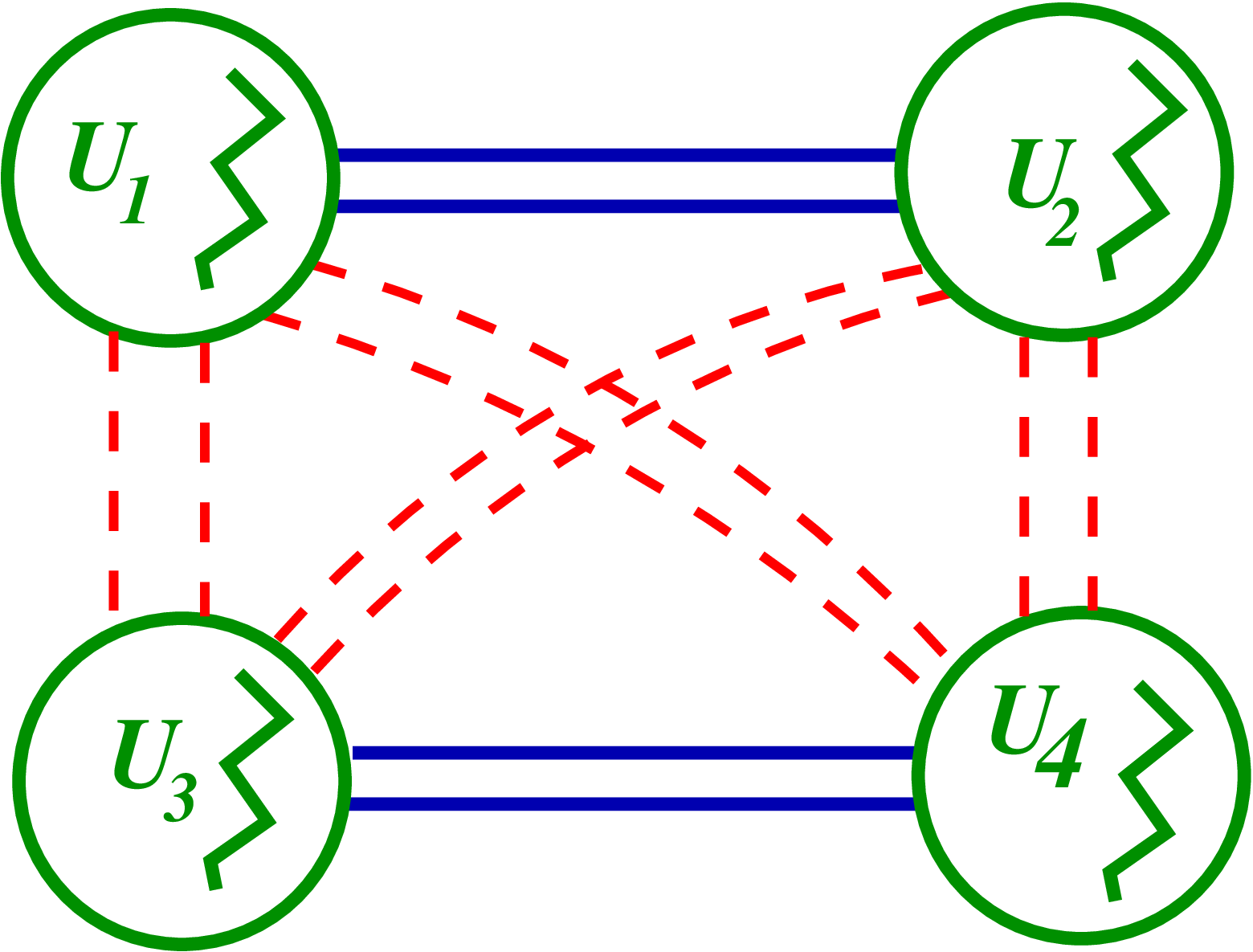}{20} 
\vskip3mm
\abrax{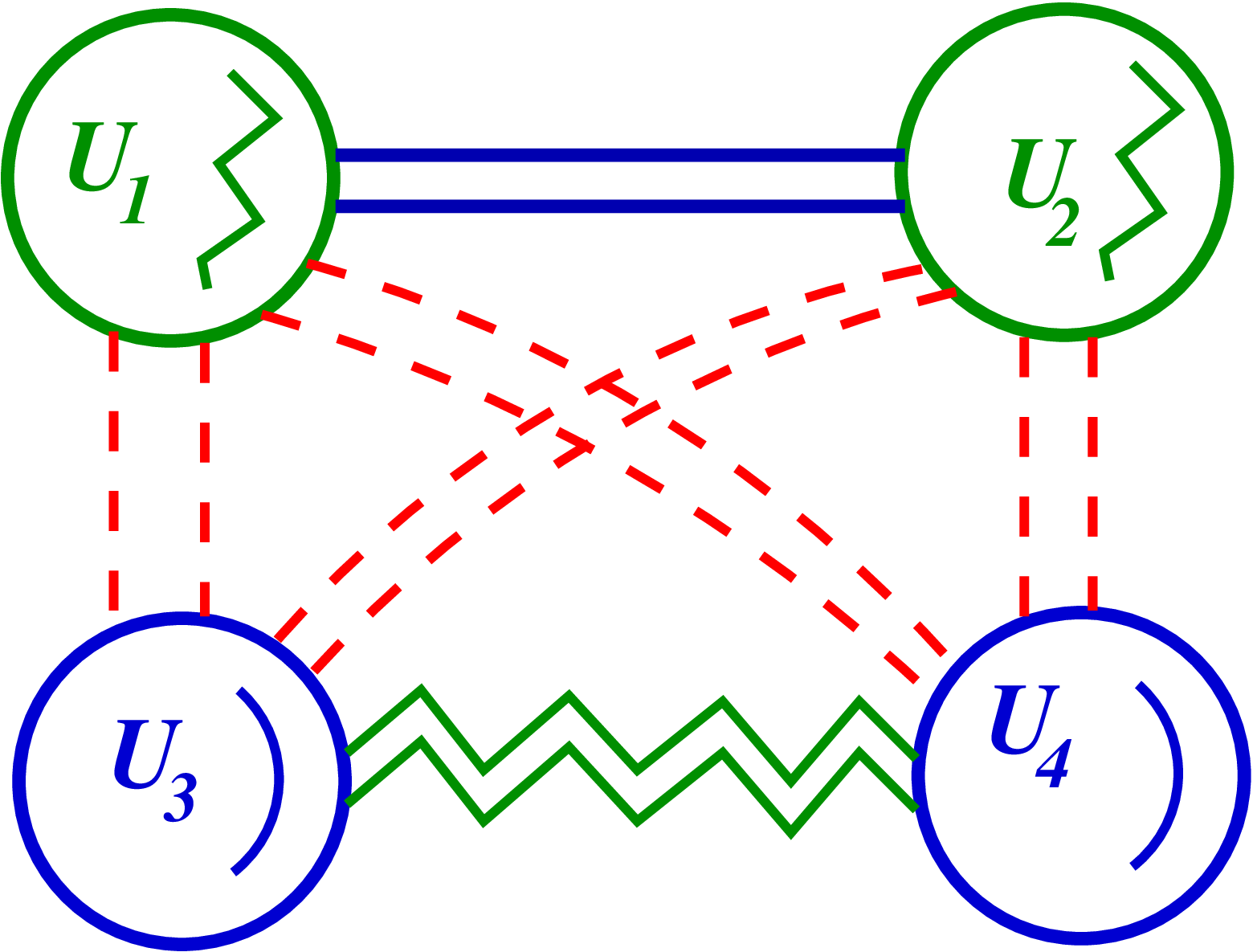}{20}
\caption{Cycle-Ramsey with three colours} 
\end{wrapfigure}

One of the many questions Bondy and Erdős \cite{BondyErd73Rams} asked was whether for an odd $n\ge3$, is it true that $R_r(C_n,C_n,\dots,C_n)=2^{r-1}(n-1)+1$? The construction suggesting/motivating this conjecture is recursive: for two colours we take two BLUE complete $K_{n-1}$ and join them completely by GREEN edges. If we have already constructions on $r-1$ colours, take two such constructions, and a new colour, say, RED, and join the two constructions completely with this new colour. (Watch out, if we use different sets of colours, there are other, similar but more complicated constructions, colour-connections between them.  As $r$ increases, the number of non-equivalent constructions increases. For $r=3$ we have two similar, yet different constructions.)

\begin{conjecture}[Bondy and Erdős] Let $n$ be odd.
If we $r$-colour $K_N$ for $N=2^{r-1}(n-1)+1$, then we shall have a monochromatic $C_n$.
\end{conjecture}

The following approximation of the Bondy-Erdős conjecture, for three colours, was a breakthrough in this area. 

\begin{theorem}[\L{u}czak \ev(1999) \cite{Lucz99BondyErd}] 
If $n$ is odd, then
$$R(C_n,C_n,C_n)=4n+o(n).$$ 
\end{theorem}

The proof was highly non-trivial. Applying stability methods, the $o(n)$ was eliminated:

\begin{theorem}[Kohayakawa, Simonovits, and Skokan \cite{KohaSimSkok05Graco}] If $n$ is a sufficiently large odd integer, then\footnote{Watch out, mostly it does not matter, but here, in case of sharp Ramsey results one has to distinguish the lower and upper Ramsey numbers. The upper one is the smallest one for which there is no good colouring, here $4k-3$. The lower Ramsey number is $R(L_1,\dots,L_r)-1$.} $$R(C_n,C_n,C_n)=4n-3.$$
\end{theorem}

Recently Jenssen and Skokan proved the corresponding general conjecture, \cite{JenssSkok16Ramsey}, for odd cycles, at least if $n$ is large. They use both the Regularity Lemma and the Stability method, adding some non-linear optimization tools to the usual methods, too.

The situation with the paths and the even cycles is different. We could say that the reason is that the above colouring contains many monochromatic even cycles. This is why for three colours the Ramsey number is half of the previous one. 
Figaj and Łuczak proved

\begin{theorem}[Figaj, Łuczak \cite{FigajLucz07Even}]
 If $\alpha_1\ge\alpha_2,\alpha_3>0$, then for $m_i:=2\lfloor\alpha_in\rfloor$ ($i=1,2,3$)
$$R(C_{m_1},C_{m_2},C_{m_3})=(2\alpha_1+\alpha_2+\alpha_3)n+o(n).$$
\end{theorem}

\begin{corollary} If $m_1\ge m_2\ge m_3$ are even, then
$$R(P_{m_1},P_{m_2},P_{m_3})=m_1+\half(m_2+m_3)+o(m_1).$$
\end{corollary}

In a recent paper Figaj and Łuczak \cite{FigLucz18RamsCyc} determined the asymptotic value of $R(C_\ell,C_m,C_n)$, when the parities are arbitrary and the asymptotic values of $\ell/n=\alpha$ and $m/n=\beta$ are given, as $n\to\infty$.
We mention here just the following theorem, that proves a conjecture of Faudree and Schelp.

\begin{theorem}[Gyárfás, Ruszinkó, G.N. Sárközy, and Szemerédi \cite{GyarfRuszSarkoSzem08-09}]
 There exists an $n_0$ such that if $n>n_0$ is even then $R(P_n,P_n,P_n)=2n-2$, if $n$ is odd then $R(P_n,P_n,P_n)=2n-1$. \end{theorem}

The results listed above use the Regularity Lemma, and one thinks they could be proved without it, too. The difference between paths and even cycles is not that large:

\begin{theorem}[Benevides and Skokan, \ev(2009) \cite{BenevSkokan09}]
If $n$ is even, and sufficiently large, then $R(C_n,C_n,C_n)=2n$.
\end{theorem}

These proofs use (among others) the Coloured Regularity Lemma. Hence all the assertions not using $o(..)$ are proved only for very large values of $n$.
This area is again large and ramified.
%% for some further related results see \wck \cite{xxxx}\dots. 
We shall return to some corresponding Ramsey hypergraph results in
Section~\ref{HyperRamseyS}.

\subsection{Ramsey Theory, general case}

Burr and Erdős conjectured \cite{BurrErd73Ramsey} that the 2-colour Ramsey number $R(H_n,H_n)$ is linear in $n$ for bounded degree graphs.\footnote{Here by linear we mean $O(n)$.} First some weaker bound was found by József Beck, but then the conjecture was proved by

\begin{theorem}[Chvátal-Rödl-Szemerédi-Trotter \ev(1983) \cite{ChvatRodlSzemTrotter83}]\label{ChvatRodlSzemTrotTh} For any $\Delta>0$ there exists a constant $\Gamma=\Gamma(\Delta)$ such that for any $H_n$ with $\maxdeg(H_n)<\Delta$, we have $R(H_n,H_n)<\Gamma n$.
\end{theorem}

The proof was based on the Regularity Lemma, applied to a graph $G_N$ defined by the RED edges of a RED-BLUE-coloured $K_N$. Then Turán's theorem and Ramsey's theorem were applied to the Coloured Cluster Graph $H_\nu$. The
RED-BLUE colouring of $K_N$ defines a RED-BLUE colouring of $H_\nu$: the cluster-edge $V_iV_j$ gets the colour which has the majority. Then the proof is completed by building up a RED copy of $H_n$ in the RED-BLUE $K_N$, if $H_\nu$
contained a sufficiently large RED complete subgraph. The next, stronger conjecture is still open.

\begin{conjecture}[Burr-Erdős]  If condition $\maxdeg(H_n)<\Delta$ is replaced by the weaker condition that for any $H^*\subset H_n$ we have $\mindeg(H^*)<\Delta$, then $R(H_n,H_n)=O(n)$.
\end{conjecture}

\begin{remark} 
  The actual bound of \cite{ChvatRodlSzemTrotter83} on the multiplicative constant in Theorem~\ref{ChvatRodlSzemTrotTh} was fairly weak. This was improved in several steps, first by N. Eaton \cite{Eaton98Ramsey}, then by Graham, Rödl, and Ruciński \cite{GrahRodlRuc00Lin}, and finally by Conlon \cite{Con09Hyper}, and Fox and Sudakov \cite{FoxSudak09}.
\end{remark}

\subsection{Ramsey-Turán theory}\label{RamsTurS}
 
When Turán, around 1969 \cite{Tur69Calga}, started a series of applications of his graph theorem in other fields of Mathematics, in some sense it turned out that $T_{n,p}$ is
too regular, has too simple structure: from the point of view of applications perhaps it would be better to have a version of his theorem for graphs with
less regular, more complicated structure, but providing better estimates.\footnote{Actually, Turán's corresponding results, or the Erdős-Sós type
 Ramsey-Turán theorems were not used in ``applications'', however the Ajtai-Komlós-Szemerédi type results are also in this category and, as we have seen in
 Sections \ref{Nibble}-\ref{GowersCollin}, they were used in several beautiful and important results.} This led Vera Sós, in \cite{Sos69Calgary} to ask
more general questions than the original Turán question, and that led her to the Turán-Ramsey theory. One of her questions was

\begin{problem}[Sós \ev(1969) \cite{Sos69Calgary}]
  Fix $r$ sample graphs $L_1,\dots,L_r$. What is the maximum of $e(G_n)$ if the edges of $G_n$ can be $r$-coloured so that the $i\th$ colour does not contain $L_i$.
\end{problem}

The answer was given by

\begin{theorem}[Burr, Erdős and Lovász \ev(1976) \cite{BurrErdLov76Ramsey}] For a family $\{L_1,\dots,L_r\}$ of fixed graphs,
let $R$ be the smallest integer for which there exists an $m$ such that if $T_{m,R}$ is arbitrarily $r$-coloured, then for some $i$ the $i\th$ colour contains an $L_i$. Then $$\RT(n;L_1,\dots,L_r)=\ext(n,K_{R-1})+o(n^2).$$
\end{theorem}

Actually, this easily follows from the Erdős-Stone theorem. (Further, obviously, we can choose $m=\sum v(L_i)$.)  Later it turned out that the really interesting problem is to determine $\RT(n;L_1,\dots,L_{r-1},o(n))$, which in the simplest case $r=2$ means the following.\footnote{One has to be cautious with this notation, when we write $o(n)$ instead of a function $f(n)$.}

\begin{problem}[Erdős-Sós \ev(1970) \cite{ErdSos70RamsTur}] Consider a graph sequence $(G_n)$. Estimate $e(G_n)$ if $L\not\subset G_n$, and the independence number, $\alpha(G_n)=o(n)$. 
\end{problem}

The case of odd complete graphs was solved by

\begin{theorem}[Erdős and Sós \ev(1970) \cite{ErdSos70RamsTur}]\label{ErdSosRT-Th} Let $(G_n)$ be a graph sequence. If $K_{2\ell+1}\not\subset G_n$ and $\alpha(G_n)=o(n)$, then $$e(G_n)\le\Tur n\ell+o(n^2),$$ and this is sharp.
\end{theorem}

Here the sharpness, i.e. the lower bound is easy: let $m=\lfloor n/\ell\rfloor$ and embed into each class of a $\Tn n,\ell$ a graph $G_m$ not containing $K_3$, with $\alpha(G_m)=o(m)$. By the ``probabilistic constructions'' of Erdős
\cite{Erd61GraphProbII} we know the existence of such graphs.\footnote{These are the graphs we considered in connection with $R(n,3)$ in Section~\ref{DegenS}.  There are many such graphs obtained by various, more involved constructions.} The
obtained graph $S_n$ provides the lower bound in Theorem \ref{ErdSosRT-Th}: it does not contain $K_{2\ell+1}$ and $\alpha(S_n)=o(n)$.

The problem of estimating $\RT(n,K_4,o(n))$ turned out to be much more difficult and this was among the first ones solved (basically!) by the Regularity Lemma, where the upper bound was given by Szemerédi \cite{Szem72K4} and the sharpness of this upper bound was proved by Bollobás and Erdős \cite{BolloErd76RT}.  It turned out that

\begin{theorem}[Szemerédi \ev(1972) \cite{Szem72K4}, Bollobás and Erdős \ev(1976) \cite{BolloErd76RT}]\label{SzemK4Th}
 $$\RT(n,K_4,o(n))=\left(\reci8 +o(1)\right)n^2.$$
\end{theorem}

The proof of the upper bound used (a simpler (?) form of) the Regularity Lemma. The first graph theoretical successes of the Regularity Lemma were on the $f(n,6,3)$-problem, see Theorem \ref{RuzsaSzemTh}, and the Ramsey-Turán theory:
Theorem~\ref{SzemK4Th}, and its generalization by Erdős, Hajnal, Sós, and Szemerédi \cite{ErdHajSosSzem83}, see also \cite{ErdHajSimSosSzem93Simple}, and \cite{ErdHajSimSosSzem94Kp} by Erdős, Hajnal, Simonovits, Sós, and Szemerédi. Since
those days this area has gone through a very fast development, Simonovits and Sós have a longer survey on Ramsey-Turán problems \cite{SimSosV01RT}. However, even since the publication of \cite{SimSosV01RT}, many beautiful new results have
been proved, e.g., by Mubayi and Sós \cite{MubayiSos06}, Schelp \cite{Schelp12RT}, Fox, Loh, Zhao \cite{FoxLohYufei15RT}, Sudakov \cite{Sudak03RT}, Balogh and Lenz \cite{BalLenz12RT} and \cite{BalLenz13RT2}, and many further ones.

In the proof of the upper bound of Theorem \ref{SzemK4Th} one needed something more than what was used in the earlier arguments: Until this point in most applications of the Regularity Lemmas the densities were near to 0 or near to 1. Here the densities in the Regular Partitions turned out to be near 0 or $\half$. Fixing the appropriate $\eps$ and $\tau=\sqrt[4]{\eps}$, one important step was that the Cluster Graph $H_\nu$ was triangle-free: $K_3\not\subset H_\nu$, and another was that in the
Regular Partition, $d\UiUj \le\half+o(1)$. This implies the upper bound in Theorem \ref{SzemK4Th}, namely
$$e(G_n)\le\ext(\nu,K_3)\cdot\left(\half+o(1)\right)\left({n\over \nu}\right)^2+\text{negligible terms}=\left(\reci8 +o(1)\right)n^2.$$

As to the lower bound in Theorem \ref{SzemK4Th}, the Erdős-Bollobás construction used a high dimensional isoperimetric inequality, similar to a related construction of Erdős and Rogers \cite{ErdRogers62}. One remaining important open question~is

\begin{problem}[Erdős]\label{ErdK222} 
Is $\RT(n,K(2,2,2),o(n))=o(n^2)$ or not? \end{problem}
 Rödl proved that

\begin{theorem}[Rödl \cite{Rodl85RT}]
  There exist a graph sequence $(G_n)$ with $\alpha(G_n)=o(n)$ and $K_4\not\subset G_n$, $K(3,3,3)\not\subset G_n$ for which $e(G_n)\ge \reci8 n^2+o(n^2)$.\footnote{Actually, Rödl proved a slightly stronger theorem, answering a question of Erdős, but the original one, Problem \ref{ErdK222}, is still open.}
\end{theorem}

\subsub Phase transition.// It can happen that $f(n)$ is a ``small'' function, much smaller than ``just'' $o(n)$, say $f(n)=o(\sqrt n)$, and then we see kind of a phase transition: one can prove better estimates on $\RT(n,L,f(n))$. In other words, sometimes when $f(n)$ is replaced by a slightly smaller $g(n)$, the Ramsey-Turán function noticeably
drops. Such result can be found, e.g., in Sudakov \cite{Sudak03RT}, in Balogh, Hu, and Simonovits \cite{BalHuSim15}, or in Bennett and Dudek \cite{BennettDudek17RT},\dots Bennett and Dudek also list several related results and papers. Some roots of this phenomenon go back to \cite{ErdHajSosSzem83,ErdHajSimSosSzem93Kp}.

We close this part with a two remarks on \Sc phase~transition results.
Perhaps the simplest approach is to investigate $RT(L_1,\dots,L_r,f(n))$,
when we know of $f(n)$ only that $f(n)=o(n)$, however, $f(n)/n\to0$ can be arbitrary slow. There are two directions from here:

(a) when we know of $f(n)$ that $f(n)/n\to0$ relatively fast, so fast that $\lim \reci n^2 RT(L_1,\dots,L_r,f(n))$ becomes smaller than for a slightly larger  $g(n)$.

(b) The other direction is when we investigate the $\delta$-dependence of $$\liminf \reci n^2 RT(L_1,\dots,L_r,\delta n).$$

There are some very new interesting results in both areas, we refer the reader to the above mentioned \cite{Sudak03RT} and \cite{BalHuSim15}, and also to papers of Lüders and Reiher \cite{LudersReih18A-RT} and of Kim, Kim
and Liu \cite{KimKimLiu18A-RT} settling several phase-transition problems earlier unsolved. They determined $\RT(n,K_3,K_s,\delta n)$ for $s=3,4,5$ and
sufficiently small $\delta$, confirming -- among others -- a conjecture of Erdős and Sós from 1979, and settling some conjectures of Balogh, Hu, and
Simonovits, according to which $\RT(n,K_8,o(\sqrt{n\log n}))=\frac{n^2}{4}+o(n^2)$.

\subsection{Ruzsa-Szemerédi theorem, Removal Lemma}\label{RemovalS1}

We start with a slightly simplified version of a result of G. Dirac \cite{Dirac63Turan}, a generalization of Turán's theorem:

\begin{theorem}[Dirac \ev(1963) \cite{Dirac63Turan}]
Let $p\ge2$, $q\in[1,p]$ and $n\ge p+q$ be integers. If $G_n$ is a graph with $e(G_n)>e(T_{n,p})$ edges, then $G_n$ contains a subgraph of $p+q$ vertices and ${{p+q}\choose 2}-(q-1)$ edges. 
\end{theorem}

Some related results were also independently obtained by Erdős, see \cite{Erd89Dir}, and also \cite{Erd65Smole}. Erdős asked the following problem. 

\begin{problem}[\cite{Erd65Smole}]
Let $\cL_{k,\ell}$ be the family of graphs $L$ with $k$ vertices and $\ell$ edges. Determine (or estimate) $f(n,k,\ell):=\ext(n,\cL_{k,\ell})$.
\end{problem}

Several results were proved on this problem by Erdős \cite{Erd65Smole}, Simonovits \cite{Sim69Kandid}, Griggs, Simonovits, and Rubin G. Thomas \cite{GriggsSimThom}, Simonovits \cite{Sim99Stirin} and others. Extending these types of results to 3-uniform hypergraphs in \cite{BrownErdSos73a} and then to $r$-uniform hypergraphs in \cite{BrownErdSos73b}, W. G. Brown, P. Erdős and V. T. Sós started a systematic investigation of $f_r(n,k,\ell)$ defined as the maximum number of
hyperedges an $r$-uniform $n$-vertex hypergraph ${\cH}_n^{(r)}$ can have without containing some $k$-vertex subgraphs with at least $\ell$ hyperedges.

Below mostly we shall restrict ourselves to the case of 3-uniform hypergraphs, i.e. $r=3$. Several subcases were solved in \cite{BrownErdSos73a}, but the problem of $f(n,6,3)$, i.e. the case when no 6 vertices contain 3 triplets, turned out to be very difficult. One can easily show that $f(n,6,3)\le\reci6 n^2$. Ruzsa and Szemerédi proved that

\begin{theorem}[Ruzsa-Szemerédi \ev(1978) \cite{RuzsaSzem78}]\label{RuzsaSzemTh}
\beq{RuzsaSzemFo} cn\cdot r_3(n)<f(n,6,3)=o(n^2).\eeq 
\end{theorem}

The crucial tool was a consequence of the Regularity Lemma:

\begin{theorem}[Ruzsa-Szemerédi Triangle Removal Lemma]\label{RuzsaSzemRemTh}
 If $(G_n)$ is a graph sequence with $o(n^3)$ triangles, then we can delete $o(n^2)$ edges from $G_n$ to get a triangle-free graph.
\end{theorem}

Clearly, \eqref{RuzsaSzemFo} implies Roth theorem that $r_3(n)=o(n)$. There is a more general form of Theorem \ref{RuzsaSzemRemTh}, the so called \Sc Removal~Lemma:

\begin{theorem}[Removal Lemma]\label{RemovalLem} 
If $v(L)=h$, then for every $\eps>0$ there is a $\delta>0$ for which,
 if a $G_n$ contains at most $\delta n^h$ copies of $L$, then one can delete $\eps n^2$ edges to destroy all the copies of $L$ in $G_n$.
\end{theorem}

\begin{remark}[Induced matchings] Let us call some edges of a graph $G_n$ \Sc Strongly~Independent, if there are no edges joining two of them.  The Ruzsa-Szemerédi theorem can be formulated also without using hypergraphs. Indeed, for each $x\in V(\cH)$, the \Sc link of $x$, i.e. the pairs $uv$ forming a hyperedge in $\cH$ with $x$, form a so-called \Sc induced~matching: not only its edges are independent but they are pairwise \Sc strongly~independent.\footnote{Here we assume
    that the mindegree is at least 3.}  So the Ruzsa-Szemerédi theorem has the following alternative form:

\begin{theorem}[Ruzsa-Szemerédi \cite{RuzsaSzem78}]
 If $E(G_n)$ can be covered by $n$ induced matchings, then $e(G_n)=o(n^2)$.
\end{theorem}
 
\end{remark}

\Sc For~some~further~applications\footnote{Some applications of the Ruzsa-Szemerédi theorem are given in Subsection \ref{RuzsaSzemAppliS}.}  and results connected to Induced Matchings see also the papers of Alon, Moitra, and Sudakov
\cite{AlonMoiSudak13}, and also of Birk, Linial, and Meshulam \cite{BirkLiniMesh93}.  Alon, Moitra, and Sudakov describe several constructions and their applications (in theoretical computer science) connected to the
Ruzsa-Szemerédi theorem, more precisely, to dense graphs that can be covered by a given number of induced matchings.

\begin{remark} Among the original questions of Brown, Erdős, and Sós, the estimate of $f(n,7,4)$ is still unsolved. For some special graphs connected to groups, Solymosi \cite{Solymo15_7_4} has a solution. Actually, very recently Nenadov, Schreiber, and Sudakov \cite{NenadSudakTyomkin19A-Groups} extended this result.
\end{remark}

\begin{remark} 
 (a) Though Theorem \ref{RemovalLem} was not explicitly formulated in \cite{RuzsaSzem78}, implicitly it was there. Later it was explicitly formulated,
 e.g., in the paper of Erdős, Frankl, and Rödl \cite{ErdFranklRodl86Asym} and Füredi \cite{Fure92DiamCrit}, \cite{Fure94Zurich}, and soon it became a central research topic, partly because it is often applicable and partly because
 -- though for ordinary graphs the Removal Lemma easily follows from the Regularity Lemma, for hypergraphs everything is much more involved.

(b) The Ruzsa-Szemerédi theorem was extended by Erdős, Frankl, and Rödl \cite{ErdFranklRodl86Asym} to $r$-uniform hypergraphs. They proved that
$$f_r(n,3r-3,3)=o(n^2).$$ 

(c) Actually, $f(n,6,3)$ is a fixed function of $n$, and -- though \cite{RuzsaSzem78} asserts only that $f(n,6,3)=o(n^2)$, -- it has some better asymptotics
and they are very interesting and important, since these results have many applications. J. Fox gave a new, more effective proof of the Removal Lemma
\cite{Fox11NewRemo}, not using the Regularity Lemma, leading also to better estimates on $f(n,6,3)$, according to which
\beq{FoxFo}f(n,6,3)<{n^2\over {2^{\log^*n}}}.\eeq

(d) A more detailed description of this topic can be found in two surveys of Conlon and Fox \cite{ConFox12Remo,ConFox13RemoveBCC}. We skip all the details connected to ``Tower'' functions and ``Wowzer'' functions,\footnote{with the
  exception of the next theorem.} (which show that the Regularity-Lemma methods are very inefficient). Conlon and Fox discuss the estimates connecting the various parameters, in details.

(e) Some related results and generalizations were proved by G.N. Sárközy and Selkow \cite{SarkoSelko04RuzsaSzem, SarkoSelko05RuzsaSzem}, and by Alon and Shapira \cite{AlonShapira06Brown}. The emphasis in \cite{AlonShapira06Brown} is on
obtaining a new lower bound to generalize Theorem \ref{RuzsaSzemTh}.\footnote{Actually, this assertion is somewhat more involved, see the introduction of \cite{AlonShapira06Brown}.}

(f) Ruzsa-Szemerédi removal lemma was applied also by Balogh and Petříč\-ková \cite{BalPetri14Max}, to give a sharp estimate on the number of maximal triangle-free graphs~$G_n$.
\end{remark} 

We close this part with quoting a result of Jacob Fox. Let $TF(t)$ be the ``tower'' of 2's of height $t$: $TF(1)=2$ and $TF(t+1)=TF(2^{TF(t)})$. 

\begin{theorem}[J. Fox \cite{Fox11NewRemo}]\label{Fox11RemoTh}
 Fix a sample graph $L_h$ (on $h$ vertices) and an $\eps>0$. Let $\de:=1/TF(\lceil5h^4\log{1/\eps}\rceil)$. Every $G_n$ with at most $\de n^h$ copies of
 $L_h$ can be made $L_h$-free by removing at most $\eps n^2$ edges.
\end{theorem}

\begin{remark}
The surveys of Conlon and Fox, e.g., \cite{ConFox13RemoveBCC} also discuss the connection of property testing if $G_n$ contains an induced copy of $H$, and its connection to the ``Induced Removal Lemma''.
\end{remark}

The interested reader is recommended to read \cite{Fox11NewRemo} or the survey of Conlon and Fox \cite{ConFox13RemoveBCC}. 

\minisepa

\begin{remark} 
(a) Mostly we skip the references to graph limits, however, here we should mention the Elek-Szegedy Ultra-product approach to graph limits, and within that to the Removal Lemma, see \cite{ElekSzegedy12}. %\Margo{Elhagyva ennek a 07=es valtozata}

(b) T. Tao also has a variant of the Hypergraph removal lemma \cite{Tao06HyperRemo} which he uses to prove a Szemerédi type theorem on the Gaussian primes \cite{Tao06Gaussian}.
\end{remark}

\begin{remark}[Szemerédi theorem and the Clique-union Lemma] 
  V. Rödl, as an invited speaker of the ICM 2014, in Seoul, spoke about the relations described above. Generalizing the above approach, Peter Frankl and he tried to prove $r_k(n)=o(n)$, using this combinatorial approach. The reader is suggested to look up the movie about his lecture.  (ICM 2014 Video Series, Aug 21)
\end{remark}

\subsub Some generalizations.// The Ruzsa-Szemerédi theorem has two parts, and both have several important and interesting generalizations, yet, in this paper we mostly skip the results connected to the lower bounds. Several related results can be found in the papers of Erdős, Frankl, and Rödl \cite{ErdFranklRodl86Asym}, or e.g., in the paper of Füredi and Ruszinkó on excluding the grids \cite{FureRuszi13Grid}, and also in G.N. Sárközy and Selkow \cite{SarkoSelko05RuzsaSzem}, and Alon and Shapira \cite{AlonShapira06Brown}.

\begin{remark}
Solymosi tried to formulate a non-trivial Removal Lemma for bipartite excluded graphs, however Timmons and Verstraëte \cite{TimVerstra15SparseRemo} provided infinitely many ``counterexamples''.
\end{remark}

\subsection{Applications of Ruzsa-Szemerédi Theorem}\label{RuzsaSzemAppliS}

As we have mentioned, one of the early successes of a simpler form of the Regularity Lemma was the answer to the question of Brown, Erdős, and Sós,
according to which $f(n,6,3)=o(n^2)$, (and the proof of the triangle removal lemma, implying this). Fox writes in \cite{Fox11NewRemo}: ``The graph removal lemma has many
applications in graph theory, additive combinatorics, discrete geometry, and theoretical computer science.'' Here we mention two graph theoretical
applications, one of which is strongly connected to Turing Machines.

\subsubsection{Füredi theorem on diameter critical graphs}\label{DiamCritical} 

Call a graph $G$ \Sc diameter-$d$-Critical, if it has diameter $d$ and deleting any edge the diameter becomes at least $d+1$ (or $G$ gets disconnected). Such graphs are, e.g., $C_k$ or $T_{n,p}$. Restrict ourselves to $d=2$: consider
diameter-critical graphs of diameter 2. As Füredi describes in \cite{Fure92DiamCrit}, Plesník observed \cite{Plesnik75Crit} that for all known diameter-2-critical graphs $G_n$, $e(G_n)\le \turtwo n$. Independently, Simon and Murty
conjectured \cite{CaccHaggkv79DiamCrit} that

\begin{conjecture}
If $G_n$ is a diameter-2-critical graph, then $e(G_n)\le \turtwo n$. 
\end{conjecture}

This seemed to be a beautiful but very difficult conjecture. Füredi \cite{Fure92DiamCrit} proved it for large $n$, using the Ruzsa-Szemerédi theorem:

\begin{theorem}[Füredi \ev(1992) \cite{Fure92DiamCrit}]\label{Fure92DiamCritTh} There exists an $n_0$ such that if $n>n_0$, then the Murty-Simon conjecture holds.
\end{theorem}

\begin{remark} In the proof of Theorem \ref{Fure92DiamCritTh} we get a very large $n_0$.  In some sense Theorem \ref{Fure92DiamCritTh} settles the conjecture, at least for most of us. Yet a lot of work has been done to prove it for reasonable values of $n$, too. Plesník proved, instead of $e(G_n)\le \reci4 n^2$, that $e(G_n)\le {3\over8}n^2$, Caccetta and Häggkvist \cite{CaccHaggkv79DiamCrit} improved this to $e(G_n)\le 0.27n^2$, and G. Fan \cite{Fan87Criti} proved for $n>25$ that $e(G_n)\le
  0.2532n^2$. \footnote{Watch out, some of these papers, e.g., \cite{Fan87Criti} are from before Füredi's result, some others are from after it.}

\end{remark} 

There are many related results proved and some conjectures on the $d$-diameter-critical graphs, see e.g. a survey of Haynes, Henning, van der Merwe, and Yeo, \cite{HayHenMerYeo15Murty}, or Po-Shen Loh and Jie Ma \cite{LohMa16DiamCrit}. This later one disproves a Caccetta-Häggkvist conjecture on the average edge-degree of diameter-critical graphs (and contains some further nice results as well).

\subsubsection{Triangle Removal Lemma in Dual Anti-Ramsey problems}

Burr, Erdős, Graham and Sós started investigating the following ``dual Anti-Ramsey'' problem \cite{BurrErdGraSos89}, (see also \cite{BurrErdFraGraSos88}):

\begin{quote}{Given a sample graph $L$, $n$, and $E$, what is the minimum number of colours $t=\chi_S(n,E,L)$ such that any graph $G_n$ with $E$ edges can be edge-coloured with $t$ colours so that all the copies $L\subset G_n$ be Rainbow coloured (i.e. its edges be coloured without colour-repetition? \footnote{The corresponding extremal value will be denoted by $\chiS(n,E,L)$. Here $S$ stands for ``strong'' in $\chiS$. It is the \Sc strong~chromatic~number of the $v(L)$-uniform hypergraph whose hyperedges are the $v(L)$-sets of vertices of the copies of $L\subset G_n$.}}
\end{quote}

They observed that for $L=\P4$ this question can be solved using $f(n,6,3)=o(n^2)$. Actually, their solution for $L=\P4$ is reducing the problem
to the problem of estimating $f(n,6,3)$. For $L=C_5$ (which seemed one of the most interesting cases), they proved:

\Proclaim Theorem 4.1 of \cite{BurrErdGraSos89}. There exists an $n_0$ such that if $n>n_0$ and $e=\turtwo n+1$, then \beq{Fo2} c_1n\le \chiS(n,e,C_5)\le
\fele n+3. \eeq

Erdős and Simonovits \cite{ErdSim18Begs}, using the Lovász-Simonovits Stability theorem, (see Subsection \ref{LovSimStabS}) proved that the upper bound (obtained from a simple construction) is
sharp.

\begin{theorem}[Erdős and Simonovits]\label{MainRes}
 There exists a threshold $n_0$ such that if $n>n_0$, and a graph $G_n$ has $E=\turtwo n+1$ edges and we colour its edges so that every $C_5$ is
 $5$--coloured, then we have to use at least $\fele n +3$ colours.
\end{theorem}

To apply the Lovász-Simonovits Stability, they needed the result of \cite{BurrErdGraSos89} on $\P4$. So, again, the removal lemma was one of the crucial tools. (The application of the Lovász-Simonovits Stability can be replaced here by a second application of the Removal Lemma and the Erdős-Simonovits Stability, however that approach would be less elementary and effective.) They also proved several further results, however, we skip them. 
Simonovits has also some further results in this area, e.g., on the problem of $C_5$ when $e(G_n)=\turtwo n +k$, and $k$ is any fixed number, or tending to $\infty$ very slowly, and also on the problem of $C_7$.
Another conjecture of Burr, Erdős, Graham, and Sós was (almost completely) proved by G.N. Sárközy and Selkow \cite{SarkoSelko06AntiRamsey}.

\subsection{Hypergraph Removal Lemma?}

There are several ways to approach the Removal Lemma and the Hypergraph Removal Lemma.  Rödl and Schacht \cite{RodlSchacht09Removal} describe a hypergraph generalization of the Removal Lemma. In the nice introduction of
\cite{RodlSchacht09Removal} they also write

\begin{quote}
{\dots the result of Alon and Shapira \cite{AlonShapira05Test} is a generalization which extends all the previous results of this type where the triangle is replaced by a possibly infinite family $\cF$ of graphs and the containment is induced.\dots }
\end{quote}

We have promised to avoid some areas that are much more difficult/technical than the others. Unfortunately, it is not easy to decide what is ``technical''. Komlós and Simonovits \cite{KomSim96SzemRegu} tried to provide an easy introduction to the Regularity Lemma. If one knows enough non-discrete mathematics, e.g., compactness, metric spaces, then
the Lovász-Szegedy paper \cite{LovSzege07Anal} is an easy introduction to the area of connections of Regularity Lemma and other, related areas.

The situation with hypergraphs is different. That area seems inherently difficult to ``learn''. The paper of Rödl, Nagle, Skokan, Schacht, and Kohayakawa \cite{RodlNagleSkokSchKoha05PNAS} and the companion paper of Solymosi
\cite{Solymo05PNAS} helps a lot to bridge this difficulty. The difficulties come from two sources. The first one is that there are useful and not so useful versions of the Hypergraph Regularity Lemma, and the useful ones are difficult to
formulate. Layers come in, e.g., for 3-uniform hypergraphs: beside considering the vertices and hyperedges, we have to consider some auxiliary graphs.

In \cite{RodlNagleSkokSchKoha05PNAS} twelve theorems are formulated, the last two are the hypergraph regularity and the hypergraph counting lemmas.
It is nice that the hypergraph removal lemma keeps its simple form. 

\begin{theorem}[See Theorem 5 of \cite{RodlNagleSkokSchKoha05PNAS}] For any fixed integers $\ell\ge k\ge2$ and $\eps>0$, there exists a $\xi=\xi(\ell,k,\eps)$ and an $n_0(\ell,k,\eps)$, such that if $\Fhyp \ell k$ is a $k$-uniform hypergraph on $\ell$ vertices and $\Hhyp nk$ in a $k$-uniform hypergraph of $n$ vertices, with fewer than $\xi n^\ell$ copies of $\Fhyp \ell k$, then it may be transformed into an $\Fhyp \ell k$-free hypergraph by deleting $\eps n^k$ hyperedges.

\end{theorem}

Solymosi remarked that the appropriate hypergraph removal lemma implies the higher dimensional version of the Szemerédi theorem. We close this part by quoting Solymosi \cite{Solymo05PNAS}:

\begin{quote}{
``There is a test to decide whether a hypergraph regularity is useful or not. Does it imply the Removal Lemma? If the answer is yes, then it is a correct concept of regularity indeed. On the contrary, applications of the hypergraph
regularity could go beyond the Removal Lemma. There are already examples for which the hypergraph regularity method, combined with ergodic theory, analysis,
and number theory, are used efficiently to solve difficult problems in mathematics.'' }
\end{quote}

The hypergraph removal lemma was approached from several directions. Among others, Tao considered it in \cite{Tao06HyperRemo}, Elek and B. Szegedy \cite{ElekSzege07Remov} approached it from the direction of Non-Standard Analysis, Rödl and Schacht from their general hypergraph regularity theory.

\subsection{The Universe of Random Graphs}\label{RandomUnivers}

The following result answers a question of Erdős:

\begin{theorem}[Babai, Simonovits and Spencer \ev(1990) \cite{BabaiSimSpenc90}]\label{BabaiSimSpencTh} There exists a $p_0<\half $ such that a random graph $R_{n,p}$ with edge-probability $p>p_0$ almost surely has the following property $\cB_L$, for $L=K_3$: all its triangle-free subgraphs with maximum number of edges are bipartite\footnote{where ``almost surely'' means that its probability tends to 1 as $n\to\infty$.}.
\end{theorem}

Actually, they proved much more general results. Consider the following assertion, depending on $L$ and $p$.

\begin{quote}{ ($\cB_{L,p}$) All the $L$-free subgraphs $F_n\subset R_{n,p}$ having maximum number \dori ~~ of edges are $\chi(L)$-chromatic, almost surely, if $n>n_0(L,p)$.}
\end{quote}

They proved ($\cB_{L,p}$) for all cases when $L$ has a critical edge \footnote{see Meta-Theorem \ref{MetaCritical}.} and $p$ is nearly $\half$. They also proved several related weaker results in the general case, when $L$ was arbitrary, and $p>0$. They could not extend their results to sparse graphs, primarily because that time the sparse Regularity Lemma did not exist. Soon it was ``invented'' by Kohayakawa.\footnote{Rödl also knew it, but it seems that he had not published it.} The results of Babai, Simonovits and Spencer were generalized, first by Brightwell, Panagiotou and Steger \cite{BrightPanaSteger12}, and then, in various ways, by others. So the first breakthrough towards sparse graphs was

\begin{theorem}[Brightwell, Panagiotou and Steger \ev(2012) \cite{BrightPanaSteger12}] 
  There exists a constant $c>0$ for which choosing a random graph $R_{n,p}$ where each edge is taken independently, with probability $p=n^{-(1/2)+\eps}$, the largest triangle-free subgraph $F_n$ of $R_{n,p}$ is bipartite, with probability tending to 1.
\end{theorem}

They remarked that the conclusion cannot hold when $p=\reci10 {\log n\over\sqrt n}$, since these $R_{n,p}$ contain, almost surely, such an induced $C_5$ {\em whose edges are not contained in triangles of this $R_{n,p}$}. All the edges of $R_{n,p}$ not covered by some $K_3\subset R_{n,p}$ must belong to $F_n$. Therefore now $C_5\subset F_n$: $F_n$ is
not bipartite. The proof of \cite{BrightPanaSteger12}, similarly to that of the original proof of Babai, Simonovits and Spencer, uses a stability argument, however, instead of the original Regularity Lemma it uses the Sparse Regularity Lemma. In some sense the ``final'' result was found by DeMarco and Jeff Kahn
\cite{DeMarcoKahn15Man,DeMarcoKahn15ArxTur}. They proved (among others) that

\begin{theorem}[DeMarco and Kahn \ev(2015) \cite{DeMarcoKahn15ArxTur}]
For each $r$ there exists a constant $C=C_r>0$ for which choosing a random graph $R_{n,p}$ where each edge is taken independently, with probability 
$$p>Cn^{-{2\over(r+1)}} \log^{2\over{(r+1)(r-2)}}n,$$ the largest $K_r$-free subgraph $T_n$ of $R_{n,p}$ is almost surely $r-1$-partite.
\end{theorem}

A hypergraph analog was proved by Balogh, Butterfield, Hu, and Lenz \cite{BalButtHuLenz16Mant}. (They again used the stability approach.)

\subsection{Embedding large trees}\label{LargeTrees}

There are many results where one fixes a sample graph $L$ and tries to embed it into a Random Graph $R_{n,p}$. (See e.g. Erdős-Rényi \cite{ErdRenyi60Evol}.) Here we try to embed a fixed tree $T_m$ of
$m\approx(1-\alpha)n$ vertices into a (random) graph $R_n$. However, we use a slightly different language.

A relatively new notion of \Sc resilience was introduced by Sudakov and Vu \cite{SudakVu08Resili}. Fix a graph property $\cP$.  The resilience of $G_n$ is its ``edit'' distance\footnote{The ``edit'' distance is the same used in
  \cite{Sim66Tihany}: the minimum number of edges to be changed to get from $G_n$ a graph isomorphic to $H_n$.}  to graphs not having property $\cP$.\footnote{Though we formulate a theorem on the local resilience of graphs for some graph
  property, we shall not define here the notion of local and global resilience: we refer the reader to the papers of Sudakov and Vu \cite{SudakVu08Resili}, or Balogh, Csaba, and Samotij \cite{BalCsabaSamo11Trees}, or suggest to read only
  Theorem~\ref{BaloCsabSamoT}.}  Balogh, Csaba, and Samotij \cite{BalCsabaSamo11Trees} proved

\begin{theorem}
  Let $\alpha$ and $\gamma$ be (small) positive constants and assume that $\Delta\ge 2$. There exists a constant $C>0$ (depending on $\alpha,\gamma$, and $\Delta$) such that for all $p=p(n)\ge C/n$, the local resilience of $R_{n,p}$ with
  respect to the property of containing all trees $T_m$ of $m:=\lfloor(1-\alpha)n\rfloor$ vertices and maximum degree $\maxdeg(T_m)\le \Delta$ is almost surely greater than $\half-\gamma$.
\end{theorem}

As a subcase, this contains

\begin{theorem}[Balogh, Csaba, and Samotij \cite{BalCsabaSamo11Trees}]\label{BaloCsabSamoT}
Let $\alpha$ be a positive constant, and assume that $\Delta\ge 2$. There exists a constant $C>0$ (depending on $\alpha$, and $\Delta$) such that for all $p=p(n)\ge C/n$, $R_{n,p}$ contains all trees $T_m$ of $m\le (1-\alpha ) n$ vertices and maximum degree $\maxdeg(T_m)\le \Delta$, almost surely, as $n\to\infty$.
\end{theorem}

The proof uses a sparse Regularity Lemma and a theorem of Penny Haxell \cite{Haxell00TreeEmb} on embedding bounded degree trees into ``expanding'' graphs.

\subsection{Extremal subgraphs of Pseudo-Random graphs}\label{PseudRandomUniv}

Another direction of research is when the Universe consists of more general objects, say of pseudo-random graphs. These are natural directions:

(a) whenever we can prove a result for complete graphs $K_p$, we have a hope to extend it to any $L$ with critical edges, and 

(b) whenever we know something for Random Graphs, there is a chance that it can be extended to random-looking objects (say to quasi-random graphs, or Pseudo-Random graphs, or to expanders graphs.\footnote{We have defined only the quasi-random graphs here, for pseudo-random graphs see e.g. \cite{Thomason87Poznan,KriveSudak06Pseudo}, for expanders see e.g. \cite{Alon86Eigen}.}.)

We mention here a few such papers: Thomason \cite{Thomason87Poznan,Thomason87Pseudo}, Krivelevich and Sudakov \cite{KriveSudak06Pseudo} are nice and detailed surveys on Pseudorandom graphs. Aigner-Horev, Hàn, and Schacht
\cite{AigHanSchacht14Odd}, and D. Conlon, J. Fox, and Yufei Zhao \cite{ConFoxZhao14SparsePseudo} also are two more recent nice papers (surveys?)  on this topic. We recommend these papers, and also Kühn and Osthus
\cite{KuhnOsthus13HamilExpan}. For hypergraphs see, e.g., Haviland and Thomason, \cite{HaviThomason89PseudoHyperC,HaviThomason92PseudoHyper}, or Kohayakawa, Mota and Schacht \cite{KohaMotaSchacht17}.

\subsection{Extremal results on ``slightly randomized'' graphs}\label{PseudRandomUniv2}

Bohman, Frieze, and Martin \cite{BohFrieMartin16Hamil} proved Hamiltonicity for graphs $R_n$ obtained from a non-random graph $G_n^0$ whose minimum degree was $d<n/2$ where $R_n$ is obtained by adding to $G_n^0$ a random binomial graph
$R_{n,p}$ whose edge-probability is a small $p=p_n>0$.  Dudek, Reiher, Ruciński, and Schacht continued this line, using this model, and proving in \cite{DudekReihRucSchacht18A} that the stronger conclusion of Pósa-Seymour conjecture holds almost surely for the obtained graph, under their conditions: instead of a Hamiltonian cycle they obtained a power of it.

\subsection{Algorithmic aspects}\label{Algorithmic}

Whenever we use an existence-proof in combinatorics, it is natural to ask if we can turn it into an algorithm. This was the case with the Lovász Local Lemma, and this is the case with the Regularity Lemma, and also with the Blow-up Lemma.

The algorithmic aspects of the Lovász Local Lemma were investigated first by József Beck \cite{Beck91LovLoc}, and later by R.A. Moser and G. Tardos \cite{Moser08LovLoc}, \cite{MosTard10LovLoc}, and many others.

\minisepa

Alon, Duke, Lefmann, Rödl, and Yuster in \cite{AlonDukeYuster94Algo} showed that one can find the Regular Partition, and the Cluster Graph of a $G_n$ fairly efficiently. Actually, they proved two theorems: (a) to decide if a partition is
$\eps$-regular is difficult, but (b) to find an $\eps$-regular partition of a given graph is easy. More precisely,

\begin{theorem}[\cite{AlonDukeYuster94Algo}]\label{AlonDukeYus1Th}
Given a graph $G_n$ and an $\eps>0$,
and a partition $(V_0,\dots,V_k)$,\footnote{Here we have $k+1$ classes, since originally there was also an exceptional class $V_0$, different from the others. This $V_0$ can be forgotten: its vertices can be distributed in the other classes.} it is 
 Co-NP-complete to decide if this partition is $\eps$-regular.
\end{theorem}

\begin{theorem}[\cite{AlonDukeYuster94Algo}]\label{AlonDukeYus2Th}
  For every $\eps>0$ and every $t>0$, there exists a $Q(\eps,t)$ such that for every $G_n$ with $n>Q(\eps,n)$ one can find an $\eps$-Regular Partition (described in Thm \ref{SzemReguTh}) in $O(M(n))$ steps, where $M(n)$ is the ``number of steps'' needed to multiply two $0-1$ matrices over the integers.\footnote{The theorem also has a version on parallel computation.} \footnote{Here we do not define the ``steps'' and ignore again the difference caused by neglecting $V_0$ in Theorem \ref{SzemReguTh}. }
\end{theorem} 

The algorithmic problem with the Regularity Lemma is that we may have too many clusters. Therefore a direct way to transform it into an efficient algorithm may be hopeless. The Frieze-Kannan version often solves this problem, see Subsection \ref{Versions-RL}, or \cite{Kannan94Focs,FriezeKannan96Survey,FriezeKannan99Weak}, or Lovász-Szegedy \cite{LovSzege07Anal}.

The above methods were needed and extended to hypergraphs, see e.g. Nagle, Rödl, and Schacht \cite{NagleRodlSchacht18AlgorRegLemma}.

\begin{remark}
We needed this short section {\em here}, since it helps to understand the next part better: otherwise it would come later.
\end{remark}

\subsection{Regularity Lemma for the Analyst}

As soon as the theory of Graph Limits turned into a fast-developing research area, it became interesting, what happens with the regularity lemmas in this area. As we have mentioned, Lovász wrote first a long survey \cite{Lov09LimitSurvey},
then a thick book on graph limits \cite{Lov12LimitBook}, so we shall not discuss it here, but mention just one aspect. The paper of Lovász and Szegedy \cite{LovSzege07Anal} with the title ``Szemerédi Lemma for the analyst'' not only
described the Regularity Lemma in terms of the Mathematical Analysis, but also described the Weak Regularity Lemma (i.e. the Frieze-Kannan version \cite{FriezeKannan99Weak}), and the Strong Regularity Lemma \cite{AlonFischKrivSzege99} of
Alon, Fischer, Krivelevich, and M. Szegedy, and connected the Regularity Lemma to $\eps$-nets in metric spaces, and to Compactness. We quote part of the Introduction of their paper.

\begin{quote}{
 "Szemerédi's regularity lemma was first used in his celebrated proof of the Erdős-Turán Conjecture on arithmetic progressions in dense sets of
 integers.\footnote{As we have mentioned, this is not quite true. It was invented to prove a conjecture of Bollobás, Erdős, and Simonovits on the
 parametrized Erdős-Stone theorem, and was first used in the paper of Chvátal and Szemerédi \cite{ChvatSzem81ErdStone}. A weaker, bipartite, asymmetric
 version of it was used to prove that $r_k(n)=o(n)$.} Since then, the lemma has emerged as a fundamental tool in Graph Theory: it has many applications
 in Extremal Graph Theory, in the area of `Property Testing' in computer science, combinatorial Number Theory, etc. \dots Tao described the lemma as a
 result in probability. Our goal is to point out that Szemerédi's lemma can be thought of as a result in analysis. We show three different analytic
 interpretations. The first one is a general statement about approximating elements in Hilbert spaces which implies many different versions of the
 Regularity Lemma, and also potentially other approximation results. The second one presents the Regularity Lemma as the compactness of an important metric
 space of 2-variable functions. \dots The third analytic interpretation shows the connection between a weak version of the regularity lemma and covering by
 small diameter sets, i.e., dimensionality. \dots We describe two applications of this third version: \dots and an algorithm that constructs the weak
 Szemerédi partition as Voronoi cells in a metric space.''}
\end{quote}

Actually, it is surprising that such a short paper can describe such an involved situation in such a compact way, also including the proofs.

\subsection{Versions of the Regularity Lemma}\label{Versions-RL}

There are many versions of the Regularity Lemma. Below we list some of them, with very short descriptions. Many of them are difficult to invent but have you invented them, their proofs are often very similar to the original proof. (In case of Hypergraph Regularity Lemmas the situation is completely different.)

We have already described the Regularity Lemma. 

\subsub Coloured version.// There is an easy generalization of the Regularity Lemma in which the edges of a graph $G_n$ are $r$-coloured, for some fixed $r$, so we have $r$ edge-disjoint graphs, and we wish to find a vertex-partition which satisfies the Regularity Lemma in all the colours, simultaneously. This is possible and often needed, e.g., in the Erdős-Hajnal-Sós-Szemerédi extension of Theorem \ref{SzemK4Th} (in \cite{ErdHajSosSzem83}), and more generally, this was used in Ramsey type theorems and Ramsey-Turán type theorems, and later in many similar cases.

\subsub Weak Regularity Lemma.// There is an important weakening of the Regularity Lemma, namely the Frieze-Kannan Weak Regularity Lemma \cite{FriezeKannan99Weak}, see also \cite{FriezeKannan96Survey,FriezeKannan99Simple} and the
Frieze-Kannan-Vempala approach \cite{FriezeKanVempa04}. The difference between the Regularity Lemma and the Weak Regularity Lemma is that the later one ensures $\eps$-regularity only for ``much larger subsets'', and (therefore) needs much
fewer clusters.  Here the Weak $\eps$-regularity means that given a partition $U_1,\dots,U_k$, for a subset $X\subset V(G_n)$ we hope to have
$$\cE(X):=\sum d(U_i,U_j)\cdot |U_i\cap X||U_j\cap X|$$ edges in $G[X]$, so we conclude that $\cE(X)$ is close to $e(X)$.

\minisepa

While the original algorithm of Frieze and Kannan is randomized, 
Dellamonica, Kalyanasundaram, Martin, Rödl, and Shapira \cite{DellaKalyMartinRodlShapi12} provided a deterministic $O(n^2)$ algorithm, analogous to Theorem \ref{AlonDukeYus2Th}, to find the Frieze-Kannan Partition. 

This regularity lemma is more connected to Statistics\footnote{Principal component analysis, see e.g. Frieze, Kannan, Vempala, and Drineas \cite{FriezeKanVempa04,DrinFriezeKanVemp04Singu}.}, and in many cases, where one can apply a regularity lemma for an existence proof, the algorithmic versions of the regularity lemmas provide algorithms in these applications too.
 
As we wrote, these algorithms are slow because of the very large number of clusters.\footnote{There are many results showing that the number of clusters must be very large. The first such result is due to Gowers \cite{Gowers97Tower}. }
but in many such cases the {\bf Algorithmic Weak Regularity Lemma} can also be used, and then it provides a much faster algorithm, basically because it requires much fewer classes in the $\eps$-regular partition. 

One should remark that the Frieze-Kannan Regularity Lemma can be iterated and then it provides a proof of the original Szemerédi Regularity Lemma.

On the other end, there is the {\bf Strong Regularity Lemma} of Alon, Fischer, Krivelevich, and M. Szegedy \cite{AlonFischKrivSzege99}.  The advantage of the Strong Regularity Lemma is that it can be applied in several cases where the
original Regularity Lemma is not enough, primarily when we are interested in induced subgraphs.

We shall not formulate this strong lemma but include an explanation of it, from Alon-Shapira \cite{AlonShapira05Monotone} (with a slight simplification). Alon and Shapira write: 

\begin{quote}{\dots
 This lemma can be considered a variant of the standard Regularity Lemma, where one can use a function that defines $\eps>0$ as a function of the size of the
 partition, rather then having to use a fixed $\eps$ as in Lemma 2.2.}
\end{quote}

\noindent
Large part of this is described in the paper of Lovász and B. Szegedy \cite{LovSzege07Anal}.

\subsection{Regularity Lemma for sparse graphs}

The Regularity Lemma can be used in many cases but has several important limitations. 

(a) Because of the large threshold $n_0$, one cannot combine it with computer programs, checking the small cases. In other words, it is a theoretical result but it cannot be used in practice.

(b) The most serious limitation is that we can apply it for embedding $H$ into $G$ only if the degrees in $H$ are (basically) bounded. 

(c) Another one is that it can be applied only to dense graphs $G_n$. This problem is partly solved by the Sparse Regularity Lemma, established by Kohayakawa \cite{Koha97SparseRegu}, and Rödl, see also \cite{KohayaRodl03Quasi}.

For sparse graph sequences, i.e. when $e(G_n)=o(n^2)$, the original Regularity Lemma is trivial but does not give any information.  Having a bipartite subgraph $H[U,V]\subseteq G_n$, consider the following ``rescaled'' density:
\beq{Modified}d_{H,p}(U,V):={e(U,V)\over p\cdot|U||V|}.\eeq If $p>0$ is very small, e.g., $p:=n^{-2/3}$, then the condition
$$|d_{H,p}(X,Y)-d_{H,p}(U,V)|<\eps$$ does not say anything without $1/p$ in \eqref{Modified}, but with $1/p$ it is a reasonable and strong restriction for sparse graphs. The Sparse Regularity lemma says that for ``nice'' graphs $H_n$ there is a partition $V_1,\dots,V_\nu$, described in the Regularity Lemma, even if we use this stronger regularity requirement \eqref{Modified}. Which are the ``nice'' graphs?

\medskip

(d) Sparse regularity lemmas are well applicable when the graph $G_n$ to be approximated by generalized quasi-random graphs does not contain subgraphs whose
density is much above the edge-density of $G_n$, e.g, for some bound $b$ and the edge density $p=e(G_n)/{n\choose 2}$, \beq{SparseRegCondit} E(V_1,V_2)<bp
|V_1||V_2|\Text{if}|V_1|,|V_2|>\eta n.\eeq This is the situation, e.g., when we consider (non-random) subgraphs of random graphs, see
Section~\ref{RandomUnivers}.

The sparse Regularity Lemma sounds the same as the original Regularity Lemma, with two differences: we use the modified uniformity: \eqref{Modified}, instead of \eqref{OriginalUnif}, and the extra condition \eqref{SparseRegCondit}, for some fixed $b$. %\MargoS{Ujra ellenorizni.}

One of the latest developments in this area is a sparse version of the Regularity Lemma due to Alex Scott \cite{Scott11SparseRegu}. Scott succeeded in
eliminating the extra condition on the sparse graph that it had no (relatively) high density subgraphs. However, this had some price, discussed by Scott in
Section 4 of \cite{Scott11SparseRegu}. This new sparse Regularity Lemma was used in several cases, e.g. by P. Allen, P. Keevash, B. Sudakov, and
J. Verstraëte, in \cite{AllenKeevSudakVers14}.
		
\begin{remark} 
  As we wrote, the Sparse Regularity Lemma can be used basically if $G_n$ does not contain subgraphs much denser than the whole graph. Gerke and Steger wrote an important survey about it and about its applicability \cite{GerkeSteg05SparseRegu}, see also \cite{GerkeKohaRodlSteg07}. Its applicability is discussed (among others) in the paper of Conlon, Gowers, Samotij, and Schacht \cite{ConGowSamoSchacht14}, with its connection to the Kohayakawa-Łuczak-Rödl conjecture \cite{KohaLuczRodl97K4free}. We also warmly recommended the paper of Conlon, Fox, and Yufei Zhao \cite{ConFoxZhao14SparsePseudo}.
\end{remark}

\subsection{Quasi-Random Graph Sequences}

Quasi-random sequences are very important, e.g., in Theoretical Computer Science, and also very interesting, for their own sake. In Graph Theory they emerged in connection with some Ramsey problems. The first detailed, pioneering results in
this direction are due to A. Thomason, see e.g. \cite{Thomason87Poznan,Thomason87Pseudo}. He was motivated by some Ramsey Problems.

 Chung, Graham and Wilson \cite{ChungGrahWilson89} developed a theory in which {\em six properties of random graphs} were formulated which are equivalent
 for any infinite sequence $(G_n)$ of graphs.  The graphs having these properties are called {\em quasi-random}.

 Quasi-randomness exists in other universes as well, e.g., there exist quasi-random subsets of integers, groups, see Gowers \cite{GowersT08QuasiGr},
 tournaments, see Chung and Graham \cite{ChungGrah92QuasiZ,ChungGrah91Tour}, of real numbers, digraphs \cite{FerberKroneLong15RandHamil}, of hypergraphs,
 e.g., Chung \cite{Chung90QuasiHyper,Chung91Regu,Chung12QuasiRevis}, Rödl and Kohayakawa \cite{KohayaRodl03Quasi}, permutations, see J. Cooper
 \cite{Cooper04Permu},Král and Pikhurko \cite{KralPikhu13QuasiPermu}, and in many other settings\dots

Quasi-randomness and the Regularity Lemmas are very strongly connected. This was first established in a paper of Simonovits and Sós \cite{SimSos91Parti}:

\begin{theorem}[Simonovits, Sós \ev(1991) \cite{SimSos91Parti}]\label{SimSos91Quasi} A sequence of graphs $(G_n)$ is $p$-quasi-random iff for every $\kappa>0$ and $\eps>0$, there exist two thresholds $k_0(\eps,\kappa)$ and
 $n_0(\eps,\kappa)$ such that for $n>n_0(\eps,\kappa)$ $G_n$ has an $\eps$-Regular Partition where all the pairs $(V_i,V_j)$ are $\eps$-regular with densities between $p-\eps<d(V_i,V_j)<p+\eps$ and $\kappa<k<k_0(\eps,\kappa)$.
\end{theorem}

\begin{remark} Actually, a slightly stronger theorem holds. On the one hand, we may allow $\eps{k\choose2}$ exceptional pairs $(V_i,V_j)$ to ensure $p$-quasi-randomness. On the other hand, if $(G_n)$ is $p$-quasi-random, we can find a partition where there are no exceptional pairs.
\end{remark}

The corresponding generalization for sparse graphs was proved by Kohayakawa and Rödl. For a longer and detailed survey see their paper \cite{KohayaRodl03Quasi}. 

\Sc For~further~related~results see, e.g., Simonovits and Sós \cite{SimSos97HeredNNI,SimSos03Hered}, Skokan and Thoma \cite{SkokThoma04Bipa}, Yuster \cite{Yuster08Quasi}, Shapira and Yuster \cite{ShapiraYuster10Ind}, Gowers on Counting Lemma \cite{GowersT06-3Unif}, on quasi-random groups \cite{GowersT08QuasiGr}, and also Janson \cite{Janson11Quasi}, Janson and Sós \cite{JansonSos15Quasi} on the connections to Graph Limits.

A generalization of the notion of Quasi-random graphs was investigated by Lovász and Sós \cite{LovSos08GenQuasi} which corresponds to generalized random matrix-graphs.

\subsection{Blow-up Lemma}

Several results exist about embedding {\em spanning} subgraphs into dense graphs. Many of the proofs use a relatively new and very powerful tool, called \Sc Blow-up Lemma. Here we describe it in a fairly concise way. The Blow-up Lemma is
mostly used to embed a {\em bounded degree} graph $H$ into a graph $G$ as a {\em spanning subgraph}. The reader is also referred to the excellent ``early'' survey of Komlós \cite{Kom99BlowSurvey} (explaining a lot of important background
details about the Regularity Lemma and the Blow-Up lemma, and how to use them) or to the surveys of Komlós and Simonovits, \cite{KomSim96SzemRegu}, Komlós, Simonovits, Shokoufandeh, and Szemerédi, \cite{KomShokoSimSzem00}. The ``Doctor of Sciences'' Thesis of Sárközy \cite{Sarko07Doktori} is also an excellent source in this area. For some newer results on the topic see e.g., Rödl and Ruciński \cite{RodlRuc99Blow}, Keevash \cite{Keev11HyperBlow} who extended the method to hypergraphs, Sárközy \cite{Sarko14QuantiBlow}, Böttcher, Kohayakawa, Taraz, and Würfl \cite{BottKohaTarazWurfl15Blow} extended it to $d$-degenerate graphs.\footnote{They call it $d$-arrangable.} A long survey of Allen, Böttcher, Han,
Kohayakawa, and Person \cite{AllenBottHanKohaPers16A-Blow} discusses several features of the Blow-up Lemma applied to random and pseudo-random graphs.  Recently Allen, Böttcher, Hàn, Kohayakawa, and Person extended the Blow-up lemma to
sparse graphs \cite{AllenBottHanKohaPers16A-Blow}. 

\minisepa

We start with a definition.

\begin{definition}[$(\eps,\de)$-super-regular pair]
Let $G$ be a graph, $U,W\subseteq V(G)$ be two disjoint vertex sets, $|U|=|W|$. The vertex-set pair $(U,W)$ is $(\eps,\delta)$-super-regular if it is $\eps$-regular and $\mindeg(G[U,W])\ge \de |U|$. 
\end{definition}

The Blow-up Lemma asserts that $(\eps,\de)$-regular pairs behave as complete bipartite graphs from the point of view of embedding {\em bounded degree} subgraphs. In other words, for every large $\Delta>0$ and small $\delta>0$ there exist an $\eps>0$ such that if in the Cluster graph $H_\nu$ the min-degree condition also holds and we replace the $(\eps,\delta)$-regular pairs $\UiUj$ by complete graphs and then we can embed (the bounded degree) $H_n$ into this new graph $\Ti G_n$, then we can embed it into the original $G_n$ as well. The low degree vertices of $G_n$ could cause problems. Therefore, for embedding spanning subgraphs, one needs all degrees of the host graph large. That's why using regular pairs is not sufficient any more, we need {\em super-regular pairs}. Again, the Blow-up Lemma plays crucial role in embedding {\em spanning} graphs $H_n$ into $G_n$.

The difficulty is in embedding the ``last few'' vertices. The original proof of the Blow-up Lemma starts with a probabilistic greedy algorithm, and then uses a König-Hall argument to complete the embedding. The proof is quite involved. 

\begin{theorem}[Blow-up Lemma, Komlós-Sárközy-Szemerédi 1994 \cite{KomSarkoSzem97BlowUp}]
 Given a graph $H_\nu$ of order $\nu$ and two positive parameters $\delta,\Delta$, there exists an $\eps>0$ such that if $n_1,n_2,\dots,n_r$ are arbitrary positive integers and we replace the vertices of $H_\nu$ with pairwise disjoint sets
 $V_1,V_2,\ldots,V_\nu$ of sizes $n_1,n_2,\dots,n_\nu$, and construct two graphs on the same vertex-set $V=\bigcup V_i$ so that 

(i) the first graph $H_\nu(n_1,\dots,n_\nu)$ is obtained by replacing each edge $\{v_i,v_j\}$ of $H_\nu$ with the complete bipartite
 graph $K(n_i,n_j)$ between the corresponding vertex-sets $V_i$ and $V_j$, 

(ii) and second, much sparser graph $H^*_\nu(n_1,\dots,n_\nu)$ is obtained by replacing each $\{v_i,v_j\}$ with an $(\eps,\delta)$-super-regular pair between $V_i$ and $V_j$,

\noindent
then if a graph $L$ with $\maxdeg(L)\le\Delta$ is embeddable into $H_\nu(n_1,\dots,n_\nu)$ then it is embeddable into the much sparser $H^*_\nu(n_1,\dots,n_\nu)$ as well.
\end{theorem}

The Blow-up Lemma has several different proofs, e.g., Komlós, Sárközy, and Szemerédi first gave a randomized embedding \cite{KomSarkoSzem97BlowUp}, and then they gave a derandomized version \cite{KomSarkoSzem98BlowAlgo} as well.
Other proofs were given, e.g., by Rödl, Ruciński \cite{RodlRuc99Blow}, 
Rödl, Ruciński and Taraz \cite{RodlRucWagner98},
and Rödl, Ruciński and Wagner \cite{RodlRucWagner98}.
%%We have also mentioned the Komlós survey \cite{Kom99BlowSurvey} on this topic, from 1999, however, many related result were proved only after it.

\begin{remark} 
  G.N. Sárközy gave a very detailed version of the proof of the Blow-up lemma, \cite{Sarko14QuantiBlow}, where he calculated all the related details very carefully in order that he and Grinshpun could use it in their later work
  \cite{GrinshpunSarko16Mono}, Theorems \ref{GrinSarkoMono},\ref{GrinSarkoMonoB}: without this they could prove only a weaker result.
\end{remark}

\subsection{Regularity lemma in Geometry}

Until now we restricted our consideration to applications of the Regularity Lemma to embed sparse graphs into dense graphs. The Regularity Lemma has several
applications in other fields as well, and beside this, versions tailored to these other fields. Here we mention only a few such versions, very briefly.

A theorem of Green and Tao \cite{GreenTao10Arith} is a good example of this. They prove a so called {\sc arithmetic regularity lemma} that can be applied in
{\sc Additive Combinatorics}, in several problems similar to Szemerédi theorem on arithmetic progressions. However, here we wish to speak of Geometry.

There are several cases where we restrict our consideration to some particular graphs, e.g., to graphs coming from geometry. In this case one may hope for much better estimates in some cases than for arbitrary graphs. A whole theory was built up around such problems, see, e.g., Erdős \cite{Erd1946-Unit,Erd1967-13}, a survey Szemerédi \cite{Szemer16Unit}, Szemerédi and Trotter \cite{SzemTrott83Geom}, Pach and Sharir \cite{PachSharir98Incid} \dots or the book of Pach and Agarwal \cite{PachAgar95Book}.  We remark here only that the Regularity Lemma also has some strengthened form, see, e.g., the improvement of some results of Fox, Gromov, Lafforgue, Naor, and Pach \cite{FoxGromLaffNaorPach} by Fox, Pach and Suk \cite{FoxPachSuk15}\dots

So the Regularity Lemma can be applied in Geometry in many cases and the fact that we apply it to geometric situation implies that in many cases the
connection between the clusters will be a (basically) complete connection, or very few edges, instead of having randomlike connections. This is not so
surprising, since many of the geometric relations are described by polynomials, or analytic functions, (?) and in these cases, if we have ``many'' 0's of a
polynomial (of several variables) then the corresponding polynomial must vanish everywhere.\footnote{One form of this is expressed in the Combinatorial Nullstellensatz of Alon \cite{Alon99NullStel}.}
 It seems that in most cases where Regularity Lemma is used in Geometry, its use can be eliminated. 

 The interested reader is suggested to read, e.g., the survey of J. Pach \cite{Pach14ICM}, or his book with Agarwal \cite{PachAgar95Book}, the papers of
 Alon, Pach, Pinchasi, Radoičić, and Sharir \cite{AlonPachPinchaRadShar05}, on \Sc semi-algebraic sets, or the paper of Fox, Gromov, Lafforgue, Naor, and
 Pach \cite{FoxGromLaffNaorPach}, or the paper of Pach and Solymosi, \cite{PachSolymo01Crossing}, or of J. Pach \cite{Pach98Tverberg}. The result of
 Karasev, Kynčl, Paták, Patáková, and Martin Tancer \cite{KarasevKyncl15} is also related to this topic.

\Section{With or without Regularity Lemma?}{RemovingReguLem}

Regularity Lemma is one of the most effective, most efficient lemmas in Extremal Graph Theory. In several cases we can prove a result with the Regularity
Lemma, but later we find out that it can easily be proven without the Regularity Lemma as well. So it is natural to ask if it is worth getting rid of the
application of Regularity Lemma, if we can. The answer is not so simple. We should mention a disadvantage and two advantages of using the Regularity Lemma.

The {\em disadvantage} is that it can be applied only to very large graphs. This means that ``in practice it is of no use''. This bothers some people, e.g., those who wish to find out, -- sometimes with the help of computers -- the truth
for the smaller values as well. Many of us feel this unimportant, some others definitely prefer eliminating the use of Regularity Lemma if it is possible. First we quote an opinion against the usage of the Regularity Lemma. In
Subsection~\ref{DiamCritical} we discussed the beautiful conjecture of Murty and Simon on the maximum number of edges a diameter-2-critical graph can have. We have mentioned that Füredi \cite{Fure92DiamCrit} solved this problem in 1992,
using the Ruzsa-Szemerédi theorem that $f(n,6,3)=o(n^2)$. Much later, in 2015, a survey \cite{HayHenMerYeo15Murty} on the topic described this area very active and wrote:

\begin{quote}{ ``The most significant contribution to date is an astonishing asymptotic result due to Füredi (1992) who proved that the conjecture is true for large $n$, that is, for $n > n_0$ where $n_0$ is a tower of 2’s of height about
    $10^{14}$. As remarked by Madden (1999), ` $n_0$ is an inconceivably (and inconveniently) large number: it is a tower of 2’s of height approximately $10^{14}$.' Since, for practical purposes, we are usually interested in graphs which
    are smaller than this, further investigation is warranted\dots ''}
\end{quote}

First of all, it is not clear if it is correct to call Füredi theorem an ``asymptotic result''. The {\em advantage} of applying the Regularity Lemma is that it
often provides a proof where we have no other proofs, and, in other cases, a much more transparent proof than the proof without it. (Thus for example, the
beautiful theorem of Erdős, Kleitman, and Rothschild \cite{ErdKleiRoth76Lincei} was proved originally without the Regularity Lemma, but for many of us the
Regularity Lemma provides a more transparent proof.)

Perhaps one of the first cases where we met a situation where the Regularity Lemma could be eliminated was a paper of Bollobás, Erdős, Simonovits, and
Szemerédi, \cite{BolloErdSimSzem77} discussing several distinct extremal problems, and one of them was which could be considered as one germ of \Sc Property~Testing:

\begin{problem}[Erdős]
 Is it true that there exists a constant $M_\eps$ such that 
if one cannot delete $\eps n^2$ edges from $G_n$ to make it bipartite then
 it has an odd cycle
 $C_\ell$ with $\ell<M_\eps$?
\end{problem}

The answer was

\begin{theorem}[Bollobás, Erdős, Simonovits, and Szemerédi \cite{BolloErdSimSzem77}]\label{BolloErdSimSzemProp} YES, 
  if one cannot delete $\eps n^2$ edges from $G_n$ to make it bipartite then one can find an odd cycle $C_\ell$ with $\ell<{1\over\eps}$.
\end{theorem}

In \cite{BolloErdSimSzem77} we gave two distinct proofs of this theorem, one using the Regularity Lemma and another one, without the Regularity Lemma.  Interestingly enough, these questions later became very central and important in the
theory of Property Testing, but there, in the works of Alon, Krivelevich, Shapira, and others, (see e.g., \cite{AlonFischKrivSzege99}) it turned out that such property testing results depend primarily on whether one can apply the Regularity Lemma or not. Just to illustrate this, we mention from the many similar results the paper of N. Alon, E.~Fischer, I. Newman and A. Shapira \cite{AlonFischNewShap06Test}, the title of which is ``A combinatorial characterization of the
testable graph properties: It's all about regularity''.\footnote{Two remarks should be made here: (a) Originally Property Testing was somewhat different, see e.g. O. Goldreich, S. Goldwasser and D. Ron \cite{GoldGoldRon98}. (b) The theory of graph limits also has a part investigating property testing, see e.g., \cite{BorgsChayLovSosVesz06Testing,BorgsChayLovSosVesz08Conv},\dots, \cite{LovSzege10PropTesting}.}

\begin{remark} Duke and Rödl \cite{DukeRodl85} extended Theorem \ref{BolloErdSimSzemProp} to higher chromatic number, answering another question of Erdős, see also the ICM lecture of Rödl \cite{Rodl14Seoul}, and also the result of Alon
  and Shapira on property testing \cite{AlonShapira05Monotone}. (Rödl: ``Further refinement was given by Austin and Tao'' \cite{AustinTao10TestHeredit}.)
\end{remark}

So we emphasise again that there are several cases where certain results can be proved with and without the Regularity Lemma, and the proof with the
Regularity Lemma may be much more transparent, however, the obtained constants are much worse. Before proceeding we formulate a meta-conjecture about the
``elimination''.

\begin{MetaConj}[Simonovits]\sl 
  The use of the Regularity Lemma can be eliminated from those proofs where \dori (i) the conjectured extremal structures are ``generalized random graphs'' with a fixed number of classes and densities 0 and 1, and 
\dori (ii) at least one of the densities is 0 and one of them is 1.
\end{MetaConj}

One has to be very careful with this -- otherwise informative - Meta-Conjecture:

(a) First of all, mathematically it is not quite well defined, what do we mean by ``eliminating the Regularity Lemma''.

(b) Further, without (ii) the Ruzsa-Szemerédi Theorem could be regarded as a ``counter-example''.\footnote{Actually, J. Fox eliminated the application of the Regularity Lemma from the proof of the Triangle Removal Lemma, see Section
  \ref{RemovalS1} or \cite{Fox11NewRemo}.}

(c) It is not well defined if using graph limits we regard as elimination of the Regularity Lemma or not? 

(d) The first two-three results which we like proving nowadays to illustrate the usage of the Regularity Lemma, e.g. the Ramsey-Turán estimate for $K_4$ (Thm \ref{SzemK4Th}) and Ruzsa-Szemerédi Theorem, originally were proved using some
weaker forms of the Regularity Lemma.

In several cases originally the Regularity Lemma was used to obtain some results, but then it was {\em easily} eliminated. Such examples are Erdős and Simonovits \cite{ErdSim80Geom}, or Pach and Solymosi \cite{PachSolymo01Crossing}, on
Geometric Graphs, or results in the paper of Erdős, S.B. Rao, Simonovits, and Sós, \cite{ErdRaoSimSos92}. Often in the published versions we do not even find the traces of the original proof with Regularity Lemma, anymore \dots

One interesting case of this discussion is the proof of

\begin{conjecture}[Lehel's conjecture \cite{Ayel79Lehel}]\label{LehConj}
If we 2-colour the edges of $K_n$, then $V(K_n)$ can be covered by two vertex-disjoint monochromatic cycles of distinct colours.
\end{conjecture}

\begin{remark} 
  The first reference to Conjecture \ref{LehConj} can be found in the PhD thesis of Ayel \cite{Ayel79Lehel}. The conjecture was first proved by Łuczak, Rödl, and Szemerédi \cite{LuczRodlSzem98Partit}, using the Regularity Lemma. Of course, this worked only for very large values of $n$. The Regularity Lemma type arguments were eliminated by Peter Allen \cite{Allen08Covering}. The difference between the two proofs is that Allen covers $K_n$ by monochromatic cliques, using Ramsey's theorem, instead of using Regularity type arguments. Hence the threshold in the first proof is very large, and in Allen's proof it is ``only'' $2^{18000}$. Finally, surprisingly, Bessy and Thomassé \cite{BessyThomasse10} found a simple and short proof of the conjecture, without using the Regularity Lemma, or any deep tool, and which worked for all $n$.
\end{remark}

\minisepa

Conjecture \ref{LehConj} was extended by Gyárfás to any number of colours:

\begin{conjecture}[Gyárfás, \cite{Gyar89Covering,ErdGyarfPyb91}]\label{GyarfCycleParti} If the edges of $K_n$ are $r$-coloured, then $V(K_n)$ can be covered by $p(r)=r$ vertex-disjoint monochromatic cycles (where $K_1,K_2$ are also considered as cycles).
\end{conjecture}

\begin{remark}\label{PokrovskiyR}
  This was ``slightly'' disproved for $k\ge3$ by Pokrovskiy \cite{Pokro14Counter}. Here ``slightly'' means that in his counterexample there is one vertex which could not be covered, however, $p(r)=r+1$ is still possible.  We shall return to the Gyárfás conjecture (often called Gyárfás-Lehel conjecture) in Subsection \ref{VertexPartitionS}.
\end{remark}

\begin{remark} 
  The results of Bessy and Thomassé and of Erdős, Gyárfás, and Pyber, (from Section \ref{Absorb}) were extended to \Sc local~$r$-colouring by Conlon and Maya Stein \cite{ConStein16Mono}, where local $r$-colouring means that each vertex is adjacent only to at most $r$ distinct colours, but the total number of colours may be much larger.
\end{remark}

\subsection{Without Regularity Lemma}

There are several cases where eliminating the use of Regularity Lemma from the proof is or would be important. Above we discussed this and listed such cases. Here we mention two further cases. Fox gave a proof of the Removal Lemma without using the Regularity Lemma, in \cite{Fox11NewRemo} (slightly simplified in the beautiful survey of Conlon and Fox \cite{ConFox12Remo}). This improves several estimates in some related cases. Conlon, Fox and Sudakov \cite{ConlonFoxSudak18Hered} recently removed using the Regularity Lemma from the proof of a theorem of Simonovits and Sós \cite{SimSos97HeredNNI}, which helped to understand the situation better.

\subsection{Embedding spanning or almost-spanning trees}\label{AlmostSpanningTrees}

Originally we planned to write -- among others, -- about our results on tree embeddings: about the solutions of the Erdős-Sós conjecture and the Komlós-Sós
conjecture.  However, they are described in \cite{FureSim13Degen} and in
\cite{HladPiguSim15Electro},\cite{HladKomPigu12}-\cite{HladKomPigu17-4}, respectively, and we shall return to these topics elsewhere.

There are many results where we try to embed large or actually spanning trees into a graph $G_n$. We mention a few of them. Some of them describe a
case when $G_n$ is a random or randomlike graph, e.g., pseudo-random, expanding,\dots If we know something for random graphs, that often can (easily?) be
extended to these cases: quasi-random, pseudo-random, or expander graphs. The ``Resilience results'' of Subsections \ref{LargeTrees} and
\ref{PseudRandomUniv} were of this type. One of the early results on expander graphs was

\begin{theorem}[Friedman and Pippenger \cite{FriedmanPippen87Expan}]
If for every $X\subseteq V(G_n)$, with $|X|\le 2k-2$, $$|\Gamma(X)|\ge (d+1)|X|,$$ then $G_n$ contains all trees $\Tk$ with $\maxdeg(\Tk)\le d$.
\end{theorem}

\begin{remark} The Friedman-Pippenger theorem ``embeds'' only relatively small trees. It was extended by Balogh, Csaba, Pei, and Samotij \cite{BalCsabaPeiSamo10Tree}, where a result of Haxell \cite{Hax97MonoCyc} was simplified, and then
  used. This guaranteed embedding almost spanning trees into ``expanding graphs''. We skip the precise formulation and the details, because they may look technical at first sight.
 \end{remark}

 For some earlier related works see Alon-Chung \cite{AlonChung88Tolerant}, and Beck \cite{Beck83Size}.

\minisepa

Now we consider a conjecture of Bollobás \cite{Bollo78ExtreBook}, on embedding bounded degree trees, proved by

\begin{theorem}[Komlós, G.N. Sárközy, and Szemerédi \cite{KomSarkoSzem95Pack,KomSarkoSzemer01SpanningTree}]\label{KomSarkoSzem-BolloTh}
 For every $\eps>0$ and $\Delta>0$, there exists an $n_0$ for which, if $T_n$ is a tree on $n$ vertices with $\maxdeg(T_n)\le\Delta$, and
 $G_n$ is a graph on $n$ vertices with \beq{AKS-Bollo}\mindeg(G_n) \geq {n\over2}+\eps n,\eeq then $T_n\subseteq G_n$, assuming that $n>n_0(\Delta)$.
\end{theorem}

This was improved by

\begin{theorem}[Csaba, Levitt, Nagy-György, and Szemerédi \cite{CsabaLevittSzem10Trees}]

 (a) For any constant $\Delta>0$ there exists a constant $c_\Delta>0$ such that if $T_n$ is a tree on $n$ vertices with $\maxdeg(T_n)\le\Delta$, and $G_n$ is a graph on $n$ vertices with 
\beq{CsabaLevittLow}\mindeg(G_n) \geq {n\over2}+c_\De\log n,\eeq
 then $T_n\subseteq G_n$, assuming that $n>n_0(\Delta)$. 
\dori (b) There exist infinitely many graphs $G_n$ with $\mindeg(G_n)\ge \half n+\reci 17 \log n$ not containing the complete ternary tree $T_n^{[3]}$.
\end{theorem}

So the bound in \eqref{CsabaLevittLow} is tight. The proofs in \cite{KomSarkoSzem95Pack,KomSarkoSzemer01SpanningTree} used the Regularity Lemma and the Blow-up lemma, while \cite{CsabaLevittSzem10Trees} did not use them, and provided
smaller $n_0$ and a sharper theorem. B. Csaba also extended the above results to \Sc well-separable graphs:

\begin{definition}
An infinite graph sequence $(H_n)$ is \Sc well-separable, if one can delete $o(n)$ vertices of $H_n$ so that each connected component of the remaining graph has $o(n)$ vertices, as $n\to\infty$.
\end{definition}

\begin{theorem}[Csaba \cite{Csaba08Separ}]\label{CsabaWSepa}
For every $\eps,\Delta>0$, there exists an $n_0=n_0(\eps,\Delta)$ such that if $(H_n)$ is well-separable, $n>n_0$, and $\maxdeg(H_n)\le \Delta$, and
\beq{CsabaSepa}\mindeg(G_n)>\left(1 -\reci2(\chi(H)-1) +\eps\right)n,\eeq then $H_n$ can be embedded into $G_n$.
\end{theorem}

For trees (or, more generally, for bipartite graphs $H_n$) \eqref{CsabaSepa} reduces to $\mindeg(G_n)\ge \half n+\eps n$: Theorem \ref{CsabaWSepa} is a
generalization of Theorem \ref{KomSarkoSzem-BolloTh}. Another version is where we assume that $G_n$ is bipartite, see Csaba \cite{Csaba07ReguSpanning}.

\subsection{Pósa-Seymour conjecture}\label{PosaSeyRevisit}

Speaking of extremal problems, we could consider problems where we ask one of the following questions:
\begin{enumerate} \dense
\item how large $e(G_n)$ ensures a property $\cP$?
\item \label{Dir} how large $\mindeg(G_n)$ ensures $\cP$?
\item \label{Ore} which (Ore type) degree sum conditions $d(x)+d(y)\ge f_O(n,\cP)$ ensure $\cP$, where we assume this only for independent vertices $x,y$; 
\item \label{DegSequ} Given a graph $G_n$ with the degree sequence $d_1,d_2,\dots,d_n$, does it ensure $\cP$?
\end{enumerate}

In the next part we consider two of these questions: \ref{Dir}, called Dirac type problems, and \ref{Ore} called Ore type problems.

\minisepa

\BalAbraCapMed{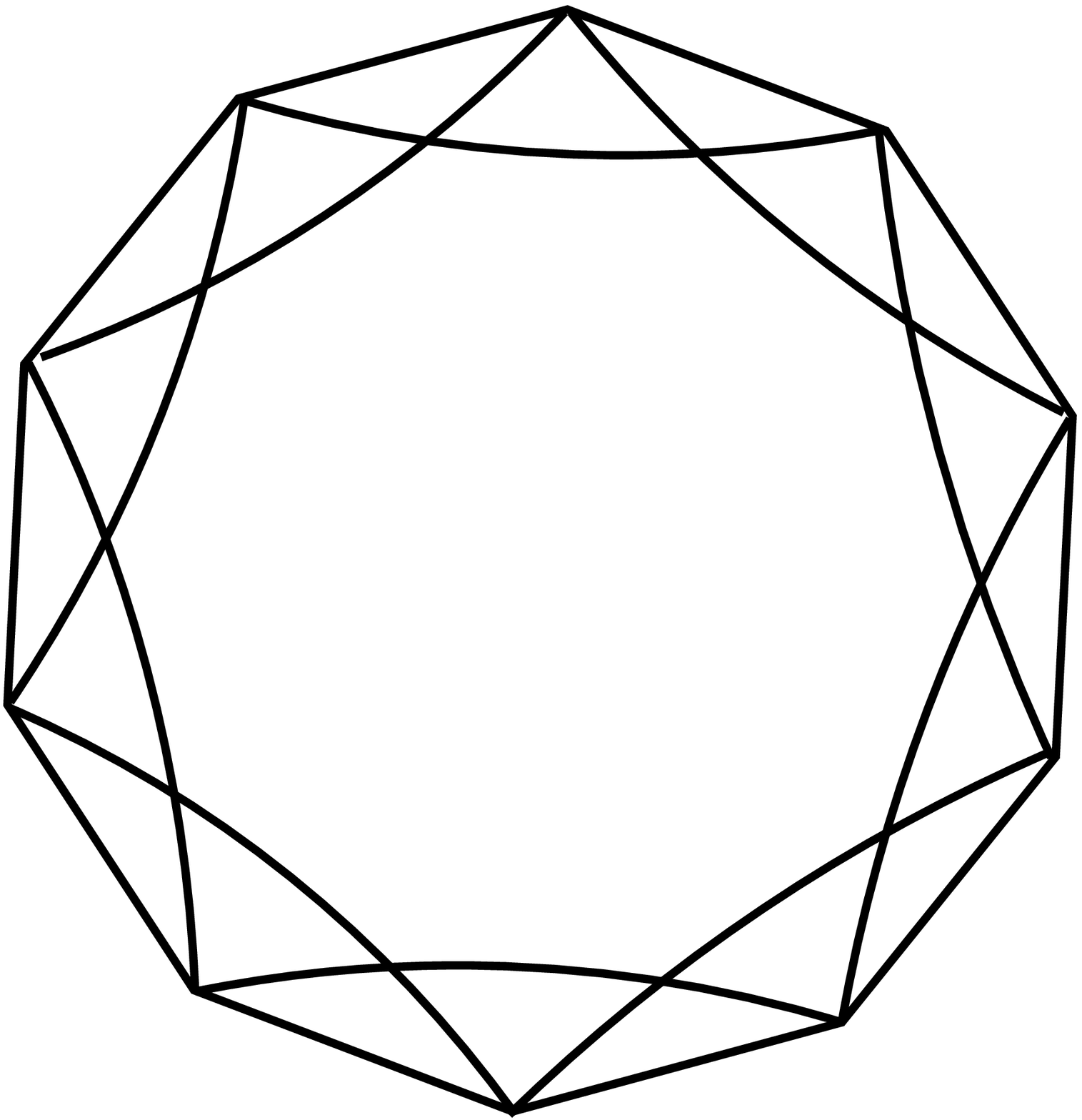}{Square of a cycle}{17}
Hamiltonicity of graphs is a central problem in graph theory. There are many results of type \ref{DegSequ}, where some conditions on the degrees ensure the Hamiltonicity. In the Introduction we formulated one of the first such results, Theorem~\ref{DiracX}: 

\Proclaim Dirac Theorem. If $\mindeg(G_n)\ge n/2$, and $n\ge3$, then $G_n$ contains a Hamiltonian cycle.

As we have mentioned, this is sharp. We shall go into two distinct directions from Dirac's Theorem: here we shall consider some generalizations for simple graphs, while in Section \ref{HyperLargeExcludedS} we shall discuss some hypergraph extensions.  A natural question analogous to Dirac's theorem was asked by Pósa (see Erdős \cite{Erd65Smole} in 1965). The reader is reminded that the $k\th$ power $L:=M^k$ of a graph $M$ is obtained from $M$ by joining all the pairs of vertices $x\ne y$ having distance $\rho_M(x,y)\le k$.

\begin{conjecture}[Pósa]\label{posa}
 If for a graph $G_n$ $\mindeg(G_n)\geq\frac23n$, then $G_n$ contains the square of a Hamiltonian %% \\ \phantom{~} \hskip 27mm 
cycle. 
\end{conjecture}

This was generalized by Seymour in 1973:

\begin{conjecture}[Seymour \cite{Seymour74Conj}]\label{sey}
Let $G_n$ be a graph on $n$ vertices. If $\mindeg(G_n)\geq\frac k{k+1}n$, then $G_n$ contains the $k^{th}$ power of a Hamiltonian cycle.
\end{conjecture}

For $k=1$, this is just Dirac's theorem, for $k=2$ the Pósa conjecture. The validity of the general conjecture implies the notoriously hard
Hajnal-Szemerédi theorem (i.e. Theorem \ref{HajSzemCorrTh}).\footnote{Actually, this was the motivation for Pósa.}

\begin{remark}
 Observe that for $\ell\ge k+1$, we have $K_{k+1}=P_{k+1}^k \subseteq P_{\ell}^k$. Hence $T_{n,k}$ does not contain $P_\ell^k$. On the other hand,
 $P_n^k\subset T_{n,k+1}$. This provides some further motivation for the above conjectures.
\end{remark}

In the earlier parts we mostly considered embedding problems where the graph $H_m$ to be embedded into the ``host graph'' $G_n$ had noticeably fewer
vertices than $G_n$, and thus one could use the Regularity Lemma. As we mentioned, when we embed {\em spanning} subgraphs, the embedding of the last few vertices may create serious difficulties and {\em this difficulty} was overcome by
using the Regularity Lemma -- Blow-up Lemma method. First in \cite{KomSarkoSzem98SeymApprox} Komlós, Sárközy and Szemerédi proved Conjecture~\ref{sey} in 
its weaker, asymptotic form:

\begin{theorem}[Pósa-Seymour conjecture - approximate form, \ev(1994) \cite{KomSarkoSzem98SeymApprox}] For any $\eps>0$ and positive integer $k$ there is an $n_k(\eps)$ such that if $n>n_k(\eps)$ and
$$\mindeg(G_n)>\left(1-\frac1{k+1}+\eps\right)n,$$ 
then $G_n$ contains the $k\th$ power of a Hamilton cycle.
\end{theorem}

Next they got rid of $\eps$ in \cite{KomSarkoSzem96Posa} and \cite{KomSarkoSzem98SeymAnnals}: they proved both conjectures for $n\geq n_k$, without the extra $\eps>0$.

\begin{theorem}[Komlós, Sárközy, and Szemerédi \cite{KomSarkoSzem98SeymAnnals}]\label{main}
For every integer $k>0$ there exists an $n_k$ such that if $n\ge n_k$, and
\begin{equation} \label{minDeg} \mindeg(G_n)\ge\left(1-\frac1{k+1}\right)n,
\end{equation}
then the graph $G_n$ contains the $k^{th}$ power of a Hamilton cycle.
\end{theorem} 

The proofs used the Regularity Lemma \cite{Szem78Regu}, the Blow-up Lemma \cite{KomSarkoSzem97BlowUp}, \cite{KomSarkoSzem98BlowAlgo} and the Hajnal-Szemerédi Theorem \cite{HajSzem69Corrad}. Since the proofs used the Regularity Lemma, the resulting $n_k$ was very large (it involved a tower function). The use of the Regularity Lemma was removed by
Levitt, Sárközy and Szemerédi in a new proof of Pósa's conjecture in \cite{LevittSarkoSzem10PosaRev}. Much later, finally, Péter Hajnal, Simao Herdade, and Szemerédi found a new proof of the Seymour conjecture \cite{HajnalHerdadeSzeme18A-Seymour} that avoids the use of the Regularity Lemma, thus resulting in a ``completely elementary'' proof and a
much smaller $n_k$.

\subsub Historical Remarks.// Partial results were obtained earlier on the Pósa-Seymour conjecture, e.g., by Jacobson (unpublished), Faudree, Gould, Jacobson and Schelp \cite{FaudGouldJacoSche91Posa}, H\"aggkvist (unpublished), Genghua Fan and H\"aggkvist \cite{FanHaggkv94}, 
%% Faudree, Gould and Jacobson \cite{FGJ}, Kezirat?
and Fan and Kierstead \cite{FanKier95Posa}. Fan and Kierstead also announced a proof of the Pósa conjecture if the Hamilton cycle is replaced by
Hamilton path. (Noga Alon observed that this already implies the Alon-Fischer theorem mentioned in Subsection \ref{TilingS}, since the
square of a Hamilton path contains all unions of cycles.) We skip the exact statements of these papers, but mention that Châu, DeBiasio, and
Kierstead \cite{ChauBiasKiers11Posa} proved Pósa Conjecture for all $n>8\times 10^9$.

\subsub Stability Remark.// A crucial lemma of the proof in \cite{HajnalHerdadeSzeme18A-Seymour} is a ``structural stability'' assertion that for some constant $\gamma>0$, either $G_n$ contains an ``almost independent'' set of size $n\over {k+1}$ or $\mindeg(G_n)<({k\over k+1}-\gamma)n.$

\subsection{Ore type results/Pancyclic graphs}\label{OreType}

As we have mentioned, the Hamiltonian problems, above all, Dirac Theorem on Hamiltonicity of graphs, led to several important research directions. Here there is a significant difference between the graph and hypergraph versions. We
shall return to the Hamiltonicity of hypergraphs in Section \ref{HyperLargeExcludedS}. As to ordinary graphs, some generalizations are the Ore-type problems, some other ones are the Pósa-Seymour type generalizations, discussed in the
previous subsection.  In an Ore type theorem we assume that any {\em two independent} vertices have large degree sums. The first such result was

\begin{theorem}[Ore, \ev(1960) \cite{Ore60Hamil}]
 If for any two independent vertices $x,y$ of $G_n$, $\deg(x)+\deg(y)\ge n$, then $G_n$ is Hamiltonian.
\end{theorem}

Bondy had an important ``meta-theorem'' according to which conditions implying Hamiltonicity imply also ``pancyclicity'', which means that $G_n$ contains cycles of all length between $3$ and $n$.\footnote{In some cases we have only weaker conclusions, e.g., in the Bondy-Simonovits theorem \cite{BondySim74C2k}, and also it may happen in some cases that we get only even cycles! See also the paper of Brandt, Faudree, and Goddard \cite{BrandtFaudGodd98} on weakly pancyclic graphs.} A beautiful illustration is

\begin{theorem}[Bondy \ev(1971) \cite{Bondy71Pan}]

(a) If $G_n$ is Hamiltonian and $e(G_n)\ge\turtwo n$, then either $G_n$ is pancyclic or $G_n=K(n/2,n/2)$.

(b) Under the Ore condition, for any $k\in[3,n]$, $G_n$ contains a $C_ k$, or $G_n=K(n/2,n/2)$.
\end{theorem}

Of course, (b) follows from (a) and Ore theorem.  Generalizations of these theorems can be found in Broersma, Jan van Heuvel, and Veldman,
\cite{BroersmaHeuVeld93Ore}, and in general, there are very many ``pancyclicity'' theorems, see e.g. Erdős \cite{Erd74Colu}, Keevash, Lee, and
Sudakov \cite{LeeSudak12Pancyc,KeevSudak10PancycHamil}, Stacho \cite{Stacho99LocPancyc}, Brandt, Faudree, and Goddard \cite{BrandtFaudGodd98}, and many
others. (In some sense the Bondy-Simonovits theorem in \cite{BondySim74C2k} is also a ``weak pancyclic theorem''.)
There are also very many results on {\sc pancyclic digraphs}, see e.g. Häggkvist and Thomassen \cite{HaggkvThoma76PancycDig}, Krivelevich, Lee, and Sudakov
\cite{KriveLeeSudak10ResiPancyc}.

{\sc For~some~related results} see e.g., 
Bollobás and Thomason \cite{BolloThom97WeakPancyc,BolloThom99WeaklyPancyc} Brandt, Faudree, and Goddard \cite{BrandtFaudGodd98}, Favaron, Flandrin, Hao Li,
and F. Tian \cite{FavaFlandLiTian99Ore}, L. Stacho \cite{Stacho99LocPancyc}, Barát, Gyárfás, Lehel, and G.N. Sárközy \cite{BarGyarfLehSarko16Ore}, Barát
and Sárközy \cite{BarSarko16Ore}, Kierstead and Kostochka \cite{KierKost08Ore}, Kostochka\footnote{Kierstead, Kostochka and others have several results
  where Ore type conditions imply Hajnal-Szemerédi type theorem \cite{KierKost08Ore,KierKostoMolYea17Ore}, or a Brooks type theorem
  \cite{KierKosto09OreBrooks}.} and Yu \cite{KostoYu07OrePack} DeBiasio, Faizullah, and Khan \cite{DeBiasFaizKhan}, and many others. 

\subsub Weakly pancyclic graphs.// There are cases, when we cannot hope for all cycles between 3 and $n$. If a graph $G_n$ contains all cycles between $\girth(G_n)$ and its circumference $\circ(G_n)$\footnote{The circumference is the length of the longest cycle. Here we exclude the trees.}, then we call it {\sc weakly pancyclic}, see e.g. Bollobás and Thomason \cite{BolloThom97WeakPancyc,BolloThom99WeaklyPancyc}.  We mention a theorem of \cite{BolloThom97WeakPancyc}, on the girth, answering some questions of Erdős.

\begin{theorem}[Bollobás and Thomason]
 Let $G_n$ be a graph with (at least) two distinct Hamiltonian cycles. Then $n\geq\lfloor(g(G)+1)^2/4\rfloor$, and therefore $\girth(G)\leq\sqrt{4n+1}-1$.
\end{theorem}

There is a pancyclicity defined for bipartite graphs, which considers only even cycles. The original Bondy-Simonovits theorem was also about such pancyclicity.

\begin{theorem}[Bondy and Simonovits \ev(1974) \cite{BondySim74C2k}]
If $e(G_n)>100kn^{1+(1/k)}$, then $G_n$ contains cycles of all lengths $2\ell$, for $\ell=k,\dots, \lfloor e(G_n)/(100 n)\rfloor$.
\end{theorem}

\subsection{Absorbing Method}\label{Absorb}

\gdef\TC#1{\mathbb{TC}_{2#1}}

As we mentioned, when we try to embed a spanning subgraph $\Hm$ into a graph $G_n$, i.e. $m=n$, some difficulties may occur at embedding the ``last few
vertices''. This problem is often solved by using the Blow-up Lemma, or the Absorbing Method, some ``Connecting Lemma'', or by some Stability Argument. Mostly
we combine more than one of them.  The stability argument in most of these papers has the form that we distinguish the ``nearly extremal'' and the
``far-from-extremal'' cases.  Then we handle these two cases separately: mostly we can handle the far-from extremal structure ``easily''.

We mention among papers combining the Stability and the Absorbing methods the newest results of P. Hajnal, Herdade, and Szemerédi \cite{HajnalHerdadeSzeme18A-Seymour} on the Pósa-Seymour Conjecture, or Balogh, Lo, and Molla, Mycroft, and
Sharifzadeh \cite{BalLoMol17Til,BalMcDowMollaMycroft18,BalMollaShari16Tria}, on (Ramsey-Turán-) tiling, (see Subsection \ref{RT-Match}). Below we describe the \Sc Absorbing~Method.

In the Absorbing Method, depending on the problem, we ``invent'' an \Sc Absorbing~Structure, e.g., in the Erdős-Gyárfás-Pyber theorem \cite{ErdGyarfPyb91}, the Triangle-Cycle $\TC\ell$, (defined below) see Figure \ref{TriaCycle}. Then we start with
choosing a special subset of vertices, $\cA\subseteq V(G_n)$, defining this special, absorbing
\BalAbraCapMed{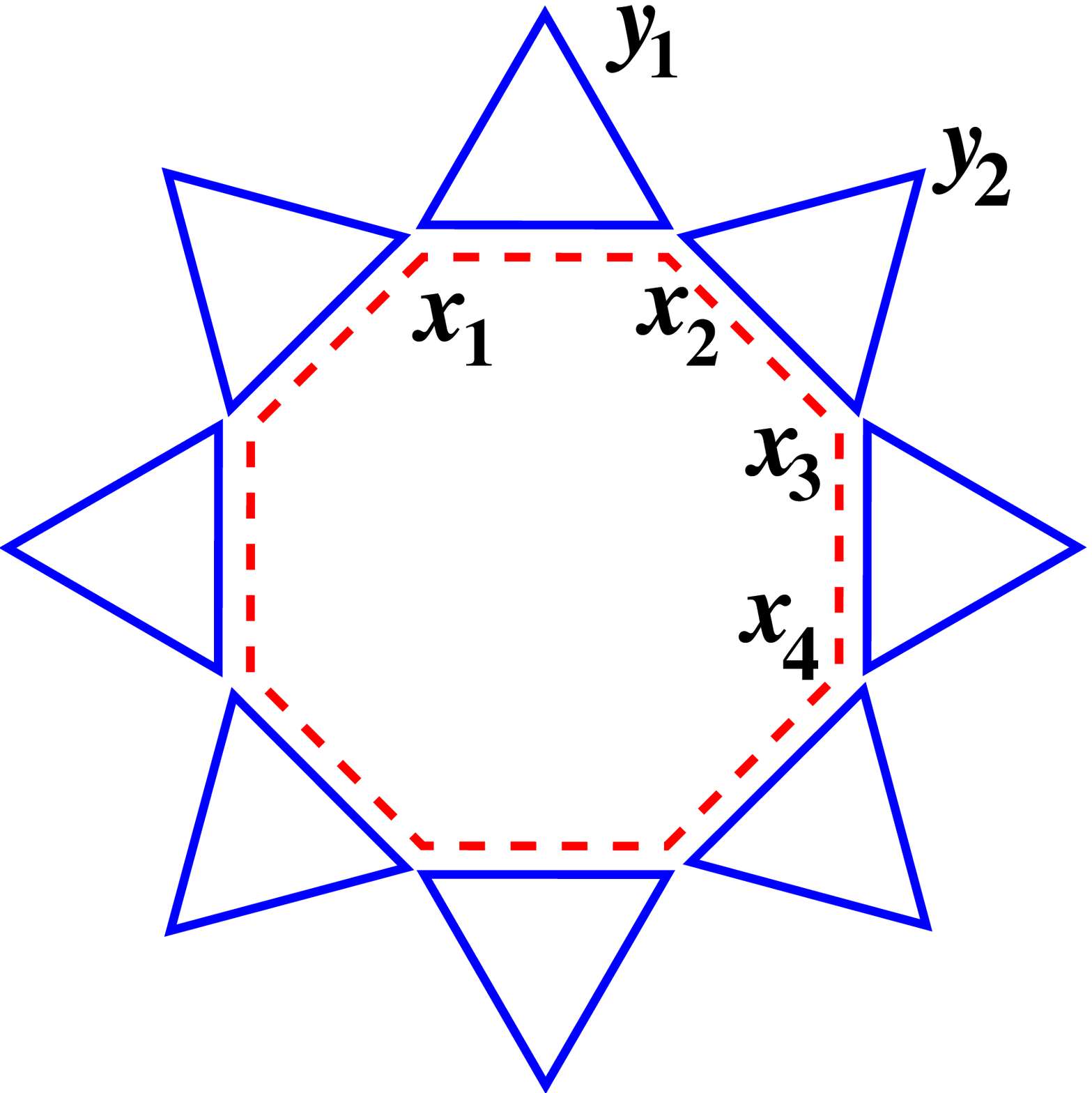}{Tri\-angle-Cycle $\TC\ell$}{30} 
 substructure $G[\cA]$ in 
our graph/hyper\-graph, e.g., in \cite{ErdGyarfPyb91} an $\cA$ spanning a Triangle-Cycle.  Next we put aside  $\cA$ and start building up the
 whole spanning structure in $G_m=G_n-\cA$ as we would do this if $m$ were ``noticeably smaller'' than $n$.  If the ``Absorbing
Structure'' $G[\cA]$ is chosen appropriately, then we will be able to add the last few unembedded vertices of $\cA$ at the end: $\cA$ will absorb/pick up these remaining, uncovered vertices.

We illustrate this, using a proof-sketch of the Erdős-Gyárfás-Pyber theorem.  We shall return to some related newer, sharper results in
Subsections~\ref{RT-Match} and \ref{VertexPartitionS}.\footnote{A stronger statement is Theorem \ref{GyarRuszSarkoSzemTh}.}

\smallskip

\begin{theorem}[Erdős, Gyárfás, and Pyber \cite{ErdGyarfPyb91}]\label{ErdGyarfPybTh} In any $r$-colouring of the edges of $K_n$, we can cover $V(G_n)$ by
  $p(r)=O(cr^2\log r)$ vertex-disjoint monochromatic cycles.
\end{theorem} 

\begin{remark}
The analogous result for bipartite graphs was proved by Haxell \cite{Hax97MonoCyc}.
\end{remark}

The basic structure of the proof is as follows. First we define a ``Triangle-Cycle'' $\TC\ell$ (Figure \ref{TriaCycle}). Its vertices are $x_1,\dots,x_\ell$ and $y_1,\dots,y_\ell$; and its $3\ell$ edges are $x_ix_{i+1}$ (where
$x_{\ell+1}:=x_1$), $y_ix_i$, and $y_ix_{i+1}$, for $i=1,\dots,\ell$.

\begin{enumerate}\dense
\item First, for some $c_1>0$, we find a {\em monochromatic} Triangle-Cycle $\TC\ell$ in $G_n$, with $\ell>c_1 n$.
\item Next we cover $G_m=G_n-\TC\ell$ with $cr^2\log r$ vertex-disjoint monochromatic cycles, also allowing to use some vertices from $Y=\{y_1,\dots,y_\ell\}$.
\item Finally, we can cover the remaining uncovered vertices with one more monochromatic cycle, since our triangle-cycle $\TC\ell$ has the nice property that deleting any subset $Y'\subseteq Y$, the remaining $\TC\ell-Y'$ is still
  Hamiltonian.
\end{enumerate}

The \Sc Absorbing~Method was used in the paper of Rödl, Ruciński, and Szemerédi \cite{RodlRucSzeme06Match}, to find a matching in a hypergraph. This seems to be the breakthrough point:  soon this method became very popular, both
for graphs and hypergraphs. (In Tables \ref{GraphAbsorb}-\ref{HypergraphsB} we list several graph- and  hypergraph applications.) Replacing the regularity method by the Absorption Method is discussed, e.g., in Szemerédi \cite{Szemer13LazyPisa}, Levitt, Sárközy and Szemerédi \cite{LevittSarkoSzem10PosaRev}.

In several cases one uses the \Sc Absorbing~Method to get sharp results for $n>n_0$ after having already a weaker, asymptotic result.  Thus, DeBiasio and Nelsen \cite{DeBiasNels17Monochr}, improving a result of Balogh, Barát, Gerbner, Gyárfás and Sárközy \cite{BalBarGerGyarSarko14} proved a conjecture from \cite{BalBarGerGyarSarko14}:

\begin{theorem}[DeBiasio and Nelsen \cite{DeBiasNels17Monochr}]
 For any $\gamma>0$ there exists an $n_0(\gamma)$ such that if
$n>n_0(\gamma)$ and $\mindeg(G_n)>(3/4+\gamma)n$ and $E(G_n)$ is 2-coloured, then $G_n$ contains two vertex-disjoint monochromatic cycles covering $V(G_n)$.
\end{theorem}

Similarly, the 3-uniform hypergraph tiling results of Czygrinow, DeBiasio, and Nagle \cite{CzygBiasNagle14Tiling} are the sharp versions of some earlier results of Kühn and Osthus \cite{KuhnOsthus06LooseHam} on hypergraph tiling.
Let us repeat that the essence of the Absorption Method is to construct certain ``advantageous configurations'', substructures $G[\cA]$, in $G_n$, called \Sc Absorbing~Structure, covering a large part (say $cn$
vertices) of the host-graph. Next -- using the standard methods, -- we cover all the vertices of $G_n-G[\cA]$ with the given configurations and finally we can expand the embedded configuration into a spanning configuration, using the
particular properties of this \Sc Absorbing~Structure. Often the large \Sc Absorbing~Structure consists of many small substructures, and we gain on each of them an uncovered vertex, obtaining at the end a spanning subgraph, as
wanted.

We refer the interested reader to the Rödl-Ruciński-Szemerédi paper \cite{RodlRucSzemer09Match}, using the ``Absorbing Method'', and to the Rödl-Ruciński survey \cite{RodlRuc10HyperSurv}, however here we mostly avoid the hypergraphs: we
return to them in Section \ref{HyperLargeExcludedS}.\footnote{In several cases we must distinguish subcases also by some divisibility conditions: not only the proofs but the results also strongly depend on some divisibility conditions.}.

\minisepa

In Table \ref{GraphAbsorb} we list just a few successful graph-applications of the \Sc Absorbing~Method.  In some other cases, later, we shall just ``point
out'' that the \Sc Absorbing~Method was successful, when we discuss the corresponding results, e.g., in Sections \ref{PosaSeyRevisit}, \ref{RT-Match},
\ref{HyperLargeExcludedS},\dots We collected some papers using the \Sc Absorbing~Method for hypergraphs in Table \ref{Hypergraphs} (primarily on hypergraph
matching) and in Table \ref{HypergraphsB}, in Section \ref{HyperLargeExcludedS} (primarily on Hamiltonian hypergraphs).

\bigskip

\newcount\tabco

\gdef\tabu#1{\global\advance\tabco by1 %%\the\tabco 
~#1}

\begin{table}[h]
\noindent
{\margofont
\begin{tabular}{|l|l|l|l|l|l|}
\hline
Authors& Year &About what? (or title)&Methods&Where\\
\hline\hline
Erdős, Gyárfás,& 1991 &cycle partition &Absorbing& \cite{ErdGyarfPyb91} \\
Pyber& &Perhaps the first absorbing?&& JCTB \\
\hline 
Levitt, Sárközy, & 2010 &Pósa, How to avoid &Absorbing&DM  \\ 
Szemerédi  &&Regularity Lemma & Connecting&\cite{LevittSarkoSzem10PosaRev}\\
& &\hfill + Stability&Reservoir &\\
\hline
Keevash&2014 &Existence of designs&Absorption&Arxiv \\
&2018& One of the most celebrated&Nibble& \cite{Keev14ExistDesign}
  \\
&&results of these years& \dots& \cite{Keev18ExistDesign} \\
\hline
Ferber, Nenadov,&2014 &Robust Hamiltonicity &Connecting&Arxiv\\
Noever, Peter&2014 & of random directed&Absorbing&Arxiv\\
Škorić & & graphs (resilience)& &\cite{FerbNenaNoeSkor15Robust}\\
\hline
Balogh, Molla, &2016&Triangle factor + small&Absorption&\tiny RSA 
\\
Sharifzadeh &&stable sets,\hfill Weighted graphs&&\cite{BalMollaShari16Tria}\\
\hline
Barber, Kühn,&2016&Edge decomposition of &Iterative &Advances \\
Lo, Osthus &&\hfill graphs with high mindeg&Absorption&\cite{BarbKuhnLoOst16} \\
\hline
Balogh-Lo-Molla&2017 &Digraph packing&Stability &JCTB \\
        &&in-out-degree $\ge{7n/18}$      &Absorption & \cite{BalLoMol17Til} \\
\hline
DeBiasio, Nelsen&2017 &Strengthening Lehel conj.&Absorbing&
\tiny JCTB \\ 
&&&&\cite{DeBiasNels17Monochr} \\
\hline
Glock, Kühn, Lo, &2018 &Existence of designs&Iterative&Arxiv \\
Osthus        &&                         &Absorption&\cite{GlockKuhnLoOsthus17A-FDesign}\\
 &         &                         &Connection&   \\
\hline
Montgomery&2018&Embedding bounded degree &Iterative &Arxiv
\\
&&trees into Random Graphs &Absorption?&\cite{Montgo14BoundedDeg} \\
&&until $p=\Delta\log^5n/n$&& \\
\hline
P. Hajnal, Her-&2018 &Pósa-Seymour&Absorption&Arxiv \\
dade, Szemerédi&& without Regularity Lemma&Connection&  \cite{HajnalHerdadeSzeme18A-Seymour} \\
\hline
\end{tabular}
}
\caption{Absorbing method for graphs}
\label{GraphAbsorb}
\end{table}

These three tables are self-explanatory, however, they contain just a short list of the applications. We could include several further results, like 
the results of Lo and Markström on multipartite Hajnal-Szemerédi results \cite{LoMarks13MultiHajSzem,LoMarkstr15Fact},   on graphs and hypergraphs, or \cite{LoMarks14PerfMatch}\dots

\subsection{Connecting Lemma, Stability, Reservoir}\label{ConnectionS}

In Table \ref{GraphAbsorb} we see several papers using the Absorbing and the Stability Methods. We illustrate this on the example of the new proof of the Pósa-Seymour conjecture \cite{HajnalHerdadeSzeme18A-Seymour}, by Péter Hajnal,
Herdade, and Szemerédi. It has two subcases. A small $\alpha>0$ is fixed and Case 1 (the ``non-extremal'' one) is when each $X\subset G_n$ of $\lfloor {n\over k+1}\rfloor$ vertices has $e(X)\ge \alpha n^2$ edges.  The remaining situation is Case~2, where we use the stability:

\begin{theorem}[Hajnal-Herdade-Szemerédi: Pósa-Seymour, stability \cite{HajnalHerdadeSzeme18A-Seymour}]

  Given an integer $k\ge2$ and an $\alpha>0$, there exists an $\eta=\eta(\alpha,k)>0$ such that in Case 1 (called $\alpha$-non-extremal), if $\mindeg(G_n)\ge (1-\reci k+1
  -\eta)k$, then $G_n$ contains the $k\th$ power of a Hamiltonian cycle.
\end{theorem}

Whenever we use the Absorbing method, mostly we use some other tools as well, tailored specifically to the problem in consideration.  The Connecting Lemma and the Reservoir method are combined with the Absorbing Method, e.g., in the earlier paper of Levitt, Sárközy, and Szemerédi \cite{LevittSarkoSzem10PosaRev} on how to eliminate the use of the Regularity Lemma and the Blow-up Lemma in the proof of the Pósa conjecture. The Regularity Lemma, and the Blow-Up Lemma are eliminated in the tree-embedding paper of Csaba, Levitt, Nagy-György, and Szemerédi \cite{CsabaLevittSzem10Trees}, using the stability and some ``elementary embedding methods''.

The very recent new proof of the Seymour conjecture, by P. Hajnal, Herdade, and Szemerédi \cite{HajnalHerdadeSzeme18A-Seymour} is much more involved and much longer than the original proof. Here the authors use a ``Connecting Lemma'',
asserting that certain parts of $G_n$ can be connected in many advantageous ways.  This means that we have a $G_n$ (with large minimum degree) and wish to find a $C_n^{k+1}\subseteq G_n$. We cover most of the vertices by Turán graphs
$T_{m,k+1}$ and $T_{m,k+2}$, where $m\to\infty$. Inside these ``blocks'' we can easily connect some vertices by $(k+1)\th$ power of a path covering this block, and we must connect these vertices from the various blocks so that altogether we
get a $(k+1)\th$ power of a Hamiltonian cycle. The Connecting Lemma does this.

\subsection{Ramsey-Turán Matching}\label{RT-Match}

This subsection is about the fifth line of Table \ref{GraphAbsorb}, about \cite{BalMollaShari16Tria}. The question is:

\begin{quote}{Does there a new, interesting phenomenon appear when we wish to ensure an almost perfect  tiling, or a Hamiltonian cycle, or some other (almost) spanning configuration in a graph $G_n$ and assume that \beq{SmallStab}\alpha(G_n)=o(n).\eeq }
\end{quote}

We have to decide if we wish to use that $e(G_n)$ is large or that $\mindeg(G_n)$ is large.  We know that for ordinary non-degenerate extremal graph problems the edge-extremal and the degree-extremal problems do not differ too much: if
$\dext(n,\cL)$ is the maximum integer $\Delta$ for which, if $G_n$ is $\cL$-free, then $\mindeg(G_n)\le\Delta$, then \beq{DegEx}\dext(n,\cL)\approx {2\over n}\ext(n,\cL),\eeq and asymptotically $T_{n,p}$ is the degree-extremal graph (where
$p$ is defined by \eqref{subchrom}). We also saw that if we assume \eqref{SmallStab} then in some cases \eqref{DegEx} can be noticeably improved, see
e.g. Theorem \ref{ErdSosRT-Th}. In other cases \eqref{SmallStab} changes the maximum only in a negligible way. However, one can also ask what happens if
we assume \eqref{SmallStab} in cases when we wish to ensure a spanning (or an almost-spanning) subgraph, e.g., a 1-factor, or a Hamilton cycle. Anyway, to ensure an (almost) 1-factor, or a Hamilton cycle we had better to have lower bounds
on $\mindeg(G_n)$ than on $e(G_n)$. Without too much explanation, we formulate two related results.
 
Balogh, McDowell, Molla, and Mycroft studied the minimum degree necessary to guarantee the existence of perfect and almost-perfect triangle-tilings in 
 $G_n$ with $\alpha(G_n)=o(n)$. Among others, they proved

\begin{theorem}[Balogh, McDowell, Molla, and Mycroft \ev(2018) \cite{BalMcDowMollaMycroft18}] Fix an $\eps>0$. If $(G_n)$ is a graph sequence with $\alpha(G_n)=o(n)$ and $\mindeg(G_n)\ge n/3+\eps n$, then $G_n$ has a triangle-tiling covering all but at most four vertices, if $n>n_0(\eps)$.
\end{theorem}

Of course, without the extra condition $\alpha(G_n)=o(n)$ we get back to the Corradi-Hajnal theorem, where $\mindeg(G_n)\ge{2\over3}n$ is needed.  The case when we do not ``tolerate'' the four exceptional vertices is described by \cite{BalMollaShari16Tria}:

\begin{theorem}[Balogh, Molla, Sharifzadeh \cite{BalMollaShari16Tria}]
For every $\eps>0$, there exists a $\ga>0$ such that if $3|n$ and 
$$\mindeg(G_n)\ge\left(\half+\eps\right) n,\Text{and}\alpha(G_n)<\ga n,$$ 
then $G_n$ has a $\K3$-factor.\footnote{The paper has an Appendix written by Reiher and Schacht, about a version of this problem, also using the Absorption technique. In this version they replace the condition that any linear sized vertex-set contains an edge by a condition that any linear sized set contains ``many edges''.}
\end{theorem}

The proofs in \cite{BalMollaShari16Tria} also use the Stability method and the Absorbing technique of Rödl, Ruciński, and Szemerédi \cite{RodlRucSzem06Dirac3CPC}, discussed in Subsection~\ref{Absorb}.

\Section{Colouring, Covering and Packing, Classification}{ColouringProblems}

\subsub The setup. // In this section we consider an \Sc $r$-edge-colouring of a $K_n$, or of a random graph $R_{n,p}$, or of any $G_n$ satisfying some
conditions. We have $r$ families of potential subgraphs, $\cL_i$, ($i=1,\dots,r$), and try to cover $K_n$ with as few monochromatic subgraphs
$L_i\in\cL_i$ in the $i\th$ colour as possible. However, there are several types of problems to be considered:

(A) sometimes we wish to cover all or almost all the {\em edges} of the coloured graph, 

(B) in other cases we wish to cover all the {\em vertices} or almost all the vertices with the vertices of our monochromatic subgraphs.

These are quite different problems and in the next subsection we shall list several versions of these problems, and in some sense, classify them.

There are several problems/results in Extremal Graph Theory simple to formulate, and when we combine some of them, we get very interesting new problems. However, sometimes it is difficult to ``classify'' these problems.  The reader could
ask: ``Why to classify them?'' The answer is that without some classification one may end up with a chaotic picture about the whole field.  Mostly,

\medskip

{\leftskip=5mm \rightskip=5mm

\noindent
we have a ``host graph'' $G_n$ satisfying some conditions, e.g., it may be a complete graph $K_n$, or a random graph $R_{n,p}$, or a pseudo-random graph
  with many edges, or with large minimum degree. There is also a family $\cL$ of ``sample'' graphs\footnote{Here we took $\cL_i=\cL$.}. Now,
  $E(G_n)$ is $r$-coloured, and

}

\smallskip

\begin{enumerate}\dense
\item either we wish to \Sc partition~the~vertices of $G_n$ into a few classes, $U_1,\dots,U_t$ so that each $G[U_i]$ spans a monochromatic $L_i\in\cL$, or
\item we wish to \Sc cover~the~edges, $E(G_n)$, by a few copies of monochromatic sample graphs, $L_i\in\cL$, or
\item we wish to approximate the situation (a) or (b), allowing a few uncovered vertices.
\end{enumerate}

\subsub Historical remarks.// The case of one colour and vertex-disjoint packing goes back to several early papers in Extremal Graph Theory: several proofs, e.g., Erdős \cite{Erd62Radem}, or later Simonovits \cite{Sim66Tihany}, used that
large part of the considered $G_n$ can be covered by vertex-disjoint copies of some $L$.\footnote{In some sense, this is used also in the original proof of Erdős-Stone theorem \cite{ErdStone46}.} Also this was used in the new proof of the
Pósa-Seymour conjecture, in~\cite{HajnalHerdadeSzeme18A-Seymour}.

Mostly we define $K_1$ and $K_2$  also as monochromatic graphs from $\cL$: this is needed to ensure the
 existence of suitable colourings. The graphs in $\cL$ may have more restricted or less restricted structure, e.g., they may be (A) all the
connected graphs, or (B) all the trees, or (C) the sets of independent edges, or (D) all the cycles, or (E) the $k\th$ powers of the cycles, or (F) the Hamiltonian graphs, (G) on the other hand, they may be copies of the same fixed
graph~$L$.

We start with Covering Problems where we have only one colour. We could continue with Ramsey problems, when we wish to find just one monochromatic subgraph of a particular type, however, we have already written about Ramsey problems, and we shall not consider them here.

\subsection{The One-colour Covering Problem}

In this section we consider edge-coverings. Here is an early problem of Gallai, corresponding to the simplest case, to the monochromatic $G_n$ i.e. $r=1$.

\begin{questio} Given a sample graph $L$ and a graph $G_n$, how many copies of $L$ and edges are enough to cover $E(G_n)$ (in the worst case)? In principle, we may require that the selected copies of $L$ be (i) edge-disjoint, or (ii)
  vertex-disjoint, or (iii) we may allow them to overlap.
\end{questio}

Of course, the $L$-free graphs need $e(L)$ edges, so the $L$-extremal graphs need $\ext(n,L)$ copies of $L$ or edges. One feels that if $e(G_n)$ is much larger than $\ext(n,L)$, then one may use many copies of $L$ and only a few edges. This motivates many results in this area. Here one of the first important results is

\begin{theorem}[Erdős, Goodman, and Pósa \ev(1966) \cite{ErdGoodPosa66}]\label{EGP66} Any graph $G_n$ can be edge-covered by $\turtwo n$ complete subgraphs of $G_n$.
\end{theorem}

\begin{remark}
 
(a) As it is stated in \cite{ErdGoodPosa66}, to cover the edges, we may restrict ourselves to edges and triangles in Theorem \ref{EGP66}.

(b) Theorem \ref{EGP66} is sharp, as shown by $\Tn n,2$.

(c) Theorem \ref{EGP66} was also proved by Lovász, as remarked in \cite{ErdGoodPosa66}.

\end{remark}

The extremal number for cycles is $n-1$.  Erdős and Gallai conjectured 
(see \cite{ErdGoodPosa66}) and Pyber proved

\begin{theorem}[Pyber \ev(1985) \cite{Pyber85Gallai}]
Every graph on $n$ vertices can be edge-covered by $n-1$ cycles and edges.
\end{theorem}

Note that for trees this is sharp.  Pyber in \cite{Pyber85Gallai} proved also some stronger results, and mentioned that the crucial tool in his proof was
the following

\begin{theorem}[Lovász \ev(1966) \cite{Lov66TihanyCover}] A graph on $n$ vertices can be edge-covered by $\fele n $ disjoint paths and cycles.
\end{theorem}

Lovász proved also several related theorems, e.g.,

\begin{theorem}[Lovász \ev(1966) \cite{Lov66TihanyCover}] Any graph $G_n$ can be edge-covered by ${2\over3}n$ double-stars (trees of diameter $\le3$).
\end{theorem}

Another result of this (early) Lovász paper is a generalization of the Erdős-Goodman-Pósa theorem.\footnote{For some related result for random or quasi-random graphs see \cite{KorKriveSudak15DecoRandom},\cite{KrivelSamot12},
  \cite{GlockKuhnLoOsthus16PathDecomp},\cite{GlockKuhnOsthus16Decomp}, and many others.}

\subsection{Embedding monochromatic trees and cycles}

In the following parts we consider edge-coloured graphs $G_n$ and try to partition $V(G_n)$ into a few subsets $V_i$ spanning some monochromatic Hamiltonian cycles,\footnote{Here ``$G[V_i]$ is Hamiltonian'' means that it has a spanning
  cycle. } or monochromatic powers of Hamiltonian cycles, or spanning trees.

These problems and results are related to Lehel's Conjecture \ref{LehConj} about partitioning the vertices of edge-coloured graphs into two monochromatic cycles, and the Gyárfás Conjecture \ref{Gyarf78Pack}, (see \cite{GyarfLehel78Pack}),
about partitioning the vertices of edge-coloured graphs into given monochromatic trees. Here the colourings are always edge-colourings.\footnote{We used vertex-colouring in connection with colouring properties of excluded subgraphs, or
  equipartitions in Hajnal-Szemerédi theorem, \dots}

\subsubsection{Gyárfás tree-conjecture}

The topic of graph packing has at least two larger parts: the vertex-packing and the edge-packing. Here we are interested in packing some given graphs $H_i$ into a $G_n$ in an edge-disjoint way. This problem is interesting on its own and
also was originally motivated by Theoretical Computer Science, more precisely, by \Sc computational~complexity, see e.g. the paper of Bollobás and Eldridge \cite{BolloEldri78Packing}. Given some graphs $H_1,\dots,H_\ell$, their \Sc packing
means finding the automorphisms $\pi_1,\dots,\pi_\ell$ which map them into $K_n$ in an edge-disjoint way.

For a family $H_1,\dots,H_m$ of graphs, we say that they \Sc pack into $G$, if they have edge-disjoint embeddings into $G$. 

\begin{conjecture}[Gyárfás \ev(1978) \cite{GyarfLehel78Pack}]\label{Gyarf78Pack}
Let for $i=1,2,\dots,n$, $T_i$ be an $i$-vertex tree. Then $K_n$ can be decomposed into these trees: $\{T_i\}$ pack into~$K_n$.
\end{conjecture}

An asymptotic weakening of the conjecture was proved by Böttcher, Hladký, Piguet, and Taraz \cite{BottHladPigTar16Pack}, for bounded degree trees.\footnote{In fact, one can allow the first $o(n)$ trees to have arbitrary degrees.} 
This was improved to a sharp version:

\begin{theorem}[Joos, J. Kim, Kühn, and Osthus \ev(2016) \cite{JoosKimKuhnOsthus16Pack}] For any $\Delta>0$ there exists an $n_\Delta$ such that for $n>n_\Delta$, if $T_1,\dots,T_n$ are trees with $\maxdeg(T_i)<\Delta$ and $v(T_i)=i$
  ($i=1,\dots,n$), then $E(K_n)$ has a decomposition\footnote{i.e. $T_1,\dots,T_n$ pack into $K_n$.}  into $T_1,\dots,T_n$.
\end{theorem}

Hence the tree packing conjecture of Gyárfás holds for all bounded degree trees, and sufficiently large $n$. Beside several further results,
\cite{JoosKimKuhnOsthus16Pack} also contains the following fairly general result.

\begin{theorem} Let $\de>0$ and $\Delta$ be fixed.
Let $\cF$ be a family of trees $T^m$ of the following properties\footnote{We used superscript since here $v(T^m)$ is not necessarily $m$.}:

(i) For each $T^m\in\cF$, $v(T^m)\le n$, and $\maxdeg(T^m)\le \Delta$.

(ii) For at least $(\half+\de)n$ values of $m$ $v(T^m)\in[\de n, (1-\de)n]$.

(iii) $\sum_m e(T^m)={n\choose 2}$. 
\\
Then $K_n$ can be decomposed into the trees $T^m$.
\end{theorem}

\begin{remark} 
  There is a vast literature on this topic and we recommend the reader to read the introduction of \cite{JoosKimKuhnOsthus16Pack}, and Gyárfás \cite{Gyarf16VertexCoverSurv}. For results where we consider fewer trees, see e.g., Balogh and
  Palmer \cite{BalPalmer13Pack}.  (The main tool of the proof in \cite{BalPalmer13Pack} is the Komlós-Sárközy-Szemerédi theorem \cite{KomSarkoSzemer01SpanningTree} on embedding spanning trees.)  Several related but slightly different
  problems are discussed, e.g., in the survey of Kano and Li \cite{KanoLi08MonoSurv}, e.g., it describes several Anti-Ramsey decomposition problems as well.
\end{remark}

\minisepa

One could ask what happens if we wish to embed some trees into a non-complete graph $G_n$. A nice result on this is

\begin{theorem}[Gyárfás \cite{Gyarf14PackingTrees}]
 If a sequence $T_1,T_2,\dots,T_{n-1}$ of trees can be packed into $K_n$ then they can be packed also into any $n$-chromatic graph.
\end{theorem}

\subsection{Bollobás-Eldridge conjecture}

The packing problems (strongly connected to Theoretical Computer Science) were discussed roughly the same time by Bollobás and Eld\-ridge~\cite{BolloEldri78Packing}, Catlin \cite{Catlin74}, and Sauer and Spencer \cite{SauerSpencer78}.  One of the most important conjectures in the field of graph-packing is

\begin{conjecture}[Bollobás-Eldridge \cite{BolloEldri78Packing}] Let $H_1$ and $H_2$ be two $n$-vertex graphs. If 
$$(\maxdeg(H_1)+1)(\maxdeg(H_2)+1)\le n+1,$$
 then there is a packing of $H_1$ and $H_2$, i.e., there are two edge-disjoint subgraphs of a complete graph $K_n$ isomorphic to $H_1$ and $H_2$, respectively. \end{conjecture} 

The complementary form of this problem is

\begin{conjecture}[Bollobás--Eldridge]
 Let $\maxdeg(H_n)\le k$. If $G_n$ is a simple graph, with $$\mindeg(G_n)\ge {{kn-1}\over {k+1}}, $$ then it contains $H_n$.
\end{conjecture}

Aigner and Brandt \cite{AigBrandt93}, and Alon and Fischer \cite{AlonFisch96Fact} proved the conjecture for $\maxdeg(H_n)=2$, i.e. when $H_n$ is the union
of cycles.  Csaba, Shokoufandeh, and Szemerédi proved this for $\maxdeg(H_n)=3$.

\begin{theorem}[Csaba, Shokoufandeh, and Szemerédi \cite{CsabaShokoSzem99}]
 If $G_n$ is a simple graph, with $\mindeg(G_n)\ge \reci4 (3n-1),$ then it contains any $H_n$ for which $\maxdeg(H_n)\le 3$, if $n$ is sufficiently large.
\end{theorem}

Csaba \cite{Csaba03BolloEldridge4} also proved this conjecture for $\maxdeg(H_n)=4$, and in \cite{Csaba07BolloEldriBipa} for bipartite $H_n$, (where $G_n$ is not necessarily bipartite) and $\maxdeg(H_n)\le\Delta$.

Improving some results of Sauer and Spencer \cite{SauerSpencer78} and of Catlin, P. Hajnal and M. Szegedy \cite{HajnalPSzegedy92Packing} considered some bipartite packing problems where they proved an asymmetric version: for the first graph they considered the maximum degree but for the second one only the average degree, thus improving the previous results for bipartite graphs.

The reader interested in more details will find a lot of information, e.g., in the survey of Kierstead, Kostochka, and Yu \cite{KierKostochYu09Pack}.

\subsection{Vertex-partitioning into monochromatic subgraphs of given type}\label{VertexPartitionS}

This area has two parts: one about ordinary graphs and the second one about hypergraphs. The hypergraph results can be found in Subsection \ref{HyperCover}.
Here we consider edge-colourings. Gerencsér and Gyárfás \cite{GerenGyarf67} proved that the vertices of any 2-coloured $K_n$ can be partitioned into the
vertex sets of two monochromatic paths of distinct colours. Here we consider problems where the {\em vertex set} of an $r$-coloured graph has to be
partitioned into the vertex sets $U_i$ of some {\em given types of monochromatic subgraphs}.
%% \footnote{We do not require that $G[V_i]$ be monochromatic, only that it contains some given monochromatic $L_i$.}
 We mention just a few such theorems, to give the flavour of these results. 

\begin{problem}
  Given a graph $G_n$ and a family $\cL$ of graphs. What is the minimum integer $t$ for which every $r$-colouring of $E(G_n)$ has a vertex partition $V(G_n)=\pcup U_i$ into $t$ vertex sets so that each $G[U_i]$ contains a monochromatic
  spanning $H_i\in\cL$.\footnote{More generally, we may fix for each colour $i$ a family $\cL_i$ and may try to cover $V(G_n)$ by vertex disjoint subgraphs $H_i\in\cL_i$.}
\end{problem}

Gyárfás, Sárközy, and Selkow \cite{GyarfSarkoSelko15Cov} discussed a natural but much more general family of problems:

\begin{problem}
  Given a family $\cL$ of graphs, (trees, connected subgraphs, matchings, cycles,\dots) and two integers, $r\ge s\ge1$. At least how many vertices of an $r$-edge-coloured $K_n$ can be covered by $s$ monochromatic subgraphs $H_i\in\cL$ in
  this $K_n$?\footnote{Formally we have here two problems, one when we $r$-colour $E(K_n)$, the other when we $r$-colour $E(G_n)$, however, the difference ``disappears'' if $r$ is large. Further, we may also ask for the largest subgraph
    $H\subseteq K_n$ that is coloured by at most $t$ colours, which is different from asking for the largest number of edges covered by $t$ monochromatic $H_i\in\cL$.}
\end{problem}

The simplest case is $r=1$, where we use only one colour.  Another simple subcase is when we RED-BLUE-colour the edges of a $K_n$ and study how many monochromatic cycles are needed (in the worst case) to cover the vertices.  We can not always partition the edges into monochromatic cycles, unless we agree that the vertices and the edges are also regarded as cycles.\footnote{If for a fixed $x$, all edges $xy$ are Red, and the other edges are Blue, then we need this.}  In this case we can always partition $E(G_n)$ into $e(G_n)$ monochromatic ``cycles'': real cycles, edges and vertices. The next part contains some repetition from Section \ref{RemovingReguLem}, primarily from Subsection \ref{Absorb}.  We start with

\begin{conjecture}[Gyárfás \ev(1989) \cite{Gyar89Covering}, proved]\label{GyarfConjCov}
There exists an integer $f(r)$ independent of $n$ such that 
if $E(K_n)$ is $r$-coloured, then $V(K_n)$ can always be covered by $f(r)$ vertex-disjoint monochromatic paths.
\end{conjecture}

Actually, Gyárfás formulated three conjectures in \cite{Gyar89Covering}. In the first one he wanted to cover $V(K_n)$ by $f(r)=r$ vertex-disjoint monochromatic paths, in the second, weaker one, the vertex-disjointness was not assumed, and the third, weakest one, was Conjecture \ref{GyarfConjCov}. ($f(r)=r$ would yield the first conjecture.) Gyárfás proved the following weakening of his conjectures:

\begin{theorem}[Gyárfás \cite{Gyar89Covering}]
There exists an integer $f(r)>0$ such that in any $r$-colouring of $E(K_n)$ one can cover the vertices by $f(r)$ monochromatic paths.
\end{theorem} 

This result does not assert that the covering paths are vertex-disjoint. The proof of Gyárfás gave an explicit $f(r)\approx r^4$. This was improved by
Erdős, Gyárfás, and Pyber in Theorem \ref{ErdGyarfPybTh}: In any $r$-edge-colouring of $K_n$, one can cover $V(G_n)$ by $p(r)=O(r^2\log r)$ {\em vertex-disjoint monochromatic cycles}. This was further improved by

\begin{theorem}[Gyárfás-Ruszinkó-Sárközy-Szemerédi \cite{GyarfRuszSarkoSzemer06CycleParti}]\label{GyarRuszSarkoSzemTh}
If $E(K_n)$ is $r$-coloured, then $V(K_n)$ can be partitioned into $O(r\log r)$ vertex-sets of monochromatic cycles.
\end{theorem}

We have to point out two things. 

(a) While the problems with cycles and paths are basically of the same difficulty in the Erdős-Gallai theorems, here the problems on covering {\em with cycles} are much more difficult.

(b) Covering with {\em vertex-disjoint} paths/cycles is significantly more difficult than the case when vertex-disjointness is not assumed.

Here the most important conjecture {\em was}

\begin{conjecture}[Erdős, Gyárfás, and Pyber \cite{ErdGyarfPyb91}]
In Theorem \ref{GyarRuszSarkoSzemTh}\, $p(r)=r$ vertex-disjoint monochromatic cycles are enough.
\end{conjecture}

This was proved for $r=2$ by Bessy and Thomassé (see Subsection \ref{VertexPartitionS}) but ``slightly disproved'' for $r\ge3$ by Pokrovskiy \cite{Pokro14Counter}, see Remark \ref{PokrovskiyR}.  A whole theory emerged around this conjecture, and below
we provide a very short description of it, basically following Sárközy's paper \cite{Sarko17Parti}, which starts with a good ``minisurvey'' about this area.
He extends the problem to covering by powers of cycles, (as Pósa and Seymour extended Dirac's Hamiltonicity theorem). Sárközy proved

\begin{theorem}[Sárközy \ev(2017) \cite{Sarko17Parti}]
For every integer $k\ge1$ there exists a constant $c(k)$ such that
in any 2-colouring of $E(K_n)$ at least $n-c(k)$ vertices can be covered by 
$200k^2\log k$ vertex-disjoint monochromatic $k\th$ powers of cycles.
\end{theorem}

The interested reader is referred to \cite{Sarko17Parti}. We shall return to the corresponding {\em hypergraph problems} in Section \ref{HyperCover}.
We close this part with a result of Grinshpun and Sárközy \cite{GrinshpunSarko16Mono} on a conjecture of Gyárfás where the cycles are replaced by an arbitrary fixed family of bounded degree graphs. (This covers the case of the $k\th$ powers of cycles.)

 Fix a degree bound $\Delta$ and let $\cF = \{F_1, F_2, \ldots \}$ be any given family of graphs, where $F_r$ is an $r$-vertex graph of maximum degree at most $\Delta$. 

\begin{theorem}[Grinshpun and Sárközy \cite{GrinshpunSarko16Mono}]\label{GrinSarkoMono} 
There exists an absolute constant $C$ such that for any 2-colouring of $E(K_n)$, there is a vertex partition of $K_n$ into (vertex sets of) monochromatic copies of members of $\cF $ with at most $2^{C\Delta \log \Delta}$ parts. 
\end{theorem}

If $\cF $ consists of {\em bipartite graphs} then $2^{C\Delta \log \Delta}$ can be replaced by $2^{c\Delta}$, which is best possible, apart from the value of the constant $c$:

\begin{theorem}[Grinshpun and Sárközy \cite{GrinshpunSarko16Mono}]\label{GrinSarkoMonoB}  
  Let $\cF$ be a family of bipartite graphs with maximum degree $\Delta$.  There is an absolute constant $c$ such that for any 2-edge colouring of $K_n$,
  there is a vertex partition of $K_n$ into (vertex sets of) monochromatic copies of $F_r\in\cF $ with at most $2^{c\Delta}$ parts.
\end{theorem}

These results are strongly connected to some results of Conlon, Fox and Sudakov \cite{ConFoxSudak12Ramsey}.
We close this section with an open problem.

\begin{problem}
Do the Grinshpun-Sárközy theorems extend to three colours?
\end{problem}

\section{Hypergraph extremal problems, small excluded graphs}\label{HyperSmallExcludedS}

Until now we primarily concentrated on two Universes: graphs and integers. Here we include a short section on the Universe of hypergraphs. We used hypergraphs, e.g., in estimating independence numbers in Section \ref{AKS-Indep} in ``uncrowded graphs and hypergraphs'', to improve a Ramsey number estimate, $R(3,k)$, to disprove the Heilbronn conjecture, see Subsection \ref{HeilbronnS}, and in case of the infinite Sidon sequences, and several further results, in Subsections \ref{AKS-Indep}-\ref{PippeSpencS}. The parts in Section~\ref{RuzsaSzemAppliS} connected to the Removal Lemma were also in some sense hypergraph results. Here we ``start again'', but go into other directions. We refer the readers interested in more details to the surveys of Sidorenko, \cite{Sidor95Know}, Füredi \cite{Fure88Match,Fure91Sur}, G.O.H. Katona \cite{Katona74ExtreHyper}, Keevash \cite{Keev11HypergBCC}, of Kühn and Osthus \cite{KuhnOsthus09Embedding}, and of Rödl and
Ruciński~\cite{RodlRuc10HyperSurv}.

Below, to emphasize that we consider hypergraphs, occasionally we shall use  a different typesetting, e.g., for $r$-uniform hypergraphs, for the excluded hypergraphs we may use $\Lhyp{}r$, instead of $L$, and $\cL^{(r)}$ for the family of excluded hypergraphs.

\BalAbraCapMed{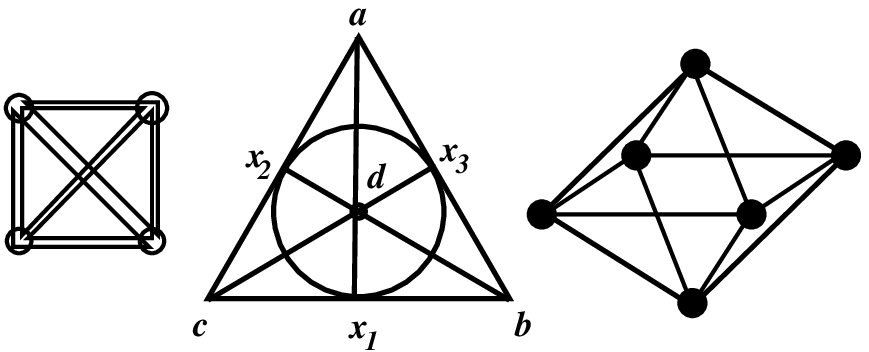}{Small excluded hypergraphs: Complete hypergraph, Fano hypergraph and the Octahedron hypergraph}{50} We start with the case of ``Small excluded subhypergraphs''. Here ``small'' means that the excluded hypergraphs
have bounded number of vertices. In Section \ref{HyperLargeExcludedS} we shall consider the case of 1-factors and the problem of ensuring Hamiltonian cycles and other ``spanning'' or ``almost spanning'' configurations. For the sake of
simplicity, we mostly (but not always) restrict ourselves to 3-uniform hypergraphs. The extremal number for $r$-uniform hypergraphs will be denoted also by $\ext(n,\cL)$, or, to emphasise that we consider $r$-uniform hypergraphs, we may
write $\ext_r(n,\cL)$, or $\ext_r(n,\cL^{(r)})$.\footnote{Here we use $L$ for an excluded graph, $\LL$ for a hypergraph, and $\cL$ for a family of graphs or hypergraphs.}  There are many interesting results on hypergraph extremal
problems. We mention the easy theorem of Katona, Nemetz and Simonovits \cite{KatonNemSim64}, according to which $\ext_r(n,\cL)/{n\choose r}$ is monotone decreasing, non-negative, and therefore convergent.

\smallskip

\BalAbraCapMed{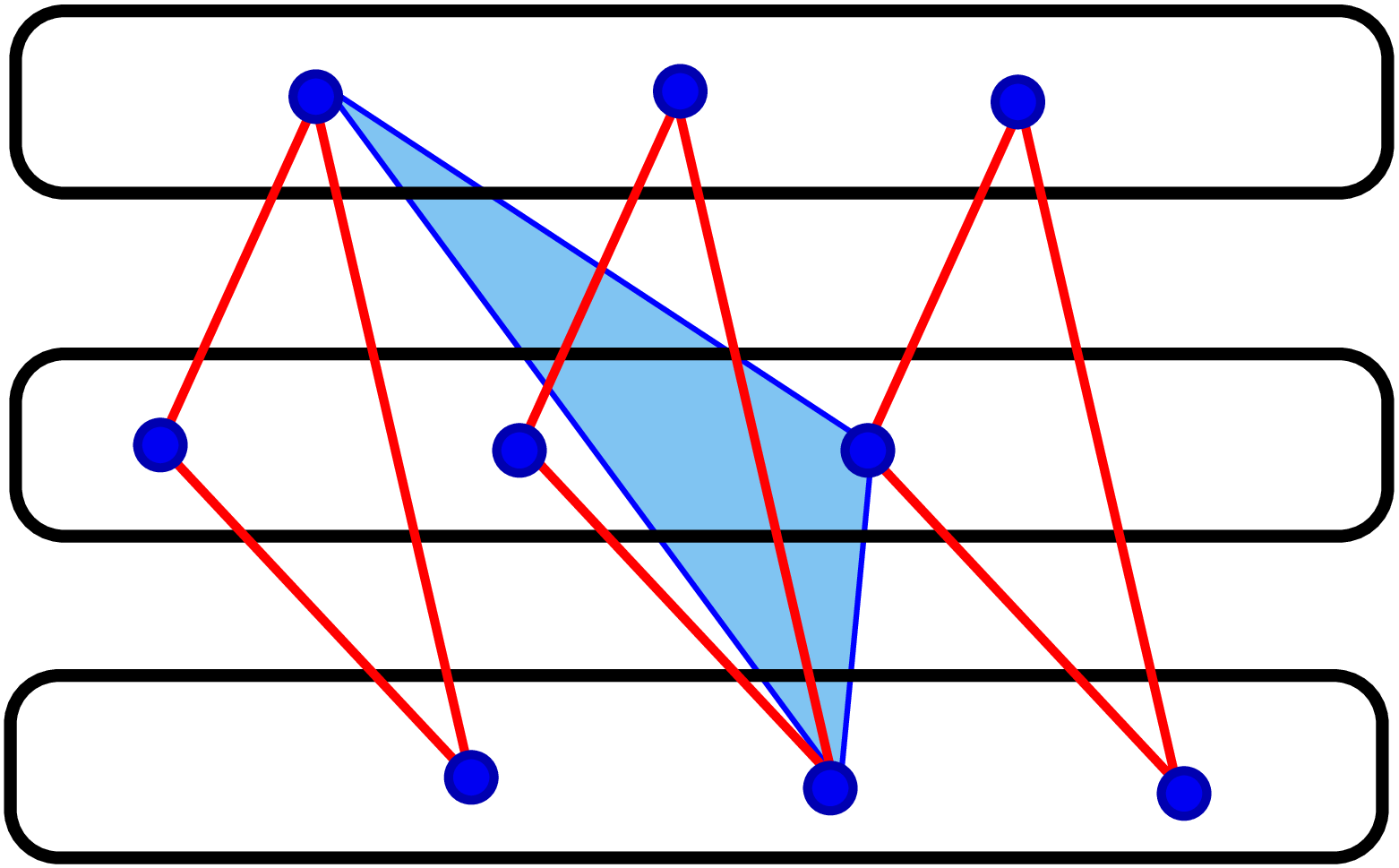}{Three-partite 3-uniform hypergraph.}{38} As for simple graphs, for hypergraphs we can also distinguish \Sc degenerate and \Sc non-degenerate extremal graph problems: an $r$-uniform problem is degenerate if $\ext_r(n,\cL)=o(n^r)$.  For graphs the Kővári-Sós-Turán theorem was the key in this.  A fairly simple result of Erdős \cite{Erd64Israel} generalizes the Kővári-Sós-Turán theorem. Let $\Khyp rr(t_1,t_2,\dots,t_r)$ denote the $r$-uniform hypergraph in which the vertex set $V$ is partitioned\\

\noindent  into $V_1,\dots,V_r$, $|V_i|=t_i$, 
and the hyperedges are the ``transversals'': $r$-tuples intersecting each $V_i$ in one vertex.

\begin{theorem}[Erdős, \ev(1964) \cite{Erd64Israel}]\label{ErdHyperKST}
If $t=t_1\le t_2\le\ldots\le t_r$, then 
$$\ext_r(n,\Khyp rr(t,t_2,\dots,t_r))=O(n^{r-1/t^{r-1}}).$$
\end{theorem}

This implies, exactly as for $r=2$, that 

\begin{corollary}\label{HyperDegen} $\ext_r(n,\cL^{(r)})=o(n^r)$ if and only if there is an $\LL\in\cL^{(r)}$ of strong chromatic number $r$.\footnote{The strong chromatic number $\chi_S(\Fhyp{}r)$ of $\Fhyp {}r$ is the minimum $\ell$ for which the vertices of $\Fhyp{}r$ can be $\ell$-coloured so that
 each hyperedge gets $r$ distinct colours.
Our condition is equivalent with that some $\Fhyp{}r\in\cL^{(r)}$ is a subgraph of $\Khyp rr(a,\dots,a)$ for some large $a$.}	
\end{corollary}

\begin{remark} As in Corollary \ref{DegenJump}, if the problem is non-degenerate, then the ``density constant jumps up'':
  \beq{JumpConst}\ext_r(n,\cL^{(r)})>\left({\reci r^r +o(1)}\right)n^r.\eeq A famous problem of Erdős was whether for 3-uniform hypergraphs 1/27 is a
  ``jumping constant'': does there exist a constant $c>0$ such that if for some $\eta>0$ $\ext_r(n,\cL^{(r)})>\left({\reci 27 +\eta+o(1)}\right)n^r$ then
  $\ext_r(n,\cL^{(r)})>\left({\reci 27 +c+o(1)}\right)n^r$ also holds.  There are many related results, here we mention only the breakthrough paper of
  Frankl and Rödl \cite{FranklRodl84HyperJump}, which however, does not decide if $1/27$ is a ``jumping constant'' or not. We mention that according to
  Pikhurko \cite{Pikhu14Dens} there are continuum many limit densities and there are among them irrational ones even for finite families of excluded
  $r$-graphs.  We also recommend to read Baber and Talbot \cite{BaberTalbot11Jump} and several papers of Y. Peng, e.g., \cite{PengZhao08}, in this area.
\end{remark} 

\subsubsection{Turán's conjecture}

Consider now 3-uniform hypergraphs: $\Hhyp{}3=(V,\cE)$. To formulate the famous hypergraph conjectures of Paul Turán in the two simplest cases, we need two
constructions. We shall call an $r$-uniform hypergraph \Sc $h$-partite if its vertices can be partitioned into $h$ classes, none of which contains
hyperedges.

\BalAbraCap{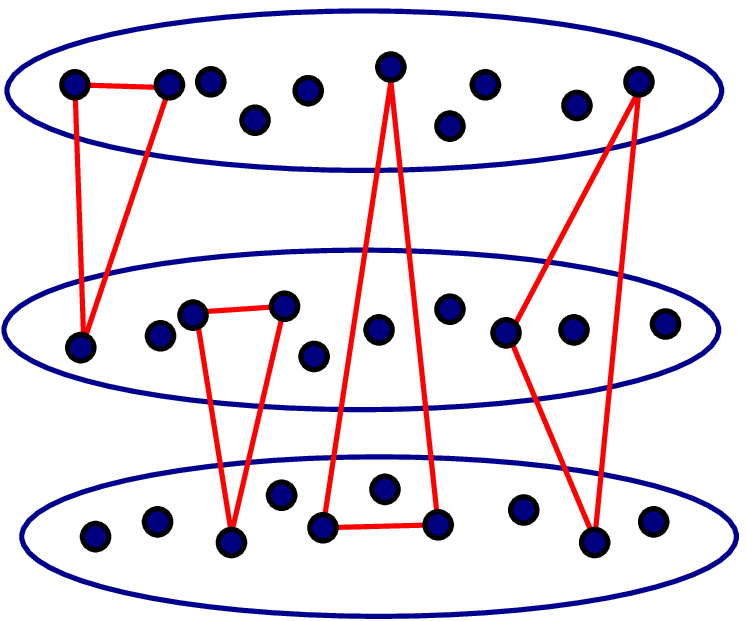}{$\Khyp43$-extremal?}  (a) For the excluded 3-uniform complete 4-graph $\Khyp43$, consider the 3-uniform hypergraph $\Hhyp n3$ obtained by partitioning $n$ vertices into 3 classes $U_1$, $U_2$ and $U_3$ as equally as possible and taking the triples of form $(x,y,z)$ where $x,y$, and $z$ belong to different classes, and the triplets $(x,y,z)$ where $x$ and $y$ belong to $U_i$ and $z$ to $U_{i+1}$, for $i=1,2,3$, where $U_4:=U_1$, see Figure \ref{turhypcon}. Turán conjectured that this is the extremal hypergraph for $\Khyp43$. This is unknown, we do not know even if this is asymptotically sharp.

Actually, first Katona, Nemetz, and Simonovits \cite{KatonNemSim64} gave examples (for $n=3m+1$) showing that if Turán's conjecture holds, then the uniqueness of extremal graphs does not always hold.

A more general construction of Brown \cite{Brown83TurHyp}, generalized by Kostochka \cite{Kosto82TurHyp}, shows that if the conjecture holds, then there are many-many extremal graph structures for the extremal hypergraph problem of $\Khyp43$ and $n=3t$. For a slightly more detailed description of this see e.g. Fon-der-Flaass \cite{Flaass98Hyper}, Razborov, \cite{Razbor11Flaass}, Simonovits \cite{Simo13TuranInflu2}.

\BalAbraCap{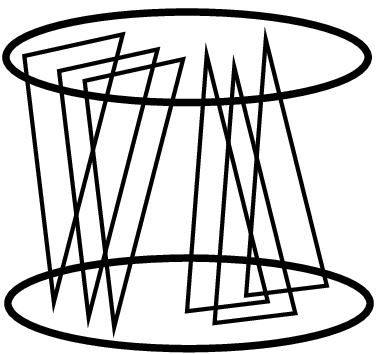}{$\Khyp53$-extremal?} (b) For the excluded complete 5-graph $\Khyp53$ Turán had a construction -- for the potential extremal hypergraph -- with four classes and another one with two classes. The one with two classes is simple, see Figure \ref{fanoext2}. It is a \Sc complete~bipartite hypergraph: we partition the vertices into $\cA$ and $\cB$ and consider all the triples intersecting both classes. V.T.~Sós observed that the construction with two classes can be obtained from the construction with four classes by moving around some triples in some simple way. Probably J. Surányi found a construction showing that Turán's conjecture for $\Khyp53$ is false for $n=9$. Kostochka
(perhaps) generalized this, founding counterexamples for every $n=4k+1$.\footnote{These constructions seem to be forgotten, ``lost'' and are not that important.}  Still Turán's conjecture may be sharp, or, at least, asymptotically sharp.\footnote{Since for hypergraphs we have at least two popular chromatic numbers, therefore the expression $r$-uniform $\ell$-partite may have at least two meanings in the related literature. }
 
Among the new achievements we mention Razborov's Flag Algebra method and his results on hypergraphs \cite{Razbor07Flag,Razbor10Config,Razbor11Flaass}. 

\subsubsection{A simple hypergraph extremal problem}\label{KatonaBolloS}

As we have often emphasized, to solve a hypergraph extremal problem is mostly hopeless, despite that recently many nice results were proved on
hypergraphs. Keevash described this situation by writing (in MathSciNet, on \cite{FureMubaPikhu08Hyper}, in 2008):

\begin{quote}{
``An important task in extremal combinatorics is to develop a theory of Turán problems for hypergraphs. At present there are very few known results, so it is interesting to see a new example that can be solved.''}
\end{quote}

Consider 3-uniform hypergraphs. The difficulties are reflected, among others, by that we do not know the extremal graph for $\cL_{4,3}^{(3)}$, i.e. when we exclude the 4-vertex hypergraph with 3 triples.\footnote{Generally, $\cL_{k,\ell}^{(r)}$ is the family of $r$-uniform hypergraphs of $k$ vertices and $\ell$ hyperedges.  As we have mentioned, the problem of $\ext(n,\cL_{k,\ell}^{(r)})$ was considered in two papers of Brown, Erdős, and Sós \cite{BrownErdSos73a,BrownErdSos73b} and turned out to be very important in this field.  Originally Erdős conjectured a relatively simple asymptotic extremal structure, for $\cL_{4,3}^{(3)}$ but his conjecture was devastated by a better construction of Frankl and Füredi \cite{FranklFure84Exact}. This construction made this problem rather hopeless.}  The next question is among the easier ones. G.O.H. Katona asked and Bollobás solved the following extremal problem.

\begin{theorem}[Bollobás \ev(1974) \cite{Bollo74Katona}]\label{Bollo74KatTh}
If a 3-graph has $3n$ vertices and $n^3+1$ triples, then there are two triples whose symmetric difference is contained in a third one. 
\end{theorem} 

$\Khyp33(n,n,n)$ -- which generalizes $T_{2n,2}$, to 3-uniform hypergraphs,-- shows the sharpness of Theorem \ref{Bollo74KatTh}.  So Theorem \ref{Bollo74KatTh} is a natural generalization of Turán's theorem: an ordinary triangle-free graph
$G_n$ is just a graph where no edge of $G_n$ contains the symmetric difference of two other edges. So the excluded hypergraph can be viewed as a hypergraph-triangle. Bollobás generalized Katona's conjecture to $r$-uniform hypergraphs. The
generalized conjecture was proved for $r=4$ by Sidorenko \cite{Sidor87Zametki}, for $r=5,6$ by Norin and Yepremyan \cite{NorinYepre17GenTria}. For related results see e.g. Sidorenko \cite{Sidor87Zametki}, Mubayi and Pikhurko
\cite{MubaPikhu07Mantel}, Pikhurko \cite{Pikhu08Exact}.
%% Expand?

\subsubsection{The Fano hypergraph extremal problem}\label{FanoExtreS}

\BalAbraCap{fanohyp2}{Fano} Here the excluded graph is the 3-uniform Fano hypergraph $\cFa$ on 7 vertices, with seven hyperedges any two of
which intersect in exactly one vertex, see Figure \ref{fanohyp2}: $\cFa$ is the simplest finite geometry.
The nice thing about the extremal problem of $\cFa$ is that it is natural, non-trivial, but has a nice solution.

\begin{conjecture}[V.~T.~Sós]
  Partition $n$ vertices into two classes $\cA$ and $\cB$ with $||\cA|-|\cB||\le 1$ and take all the triples intersecting both $\cA$ and $\cB$. The obtained 3-uniform complete bipartite hypergraph $\bH[\cA,\cB]$ is extremal for $\cFa$ (if $n$ is sufficiently large).
\end{conjecture}

Using some multigraph extremal results of Kündgen and Füredi \cite{FureKund02Weight}, first de Caen and Füredi proved 

\begin{theorem}[de Caen and Füredi \ev(2000) \cite{CaenFure00Fano}]
$$\ext(n,\cFa)={3\over 4}{n\choose 3}+O(n^2).$$
\end{theorem}

Next, applying the stability method, the sharp result was obtained, independently, by Füredi and Simonovits and by Keevash and Sudakov. Since $\chi(\cFa)=3$, 
a 3-uniform bipartite hypergraph cannot contain $\cFa$.

\begin{theorem}[Füredi-Simonovits \ev(2005) \cite{FureSim05Fano}/Keevash-Sudakov \cite{KeevSudak05Fano}]\label{FanoExt}
  If $\Hhyp n3$ is a triple system on $n> n_1$ vertices not containing $\cFa$ and of maximum number of hyperedges under this condition, then $\Hhyp n3$ is bipartite: $\chi(\Hhyp n3)=2$.
\end{theorem}

Theorem \ref{FanoExt} implies that
$${\rm ex}_3(n,{\cFa}) = {n\choose 3}-{{\lfloor n/2 \rfloor}\choose 3} -{{\lceil n/2 \rceil}\choose 3}. $$

There are two important ingredients of the proof. The first one is a multigraph extremal theorem:

\begin{theorem}[Füredi-Kündgen \cite{FureKund02Weight}]
 If $M_n$ is an arbitrary multigraph (without restriction on the edge multiplicities, except that they are non-negative) and each
 4-vertex subgraph of $M_n$ has at most 20 edges (with multiplicity), then $$e(M_n)\le 3{n\choose 2}+O(n).$$
\end{theorem}

The other ingredient of the proof of Theorem \ref{FanoExt} was that it is enough to prove the theorem for those hypergraphs where the low-degree vertices are deleted, and it is enough to prove a corresponding stability theorem.

\begin{theorem}
 There exist a $\gamma_2>0$ and an $n_2$ such that if $\cFa\not\subseteq\Hhyp n3$ and $$\deg(x)> \left({3\over 4}-\gamma_2\right){n\choose 2} \Text{for each}x\in V(\Hhyp n3),$$ then $\Hhyp n3$ is bipartite. 
 \end{theorem}

 Recently L. Bellmann and C. Reiher \cite{BellmannReiher18Fano} proved that Theorem \ref{FanoExt} holds for any $n\ge7$.  One could ask, what is the essence of their proof.  Often when we use stability arguments for extremal problems, one thinks that perhaps induction would also work. Unfortunately, for hypergraphs this mostly breaks down. As an exception, Bellmann and Reiher have found the good way to use here induction. (A similar situation was when the Lehel Conjecture was proved by Bessy and
 Thomassé \cite{BessyThomasse10}, yet for hypergraphs we do not know of such cases.)

\minisepa

Füredi, Pikhurko, and Simonovits used the stability method also in \cite{FurePikhSim05}, to prove a conjecture of Mubayi and Rödl. For further related results see, e.g., Keevash \cite{Keev05ProjeTur}, Füredi, Pikhurko, and Simonovits \cite{FurePikhSim05,FurePikhSim06}, \dots and many similar cases.

\minisepa

We should mention here a beautiful result of Person and Schacht:

\begin{theorem}[Person and Schacht \ev(2009) \cite{PersonSchacht09Fano}]\label{PersonSchacht}
Almost all Fano-free 3-uniform hypergraphs are bipartite.
\end{theorem}

\begin{remark}\label{FanoTypi}\rm 
(a) One could ask what is the connection between Theorems \ref{FanoExt} and \ref{PersonSchacht}. Here one should be careful, e.g., Prömel and Steger 
 \cite{PromStege92Berge} proved that almost all Berge graphs are perfect, which means that for most graphs the Berge Strong Perfect Graph Conjecture is true. This
 was a beautiful result, however, the actual proof of the Perfect Graph Conjecture \cite{ChudThomas06Perfect} (which came much later) was much more difficult.

 (b) Several similar results are known, where the typical structure of $\cF$-free graphs are nicely described. For graphs this was discussed in Subsection \ref{ErdKleitRothS}. Several related results were also proved for hypergraphs, e.g., Lefmann, \cite{LefmPerRodlSchacht09}

(c) Some further related results can be found in the papers of Cioabă \cite{Cioaba04FinGeom}, Keevash \cite{Keev05ProjeTur}, Balogh-Morris-Samotij-Warnke \cite{BalMorrSamoWarn16Typic}, or Balogh and Mubayi \cite{BalMubayi11Semi}, and in many further cases.

\end{remark} 

\subsection{Codegree conditions in the Fano case}\label{Codegree}

As we have already discussed, as soon as we move to hypergraphs, many notions, e.g., the notion of path, cycles and of degrees also can be defined in several different ways. Restricting ourselves to the simplest case of 3-uniform hypergraphs, let $\codeg(x,y)$ be the number of vertices $z$ for which $(x,y,z)$ is a hyperedge, and $\minco(\Hhyp n3)=\min\codeg(x,y)$. Usually $\codeg(x,y)$ is called the \Sc codegree of $x$ and $y$. Often it is more natural to have conditions on the minimum codegree than on the min-degree.

\begin{theorem}[Mubayi \cite{Mubayi05FanoCodeg}]
 For any $\alpha>0$, there exists an $n_0(\alpha)$ for which, if $n>n_0$ and $\minco(\Hhyp n3)\ge(\half +\alpha) n$, then $\Hhyp n3$ contains the Fano plane $\cFa$.
\end{theorem}

The sharpness comes from the complete bipartite 3-uniform hypergraph.  Mubayi conjectured that, for large $n$, $(\half+\alpha)$ can be replaced by $\half$. This was proved first by Keevash, and then, in a simpler way, by DeBiasio and Jiang:

\begin{theorem}[Keevash \cite{Keev09GenTur}/DeBiasio and Jiang \cite{DeBiasJiang14Fano}] There is an $n_0$ such that
if $n>n_0$ and $\minco(\Hhyp n3)>\lfloor\half n\rfloor $ then $\Hhyp n3$ contains $\cFa$.
\end{theorem}

Actually, there are more and more extremal results where we derive the existence of some substructure by knowing that the minimum codegree is high. Later we shall see several further such examples, e.g., in Section \ref{HyperLargeExcludedS}: to ensure a \Sc tight Hamiltonian cycle,\footnote{defined after Definition \ref{k-ellCycle}.} Theorem \ref{RRRSchSz} will use the minimum vertex-degree, while Theorem \ref{RodlRucSzemDirCodeg} uses that the minimum codegree is large.

\begin{remark} 
(a) DeBiasio and Jiang gave a fairly elementary proof, not using the Regularity Method. They had to assume that $n$ is sufficiently large only to use some supersaturation argument.

(b) In their paper DeBiasio and Jiang \cite{DeBiasJiang14Fano} start with a nice introduction and historical description of the situation.

(c) Both the Keevash paper and the DeBiasio-Jiang paper go beyond just proving the above formulated extremal result. 
\end{remark}

\section{Large excluded hypergraphs}\label{HyperLargeExcludedS}

While for ordinary graphs the definition of cycles is very natural, for hypergraphs there are several distinct ways to define them, leading to completely different problems and results. Below mostly we shall restrict ourselves to 3-uniform hypergraphs, within that to \Sc Loose (also called \Sc Linear) and \Sc Tight cycles, and we mostly skip Berge cycles.\footnote{Hamiltonicity for Berge cycles were discussed by Bermond, Germa, Heydemann, and Sotteau \cite{BermondGermaHeydSotte76Orsay}. Perhaps the first Tight Hamiltonicity was discussed in \cite{KatonaKier99Hamil}, however, the tight Hamilton cycles were called there as Hamilton chains. }  The reader interested in more details about Berge cycles is referred to Gyárfás, Sárközy, and Szemerédi \cite{GyarSarkoSzem10Berge,GyarfSarko11BergRams} and also to the excellent surveys of Rödl and Ruciński \cite{RodlRuc10HyperSurv}, Kühn and Osthus \cite{KuhnOsthus14ICM}.

We start this section with a short subsection on Hypergraph cycle Ramsey problems which could be regarded as medium range extremal problems, since we get in a hypergraph on $N$ vertices a (monochromatic) subgraph of $(1-c)N$ vertices.
Then we shall consider hypergraph extremal problems where the excluded configuration is a spanning subhypergraph, or at least an {\em almost} spanning one. 
To make the life of the reader (interested in more details) easier, we follow the definitions and notation of \cite{RodlRuc10HyperSurv}, as much as we could.

\subsection{Hypergraph cycles and Ramsey theorems}\label{HyperRamseyS}

We shall need the definition of the so called $(k,\ell)$-cycles, covering the matchings, the \Sc loose cycles and the \Sc tight cycles.

\begin{definition}[Tight/loose cycles]\label{HyperCycD} Consider 3-uniform hypergraphs. Let $E_1,\dots,E_\ell$ be a cyclically arranged family of triples. 
  If the consecutive ones intersect in 2 vertices, and there are no other intersections among them, then this configuration $\Cthyp \ell3$ will be called a \Sc Tight $\ell$-cycle.  If the consecutive ones intersect in 1 vertex, and there are no other intersections among them, then this configuration $\Chyp \ell3$ will be called a \Sc Loose $\ell$-cycle.
\end{definition}

\bigskip

We could have started with the more general

\begin{definition}[$(k,\ell)$-cycle]\label{k-ellCycle}
  Fix some $0\le \ell\le k-1$. In a $k$-uniform hypergraph $\Hhyp nk$ the cyclically ordered vertices $a_1,\dots,a_t$ and the hyperedges $E_1,\dots,E_t$ form a $(k,\ell)$-cycle, if the vertices of $E_i$ are cyclically consecutive (form a segment), and $|E_i\cap E_{i+1}|=\ell$, for $i=1,\dots,t$. (Here $E_{t+1}:=E_1$).

  If $\ell=1$, that we call a \Sc loose cycle, (or \Sc linear cycle) and if $\ell=k-1$ that is the \Sc tight cycle.
\end{definition}

For $\ell=0$, this reduces to a matching: a family of $t/k$ independent edges.

\subsub Divisibility.// Mostly we have some hidden divisibility conditions, e.g., speaking of a $(k,\ell)$-cycle $\Chyp m{k,\ell}$ on $m$ vertices, we assume that $m$ is divisible by $k-\ell$.

\minisepa

It was known from the beginning that the (ordinary) Ramsey number for cycles strongly depends of the parity: $R(C_n,C_n)={3\over 2}n-1$ if $n$ is even and $R(C_n,C_n)=2n-1$ if $n$ is odd (Bondy \cite{BondyErd73Rams}, Faudree and Schelp
\cite{FaudSchelp74DM}, Vera Rosta \cite{Rosta73}). So it is not too surprising that for hypergraphs the parity also strongly influences the results. (Analogous results are known for ordinary graphs and 3
or more colours, e.g., as we have mentioned, $R(n,n,n)=4n-3$ for odd $n>n_0$, see Łuczak, \cite{Lucz99BondyErd}, Kohayakawa, Simonovits, and Skokan \cite{KohaSimSkok05Graco}, but it is $2n+o(n)$ if $n$ is even, see Figaj and Łuczak
\cite{FigajLucz07Even,FigLucz18RamsCyc}.)\footnote{The breakthrough result of \L{u}czak \cite{Lucz99BondyErd} was improved by Kohayakawa, Simonovits, and Skokan \cite{KohaSimSkok05Graco} and generalized by Jenssen and Skokan
  \cite{JenssSkok16Ramsey}, see also Gyárfás, Ruszinkó, Sárközy and Szemerédi \cite{GyarfRuszSarkoSzem08-09} and Benevides and Skokan \cite{BenevSkokan09} and others.}  For hypergraphs the situation is even more involved, since the Ramsey numbers depend on the types of cycles we consider (loose, tight, Berge). Haxell, Łuczak, Peng, Rödl, Ruciński, Simonovits, and Skokan proved the following

\begin{theorem}[\ev(2006) \cite{HaxLuczPengSimSkok06}]
 Consider 3-uniform hypergraphs. If $\Chyp n3$ denote the \Sc loose $n$-vertex hypergraph cycle, then $R(\Chyp n3,\Chyp n3,\Chyp n3)=5n/4+o(n)$.
\end{theorem}

Here the upper bound is the difficult part: an easy construction shows that $R(\Chyp n3,\Chyp n3,\Chyp n3)>5n/4-c$, for an appropriate constant $c>0$.
This was generalized from 3-uniform to $k$--uniform graphs and loose cycles by Gyárfás, Sárközy, and Szemerédi \cite{GyarfSarkoSzeme08}. 
As to the tight cycles, Haxell, Łuczak, Peng, Rödl, Ruciński, and Skokan proved the following

\begin{theorem}[\ev(2006) \cite{HaxLuczPengSkok09Tight}]
  Consider 3-uniform hypergraphs.  If\, $\Cthyp n3$ denotes the \Sc tight $n$-vertex hypergraph cycle, then $R(\Cthyp n3,\Cthyp n3,\Cthyp n3)=4n/3+o(n)$ when
  $n$ is divisible by 3, and $\approx 2n$ otherwise.
\end{theorem}

On Ramsey numbers for Berge cycles see the papers of Gyárfás, Sárközy, and Szemerédi \cite{GyarfSarkoSzemer10Berge}, and Gyárfás and Sárközy
\cite{GyarfSarko11BergRams}, or Gyárfás, Lehel, G. Sárközy, and Schelp \cite{GyarfLehSarkoSchelp08BergeRams}.

\subsection{Hamilton cycles}

\noindent This research area has two roots: 

(A) For simple graphs the extremal graph problem of $k$ independent edges goes back to Erdős and Gallai \cite{ErdGallai59Path}. It may be surprising, but to
ensure $k$ independent edges, or a path $P_{2k}$ requires basically the same degree condition. On the other hand, maybe it is not so surprising.  Rödl and Ruciński \cite{RodlRuc10HyperSurv} write: ``\dots for $\ell=0$ a Hamiltonian $\ell$-cycle in a $k$-graph $H$
becomes a perfect matching in $H$. Moreover, any Hamiltonian $(k-\ell)$-cycle
contains a matching of size $\lfloor n/k\rfloor$. Hence, not surprisingly, the results for Hamiltonian cycles and perfect (or almost perfect) matchings are related.''
The Erdős-Gallai results were extended to hypergraphs, by
Erdős \cite{Erd65Annales}. Bollobás, Daykin, and Erdős \cite{BolloDayErd76} ensured $t$ independent hyperedges, Daykin and Häggkvist \cite{DayHaggkv81Hyper}
guaranteed a perfect matching.

(B) Katona and Kierstead \cite{KatonaKier99Hamil} defined the tight Hamilton cycle and tried to generalize Dirac's theorem to hypergraphs.

The area described in the next few subsections in a fairly concise way is described in much more details, e.g., in the excellent surveys of Kühn and Osthus \cite{KuhnOsthus09Embedding} and of Rödl and Ruciński \cite{RodlRuc10HyperSurv}, in the Introduction of the paper of Alon, Frankl, Huang, Rödl, Ruciński, and Sudakov \cite{AlonFranklHuaRodlRucSudak12Samu}, and in the survey of Yi Zhao \cite{YiZhao16Dirac}.

\subsubsection{The beginnings: Hamiltonicity}\label{HamilHypKezdet}

Since hypergraph problems seemed mostly too technical even for combinatorists working outside of hypergraph extremal problems, the research on hypergraph
Hamiltonicity started relatively late, with a paper of G.Y. Katona and Kierstead \cite{KatonaKier99Hamil}.\footnote{For Berge cycle Hamiltonicity there were
  earlier results, e.g., by Bermond, Germa, Heydemann, and Sotteau \cite{BermondGermaHeydSotte76Orsay}.} For $k$-uniform hypergraphs, for a subset
$\cS\subseteq V(\ckH n)$ with $\ell:=|\cS|$ we define the $\ell$-degree $\de_{\ell}(\cS)$ as the number of $k$-edges of $\ckH n$ containing $\cS$ and
 $\de_\ell(\ckH n)$ is the minimum of $\de_\ell(\cS)$ for the $\ell$-tuples $\cS\subseteq V(\ckH n)$. G.Y.~Katona and Kierstead considered
the tight cycles\footnote{They called it chain.} and paths and proved

\begin{theorem}[G.Y. Katona, H. Kierstead \ev(1999) \cite{KatonaKier99Hamil}]
If $\ckH n=(V,\cE)$ is a $k$-uniform hypergraph and
$$\de_{k-1}(\ckH n)\ge \left(1-\reci2k \right)n+4-k-{5\over2k},$$ then it contains a \Sc tight Hamiltonian cycle $\Cthyp nk$.
\end{theorem}

This is far from being sharp:

\begin{conjecture}[G.Y. Katona, H. Kierstead \cite{KatonaKier99Hamil}]
If $\de_{k-1}(\ckH n)> \left\lfloor{n-k+3\over 2}\right\rfloor,$
then $\ckH n$ has a \Sc tight Hamilton cycle.
\end{conjecture} 

Katona and Kierstead also provided the construction supporting this:

\begin{theorem}[Katona and Kierstead \cite{KatonaKier99Hamil}]
  For any integers $k\ge2$ and $n> k^2$ there exists a $k$-uniform hypergraph $\Hhyp nk$ without \Sc tight Hamilton cycles for which $\de_{k-1}(\ckH
  n)=\left\lfloor{n-k+3\over 2}\right\rfloor$.
\end{theorem}

Rödl and Ruciński write in \cite{RodlRuc10HyperSurv}: 

\begin{quote}{
``In 1952 Dirac \cite{Dirac52Abstr} proved a celebrated theorem\dots ~In 1999, Katona and Kierstead initiated a new stream of research to studying similar questions for hypergraphs, and subsequently, for perfect matchings\dots''}
\end{quote}

\subsection{Problems: Hypergraph Hamiltonicity}\label{ProbHam}

We restrict our attention primarily to 3-uniform hypergraphs.  As Rödl, Ruciński, Schacht, and Szemerédi \cite{RodlRucSchachtSzem17Hamil} describe, here there are (at least) six different questions:

(i) how large \Sc vertex-degree ensures a \Sc tight Hamiltonian cycle, 

(ii) how large \Sc vertex-degree ensures a \Sc loose Hamiltonian cycle, 

(iii) how large \Sc co-degree ensures a \Sc tight Hamiltonian cycle, 

(iv) how large \Sc co-degree ensures a \Sc loose Hamiltonian cycle, 

(v-vi) and how large degrees/co-degrees ensure a \Sc perfect~matching or an almost perfect matching?

In all these problems some divisibility questions should also be handled.

\subsub The goal.// As a first step, we would say that most of the research in this area is related to describing two functions defined for $k$-uniform $n$-vertex graphs, the thresholds for the Hamiltonicity, $h^\ell_d(k,n)$, and for the Matching, $m^r_d(k,n)$:

\begin{quote}{What is the minimum integer $t$ for which, if the $d$-degree 
$\de_d(\Hhyp nk)\ge t$
in a $k$-uniform hypergraph $\Hhyp nk $, then  $\Hhyp nk $ contains
\dori (i) a Hamilton $\ell$-cycle $\Chyp n{k,\ell}$
 \dori (ii) an almost-matching $\cM_n^{k,r}$, leaving out at most $r$ vertices,
\newline respectively.
}\end{quote}

Before going into details, we formulate two typical theorems, for 3-uniform hypergraphs, for \Sc Tight cycles.

\begin{theorem}[Reiher, Rödl, Ruciński, Schacht, and Szemerédi \ev(2016)
\cite{ReiherRodlRucSchachtSzem16Hamil}] For any $\eta>0$, there exists an $n_0(\eta)$ such that if $n>n_0$ and in an $n$-vertex 3-uniform hypergraph $\bH^{(3)}_n$ the minimum degree $$\mindeg(\bH^{(3)}_n) \ge\left({5\over9}+\eta\right){n\choose2},$$ then it contains a \Sc tight Hamiltonian cycle.
The estimate $({5\over9}+o(1)){n\choose2}$ is sharp.
\end{theorem}

The codegree problem is answered by 

\Proclaim Theorem B (Rödl, Ruciński, and Szemerédi \ev(2011) \cite{RodlRucSzemer11DiracHyp}).
 There exists an $n_0$ such that if $n>n_0$ and in a 3-uniform hypergraph $\Hhyp n3$, for any $x\ne y$,
 $$\codeg(x,y) \ge \fele{n},$$
 then $\Hhyp n3$ contains a \Sc tight Hamiltonian cycle. This is sharp: for any $n>4$, there exists a 3-uniform hypergraph with $\codeg(\Hhyp n3)=\fele n -1$, without containing a \Sc tight Hamiltonian cycle.

\silent
\begin{table}[h]
\noindent
{\margofont
\begin{tabular}{|l|c|c|c|c|}
%% \hline year&authors&Topic&Method&Reference\\
\hline
\multicolumn{5}{|c|}{\margofont Matching}\\
\hline
2006& Rödl, Rucinski, Szemerédi &Perfect Matching, mindegree &Absorbing& Annal \cite{RodlRucSzeme06Match}  \\
2009&&codegree&& JCT \cite{RodlRucSzemer09Match}  \\
\hline 
2009&Han, Person, Schacht &Perfect Matching, mindegree &Absorbing&SIDMA \\
&&&& \cite{HanPersSchacht09Match} \\
\hline
\multicolumn{5}{|c|}{\margofont Loose Hamiltonian}\\
\hline 
2014&Czygrinow-Molla&{\margofont Loose Hamilton} 3-unif &Absorbing&SIDMA\\
  & &codegree, perfect matching &Stability&\cite{CzygMolla14Loose} *\\
\hline
2010& Hàn, Schacht & Dirac, for loose Hamilton cycles &Absorbing& \cite{HanSchacht10Hamil}\\
\hline
\multicolumn{5}{|c|}{\margofont Tight Hamiltonian }\\
\hline
2011&Rödl, Ruciński, Szemerédi&Tight Hamiltonian, codegree & $\fele n -1$
& \cite{RodlRucSzemer11DiracHyp}\\
\hline
\end{tabular}
}
\end{table}

Table \ref{HypergraphsB} contains some of the results, discussed below, and some others.

\begin{table}[h]
\noindent
{\margofont
\begin{tabular}{|l|l|l|l|l|}
\hline
Authors& Year &About what? &Methods&\tiny Journal\\
\hline\hline
Hàn, Schacht &2010 & Dirac, loose Hamilton cycles &Absorb.&\cite{HanSchacht10Hamil}\\
\hline
Rödl, Ruciński& 2010 &Dirac type questions, Survey &&\cite{RodlRuc10HyperSurv}\\
\hline
Rödl, Ruciński, & 2011 &Dirac 3-hyper, approx &&Advances\\
Szemerédi& & && \cite{RodlRucSzemer11DiracHyp} \\
\hline
Glebov-Person-Weps&2012&Hypergraph Hamiltonicity & ????.&\cite{GlebovPersWeps12HyperHamil}  \\
\hline 
Czygrinow-Molla&2014 & Loose Hamilton 3-unif &Absorb.&SIDMA\\
  & &codegree, perfect matching &Stabil.&\cite{CzygMolla14Loose} *\\
\hline
Czygrinow,&2014&Tiling with $K_4^{(3)}-2e$&Absorb.&JGT, \\
DeBiasio, Nagle&  &       & & \cite{CzygBiasNagle14Tiling}\\
\hline
Jie Han, Yi Zhao&2015 &Minimum codeg threshold &Absorb.&\cite{HanZhao15TilingK43} \\
\hline
Lo, Allan, &2015 &F-factors in hypergraphs via&Absorb.&GC\\
 Markström & &Absorb. & &\cite{LoMarkstr15Fact}\\
\hline
Reiher-Rödl-Ruciń-&2016 &{\margofont Tight Hamiltonian} 3-hyper &Absorb.&SIDMA\\
ski-Schacht-Szem& & &Absorb.&\cite{ReiherRodlRucSchachtSzem16Hamil}\\
\hline
Ferber, Nenadov,&&Universality of Random &Absorb.&RSA\\
 Peter&& graphs &Absorb.&\cite{FerberNenaPeter} \\
\hline
Rödl, Ruciński,& 2017 &{\margofont Hamiltonicity} of triple&&Annales \\
Schacht,& & systems&&Comb \\
 Szemerédi& & && \cite{RodlRucSchachtSzem17Hamil} \\
\hline
\end{tabular}
}
\caption{Using the absorption method (explained in Subsection \ref{Absorb}), now for hypergraph Hamiltonicity}
\label{HypergraphsB}

\end{table}

\subsection{Lower bounds, constructions}\label{HyperLower}

In all cases considered here we have some relatively simple constructions providing a hypergraph not containing the required configuration and having large
minimum degree (of a given type). The proofs of that {\em ``the constructed hypergraphs do not contain the excluded configurations''} mostly follow from a
Pigeon Hole principle argument or from some parity arguments.  A new feature is (compared to ordinary non-degenerate extremal graph problems) that here the
{\em results} and {\em constructions} often depend on parities, divisibilities, and they may be more complicated.

We start with some constructions, conjectured extremal structures (see e.g. Figure \ref{ExtreConf}).  In the corresponding papers/proofs it turns out that
these are indeed, (almost) extremal constructions. They can be found, e.g., in Han, Person, Schacht, \cite{HanPersSchacht09Match}, or earlier, in Daykin and
Häggkvist \cite{DayHaggkv81Hyper}, Kühn and Osthus \cite{KuhnOsthus06MatchHyper}, Rödl, Ruciński and Szemerédi \cite{RodlRucSzeme06Match}, and Pikhurko
\cite{Pikhu08HyperTiling}.

\begin{construction} Assume that $n$ is divisible by $k$ and $0\le t<k$. Let $\Hhyp nk$ have two vertex classes $\cA$ and $\cB$, 
$|\cA|={n\over k}-1$,
$|\cB|={{k-1}\over k}n+1$,
 and the hyperedges be all the $k$-tuples intersecting $\cA$. This hypergraph has no perfect matchings. Further,
\beq{MinDegMatch}\de_t(\Hhyp nk) =\left(1-\left({k-1\over k}\right)^{k-t}\right)\dbinom{n}{k- t}+o(n^{k- t}). \eeq
\end{construction}

\begin{construction} Assume that $n$ is divisible by $k$ and $0\le t<k$. Let $\tHhyp nk$ have two vertex classes $\cA$ and $\cB$, $|\cA|$ be the maximum odd integer $\le n/2$, $|\cB|=n-|\cA|$, and the hyperedges be all the $k$-tuples intersecting $\cA$ in a positive even number of vertices. This hypergraph has no perfect matchings. Further, \beq{MinDegMatchB}\de_{ t}(\tHhyp nk) =\half{n\choose {k- t}}+O(n^{k- t-1}). \eeq
\end{construction}

These constructions provide the lower bounds in many results in this area. Thus, e.g., for $k=3$ the obtained coefficients of $n\choose {k- t}$ are
$5\over9$ for mindegree and $\half$ for codegree, respectively, which will turn out to be sharp in Theorems \ref{CooleyMycroft17}, \ref{RRRSchSz} and \ref{RodlRucSzemDirCodeg} below.\footnote{identical with Theorems A,B, above.}
 
\subsection{Upper bounds, asymptotic and sharp}

\noindent
Most of the above questions first were solved only in asymptotic forms, then in sharp forms. But, as it is remarked in \cite{RodlRucSchachtSzem17Hamil}, one of these problems, namely Theorem~\ref{RRRSchSz} below, the first one, (i), in the above list, is more difficult than the others. There, even the asymptotic result (i.e. to prove $\approx {5\over9}{n\choose2}+o(n^2)$ ``needed'' four steps: first Glebov, Person and Weps improved the trivial ${n-1\choose 2}$ to $(1-\eps){n-1\choose 2}$, then Rödl and Ruciński \cite{RodlRuc14Hamil}
improved $(1-\eps)$ to $\reci3 (5-\sqrt5)\approx0.91$, next Rödl, Ruciński, and Szemerédi \cite{RodlRucSzemer08ApproxDirac} to 0.8, and only then, they with Reiher and Schacht, obtained the sharp $5/9=0.555$, in \cite{ReiherRodlRucSchachtSzem16Hamil}.

\minisepa

Below we write about this theory, staying mostly with the simplest cases. For $k=3$ the above codegree threshold is $\fele n $ in
\eqref{MinDegMatchB}. One could ask, what is the appropriate degree condition. Cooley and Mycroft \cite{CooleyMycroft17Hamil}, using an appropriate
Regularity Lemma of Allen, Böttcher, Cooley, and Mycroft~\cite{AllenBottCoolMycroft}, proved

\begin{theorem}[Cooley and Mycroft \ev(2017) \cite{CooleyMycroft17Hamil}]\label{CooleyMycroft17} 
For any $\eta>0$ there exists an $n_0(\eta)$ such that if $n>n_0(\eta)$, then
any 3-uniform hypergraph $\Hhyp n3$ with $\mindeg(\Hhyp 3n)\ge(\frac 59+\eta)\binom n2$ contains a \Sc tight cycle $\Cthyp m3$ with $n-m=o(n)$.
\end{theorem}

In other words, $\Hhyp n3$ contains an almost-Hamiltonian cycle: only $o(n)$ vertices are left out. The constant $5\over 9$ is asymptotically best possible. Reiher, Rödl, Ruciński, Schacht, and Szemerédi improved Theorem \ref{CooleyMycroft17}:

\begin{theorem}[\ev(2016) \cite{ReiherRodlRucSchachtSzem16Hamil}]\label{RRRSchSz} For any $\eta>0$, there exists an $n_0(\eta)$ such that if $n>n_0$ and in an $n$-vertex 3-uniform hypergraph $\bH^{(3)}_n$ the minimum degree $$\mindeg(\bH^{(3)}_n) \ge\left({5\over9}+\eta\right){n\choose2},$$ then it contains a \Sc tight Hamiltonian cycle.
\end{theorem}

The sharpness follows from the sharpness of Theorem \ref{CooleyMycroft17}. Actually, \cite{ReiherRodlRucSchachtSzem16Hamil} describes three constructions proving the sharpness of Theorem \ref{RRRSchSz}.  The difference between Theorems \ref{CooleyMycroft17} and \ref{RRRSchSz} is that Theorem~\ref{RRRSchSz} has no left-out vertices. The proof of Theorem \ref{RRRSchSz} uses the hypergraph regularity and then the ``absorption'' method, to pick up the last few vertices.\footnote{This, with a sketch of the proof, is nicely explained by Rödl and Ruciński, in \cite{RodlRuc10HyperSurv}. The introduction of \cite{ReiherRodlRucSchachtSzem16Hamil} and a survey of Yi Zhao \cite{YiZhao16Dirac} are also good descriptions of the otherwise fairly complicated situation.}

As to the codegree, we have

\begin{theorem}[Rödl, Ruciński, and Szemerédi \cite{RodlRucSzemer11DiracHyp}]\label{RodlRucSzemDirCodeg} 
 There exists an $n_0$ such that if $n>n_0$ and in a 3-uniform hypergraph $\Hhyp n3$, for any $x\ne y$,
 $$\codeg(x,y) \ge \fele{n},$$
 then $\Hhyp n3$ contains a \Sc tight Hamiltonian cycle. This is sharp: for any $n>4$, there exists a 3-uniform hypergraph with $\codeg(\Hhyp n3)=\fele n -1$, without containing a \Sc tight Hamiltonian cycle.
\end{theorem}

The sharpness follows from the constructions of Subsection \ref{HyperLower}. A similar result holds for Hamiltonian paths.

\subsubsection{Matchings} 

We define a perfect matching in a $k$-uniform hypergraph $H$ on
$n$ vertices as a set of $\lfloor n/k \rfloor$ disjoint hyperedges.
To ensure a 1-factor in hypergraphs is a fascinating and important topic and would deserve a much longer survey. Here we again restrict ourselves to just a few related results and references. First we formulate two results on the degree-extremal problem of a 1-factor, for 3-uniform and 4-uniform hypergraphs.  Both results are sharp. The following theorem was obtained independently, by Lo and Markström \cite{LoMarkstr15Fact}, Kühn, Osthus, and Treglown and by Imdadullah Khan:

\bigskip
\begin{figure}[h]
\begin{center}
\epsfig{file=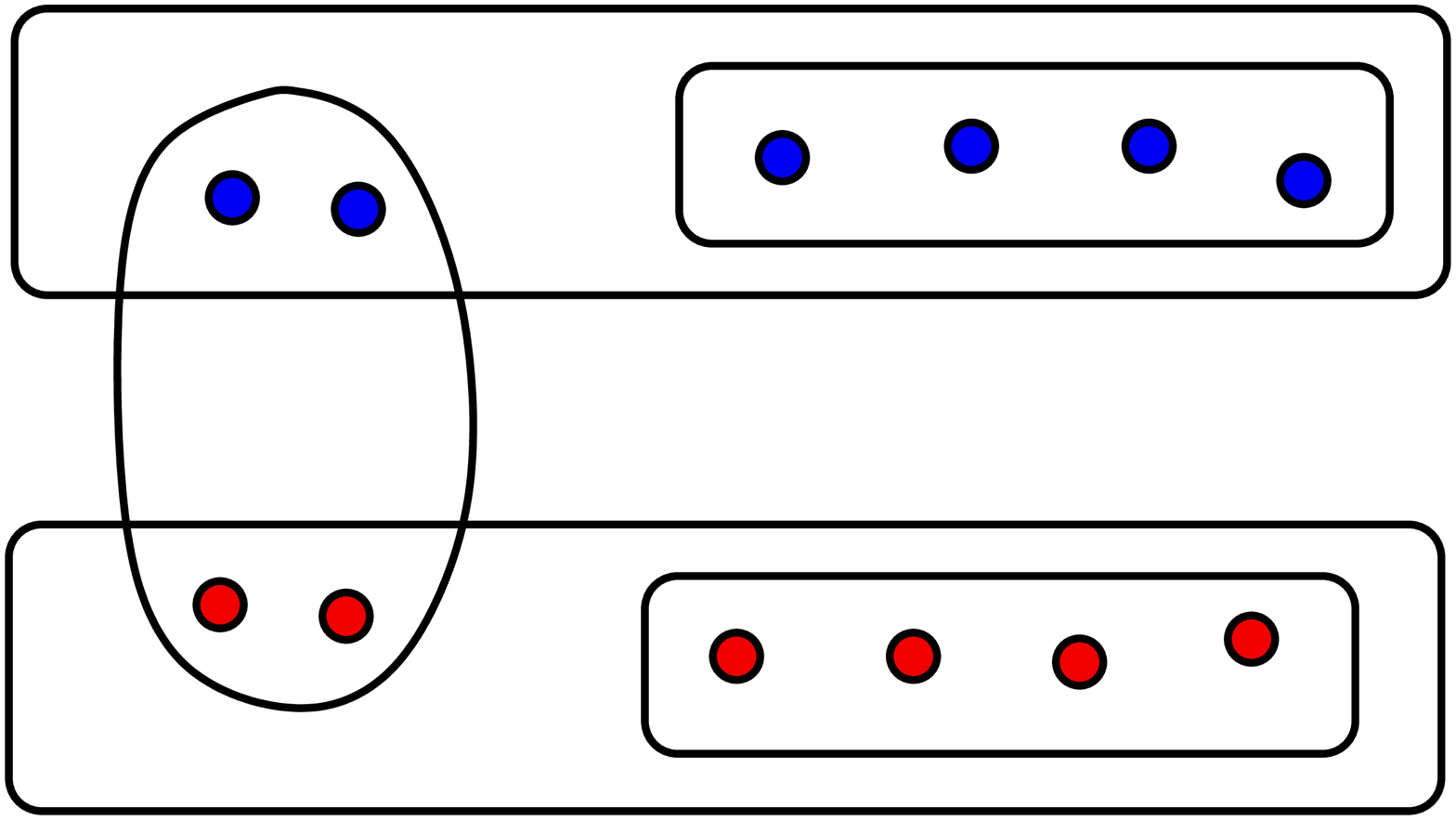,height=17mm}
\qquad \epsfig{file=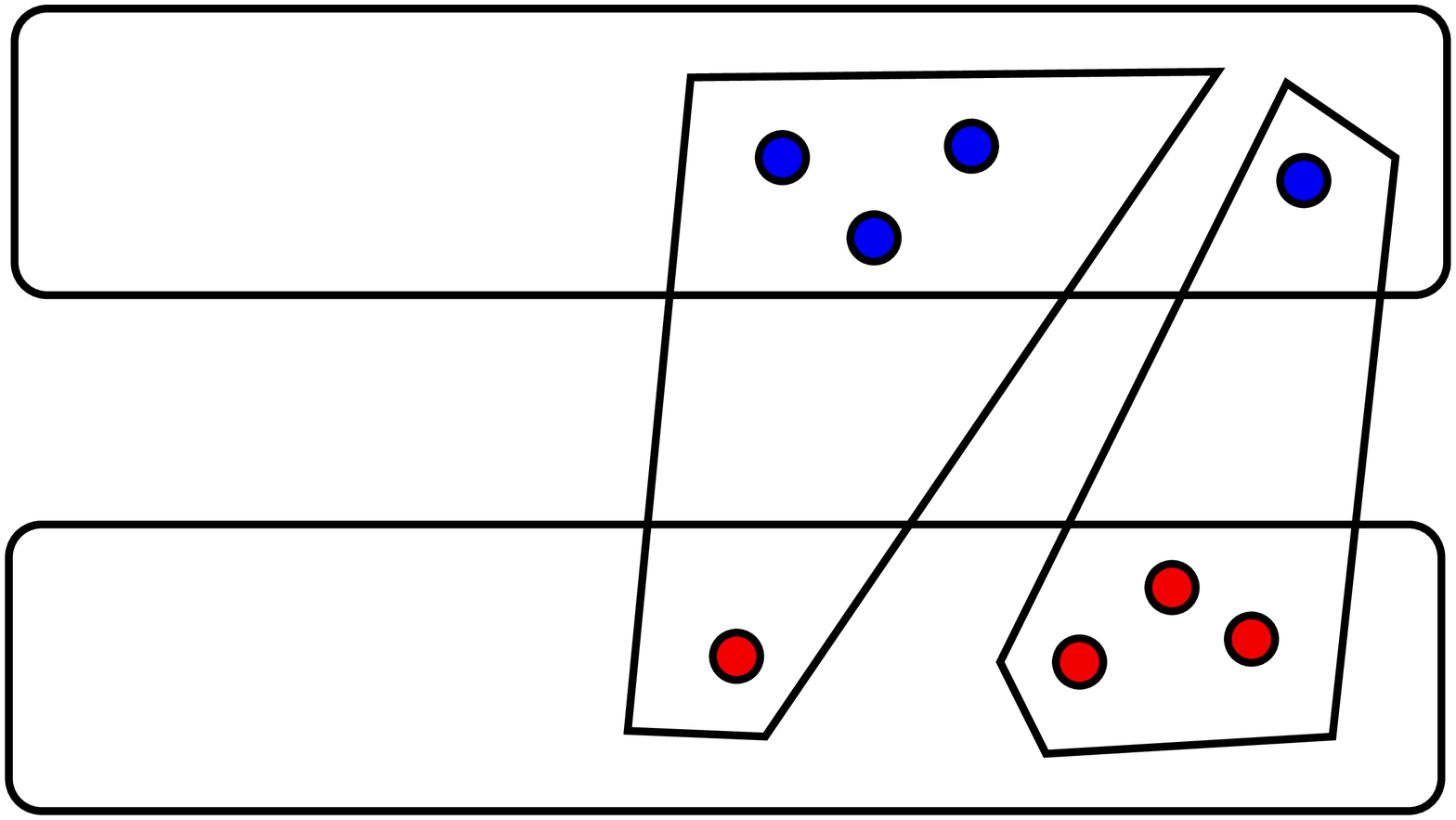,height=17mm}
\end{center}
\caption{Extremal structures for 4-uniform graphs and perfect matching in \cite{CzygKamat12CodegMatch}: (a) 4-tuples intersecting both classes in 2 vertices, (b) 4-tuples intersecting one of the classes in 1 vertex.}
\label{ExtreConf}
\end{figure}

\medskip

\begin{theorem}[Kühn, Osthus, and Treglown \cite{KuhnOsthusTreg13} / I. Khan \cite{Khan13Match3}] There is an $n_0$ such
 that if $n>n_0$ is divisible by 3 and in a 3-uniform hypergraph $\Hhyp n3$
\beq{3Match}\mindeg(\Hhyp n3)>{n-1\choose 2}-{{2n/3}\choose2}+1\approx{5\over9}{n-1\choose2},\eeq then $\Hhyp n3$ contains a 1-factor.
\end{theorem}

Observe that the LHS of \eqref{3Match} is $\approx{5\over9}{n\choose2}$, i.e. the same what we need for a \Sc tight Hamilton cycle.
Khan also proved the analog theorem for 4-uniform hypergraphs:

\begin{theorem}[Imdadullah Khan \cite{Khan16PerfectHyper4}] 
There exists a threshold $n_0$ such that if 
$n>n_0$ is divisible by 4 and in a 4-uniform hypergraph $\Hhyp n4$ 
$$\mindeg(\Hhyp n4)>{n-1\choose 3}-{{3n/4}\choose3}+1, $$ then $\Hhyp n4$ contains a 1-factor.
\end{theorem}

For related results see also Pikhurko \cite{Pikhu08HyperTiling}, Han, Person, and Schacht \cite{HanPersSchacht09Match}, Czygrinow and Kamat
\cite{CzygKamat12CodegMatch} providing sharp results and describing the earlier asymptotic results, and Alon, Frankl, Huang, Rödl, Ruciński, and Sudakov
\cite{AlonFranklHuaRodlRucSudak12Samu}.

\subsubsection{Upper bounds and the Absorbing method}

In the proofs of the results discussed in this section one often uses the \Sc Absorbing~Method, -- described in Section \ref{Absorb}, --  to ensure
\Sc spanning or \Sc almost-spanning substructures, where \Sc Almost-Spanning means a subhypergraph covering the whole hypergraph with the exception of at most $O(1)$ (or, occasionally, $o(n)$) vertices. As we wrote,
in most cases considered in these subsections we have a simple, nice conjectured extremal structure, enabling us to use stability arguments. Yet, the actual
``upper bounds'' (non-construction parts) are fairly complicated. 

In Section \ref{ProbHam} we wrote about the thresholds for the Hamiltonicity, $h^\ell_d(k,n)$, and for the Matching, $m^r_d(k,n)$. We have mentioned that sometimes we get the same results for Hamiltonicity and Perfect Matching.

\begin{conjecture}[Informal/Formal \cite{RodlRuc10HyperSurv}]
  The $d$-degree threshold $\de_d(\Hhyp nk)$ for the Hamiltonian problem and the Matching problem are roughy the same. More formally, $h_d(k,n)\approx m_d(k,n)$, if $d,k$ are fixed, $n\to\infty$.
\end{conjecture}

There are many related results but we skip them, mentioning only that, as it is emphasized in \cite{RodlRuc10HyperSurv}, not only these two quatities are often near to each other, but in several cases they are proved to be equal, see e.g. Theorem 3.1 of \cite{RodlRuc10HyperSurv}, quoted from \cite{RodlRucSzemer09Match}.

\begin{theorem}[Han, Person, and Schacht \ev(2009) \cite{HanPersSchacht09Match}]\label{HanPersSchacht09}
  For any $\gamma>0$, if $\delta_1(\Hhyp n3)\ge({5\over 9}+\gamma ){n\choose 2}$ and $n>n_0(\gamma)$, then $\Hhyp n3$ contains a 1-factor.
\end{theorem}

There is a surprising, sharp difference between the cases when $n$ is divisible by $k$ and when it is not. Consider the simplest hypergraph case, $k=3$. 

\begin{theorem} For any $\gamma>0$, if 
$\delta_2(\Hhyp nk)\ge(\half+\gamma ){n\choose 2}$ and $n>n_0(\gamma)$, then $\Hhyp nk$ contains a \Sc Tight Hamiltonian cycle.
\footnote{This implies Theorem \ref{HanPersSchacht09}.} 
If $n$ is not divisible by $3$ and $\Hhyp n3$ does not contain a 1-factor, then $\de_2(\Hhyp n3)\le {n\over k}+o(n)$.
\end{theorem}

So in case of non-divisibility, much smaller degrees ensure a perfect matching.

\minisepa

One of the first results where the Absorbing Lemma was used is the theorem of Rödl, Ruciński, and Szemerédi \cite{RodlRucSzemer09Match}, on ensuring a perfect matching in a $k$-uniform hypergraph. As we have mentioned, an appropriate construction of Bollobás, Daykin, and Erdős \cite{BolloDayErd76} and further similar constructions shows the sharpness of these results.

There is a sequence of papers of Rödl, Ruciński, and Szemerédi on this topic: 
\cite{RodlRucSzeme06Match} assumes large minimum degree, 
\cite{RodlRucSzemer09Match} assumes large codegree.
To ensure Hamiltonian cycles in hypergraphs, see 
\cite{RodlRucSzem06Dirac3CPC},
\cite{RodlRucSzemer08ApproxDirac},
\cite{RodlRucSzemer11DiracHyp}, and see also
the above authors with Schacht \cite{RodlRucSchachtSzem08Carol},
and Reiher \cite{ReiherRodlRucSchachtSzem16Hamil}.

The problem of 1-factors, or almost 1-factors is, of course, connected to (almost perfect) tilings, and there are several results in that direction, as well,
see, e.g., Markström and Ruciński \cite{MarksRuci11Perfect}, Pikhurko \cite{Pikhu08HyperTiling}, Lo and Markström \cite{LoMark13Codeg,LoMarks14PerfMatch,LoMarkstr15Fact}.
\footnote{The earlier excellent surveys of Füredi \cite{Fure88Match,Fure94Zurich} are primarily on small excluded subgraphs.}

Several results mentioned above or discussed below use the Absorption Method, as is shown in Tables \ref{HypergraphsB}, \ref{Hypergraphs}. (See also the earlier Rödl-Ruciński survey \cite{RodlRuc10HyperSurv} on these hypergraph problems and results.)

\bigskip

The first line of Table \ref{Hypergraphs} differs from the other ones by that it is connected to Graph Packing, and within that to the Sauer-Spencer theorem. The subsequent lines try to show in time order some important papers using the absorption method, on loose and tight hypergraph Hamilton cycles (Table \ref{HypergraphsB}) and matching (Table \ref{Hypergraphs})

\bigskip

These tables also contain some results on random graphs, and also some lines on Universal graphs, but we skip explaning them. For details see the papers of Alon and Capalbo, e.g., \cite{AlonCapal07Univ}, or Alon, Capalbo,
Kohayakawa, Rödl, Rucinski, and Szemerédi \cite{AlonCapalKohaRodlRucSzem01}. Further, Table \ref{GraphAbsorb} has an item on monochromatic cycle partitions as well.

%% \Margo{Rendbetenni a tablat}

\subsection{Tiling hypergraphs}\label{TilingHyperS}

Of course, a matching that (almost) covers a hypergraph $V(\Hhyp kn)$ is an (almost) tiling. One can also ask if for $k$-uniform hypergraphs do we get interesting results when we wish to tile them for a fixed $\bF$ by vertex-disjoint copies of $\bF$. The outstanding results of Keevash \cite{Keev14ExistDesign}\footnote{see also Keevash, \cite{Keev15CountingDesigns}, Barber, Kühn, Lo, and Osthus \cite{BarbKuhnLoOst16}, Barber, Kühn, Lo, Montgomery, and Osthus \cite{BarbKuhnLoMontgoOsthus17},\dots } on designs and the related results are also tiling theorems, however, here we mention only that
 Pikhurko, beside investigating codegree matching problems, in \cite{Pikhu08HyperTiling} also investigated the problem of ensuring $K^{(3)}_4$-tilings in a hypergraph with large codegree. These problems were connected to each other.

\medskip
\begin{table}[h]
\noindent
{\margofont
\begin{tabular}{|l|l|l|l|l|}
\hline
Authors& Year &About what? &Methods&Where\\
\hline\hline
Rödl, Ruciński,& 1999 &Hypergraph packing,  & Mislead-&CPC \\
Taraz& 1999 & embedding & ing!&\cite{RodlRucTaraz99}\\
\hline
Kühn and Osthus&2006& Matching, implied by min&Stability&\cite{KuhnOsthus06MatchHyper}\\
&&deg  in $r$-partite hypergraph&???& JGT \\
\hline
Pikhurko&2008& Matching, Tiling, min codeg &???& GC \cite{Pikhu08HyperTiling}\\
\hline
Han, Person,
& 2009 &Perfect Matching, mindegree &Absorb.&SIDMA \\
 Schacht&&&& \cite{HanPersSchacht09Match} \\
\hline
Czygrinow-Kamat&2012 &{\margofont perfect matching} Sharp&Absorb.&ELECT \\
 & &codegree condition, 4-unif &Stability&JC \cite{CzygKamat12CodegMatch} \\
\hline
Alon, Frankl, Huang, 
&2012 &  Large matchings in uniform 
&Fractional &JCTA\\ 
 Rödl,  Ruciński, 
&&hypergraphs and the conjec-&matching &\cite{AlonFranklHuaRodlRucSudak12Samu} \\
 Sudakov:&&ture of Erdős and Samuels,& & \\
\hline
I. Khan&2013 &Perfect matching, min degree &Absorb.&SIDMA\\
& & 3-uniform &Stability&\cite{Khan13Match3}\\
\hline
I. Khan&2016 &Perfect matching, min degree &Absorb.&JCTB\\
& & 4-uniform &Stability&\cite{Khan16PerfectHyper4}\\
\hline
\end{tabular}
}
\caption{Absorbing for hypergraphs, Matching, \dots}
\label{Hypergraphs}
\end{table}

\subsection{Vertex-partition for hypergraphs}\label{VertexPartiHyper}

The vertex-partition results of Subsection \ref{VertexPartitionS} can also be extended to hypergraphs, see Gyárfás and G.N. Sárközy \cite{GyarfSarko13MonoHyper,GyarfSarko14Loose} and of G.N. Sárközy \cite{Sarko14ImprovedLoose}. Thus, e.g., G.N. Sárközy proved the following.

\begin{theorem}[G.N. Sárközy \ev(2014) \cite{Sarko14ImprovedLoose}]
  For all integers $k,r\ge2$, there exists an $n_0=n_0(k,r)$ such that if $n>n_0(k,r)$ and $\Khyp nk$ is $r$-edge-coloured, then its vertex set can be partitioned into at most $50rk\log(rk)$ vertex disjoint \Sc loose monochromatic cycles.
\end{theorem}

For the case of two colours and loose/tight cycles Bustamante, Han, and M. Stein \cite{BustaHanStein18AlmostParti} proved some hypergraph versions of Lehel's conjecture, where, however, a few vertices remain uncovered. Bustamante, Corsten, Nóra Frankl, Pokrovskiy and Skokan \cite{BusCorsFraPokroSkok19A} proved the tight-cycle version.

\begin{theorem}
  There exists a constant $c(k,r)$ such that if $\ckH {}$ is a $k$-uniform hypergraph whose edges are coloured by $r$ colours, then $V(\ckH {})$ can be partitioned into $c(k,r)$ subsets, each defining monochromatic \Sc Tight hypergraph-cycles.
\end{theorem}

For some related results see also the surveys of Gyárfás \cite{Gyarf16VertexCoverSurv}, and Fujita, H. Liu, and Magnant \cite{FujiLiuMagnant15Mono}.

\subsection{Generalizing Gyárfás-Ruszinkó-Sárközy-Szemerédi theorem to hypergraphs}\label{HyperCover}

We have mentioned that several graph results were generalized to hypergraphs.
Theorem \ref{GyarRuszSarkoSzemTh} was generalized to hypergraphs first by Gyárfás and G.N. Sárközy \cite{GyarfSarko14Loose}: they covered the coloured $\Khyp nk$ by Berge paths, by Berge cycles, and by \Sc loose cycles. Perhaps their \Sc loose cycle result was the deepest. Its proof followed the method of Erdős, Gyárfás, and Pyber \cite{ErdGyarfPyb91} and used the linearity of Ramsey number for a ``crown'' that was a generalization of the Triangle-Cycle used in Subsection \ref{Absorb}. The \Sc loose cycle result was improved by G.N. Sárközy:

\begin{theorem}[\cite{Sarko14ImprovedLoose}] For any integers $r,k\ge2$, there exists a threshold $n_0(k,r)$ for which, if $n>n_0(k,r)$, then for any $r$-edge-colouring of a $k$-uniform complete hypergraph $\Khyp nk$, we can partition the vertices into at most $50rk\log(rk)$ classes $U_i$, each defining a monochromatic \Sc Loose $k$-cycle.
\end{theorem}

\subsection{A new type of hypergraph results, strong degree}

This subsection is motivated by a paper of Gyárfás, Győri, and Simonovits \cite{GyarfGyoriSim16}. Their original motivation was to prove the
following conjecture that is still open.

\begin{conjecture}[Gyárfás--G.N. Sárközy, \cite{GyarfSarko14Loose}]\label{GyarSar}
One can partition the vertex set of every $3$-uniform hypergraph $\Hhyp n3$ into $\alpha(\Hhyp n3)$ linear (i.e. \Sc loose) cycles, hyperedges and subsets of hyperedges.
\end{conjecture}

A theorem of Pósa \cite{Posa64Partition} is

\begin{theorem}[Pósa \ev(1964) \cite{Posa64Partition}]
For every graph $G$ one can partition  $V(G)$ into at most $\alpha(G)$ cycles, where a vertex or an edge is accepted as a cycle.
\end{theorem}

Conjecture \ref{GyarSar} would extend Pósa theorem, see \cite{Lov79CombExerc} from graphs to $3$-uniform hypergraphs. Conjecture \ref{GyarSar} was proved in
\cite{GyarfSarko14Loose} for ``weak cycles'' instead of linear cycles, where ``weak cycle'' differs from a ``loose cycle'' by that the consecutive edges
must intersect, but not necessarily in 1 vertex. Here one has to consider subsets of hyperedges also as cycles, in Conjecture \ref{GyarSar}, as shown, e.g.,
by the complete hypergraph $\KK35$.

The following weaker version of Conjecture \ref{GyarSar} was proved recently:

\begin{theorem}[Ergemlidze, Győri, and Methuku \cite{ErgemGyoriMethu18}]
One can cover the vertices of any $3$-uniform hypergraph $\Hhyp n3$ by $\alpha(\Hhyp n3)$ {\bf edge-disjoint} linear (i.e. \Sc loose) cycles, hyperedges and subsets of hyperedges.
\end{theorem}

\minisepa

Let $\rho(\Hhyp n3)$ denote the minimum number of linear cycles, hyperedges or subsets of hyperedges needed to partition $V(\Hhyp n3)$ as described in Conjecture \ref{GyarSar} and let $\chi(\Hhyp n3)$ denote the chromatic number of $\Hhyp n3$, the minimum number of colors in a vertex coloring of $\Hhyp n3$ without monochromatic edges. The following result proves that Conjecture \ref{GyarSar} is true if there are no linear cycles in $\Hhyp n3$.

\begin{theorem}[Gyárfás-Győri-Simonovits \cite{GyarfGyoriSim16}]\label{speccase}
  If $\Hhyp n3$ is a $3$-uniform hypergraph without linear cycles, then $\rho(\Hhyp n3)\le \alpha(\Hhyp n3)$. Moreover, $\chi(\Hhyp n3)\le 3$.
\end{theorem}

The family of hypergraphs without linear cycles seems to be intriguing. \cite{GyarfGyoriSim16} uses a {\bf new degree concept:} the \Sc strong~degree. Let $\Hhyp n3=(V,\cE)$ be a $3$-uniform hypergraph, for $v\in V$ the \Sc link~graph of $v$ in $\Hhyp n3$ is the graph with vertex set $V-v$ and edge set $\{(x,y): (v,x,y) \in \cE \}$. The \Sc strong~degree $d^+(v)$ for $v\in V$ is the maximum number of independent edges (i.e. the size of a maximum matching) in the link graph of $v$. 
%% The {\em underlying} graph of $\Hhyp n3$ is the ordinary graph $G_n$ the edges of which are the pairs covered by the hyperedges of $\Hhyp n3$. 
The main results of \cite{GyarfGyoriSim16} are motivated by the following trivial assertions: a graph
of minimum degree 2 contains a cycle; if $G_n$ has no cycles then $\alpha(G_n)\ge n/2$.

\begin{theorem}[Gyárfás-Győri-Simonovits \cite{GyarfGyoriSim16}]\label{min3} If $\Hhyp n3$ is a $3$-uniform hypergraph with $d^+(v)\ge 3$ for all $v\in V$,  then $\Hhyp n3$ contains a linear cycle.
\end{theorem}

Theorem \ref{min3} is ``self-improving'':

\begin{theorem}[\cite{GyarfGyoriSim16}]\label{min3S}
  Suppose that $\Hhyp n3$ is a $3$-uniform hypergraph with $d^+(v)\ge 3$ for all but at most one $v\in V$. Then $\bH$ contains a linear cycle.
\end{theorem}

The proofs are not easy.  One thinks that several interesting embedding results can be proved where one uses $\min_{v\in V(G)} d^+(v)$ instead of $\mindeg(G)$.

\bigskip

\subsub Acknowledgement.// We thank to several friends and colleagues the useful discussions about the many topics discussed in this survey. We thank above all to József Balogh and András Gyárfás, and also to Zoltán Füredi, János Pach, Jan Hladk\'y, Zoltán L. Nagy, János Pintz, Andrzej Ruciński, Imre Ruzsa, Gábor Sárközy, Andrew Thomason, and Géza Tóth.

\margonum=0
\input SimSzemRefSub3    %% References

\end{document}

%% file: SimSzemRefSub3.tex
%% 2019 febr 12

\gdef\germ#1{{\Rd\mathbb #1}}

\gdef\Xilent#1\par{\par} % Vissza-kikapcsolas

\gdef\Mutato#1{\marginbox{\fbox{#1}}}

\def\MargoF#1{}

\gdef\MR#1 {\fbox{\margofont MR#1}}
\gdef\MR#1 {}

\gdef\Megvan#1{}

\newcount\eldob

\gdef\bibiteX#1#2:#3\par{\bibitem{#1}{\global\advance \eldob by1
 \tiny #1 \fbox{Elhagy:~ \the\eldob} #2} \par }

\newcount\cit

% =================================================================

\newcount\unuse

\gdef\bibitemS#1{\global\advance\unuse by1  
\smallskip {[\the\unuse]:~[#1]:\qquad }}

%% \newpage

%% \gdef\DenseBox#1{}

\silent

\gdef\APr{}
\gdef\FGC{}
\gdef\Megvan#1{}
\gdef\Design{}
\gdef\RegLemma{}
\gdef\Tiling{}

\gdef\Cn{{\Rd C_n}}